\documentclass[a4paper,twoside]{article}

\usepackage{color}
\usepackage{makeidx}         
\makeindex   
\usepackage{graphicx}        
\usepackage[all]{xy}
\usepackage{amssymb}
\usepackage{amsmath}
\usepackage{amsthm}
\usepackage{hhline}
\usepackage{longtable}
\usepackage[a4paper, left=2cm, right=5cm, top=2cm,bottom=2cm]{geometry}
\usepackage{hyperref}

\theoremstyle{definition}
\newtheorem{theorem}{Theorem}[section]
\newtheorem{definition}[theorem]{Definition}
\newtheorem{remark}[theorem]{Remark}

\newtheorem{proposition}[theorem]{Proposition}
\newtheorem{lemma}[theorem]{Lemma}
\newtheorem{example}[theorem]{Example}%
\newtheorem{conjecture}[theorem]{Conjecture}%
\newtheorem{corollary}[theorem]{Corollary}

\usepackage[all]{xy}
\usepackage{amssymb}
\usepackage{amsmath}
\newcommand{\Span}{\mathrm{span}\,}
\newcommand{\Tor}{\mathrm{Tor}\,}
\newcommand{\Der}{\mathrm{Der}}
\newcommand{\KK}{\mathbb{K}}
\newcommand{\fsl}{\mathfrak{sl}}
\newcommand{\Eu}{\mathrm{Eu}}
\newcommand{\fgl}{\mathfrak{gl}}
\newcommand{\sq}{\mathrm{sq}}
\newcommand{\Hom}{\mathrm{Hom}}
\newcommand{\Pf}{\mathrm{Pf}\,}
\newcommand{\Jac}{\mathrm{Jac}}
\newcommand{\Grass}{\mathrm{Grass}}
\newcommand{\im}{\mathrm{im}\,}
\newcommand{\RR}{\mathbb{R}}
\newcommand{\Id}{\mathrm{Id}}
\newcommand{\PP}{\mathbb{P}}

\newcommand{\NN}{\mathbb{N}}
\newcommand{\Spec}{\mathrm{Spec}\,}
\newcommand{\CC}{\mathbb C}
\newcommand{\sk}{\mathrm{sk}}
\newcommand{\sym}{\mathrm{sym}}
\newcommand{\QQ}{\mathbb Q}
\newcommand{\ZZ}{\mathbb Z}
\newcommand{\GL}{\mathrm{GL}}
\newcommand{\SL}{\mathrm{SL}}
\newcommand{\OO}{\mathcal{O}}
\newcommand{\m}{\mathfrak{m}}
\newcommand{\height}{\mathrm{height}\,}

\newcommand{\grade}{\mathrm{grade}\,}

\newcommand{\Kosz}{\mathrm{Kosz}}
\newcommand{\rank}{\mathrm{rank}\,}

\newcommand{\D}{\operatorname{d}\!}

\begin{document}

\title{Determinantal singularities}

\author{\begin{tabular}{ll}
        Anne Fr\"uhbis-Kr\"uger & Matthias Zach\footnote{Current address: Department of Mathematics, 
	University of Kaiserslautern, Erwin-Schr\"odinger-Str., 67663 Kaiserslautern, Germany}\\ 
        Inst. f. Mathematik & Inst. f. Alg. Geom\\ 
        Universit\"at Oldenburg & Leibniz Univ. Hannover\\
        26129 Oldenburg & 30167 Hannover\\
        Germany & Germany \\
        anne.fruehbis-krueger@uol.de & zach@math.uni-hannover.de
        \end{tabular}
}

\vspace{2cm}
\maketitle

\vspace{3cm}

\abstract{
  We survey determinantal singularities, their deformations, and their topology. This class 
  of singularities generalizes the well studied case of complete intersections in several 
  different aspects, but exhibits a plethora of new phenomena such as for instance non-isolated 
  singularities which are finitely determined, 
  or smoothings with low connectivity; already the union of the coordinate 
  axes in $(\CC^3,0)$ is determinantal, but not a complete intersection. We start with
  the algebraic background and then continue by discussing the subtle
  interplay of unfoldings and deformations in this setting, including a 
  survey of the case of determinantal hypersurfaces, 
  Cohen-Macaulay codimension $2$ and Gorenstein codimension $3$ singularities, 
  and determinantal rational surface singularities. 
  We conclude with a discussion of essential
  smoothings and provide an appendix listing known classifications of simple
  determinantal singularities.
}

\newpage 

\tableofcontents

\newpage

\section{Singularities of matrices and determinantal varieties}
\label{sec:1}

Before this article can focus on geometric and topological properties 
of determinantal singularities, the first section starts by presenting 
background material from commutative algebra. This material is indispensible 
for understanding some common structural properties of this class. 
Already in the study of the much smaller class of complete intersection 
singularities, singled out by the algebraic property that their ideal 
possesses a set of generators which form a \textit{regular sequence}\footnote{A 
sequence of elements $a_1,\dots,a_n \in R$ is called \textit{weakly regular} 
on a module $M$ if multiplication by $a_{i+1}$ is 
  injective on $M/\langle a_1,\dots,a_i\rangle M$
  for every $i<n$. It is called \textit{regular} if, moreover, $M/\langle a_1,\dots,a_n\rangle M \neq 0$,
  cf. \cite[Definition 1.1.1]{BrunsHerzog93}.}, a similar 
approach to the exposition of material has at times been used.

\subsection{Determinantal Ideals} 
\label{sec:DeterminantalIdeals}

Let $R$ be a commutative ring with unity. We write $R^{m \times n}$ 
for the space of $m\times n$-matrices with entries in $R$. 

\begin{definition}
  Let $A \in R^{m\times n}$ be a matrix and $I \subset R$ the 
  ideal generated by the $t$-minors of $A$. Then we say that 
  $I$ is \textit{determinantal of type $(m,n,t)$ or} \textit{determinantal 
  of type $t$ with matrix} $A$.
\end{definition}
Note that the same ideal can be determinantal in various different 
ways, that is for different matrices. 

\begin{example}
  \label{exp:VariousDeterminantalStructures}
  The principal ideal generated by $f = xw-yz \in {\mathbb C}\{x,y,z,w\}$ 
  is determinantal of type $1$ with the $1\times 1$-matrix $(f)$ but also 
  determinantal of type $2$ with matrix 
  \[
    \begin{pmatrix}
      x & y \\
      z & w
    \end{pmatrix}.
  \]
  Another famous example
  is due to Pinkham \cite{Pinkham73}. 
  He observed that the ideal 
  $I \subseteq {\mathbb C}\{x_0,\dots,x_4\}$ given by the $2$-minors of the matrix
  \[
    \begin{pmatrix}
      x_0 & x_1 & x_2 & x_3 \\
      x_1 & x_2 & x_3 & x_4
    \end{pmatrix}
  \] 
  can also be generated by the $2$-minors of 
  \[
    \begin{pmatrix}
      x_0 & x_1 & x_2 \\
      x_1 & x_2 & x_3 \\
      x_2 & x_3 & x_4
    \end{pmatrix}.
  \] 
  Hence, $I$ is determinantal of type $(2,4,2)$ or of type $(3,3,2)$, depending 
  on the choice of the matrix. 
\end{example}
Later on, we will encounter cases of 
ideals which are determinantal in a unique way (cf. Theorem \ref{thm:HilbertBurch}),
but the examples given here show that 
in general, one really has to specify the matrix in order to describe the 
determinantal structure of a 
given ideal. 

\medskip

We now introduce some notation to keep the presentation readable in what follows. 
For a matrix $A \in R^{m\times n}$, we denote by $\langle A \rangle$ 
the ideal generated by the entries $a_{i,j}$ of $A$ in $R$. 
Any such matrix $A$ 
can be understood as a homomorphism of free modules 
$A \colon R^n \to R^m$ taking $e_i$ to  $\sum_{j=1}^m a_{i,j} \cdot f_j$,
where $\{e_i\}_{i=1}^n$ is a free basis of $R^n$ and $\{f_j\}_{j=1}^m$ of 
$R^m$. 
For any number $t$ there is a natural induced morphism on the exterior 
powers\footnote{For a discussion of the exterior algebra of a module, 
  see e.g. \cite[Appendix 2.3]{Eisenbud95}}
which we denote by 
\begin{equation}
  A^{\wedge t} \colon \bigwedge^t R^n \to \bigwedge^t R^m.
  \label{eqn:InducedHomomorphismExteriorPowers}
\end{equation}
With the same free bases for $R^m$ and $R^n$ as before,
the products $e_{i_1} \wedge e_{i_2} \wedge \dots \wedge e_{i_t}$ 
with $0< i_1 < i_2 < \dots < i_t \leq n$ and 
$f_{j_1} \wedge \dots \wedge f_{j_t}$ with $0 < j_1 < \dots < j_t \leq m$
form free bases of the free modules $\bigwedge^t R^n$ and 
$\bigwedge^t R^m$ respectively. We will write 
 ${\bf i} = (i_1<i_2<\dots< i_t)$ for the 
\textit{ordered multiindex} of length $\# {\bf i} = t$ and $e_{\bf i}$ 
for corresponding generator.
With this notation we find 
\[
  \bigwedge^t R^n \to \bigwedge^t R^m, \quad e_{\bf i} \mapsto \sum_{\# {\bf j} = t} A^{\wedge t}_{{\bf i} ,{\bf j}} f_{\bf j}
\]
where by $A^{\wedge t}_{{\bf i},{\bf j}}$ we denote the determinant of 
the $t\times t$-submatrix $A_{{\bf i},{\bf j}}$ of $A$ 
specified by the selection of columns in ${\bf i}$ and rows in ${\bf j}$. 
Thus $(A^{\wedge t}_{{\bf i},{\bf j}})$ is the matrix for the 
induced homomorphism (\ref{eqn:InducedHomomorphismExteriorPowers}) 
on the exterior powers  w.r.t. the chosen bases and the ideal of $t$-minors 
of $A$ is nothing but $\langle A^{\wedge t} \rangle$.

\begin{lemma}
  \label{lem:InvarianceOfDeterminantalIdealUnderLeftRightMultiplication}
  Let $R$ be a ring and $A \in R^{m\times n}$ a matrix. For any 
  pair of invertible matrices $P \in \GL(m;R)$ and $Q \in \GL(n;R)$ 
  and every number $t$ one has 
  \[
    \left\langle \left( P \cdot A \cdot Q^{-1} \right)^{\wedge t} \right\rangle 
    = \langle A^{\wedge t} \rangle.
  \]
\end{lemma}

\begin{proof}
  This is immediate for the case of $1$-minors, i.e. the ideal of entries of $A$. 
  For the general case note that 
  \[
    \left( P \cdot A \cdot Q^{-1} \right)^{\wedge t} = P^{\wedge t}\cdot A^{\wedge t} \cdot 
    \left( Q^{-1} \right)^{\wedge t} 
  \]
  where $P^{\wedge t} \in \GL(M;R)$ and 
  $(Q^{-1})^{\wedge t} \in \GL(N;R)$ for
  $N={ n \choose t }$ and $M={m \choose t}$. Since the $t$-minors of 
  $A$ are the entries of $A^{\wedge t}$, this reduces the problem to the 
  case of $1$-minors.
\end{proof}

\begin{remark}
  \label{rem:ReductionOfUnitsInTheMatrix}
  As a consequence of Lemma \ref{lem:InvarianceOfDeterminantalIdealUnderLeftRightMultiplication} 
  we note the following. Suppose one of the entries 
  of the matrix $A$ is a unit in $R$.
  Then there exist matrices $P$ and $Q$ as above such that 
  \[
    P \cdot A \cdot Q^{-1} = 
    \begin{pmatrix}
      1 & 0 \\
      0 & \tilde A
    \end{pmatrix}
  \]
  with $\tilde A$ of size $(m-1)\times(n-1)$. It is now easy to see that 
  the ideal of $t$-minors of $A$ coincides with the ideal of $(t-1)$-minors 
  of $\tilde A$. This allows for a \textit{reduction of the defining matrix} 
  $A$ to $\tilde A$ for the determinantal ideals 
  $\langle A^{\wedge t} \rangle = \langle \tilde A^{\wedge (t-1)} \rangle$. 
  In the particular case where $R$ is a \textit{local} ring with maximal 
  ideal $\m$ this allows us to always reduce to the case of matrices 
  $A \in {\m}^{m \times n}$ and 
  we will in the following always assume this, unless specified 
  otherwise. 
\end{remark}

\begin{remark}
  \label{rem:Pfaffians}
  In the following we will also consider \textit{Pfaffian ideals} for 
  skew symmetric matrices $A \in R^{m\times m}_\sk$. Any such matrix with 
  entries $a_{i,j} = - a_{j,i}$ can 
  naturally be interpreted as an element of the second exterior power of 
  a free module $R^m$ in some generators $e_1,\dots,e_m$, i.e.
  \[
    A = \sum_{0< i < j \leq m} a_{i,j} \cdot e_i \wedge e_j \in \bigwedge^2 R^m.
  \]
  In this setting, we can consider the exterior powers of $A$ in the usual sense 
  \begin{equation}
    A_\sk^{\wedge s} = \underbrace{ A \wedge A \wedge \dots \wedge A}_{s \textnormal{ times}}
    = \sum_{\# \mathbf i = 2s} A^{\wedge s}_{\mathbf i} \cdot e_{i_1} \wedge e_{i_2} \wedge 
    \dots \wedge e_{i_{2s}}
    \in \bigwedge^{2s} R^m
    \label{eqn:DefinitionPfaffian}
  \end{equation}
  and we write $\langle A^{\wedge s}_\sk \rangle$ for the ideal generated by 
  the coefficients $A^{\wedge s}_{\mathbf i}$ in the expansion above. 
  Note that in case 
  $A \in R^{2n \times 2n}_{\sk}$ is a skew symmetric matrix of even size, 
  we indeed recover the Pfaffian of $A$ as the coefficient in 
  the top exterior power: 
  \[
    \Pf A = A^{\wedge n}_{(1,\dots,2n)}.
  \]
  More generally, the coefficient $A^{\wedge s}_{\mathbf i}$ is nothing but 
  the Pfaffian of the skew-symmetric matrix obtained from $A$ by selecting the 
  even number of rows and columns specified by $\mathbf i$. The ideal 
  $\langle A^{\wedge s}_\sk \rangle$ thus coincides with the ideal of 
  $2s$-Pfaffians of $A$ as for example in \cite{KleppeLaksov80}. 

  The natural $\GL(m;R)$-operation on skew-symmetric matrices is given by 
  \begin{equation}
    \GL(m;R) \times R^{m\times m}_\sk \to R^{m\times m}_\sk, \quad
    (S,A) \mapsto A' = S \cdot A \cdot S^T.
    \label{eqn:MatrixActionForSkewSymmetricCase}
  \end{equation}
  Regarding $S$ as a change of basis in $R^m$ taking $e_i$ to 
  $S( e_i )= \sum_{k=1}^m s_{i,k} \cdot f_k$, it is easy to see that 
  (\ref{eqn:MatrixActionForSkewSymmetricCase}) is compatible with the interpretation 
  of $A$ as a bi-vector, given that 
  \begin{eqnarray*}
    \sum_{0<i<j\leq m} a_{i,j} \cdot e_i \wedge e_j 
    &=&\sum_{0<i<j\leq m} a_{i,j} \cdot \left(\sum_{0<k<l\leq m} (s_{i,k} s_{j,l} - s_{i,l} s_{j,k})\cdot
    f_k \wedge f_l \right) \\
    &=&\sum_{0<k<l\leq m} a'_{k,l} \cdot
    f_k \wedge f_l.
  \end{eqnarray*}
  with $a'_{k,l}$ the entries of the matrix $A'$ above. 

  We will in the following deliberately subsume the Pfaffian case under 
  determinantal ideals in general, but indicate the differences whenever necessary.
\end{remark}

\medskip

One reason, why determinantal ideals received particular interest, is that 
in general they are \textit{not} complete intersection ideals, 
but still provide sufficient additional structure for obtaining stronger 
results than in the general case, e.g. on the module of relations and on 
their free resolutions. For the remainder of this section, we will 
elaborate on the algebraic aspects in which determinantal ideals 
generalize complete intersections.

\medskip

Let $R$ be a commutative Noetherian ring and $I = \langle a_1,\dots,a_n\rangle$ 
an ideal in $I$ generated by $n$ elements. 
By Krull's principal
ideal theorem, we know that the height\footnote{Depending on a textbook,
  $\height I = \inf \{\dim A_{\mathfrak p} \mid {\mathfrak p} \textrm{ prime containing } I\}$ is also 
  called the codimension of $I$, alluding to the fact that e.g.
  for $I \subset k[x_1,\dots,x_s]$ it is indeed the codimension of the variety
  $V(I) \subset k^s$.
} 
of $I$ can at most be $n$. 
This is the \textit{expected height} or 
\textit{expected codimension} for an ideal generated by $n$ arbitrary elements. 
Observe that every such ideal 
is determinantal in a trivial way by setting $I = \langle A^{\wedge 1} \rangle$ for 
the $1\times n$-matrix $A = (a_1,\dots,a_n)$. 
For determinantal ideals in general, the following theorem establishes 
a similar bound. It is worth noting that this theorem has been proved for the 
special case of maximal minors by Macaulay \cite{Macaulay16} already in 1916, 
even before Krull established his principal ideal theorem in 1928.

\begin{theorem}(\cite[Theorem 3]{EagonNorthcott62}, cf. also \cite[Theorem 2.1]{BrunsVetter88})
  \label{thm:BoundOnHeightForDeterminantalIdeals}
  Let $R$ be a Noetherian ring and $A \in R^{m\times n}$ a matrix. If 
  $\langle A^{\wedge t} \rangle \neq R$, then $\height \langle A^{\wedge t}\rangle \leq (m-t+1)(n-t+1)$.
\end{theorem}

In analogy to the complete intersection case, we will refer to this bound on the height of a determinantal ideal as the 
\textit{expected codimension} of the \textit{determinantal ideal} $\langle A^{\wedge t}\rangle$. 
Similar bounds have been established in for determinantal ideals of symmetric 
matrices, see \cite{Kutz74}, where the expected codimension for the ideal of $s$-minors of an
$n \times n$ matrix is $\frac{1}{2} (n-s+2)(n-s+1)$. For 
skew-symmetric matrices $A \in R^{m\times m}_\sk$ the expected codimension of the 
the Pfaffian ideal $\langle A_\sk^{\wedge s}\rangle$ is 
$\frac{1}{2}(m-2s+2)(m-2s+1)$, cf. \cite[Theorem 17]{KleppeLaksov80}.

\medskip
In what follows, all varieties and singularities will usually 
be embedded in a sufficiently ``nice'' (e.g. smooth) ambient space. On the algebraic 
side this corresponds to $R$ having certain favourable properties; for 
example it can be a polynomial ring over a field, or a 
regular local ring. A reasonable and strictly less demanding assumption on $R$ is to be Cohen-Macaulay: 

\begin{definition}(cf. \cite[Section 18.2]{Eisenbud95})
  \label{def:Cohen-Macaulayness}
  A Noetherian ring $R$ is called a \textit{Cohen-Macaulay} ring, if for 
  every maximal ideal $\m$ of $R$ one has $\grade (\m) = \height (\m)$.
\end{definition}

Recall that the \textit{grade} of an ideal $I \subset R$ can be defined as the 
maximal length of a regular $R$-sequence in $I$. In particular, regular local 
rings are Cohen-Macaulay (see e.g.\cite[Section 18.5]{Eisenbud95}) and 
polynomial rings over Cohen-Macaulay rings are again Cohen-Macaulay
(see e.g. \cite[Proposition 18.9]{Eisenbud95}). In a Noetherian ring, 
the grade of an ideal is bounded from above by the height of this 
ideal\footnote{See e.g. \cite[Proposition 1.2.14]{BrunsHerzog93}}, so 
the condition in the definition of a Cohen-Macaulay ring establishes precisely
the opposite inequality for maximal ideals. It can be shown, that this already
implies
\begin{equation}
  \height I = \grade I
  \label{eqn:HeightEqualsGradeInCohenMacaulayRings}
\end{equation}
for every proper ideal $I \subset R$, cf. \cite[Theorem 18.7]{Eisenbud95}.
If, additionally, the ring $R$ is \textit{local},
then the height of any proper ideal $I\neq R$ can actually be 
understood as codimension in the sense that one has 
\begin{equation}
  \height I + \dim R/ I = \dim R,
  \label{eqn:HeightAsCodimension}
\end{equation}
see \cite[Corollary 2.1.4]{BrunsHerzog93}.

\medskip
In our context of a determinantal ideal $I$ in a Cohen-Macaulay ring $R$, 
the grade (or height) is not the main focus of our interest, it is merely one 
ingredient to acquiring more information on free resolutions of $R/I$ which are
known to provide a wealth of subtle algebraic and geometric information. For 
instance, flatness of families can be checked using the first syzygy module,
which is just the beginning of a free resolution (see \cite[I 1.91]{GreuelLossenShustin07})
, and for a Gorenstein ring $R$ free resolutions allow computation of the 
dualizing module $\omega_{R/I}$ of $R/I$ (cf.\cite[Theorem 3.3.7 (b)]{BrunsHerzog93}). 
On the other hand, explicitly computing a free resolution for a given ideal $I$ 
can be an expensive task relying on the standard basis algorithm; knowing the 
general structure of a free resolution in advance is hence a precious advantage.
For determinantal ideals, this is the case. The following theorem is the key to
understanding how this arises.\\   
Recall that a module $M$ over a noetherian ring $R$ is called \textit{perfect}, 
if its \textit{projective dimension}, i.e. the minimal length of a projective 
(or free) resolution thereof, is equal to its \textit{grade}.

\begin{theorem}[\cite{BrunsVetter88}, Theorem 3.5]
  Let $S$ be a Noetherian ring and $M$ a perfect $S$-module of grade $\mu$. 
  Let $R$ be a Noetherian $S$-algebra such that $\grade\!_R (M \otimes_S R) \geq \mu$ 
  and $M \otimes_S R \neq 0$. Then $M \otimes_S R$ is perfect of grade $\mu$ 
  and furthermore $K_\bullet \otimes_S R$ is a free resolution of $M \otimes_S R$ 
  for every free resolution $K_\bullet$ of $M$ of length $\mu$. 
  \label{thm:InheritanceOfPerfection}
\end{theorem}

So in a sufficiently nice setting a known free resolution of an 
$S$-module $M$ provides a free resolution of $M \otimes_S R$ for an $S$-algebra 
$R$. In our context the ring $S$ will be a polynomial ring over ${\mathbb Z}$ or 
${\mathbb Q}$ and the class of modules consists of those of the form 
$S/\langle Y^{\wedge s} \rangle$ for a matrix
   \[
    Y =
    \begin{pmatrix}
      y_{1,1} & \dots & y_{1,n} \\
      \vdots & & \vdots \\
      y_{m,1} & \dots & y_{m,n}
    \end{pmatrix}
    \in S[y_{i,j} : 1 \leq i \leq m, 1 \leq j \leq n].
  \]
This is possible due to a result of Eagon and Hochster 
\cite[Corollary 4]{EagonHochster71} establishing the following property of generic 
perfection for determinantal ideals:

\begin{definition}(\cite[Section 3.A]{BrunsVetter88})
  A finitely generated $\ZZ[y]$-module $M$ is called \textit{generically perfect} 
  if it is perfect and faithfully flat as a $\ZZ$-module. An ideal $I$ is 
  called generically perfect, if $\ZZ[y]/I$ is generically perfect.
\end{definition}

Various different characterizations have been given for generically perfect 
ideals and modules, see \cite{EagonNorthcott67}, or \cite[Proposition 3.2]{BrunsVetter88}. We will treat the two smallest examples, which are relevant to our setting,
as an illustration and for later reference:

\begin{example}(Complete Intersections and the Koszul Complex)
  \label{exp:KoszulComplex}
Recall that the Koszul complex in $n$ elements can be defined as follows. 
Let $y = y_1,\dots,y_n$ be a set of indeterminates over the ring $\ZZ$ and $F$
the free $\ZZ[y]$-module in $n$ generators $e_1,\dots,e_n$. Then exterior 
multiplication with the element 
$\theta = y_1\cdot e_1 + y_2 \cdot e_2 + \dots + y_n \cdot e_n \in F$ gives rise to 
an exact complex
\begin{equation}
  \Kosz \colon 
  \xymatrix{ 
    0 \ar[r] & 
    \bigwedge^0 F \ar[r]^{\theta\wedge} & 
    \bigwedge^{1} F \ar[r]^{\theta\wedge} & 
    \bigwedge^{2} F \ar[r]^{\theta\wedge} & 
    \dots 
    & & \\
    \cdots \ar[r]^{\theta\wedge} & 
    \bigwedge^{n-1} F \ar[r]^{\theta\wedge} & 
    \bigwedge^{n} F \ar[r]^-{\varepsilon} & 
    \ZZ[y]/\langle y_1,\dots,y_n \rangle \ar[r] & 
    0
  }
  \label{eqn:UniversalKoszulComplex}
\end{equation}
where $\varepsilon$ takes the generator $e_1\wedge e_2\wedge \dots \wedge e_n$ to $1 \in R$.
For an arbitrary unital commutative ring $R$ and $n$ elements $a_1,\dots,a_n \in R$,
there is a unique structure of $R$ as a $\ZZ[y]$-module substituting $a_i$ for 
the variable $y_i$. Then the Koszul complex in the elements $a_1,\dots,a_n$ can 
be written as 
\begin{equation}
  \Kosz(a_1,\dots,a_n;R) := \Kosz \otimes_{\ZZ[y]} R.
  \label{eqn:SpecializedKoszulComplex}
\end{equation}
If $\grade \langle a_1,\dots,a_n \rangle = n$, which is, in particular, the case if 
the $a_i$ form a regular sequence, the \textit{Koszul Complex} provides a free 
resolution of the quotient ring $R/\langle a_1,\dots,a_n\rangle$, seen as an 
$R$-module (for textbook references see \cite[Theorem 1.6.17]{BrunsHerzog93} or
\cite[Corollary 17.5]{Eisenbud95}.
\end{example}

\begin{example}(Perfect ideals of grade $2$ and the Hilbert-Burch theorem)
  \label{exp:HilbertBurchComplex}
  Another famous and early example of generic perfection 
  appears in the context of the Hilbert-Burch theorem, see 
  \cite{Hilbert90} and \cite{Burch68}, or \cite{Eisenbud95} for 
  a modern textbook account; the theorem is also stated explicitly as 
  Theorem \ref{thm:HilbertBurch} in Section 
  \ref{sec:DeformationsOfICMC2Singularities} below.

  Consider the $(m+1)\times m$-matrix 
  \[
    Y =
    \begin{pmatrix}
      y_{1,1} & \dots & y_{1,m} \\
      \vdots & & \vdots \\
      y_{m+1,1} & \dots & y_{m+1,m}
    \end{pmatrix}.
  \]
  The homogeneous ideal $I = \langle Y^{\wedge m} \rangle$ is a perfect ideal of 
  grade $2$ and a resolution of $\ZZ[y]/I$ is given by the complex 
  \begin{equation}
    \xymatrix{ 
      0 \ar[r] & 
      \ZZ[y]^m\ar[r]^{Y} &
      \ZZ[y]^{m+1} \ar[r]^{F} & 
      \ZZ[y] \ar[r]^\varepsilon & 
      \ZZ[y]/I \ar[r] & 
      0
    }
    \label{eqn:HilbertBurchResolution}
  \end{equation}
  where $F$ is the $1 \times (m+1)$-matrix with $i$-th entry equal to $(-1)^i$ times 
  the $i$-th minor of $Y$. This can be proved directly using the results 
  by Hilbert and Burch above, but it can also be regarded as a particular 
  case of the Eagon-Northcott complex \cite{EagonNorthcott62}. 

  Now suppose that $R$ is an arbitrary Noetherian ring and 
  \[
    A = 
    \begin{pmatrix}
      a_{1,1} & \dots & a_{1,m} \\
      \vdots & & \vdots \\
      a_{m+1,1} & \dots & a_{m+1,m}
    \end{pmatrix}
  \]
  a matrix with entries in $R$. Then the quotient ring with respect to the ideal 
  $I = \langle A^{\wedge m} \rangle$ can be identified with 
  \[
    R/I \cong \ZZ[y]/\langle Y^{\wedge m} \rangle \otimes_{\ZZ[y]} R
  \]
  via the map $\ZZ[y] \to R, y_{i,j} \mapsto a_{i,j}$. Whenever $I$ has 
  grade $2$ in $R$, Theorem \ref{thm:InheritanceOfPerfection} assures that, moreover, 
  the complex obtained from (\ref{eqn:HilbertBurchResolution}) by applying 
  $-\otimes_{\ZZ[y]} R$ is in fact a free resolution of $R/I$. 
\end{example}

\begin{remark}
  \label{rem:GenericPerfectionForDeterminantalIdealsAndResolutions}
  While generic perfection of determinantal ideals has been established by 
  Eagon and Hochster in \cite[Corollary 4]{EagonHochster71},
  explicit free resolutions of minimal length of $\langle Y^{\wedge t} \rangle$
  have only been constructed in special cases over the ring $\ZZ[y]$, but 
  are known for all values of $m,n$ and $t$ over the ring $\QQ[y]$.
  Passing from $\ZZ[y]$ to $\QQ[y]$ does not impose any severe 
  restrictions in the setting of determinantal singularities, since the main
  focus is on complex algebraic and analytic spaces: 
  All rings $R$ in question will be of 
  characteristic zero so that again any choice of elements $a_{i,j} \in R$ 
  turns the ring into a $\QQ[y]$-algebra by substitution of $y_{i,j}$ by $a_{i,j}$.
\end{remark}

To conclude this brief discussion of facts from commutative algebra, we
give a short overview (and non-exhaustive) on results concerning the 
construction of resolutions of generic determinantal ideals 
$\langle Y^{\wedge t}\rangle$ over rings $k[y]$ for $k$ being either 
$\ZZ$ or $\QQ$. For a survey on such resolutions we refer the interested 
reader to \cite{PragaczWeyman86}.

As already mentioned, in the case of maximal minors $t = m \leq n$, a 
resolution of for $k=\ZZ$ is provided by the well known
\textit{Eagon-Northcott-complex} \cite{EagonNorthcott62}. Another 
construction of this complex has been given by Buchsbaum in 
\cite{Buchsbaum79}. This complex is covered in various textbooks 
such as \cite{Eisenbud95}, or \cite{BrunsVetter88} and appears 
as a special case of the family of complexes described independently 
by Buchsbaum and Eisenbud \cite{BuchsbaumEisenbud73} and 
Kirby \cite{Kirby74}. The Hilbert-Burch theorem fits into this 
setting as the even more special case $n=m+1$.

For submaximal minors of square matrices, i.e. for $t = m-1, m = n$ 
a resolution of $\ZZ[y]/\langle Y^{\wedge t} \rangle$ 
was found by Gulliksen and Negard \cite{GulliksenNegard72}; the case of
non-square matrices, i.e. $t=m-1, m \leq n$, has been treated a decade later 
by Akin, Buchsbaum, and Weyman in \cite{AkinBuchsbaumWeyman81}.

For the general case $t \leq m \leq n$, Lascoux provided a free resolution
in \cite{Lascoux78}, but only over the rationals (or, more generally, 
a field of characteristic zero). This is due to the fact that representation 
theory and the use of Schur-functors in his article, which made the change 
of coefficients necessary. The methods by Lascoux also work for symmetric 
matrices and for skew symmetric matrices with ideals generated by Pfaffians.

Resolutions of ideals of submaximal minors for symmetric matrices were 
also described by J\'ozefiak in \cite{Jozefiak78}. 
The interest in Gorenstein rings also led to the construction of free 
resolutions. While Gorenstein rings of codimension $1$ and $2$ are known to be
complete intersections, those of codimension $3$ are submaximal Pfaffians of 
skew-symmetric matrices (see \cite{BuchsbaumEisenbud77} and for a discussion
of a generic free resolution \cite{JozefiakPragacz79},
for a slightly
different perspective see also\cite{Waldi79}). The further study of Gorenstein
rings of higher codimension led to partial results, including (non-generic) 
resolutions, and is still a topic of active research in commutative algebra.
A full description of the first syzygy module determinantal ideals for 
arbitrary $t \leq m \leq n$ over arbitrary unitary rings has been 
given by Kurano \cite{Kurano89} and by Ma \cite{Ma93} solving a conjecture 
of Sharpe.

\medskip
In what follows, we will refer to any free resolution of the quotient rings 
$k[y]/\langle Y^{\wedge t}\rangle$, for $k=\ZZ$ or $\QQ$,  
as 
\begin{equation}
  K(m,n,t) \colon 
  \xymatrix{
    0 \ar[r] &
    K_{c} \ar[r] & 
    K_{c-1} \ar[r] & 
    \dots \ar[r] & 
    K_1 \ar[r] & 
    K_0 
  }.
  \label{eqn:GenericDeterminantalComplex}
\end{equation}
where $c$ is the expected grade or codimension for the given values of $m,n$ and $t$. 
In the case of symmetric or skew-symmetric matrices, we will occasionally write 
$K^\sym(m,t)$ and $K^\sk(m,t)$ for the resolutions of the quotients by the determinantal, 
respectively the Pfaffian ideals.

\subsection{Determinantal singularities and their deformations}

Let again 
\[
  Y =
  \begin{pmatrix}
    y_{1,1} & \dots & y_{1,n} \\
    \vdots & & \vdots \\
    y_{m,1} & \dots & y_{m,n}
  \end{pmatrix}
\]
denote the matrix of $m \cdot n$ indeterminates but over a field $k$. 
The \textit{generic determinantal varieties} $M_{m,n}^s(k)$ are defined 
as the vanishing loci of the ideals $\langle Y^{\wedge s}\rangle$ of 
$s$-minors of $Y$. It is evident that 
\[
  M_{m,n}^s(k) = \left\{ \varphi \in k^{m\times n} : \rank \varphi < s \right\}.
\]
Now let $R$ be a $k$-algebra and $A \in R^{m\times n}$ a matrix with entries 
$a_{i,j}$ in $R$. 
Then $A$ gives rise to a homomorphism of rings $k[y] \to R, y_{i,j} \mapsto a_{i,j}$ 
which corresponds to a map $\Spec R \to k^{m\times n}$ on the geometric side. 
In accordance with Theorem \ref{thm:BoundOnHeightForDeterminantalIdeals} 
this suggests:

\begin{definition}
  Let $R$ be a Noetherian $k$-algebra and $A \in R^{m\times n}$ a matrix. The 
  variety defined by the ideal $I = \langle A^{\wedge s} \rangle$ is a 
  \textit{determinantal variety of type $s$ for the matrix} $A$ 
  if $\height I$ 
  is equal to the expected codimension; that is 
  $(m-s+1)(n-s+1)$ in the general case, $\frac{1}{2}(n-s+1)(n-s+2)$ for 
  the ideal of $s$-minors of symmetric matrices, and 
  $\frac{1}{2}(m-2s+2)(m-2s+1)$ for $2s$-Pfaffian ideals of 
  $m\times m$-skew-symmetric matrices.
\end{definition}
When $R$ is Cohen-Macaulay and $I = \langle A^{\wedge s} \rangle$ defines a determinantal 
variety, then due to the equality of height and grade 
(\ref{eqn:HeightEqualsGradeInCohenMacaulayRings}) 
Theorem \ref{thm:InheritanceOfPerfection} 
applies so that $K(m,n,s)\otimes_{k[y]} R$ 
is a free resolution of the module $R/I$. 

We will in the following mostly be concerned with the case $k = \CC$
and $R = \CC\{x\} = \CC\{x_1,\dots,x_p\}$, the ring of convergent power series 
at the origin in $\CC^p$. A matrix $A \in \CC\{x\}^{m\times n}$ will then 
be interpreted as a holomorphic map germ 
\[
  A \colon (\CC^p,0) \to (\CC^{m\times n},0).
\]
By abuse of notation, let $A$ denote a representative of this germ  
defined on some open neighbourhood $U$ of the origin. 
We write $X_A^s\subset U$ for the complex analytic space defined by the 
sheaf of ideals $\langle A^{\wedge s} \rangle$ in $\OO_U$. 
Then by construction one has an equality of sets 
\[
  X_A^s = A^{-1}( M_{m,n}^s ) \subset U 
\]
on the geometric side and an isomorphism 
\[
  \CC\{x\}/\langle A^{\wedge s}\rangle \cong 
  \CC[y]/\langle Y^{\wedge s} \rangle \otimes_{\CC[y]} \CC\{x\}
\]
on the algebraic side, together with its sheafification on $U$.


\begin{remark}
  \label{rem:GraphTransformationAndRelativeCompleteIntersections}
  For $A$ as above let $\Gamma_A = \{ (x, \varphi) \in U \times \CC^{m\times n} : \varphi = A(x) \}$
  be the \textit{graph} of $A$. 
  This is a subvariety of the product $U \times \CC^{m\times n}$ 
  defined by the $m\cdot n$ equations $h_{i,j} = y_{i,j} - a_{i,j}(x) \in \CC\{x\}[y]$. 
  A close inspection of the arguments of 
  Eagon and Northcott show that $\langle A^{\wedge s} \rangle$ has expected 
  codimension if and only if the equations $h_{i,j}$ form a regular sequence 
  on the coordinate ring $\CC\{x\}[y]/\langle Y^{\wedge s}\rangle$
  of the product variety $(\CC^p,0)\times (M_{m,n}^s,0)$.

  It has been shown in \cite[Proposition 2]{EagonNorthcott67} that 
  the homology of the complex $K(m,n,s) \otimes_{\ZZ[y]} \CC\{x\}$ is isomorphic 
  to the \textit{Koszul homology}
  \[
    H_i\left( K(m,n,s)\otimes_{\ZZ[y]}\CC\{x\} \right) \cong 
    H_i\left( \Kosz( h, \CC\{x\}[y]/\langle Y^{\wedge s}\rangle \right)
  \]
  of the functions $h_{i,j}$ on the $\CC\{x\}[y]/\langle Y^{\wedge s} \rangle$. 
  We saw earlier in Example \ref{exp:KoszulComplex} 
  that the Koszul homology vanishes when $h_{1,1},\dots,h_{m,n}$ is a 
  regular sequence on $\CC\{x\}[y]/\langle Y^{\wedge s} \rangle$. 
  Conversely, in a local ring, the Koszul homology 
  can be used to measure the grade of an ideal: If the right hand side 
  vanishes for all $i>0$ one has $\grade \langle h \rangle = m\cdot n$, 
  cf. \cite[Theorem 1.6.17]{BrunsHerzog93}. Since $\CC\{x\}[y]$ is a
  graded ring over a local ring, this implies that the $h_{1,1},\dots,h_{m,n}$ 
  have to be a regular sequence already, cf. \cite[Corollary 1.6.19]{BrunsHerzog93}.
  From this we see that $K(m,n,s)\otimes_{\ZZ[y]} \CC\{x\}$ is a
  free resolution of $\CC\{x\}/ \langle A^{\wedge s} \rangle$ if and 
  only if the equations $h_{i,j} = y_{i,j}-a_{i,j}(x)$ for the graph 
  $\Gamma_A$ form a regular sequence on $\CC\{x\}[y]/\langle Y^{\wedge s}\rangle$.

  This allows us to always consider $(X_A^s,0)$ as a complete intersection singularity 
  in the singular, but yet well known ambient space $(\CC^p,0) \times (M_{m,n}^s,0)$.

%
\end{remark}

The interpretation of a matrix $A \in \CC\{x\}^{m\times n}$ as a map germ 
is of particular interest when it comes to 
\textit{deformations}\footnote{
  See \cite[Definition 7.1.1]{Volume1} for a definition of deformations 
  of complex analytic germs.
}
of the germ $(X_A^s,0)$.
Recall that for hypersurface and complete intersection singularities 
$(X,0) \subset (\CC^p,0)$ of codimension $m$ and a map 
$f \colon (\CC^p,0) \to (\CC^m,0)$ defining $(X,0) = (f^{-1}(\{0\}),0)$, an 
\textit{unfolding} of $f$ gives rise to a \textit{deformation} 
of $(X,0)$, cf. \cite[Proposition 7.1.11]{Volume1}. 
The generic perfection of determinantal ideals is crucial 
for proving that more generally, an unfolding of the matrix 
$A$ gives rise to a deformation of the associated determinantal 
singularity, see Lemma \ref{lem:UnfoldingInducesDeformation}.
We briefly recall the notions involved;
for a more thorough discussion, the reader is refered to 
\cite[Sections 7.1 and 7.2]{Volume1}. 

An \textit{unfolding} on $k$ parameters of a map germ 
$f \colon (\CC^p,0) \to (\CC^n,0)$ is given by a map 
\[
 F \colon (\CC^p,0) \times (\CC^k,0) \to (\CC^n,0), \quad 
 (x,t) \mapsto F(x,t) = f_t(x)
\]
such that for $t=0$ the germ $f_0$ coincides with $f$.  
In practice, an unfolding is nothing but a perturbation 
\[
  F(x,t) = f(x) + t_1\cdot g_1(x,t) + t_2 \cdot g_2(x,t) + \dots
\]
of the original map germ $f_0$. 

An (embedded) \textit{deformation} of a complex analytic germ $(X,0) \subset (\CC^p,0)$ 
over a germ $(T,0)$ is given by a commutative diagram 
\[
  \xymatrix{
    (X,0) \ar@{^{(}->}[r] \ar[d] &
    (\mathcal X,0) \ar@{^{(}->}[r] \ar[d]^\pi & 
    (\CC^p,0) \times (T,0) \ar[dl] \\
    \{0\} \ar@{^{(}->}[r] & 
    (T,0) & 
  }
\]
where $(\mathcal X,0)\subset (\CC^p,0) \times (T,0)$ is another 
complex analytic germ such that the projection $\pi$ turns $\OO_{\mathcal X,0}$ 
into a \textit{flat} $\OO_{S,0}$-module so that the above diagram becomes a 
\textit{flat family} with special fiber $(X,0)$.

Flatness of a family is a technical 
algebraic criterion to assure that the fibers in a family vary nicely; 
for example, there can be no jumps in dimension of the fibers 
in a flat family, see \cite{Eisenbud95}. 
When $f_1,\dots,f_n \in \CC\{x\}$ are equations defining $(X,0)$ as above, then 
an arbitrary perturbation 
of these defining equations given by functions 
\[
  F_i(x,t) = f_i(x) + \sum_{j=1} t_j \cdot g_j(x,t) \in \CC\{x,t\}, 
  \quad i=1,\dots,n, 
\]
in additional parameters $t = t_1,\dots,t_k$ will in general not lead to a flat family 
\[
  (\{F_1 = \dots F_n = 0\},0) = (\mathcal X,0) \overset{\pi}{\longrightarrow} (\CC^k,0)
\]
unless the $f_i$ form a regular sequence, 
cf. Example \ref{exp:DeformationOfThreeCoordinateAxisFlatness} below.

A prominent characterization of flatness in the 
general case is the 
``Flatness by relations'', \cite[Proposition 7.1.2]{Volume1}. Applied to 
the above situation it says 
that the family $\pi \colon (\mathcal X,0) \to (\CC^k,0)$ is flat 
if and only if for every relation $r \in \CC\{x\}^n$ among the $f_i$ of the form 
\[
  r_1\cdot f_1 + r_2 \cdot f_2 + \dots + r_n \cdot f_n = 0,
\]
i.e. the \textit{first syzygies} of the $f_i$,
there exists a relation $R \in \CC\{x,t\}^n$ of the $F_i$ 
\[
  R_1\cdot F_1 + R_2 \cdot F_2 + \dots + R_n \cdot F_n = 0
\]
with $R_i$ congruent to $r_i$ modulo $\langle t \rangle$.
In general, not every perturbation $F_i$ of the $f_i$ admits such a ``lifting 
of relations''. For determinantal singularities, however, the inheritance 
of free resolutions granted by Theorem \ref{thm:InheritanceOfPerfection} 
assures that flatness by relations always applies:

\begin{lemma}
  \label{lem:UnfoldingInducesDeformation}
  Let $A \in \CC\{x_1,\dots,x_p\}^{m\times n}$ be a matrix defining a 
  determinantal singularity $(X_A^s,0) \subset (\CC^p,0)$ of type $s$. 
  Then any \textit{unfolding} 
  \[
    (\CC^p,0) \times (\CC^k,0) \to (\CC^{m\times n},0) \times (\CC^k,0), 
    \quad 
    (x,t) \mapsto (\mathbf A(x,t),t)
  \]
  of $A$ on $k$ parameters induces a \textit{deformation} 
  \[
    \xymatrix{
      (X_A^s,0) \ar@{^{(}->}[r] \ar[d] &
      (\mathcal X_{\mathbf A}^s,0) \ar@{^{(}->}[r] \ar[d]^\pi & 
      (\CC^p,0) \times (\CC^k,0) \ar[dl] \\
      \{0\} \ar@{^{(}->}[r] & 
      (\CC^k,0) & 
    }
  \]
  of the germ $(X_A^s,0)$.
\end{lemma}

\begin{proof}
  The crucial point is to verify flatness of the family $(X_{\mathbf A}^s,0) 
  \overset{\pi}{\longrightarrow} (\CC^k,0)$, where 
  $(X_{\mathbf A}^s,0) \subset (\CC^p \times \CC^k, 0)$ 
  is the complex analytic germ defined by the ideal 
  $\mathbf I = \langle \mathbf A^{\wedge s} \rangle$ in 
  the ring $\CC\{x_1,\dots,x_p,t_1,\dots,t_k\}$. 
  Let $I = \langle A^{\wedge s} \rangle \subset \CC\{x\}$ be the ideal 
  defining $X_A^s \cong \mathcal X_{\mathbf A}^s \cap \{ t = 0 \}$. 
  We have an inequality of dimensions 
  \[
    \dim \CC\{x,t\}/\mathbf I \leq \dim \CC\{t\} + \dim \CC\{x\}/ I,
  \]
  see \cite[Theorem 10.10]{Eisenbud95}, where we identify 
  $\CC\{x,t\}/\mathbf I + \langle t \rangle \cong \CC\{x\}/I$.
  Since $(X_A^s,0) \subset (\CC^p,0)$ has expected codimension 
  $c = (m-s+1)(n-s+1)$, this inequality entails that also $(\mathcal X_{\mathbf A}^s,0)$ 
  must have expected codimension $c$ in $(\CC^{p},0) \times (\CC^k,0)$. 
  Due to the generic perfection of the determinantal ideals,
  a free resolution of $\CC\{x,t\}/\mathbf I$ is given by 
  $K_\bullet(m,n,s) \otimes_{\mathbf A} \CC\{x,t\}$ according to 
  Theorem \ref{thm:InheritanceOfPerfection}. This resolution specializes to one 
  of $\CC\{x\}/I$ at $t=0$ by using the same theorem again. This implies
  in particular flatness of the family by evoking its characterization via
  the relation lifting property (see e.g. \cite[Proposition 7.1.2]{Volume1})
  in the following way: A relation among elements of $I$ (and $\mathbf I$ 
  respectively) is an element of the corresponding first syzygy module, which 
  can be read off as the image of the rightmost morphism in the given 
  free resolution. As the resolution of $\CC\{x,t\}/\mathbf I$ specializes
  to the one of $\CC\{x\}/I$, clearly any element of the first syzygy 
  module of $I$ arises in this way, i.e. lifts to one of the first syzygy 
  module of $\mathbf I$.
\end{proof}

\begin{definition}
  \label{def:DeterminantalDeformation}
  Let $(X_A^s,0)\subset (\CC^p,0)$ be a \textit{determinantal} singularity of 
  type $s$ defined by a matrix $A \in \CC\{x\}^{m\times n}$. Any 
  deformation induced from an unfolding of $A$ as 
  in Lemma \ref{lem:UnfoldingInducesDeformation} is called a 
  \textit{determinantal deformation} of $(X_A^s,0)$. 

  Conversely, we will say that a given deformation 
  $(X,0) \hookrightarrow (\mathcal X,0) \overset{\pi}{\longrightarrow} (S,0)$ 
  of an \textit{arbitrary} singularity $(X,0) \subset (\CC^p,0)$ 
  is \textit{determinantal for} $A$, if there 
  exists an integer $s$ and a matrix $A \in \CC\{x\}^{m\times n}$ such that 
  $(X,0) \cong (X_A^s,0)$ is determinantal for $A$ of type $s$
  and an unfolding $\mathbf A \colon (\CC^p,0) \times (\CC^k,0) \to (\CC^{m\times n},0)$ 
  of $A$ 
  together with a commutative diagram 
  \[
    \xymatrix{
      (X,0) \ar@{^{(}->}[r] \ar[d] & 
      (\mathcal X,0) \ar[r] \ar[d]^\pi & 
      (\mathcal X_{\mathbf A}^s,0) \ar[d] \\
      \{0\} \ar@{^{(}->}[r] & 
      (S,0) \ar[r]^{\Psi} & 
      (\CC^k,0)
    }
  \]
  where $(\mathcal X_{\mathbf A}^s,0) \to (\CC^k,0)$ is the family induced 
  from $\mathbf A$ as in Lemma \ref{lem:UnfoldingInducesDeformation}.
\end{definition}

\begin{example}
  \label{exp:DeformationOfThreeCoordinateAxisFlatness}
Consider the three coordinate axes in $({\mathbb C}^3,0)$, which form a 
determinantal singularity of type $s=2$ with matrix 
$$A=\begin{pmatrix} x & 0 & z \cr
                    0 & y & z \end{pmatrix}$$ 
The corresponding ideal is 
$\langle A^{\wedge 2} \rangle =\langle -yz,-xz,xy \rangle =: \langle f_1,f_2,f_3 \rangle$. 
If we think of these generators as a map 
\[
  f = (f_1,f_2,f_3) = A^{\wedge 2} \colon \CC^3 \to \CC^3,
\]
then according to Sard's theorem \cite{Sard42}, the set of points $v \in \CC^3$ 
for which $f^{-1}(\{v\})$ is regular of complex codimension $3$, is dense. 
Hence a generic perturbation of the $f_i$ will ``deform'' the curve $(X_A^2,0)$ 
to a collection of points.

Consider now the 1-parameter unfolding of $X^s_A$ determined by the 
perturbed matrix
$$\mathbf A = \begin{pmatrix} x & t & z \cr
                          0 & y & z \end{pmatrix}$$
and the corresponding ideal 
$\langle {\bf A}^{\wedge 2}\rangle = \langle (t-y)z, -xz,xy \rangle \subset \CC\{x,y,z,t\}$.
This is a very specific perturbation of the generators given by 
\[
  F_1 = f_1 + t\cdot z, \quad F_2 = f_2 + t \cdot 0, \quad F_3 = f_3 + t \cdot 0.
\]
According to Lemma \ref{lem:UnfoldingInducesDeformation}, the projection to the parameter $t$
\[
  ({\mathcal X}^s_{\bf A},0) \to (\CC,0), \quad (x,y,z,t) \mapsto t
\]
is indeed a deformation of $(X^s_A,0)$. This is   
because all relations among the $f_i$ arise from the following two:
\begin{eqnarray*}
x \cdot f_1 + 0 \cdot f_2 + z \cdot f_3 & = & 0\cr
0 \cdot f_1 + y \cdot f_2 + z \cdot f_3 & = & 0
\end{eqnarray*}
and these lift to relations
\begin{eqnarray*}
x \cdot F_1 + t \cdot F_2 + z \cdot F_3 & = & 0 \cr
0 \cdot F_1 + y \cdot F_2 + z \cdot F_3 & = & 0
\end{eqnarray*}
among the $F_i$.

Geometrically, this manifests itself in the fact that the fiber over $t\neq 0$ is 
indeed also of complex dimension $1$. 
\begin{figure}[h]
  \centering
  \includegraphics[trim=2cm 10cm 7cm 1cm, clip, scale=0.5]{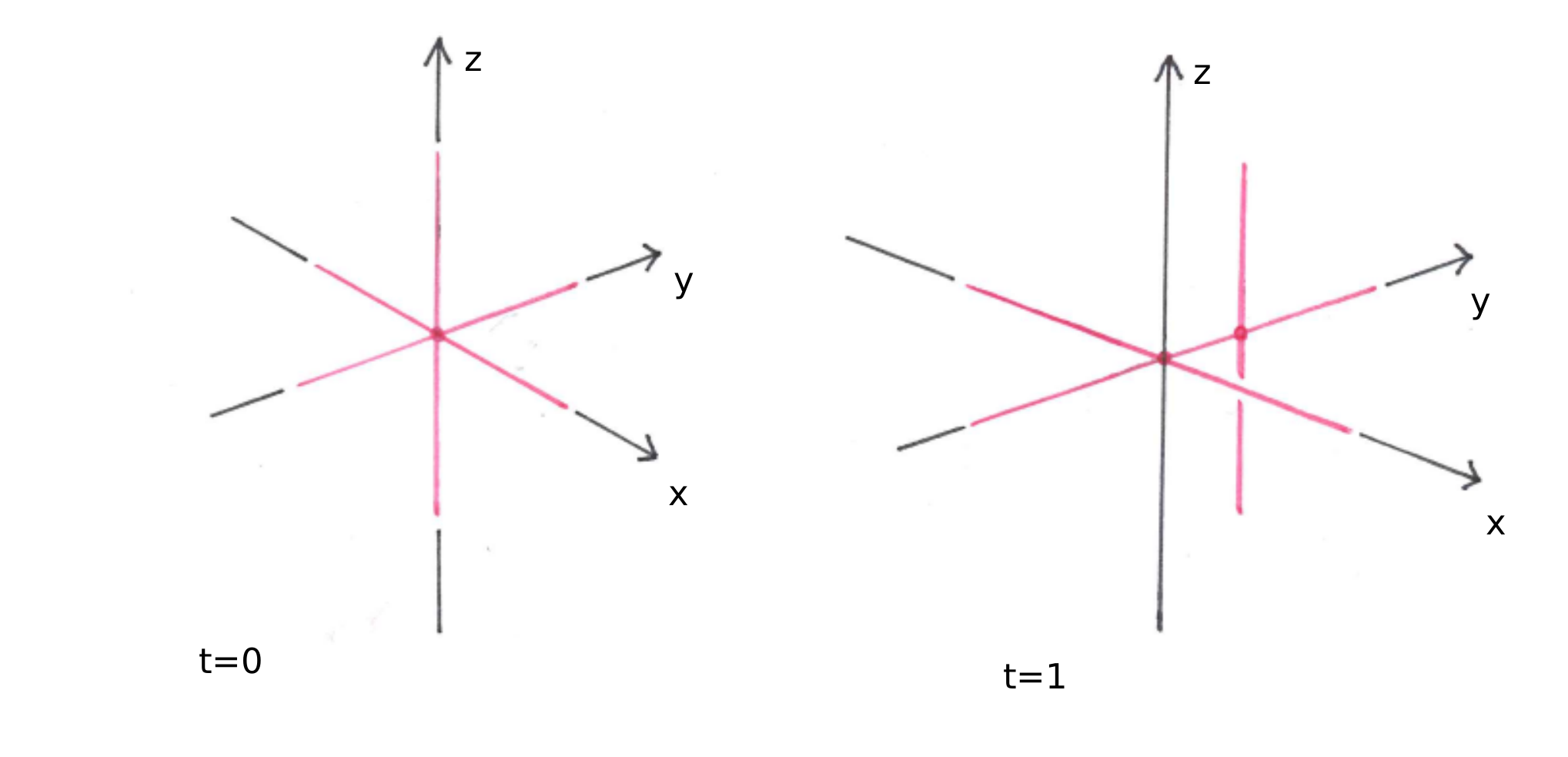}
  \caption{A deformation of the three coordinate axis in $\CC^3$}
  \label{fig:DeformationThreeCoordinateAxisNonSmoothing}
\end{figure}
In this sense, Lemma \ref{lem:UnfoldingInducesDeformation} assures that of all 
possible perturbations of the generators $f_i$ the very few ones which arise from perturbations 
of the matrix $A$ are similarly well behaved.

\end{example}

In general, not every deformation of a given determinantal 
singularity is determinantal for the specific matrix and the interplay between unfoldings of matrices 
and the deformations of the associated singularities can be quite complicated. 
This is illustrated by several examples gathered in Section 
\ref{sec:ComparisonOfUnfoldingsAndSemiUniversalDeformations}.
For some particular classes of singularities such as complete intersections,
Cohen-Macaulay singularities of codimension $2$, or Gorenstein singularities 
of codimension $3$, there are canonical choices of determinantal structures 
and the theory of unfoldings for the defining matrices indeed agrees 
with the deformation theory of the germs. These particular cases 
will be discussed in detail in Section 
\ref{sec:ComparisonOfUnfoldingsAndSemiUniversalDeformations}. 
Finally, as we will report in Section \ref{sec:FurtherReferencesAndTechniques}, 
results of Buchweitz \cite{Buchweitz81} and Svanes \cite{Svanes72} 
imply that for a given determinantal singularity $(X_A^s,0)$ which is not 
a hypersurface and which is 
\textit{unobstructed}\footnote{For the definition of obstructions in the 
  context of deformation theory see e.g. \cite[Section 7.1.5]{Volume1}.
},
in fact every deformation of $(X_A^s,0)$ is determinantal for the defining 
matrix $A$.

However, before we can 
present all these discussions in their full detail, we first need to develop 
the underlying notions of equivalence for unfoldings of matrices and the 
associated concepts of finite determinacy and versal unfoldings in 
Section \ref{sec:UnfoldingsAndEquivalenceOfMatrices}.

\subsection{Geometry of the generic determinantal varieties}

The preceding sections made it clear that the generic determinantal 
varieties $M_{m,n}^s$ play a fundamental role in the 
study of arbitrary determinantal varieties and singularities. As we shall see in 
what follows, 
this does not only apply to their algebraic, but also to their 
geometric and even their topological properties when working over the real, 
or the complex numbers.
In this section we gather some auxiliary results on the geometry 
of the generic determinantal varieties $M_{m,n}^s(\mathbb K) \subset \KK^{m\times n}$ 
where $\KK=\CC$ or $\RR$. For some of the results, the reader can 
also make the obvious translations to other fields, if needed.

\subsubsection{Resolution of singularities for $M_{m,n}^s$}

We start by observing that the generic determinantal varieties admit 
certain canonical resolutions of singularities
\begin{equation}
  \xymatrix{
    & \tilde M_{m,n}^s \ar[dl] \ar[dd]^{\tilde \nu} \ar[dr] & \\
    \hat M_{m,n}^s \ar[dr]^{\hat \nu} & & 
    \check {M}_{m,n}^s \ar[dl]_{\check \nu} \\
    & M_{m,n}^s & 
  }
  \label{eqn:ResolutionsOfSingularitiesGenericDeterminantalVariety}
\end{equation}
where 
\begin{eqnarray*}
  \hat M_{m,n}^s &=&  \left\{ 
    (\varphi, V) \in \KK^{m\times n}\times \Grass(n-s+1,n) :
    V \subset \ker \varphi
  \right\} \\
  &\cong & 
  \left\{ (\varphi, V^\perp) \in \KK^{m\times n} \times \Grass(s-1,n) : 
  \im \varphi^\vee \subset V^\perp \right\}
\end{eqnarray*}
is the \textit{Tjurina transform}
with its projection $\hat \nu \colon (\varphi,W) \mapsto \varphi$, 
\begin{eqnarray*}
  \check M_{m,n}^s &=&  \left\{ 
    (\varphi, W) \in \KK^{m\times n} \times \Grass(s-1,m) : 
    \im \varphi \subset W
  \right\}\\
  &\cong& 
  \left\{ (\varphi,W^\perp) \in \KK^{m\times n} \times \Grass(m-s+1,m) \colon 
  W^\perp \subset \ker \varphi^\vee \right\}
\end{eqnarray*}
the \textit{dual Tjurina transform} with projection $\check \nu$, and 
\[
   \begin{array}{ccll}
  \tilde M_{m,n}^s &=& \hat M_{m,n}^s \times_{M_{m,n}^s} \check M_{m,n}^s & \\
  &=& \{ (\varphi, V,W) \in 
    \KK^{m\times n} \times \Grass(n-s+1,n) & \times \, \Grass(s-1,m) :\\
  & & &  V \subset \ker \varphi, \, \im \varphi \subset W \}
\end{array}
\]
the \textit{Nash transform} of $M_{m,n}^s$ with projection $\tilde \nu$. 

Here we used the canonical isomorphisms of dual Grassmannians $\Grass(r,n) \cong \Grass(n-r,n)$ 
induced by the correspondence 
\[
  V \mapsto V^\perp = \{ f \in (\KK^n)^\vee : f|_V = 0 \}
\]
with $(\KK^n)^\vee = \Hom_\KK( \KK^n,\KK)$. 
The natural isomorphism 
\[
  \Hom(\KK^n,\KK^m) \to \Hom\left( (\KK^n)^\vee, (\KK^m)^\vee \right), \quad 
  \varphi \mapsto \varphi^\vee,
\]
which is given by transposition 
$\KK^{m\times n} \overset{\cong}{\longrightarrow} \KK^{n\times m}, \, A \mapsto A^T$
in terms of matrices, takes $M_{m,n}^s$ into $M_{n,m}^s$. 
It is now easy to see that this identification extends to the Tjurina transforms 
and their duals so that $\hat M_{m,n}^s \cong M_{m,n}^s \times_{M_{n,m}^s} \check M_{n,m}^s$ 
and $\check M_{m,n}^s \cong M_{m,n}^s\times_{M_{n,m}^s} \hat M_{n,m}^s$.

Either one of the varieties $\hat M_{m,n}^s$, $\check M_{m,n}^s$ 
and $\tilde M_{m,n}^s$ is the total space of an algebraic 
vector bundle over the respective Grassmannians. For instance, if we let 
$0 \to S \to \OO^n \to Q \to 0$ be the tautological sequence over 
$\Grass(n-s+1,n)$, then 
\[
  \hat M_{m,n}^s \cong \left| \Hom( Q, \OO^m )\right|,
\]
where we write $|E|$ for taking the associated total space of 
a vector bundle $E$.
Similar descriptions can be made for the dual Tjurina transform $\check M_{m,n}^s$ 
and the Nash transform.
In particular, all of the three spaces $\hat M_{m,n}^s, \check M_{m,n}^s$, and 
$\tilde M_{m,n}^s$ are smooth. 

\begin{remark}
  That $\tilde \nu \colon \tilde M_{m,n}^s \to M_{m,n}^s$ is really the 
  Nash modification in the sense of e.g. \cite[Definition 3.9.2]{Volume1} was 
  observed by Ebeling and Gusein-Zade in \cite{EbelingGuseinZade09}. It is a 
  well known fact that for any fixed rank $r$
  the tangent space to the stratum $V_{m,n}^r$ at a point 
  $\varphi$ is given by 
  \[
    T_\varphi V_{m,n}^r = \{ \psi \in \KK^{m\times n} : \psi ( \ker \varphi ) \subset \im \varphi \}.
  \]
  It follows that the Gauss map 
  $\gamma \colon V_{m,n}^r \to \Grass( (m+n)r - r^2, m n )$ taking 
  a point $\varphi$ to its tangent space $T_\varphi V_{m,n}^r \subset T_\varphi \KK^{m\times n}$ 
  factors through the product 
  \[
    \xymatrix{
      V_{m,n}^r \ar[r]^-\alpha \ar@/_2.0pc/[rr]_\gamma &
      \Grass(n-r,n) \times \Grass(r,m) \ar[r]^-\beta &
      \Grass( (m+n)r - r^2, m n) 
    }
  \]
  where $\alpha \colon \varphi \mapsto (\ker \varphi, \im \varphi)$ and 
  $\beta \colon (V,W) \mapsto \Hom( \KK^n/V, W) \subset \Hom(\KK^n,\KK^m)$, 
  and consequently, the Nash blowup of $M_{m,n}^{r+1}$ can be performed 
  using this product of Grassmannians.
\end{remark}

For any integer $r$ we let 
\[
  V_{m,n}^r = \{ \varphi \in \KK^{m\times n} : \rank \varphi = r \}
\]
be the set of matrices of fixed rank $r$. 
Whenever $\varphi \in M_{m,n}^s$ belongs to $V_{m,n}^{s-1}$, the spaces $V$ and $W$ 
in the definitions of $\hat M_{m,n}^s$, $\check M_{m,n}^s$ and $\tilde M_{m,n}^s$ 
above are 
uniquely determined by the kernel and the image of $\varphi$ and vary 
algebraically with $\varphi$. Hence, either one of the three projections 
$\hat \nu$, $\check \nu$, and $\tilde \nu$ is a local isomorphism over the dense, open subset 
$V_{m,n}^{s-1} \subset M_{m,n}^s$. It is easy to see that $M_{m,n}^s$ is singular along 
the complement $M_{m,n}^{s-1}$ of $V_{m,n}^{s-1}$ in $M_{m,n}^s$ so that 
indeed either one of the three maps $\hat \nu$, $\check \nu$, and $\tilde \nu$ provides a 
resolution 
of singularities\footnote{See for instance \cite[Chapter 3]{Volume1} for a definition 
of resolution of singularities.}

\subsubsection{The rank stratification}

The sets $V_{m,n}^r$ of matrices of a fixed rank $r$ form a 
complex algebraic \textit{stratification} of $\KK^{m\times n}$, that is 
a decomposition as a disjoint union of locally closed, complex algebraic submanifolds 
\begin{equation}
  \KK^{m\times n} = \bigcup_{r=0}^{\min\{m,n\}} V_{m,n}^r.
  \label{eqn:TheRankStratification}
\end{equation}
In particular, every $M_{m,n}^s = \bigcup_{r<s} V_{m,n}^r$ is a union of strata.
We will in the following refer to this particular stratification of the 
space of matrices as the \textit{rank stratification}.
For an account on stratification theory see, for instance, \cite[Chapter 4]{Volume1}, 
or \cite{GoreskyMacPherson88}.

Left- and right-multiplication by invertible matrices does not change 
the rank of a matrix, so the action of the group 
$G = \GL(n;\KK) \times \GL(m;\KK)$ given by 
\begin{equation}
  G \times \KK^{m\times n} \to \KK^{m\times n},
  \quad 
  ((P,Q),\varphi) \mapsto P \cdot \varphi \cdot Q^{-1}
  \label{eqn:LieGroupActionOnTheSpaceOfMatrices}
\end{equation}
preserves the rank stratification.
Using this action, one can even construct local analytic trivializations:
Let $\varphi \in \KK^{m\times n}$ be an arbitrary matrix and $r$ its rank so that 
$\varphi \in V_{m,n}^r$. By virtue of the $G$-action, we may assume that $\varphi$ is 
of the block form 
\[
  \varphi = 
  \left( 
  \begin{array}[h]{c|c}
    \mathbf 1_r & 0 \\
    \hline
    0 & 0
  \end{array}
  \right)
  \in \KK^{m\times n}
\]
where by $\mathbf 1_r$ we denote the $r\times r$ unit matrix.
It is now easy to see that the map 
\[
  \Phi \colon (\GL(r;\KK),\mathbf 1_r) \times (\KK^{(m-r)\times r},0) \times 
  (\KK^{r \times (n-r)},0) \times (\KK^{(m-r)\times(n-r)},0) \to \KK^{m\times n}, 
\]
taking a tuple $(A,P,Q,N)$ to the block matrix 
\[
  \left( 
  \begin{array}[h]{c|c}
    \mathbf 1_r & 0 \\
    \hline
    P \cdot A^{-1} & \mathbf 1_{m-r}
  \end{array}
  \right)\cdot
  \left( 
  \begin{array}[h]{c|c}
    A & 0 \\
    \hline
    0 & N
  \end{array}
  \right)\cdot
  \left( 
  \begin{array}[h]{c|c}
    \mathbf 1_r & A^{-1} \cdot Q \\
    \hline
    0 & \mathbf 1_{n-r}
  \end{array}
  \right)
  =
  \left( 
  \begin{array}[h]{c|c}
    A & Q \\
    \hline
    P & PA^{-1}Q + N
  \end{array}
  \right)
\]
yields a local isomorphism 
$(\KK^{(m+n)r -r^2},0) \overset{\cong}{\longrightarrow} (V_{m,n}^r,\varphi)$ 
by setting $N=0$ and restricting to the first three factors. The entries of 
$N \in \KK^{(m-r)\times (n-r)}$ 
can be understood as normal coordinates to that stratum and 
using Lemma \ref{lem:InvarianceOfDeterminantalIdealUnderLeftRightMultiplication}
we find that 
\begin{equation}
  \Phi \colon 
  (\KK^{(m+n)r - r^2},0) \times 
  (M_{m-r,n-r}^{s-r},0) 
  \overset{\cong}{\longrightarrow} 
  (M_{m,n}^s,\varphi) 
  \label{eqn:TrivializationOfTheRankStratification}
\end{equation}
is a stratification preserving isomorphism for every $s>r$. 
In particular, this local analytic triviality implies that the rank 
stratification satisfies Whitney's conditions $(a)$ and $(b)$, cf. 
\cite[Chapter 4.2]{Volume1}, so that is in fact a \textit{Whitney stratification} 
of the space of matrices.

\subsubsection{Logarithmic vector fields in $\CC^{m\times n}$}
\label{sec:LogarithmicVectorFieldsRankStratification}

The $G$-action (\ref{eqn:LieGroupActionOnTheSpaceOfMatrices}) 
on the space of matrices can also be exploited to construct 
so-called logarithmic vector fields, see \cite{Damon87} and \cite{Saito80}, 
which will be important for the theory of unfoldings later on. 

Let $U \subset \CC^n$ be an open domain and $V\hookrightarrow U$ the 
closed embedding of a complex analytic subvariety. The module 
of logarithmic vector fields\footnote{We use the notation $\Der(-\log(V))$ 
  rather than $\Der\log(V)$ following \cite{Goryunov21}.
}
for $V$ at a point $q \in U$ is given by 
\[
  \Der(-\log(V)) = \{ \xi \in T_{\CC^n,q} \colon \xi(I) \subset I \} \subset T_{\CC^n,q} 
\]
where $I = I(V)$ denotes the ideal of functions vanishing on $V$.
These form a coherent sheaf on $U$ which coincides with the tangent sheaf 
$T_{U}$ outside $V$.
In the real case, a similar construction can be made whenever 
the ideal $I$ is coherent, see \cite{Damon87}. We will restrict our exposition 
to the complex analytic setup.

The \textit{logarithmic tangent space} 
to $V$ at a point $q\in V$ is defined as the subspace 
\begin{equation}
  T^{\log}_q V = \Der(-\log(V))/(\Der(-\log(V)) \cap \m_q T_{\CC^n,q}) 
  \subset T_{\CC^n,q}/\m_q T_{\CC^n,q} = T_q \CC^n
  \label{eqn:DefinitiionLogarithmicTangentSpace}
\end{equation}
of vectors $v \in T_q \CC^n$ which extend locally to logarithmic vector fields 
for $V$ in $\CC^n$. 
It is easy to check that whenever $q$ is a smooth point of $V$ then one has 
$T_q V = T^{\log}_q V$. At singular points very little is known about 
$T_q^{\log} V$ in general, but if we endow $V$ with its \textit{canonical} 
Whitney stratification as in \cite{LeTeissier81},
$V = \bigcup_{i=0}^m V^i$,
then Damon and Mond have shown in \cite[Proposition 3.11]{DamonMond91} that 
the logarithmic tangent space to $V$ at a point $q$ 
\begin{equation}
  T_q^{\log} V \subset T_q V^i
  \label{eqn:LogarithmicTangentSpaceInStratumTangentSpace}
\end{equation}
is contained in the tangent space of the stratum $V^i$ containing $q$ 
so that logarithmic vector fields are always tangent to the strata of 
the \textit{canonical} Whitney stratification.
However, 
one need not have equality in (\ref{eqn:LogarithmicTangentSpaceInStratumTangentSpace}) 
and those strata $V^i$ for which equality holds at all points are called 
\textit{holonomic strata}.

The purpose of this section is to show:

\begin{lemma}
  The strata $V_{m,n}^r \subset M_{m,n}^s \subset \CC^{m\times n}$ of 
  the rank stratification of the generic determinantal varieties are 
  holonomic for all values of $r<s\leq \min\{m,n\}$.
  \label{lem:HolonomicityOfTheRankStratification}
\end{lemma}

To this end, let $\fgl_m \oplus \fgl_n \cong \CC^{m\times m} \oplus \CC^{n\times n}$ 
be the Lie-algebra of $G = \GL(m,\CC) \times \GL(n;\CC)$. For any vector 
$(L,R) \in \fgl_m \oplus \fgl_n$ the exponential map gives rise to a 
holomorphic $1$-parameter family
$t \mapsto \exp(t \cdot (L,R))$ in $G$. This induces a family of stratification 
preserving automorphisms 
\[
  \gamma \colon (\CC,0) \times \CC^{m\times n} \to \CC^{m\times n}, \quad 
  (t,\varphi) \mapsto \gamma_t(\varphi) = \exp(t\cdot(L,R))*\varphi.
\]
via the $G$-action on $\CC^{m\times n}$ and 
we denote by 
\begin{equation}
  \xi_{(L,R)}(\varphi) = \frac{\D}{\D t}|_{t=0} \exp(t \cdot (L,R)) * \varphi
  = L \cdot \varphi - \varphi \cdot R
  \label{eqn:LogarithmicVectorFieldsOfExponentialActions}
\end{equation}
the vector field on $\CC^{m\times n}$ generated by that action. Here, 
we deliberately identified $T_\varphi \CC^{m\times n} \cong \CC^{m\times n}$ 
with the space of matrices again.

\begin{proof}(of Lemma \ref{lem:HolonomicityOfTheRankStratification})
  We first remark that the rank stratification $M_{m,n}^s = \bigcup_{r<s} V_{m,n}^r$ 
  always coincides with the canonical Whitney stratification constructed by 
  L\^e and Teissier. 
  
  The vector fields $\xi = \xi_{(L,R)}$ constructed above are logarithmic for $M_{m,n}^s$: 
  If we let $\exp(t\cdot (L,R)) = (P_t,Q_t)$ be the associated $1$-parameter family in $G$ 
  and $\gamma \colon (t,\varphi) \mapsto \gamma_t(\varphi)$ the induced 
  family of automorphisms of $\CC^{m\times n}$, then 
  \begin{eqnarray*}
    \frac{\D }{\D t} (Y\circ \gamma_t)^{\wedge s} &=& 
    \frac{\D }{\D t} \left( P_t \cdot Y \cdot Q_t^{-1} \right)^{\wedge s} \\
    &=& 
    \frac{\D }{\D t} \left( P_t^{\wedge s} \cdot Y^{\wedge s} \cdot 
    \left(Q_t^{-1}\right)^{\wedge s} \right) \\
    &=& \left( \frac{\D }{\D t} P_t \right)^{\wedge s} \cdot Y^{\wedge s} + 
    Y^{\wedge s} \cdot \left( \frac{\D }{\D t} Q_t^{-1} \right)^{\wedge s}
  \end{eqnarray*}
  and therefore, setting $t=0$,
  clearly $\xi(Y^{\wedge s}_{\bf i, \bf j}) \in \langle Y^{\wedge s} \rangle$ 
  for every generator $Y^{\wedge s}_{\bf i, \bf j}$ of $\langle Y^{\wedge s}\rangle$.
  The claim now follows from the observation that, 
  since every stratum $V_{m,n}^r$ is a $G$-orbit in $\CC^{m\times n}$, 
  the linear map 
  \[
    \fgl_m \oplus \fgl_n \to T_\varphi V_{m,n}^r, \quad 
    (L,R) \mapsto \xi_{(L,R)}(\varphi)
  \]
  is surjective at every point $\varphi \in V_{m,n}^r$.
\end{proof}

\begin{remark}
  \label{rem:LogarithmicVectorFields}
  Note that even though $T_\varphi^{\log} M_{m,n}^s = T_\varphi V_{m,n}^r$ 
  for every $\varphi \in V_{m,n}^r \subset M_{m,n}^s$, the proof of 
  Lemma \ref{lem:HolonomicityOfTheRankStratification} does \textit{not} imply that 
  the module $\Der(-\log(M_{m,n}^s))$ is generated by the 
  vector fields $\xi_{(L,R)}$ in (\ref{eqn:LogarithmicVectorFieldsOfExponentialActions}). 
  This is true for the varieties $M_{m,m}^m$ of degenerate square, 
  symmetric, and skew-symmetric matrices, cf. \cite{Bruce03}, \cite{BruceTari04}, 
  and \cite{GoryunovMond05}, but false for the generic determinantal varieties 
  $M_{m,n}^s$ defined by non-maximal minors for $s<\min\{m,n\}$: Considering 
  $V_{m,n}^{s-1}$ as a stratum of $M_{m,n}^s$ only, we may extend any tangent vector 
  field $\zeta$ to $V_{m,n}^{s-1}$ in an arbitrary way to a vector field on a 
  neighborhood in $\CC^{m\times n}$ 
  while an extension as a linear combination of the $\xi_{(L,R)}$'s will necessarily 
  be tangent to the orbits $V_{m,n}^t$ for $t\geq s$, as well.
  Thus, the vector fields $\xi_{(L,R)}$ only generate a proper submodule of 
  $\Der(-\log(M_{m,n}^s))$ for non-maximal $s$.
\end{remark}

\section{Unfoldings and equivalence of matrices}
\label{sec:UnfoldingsAndEquivalenceOfMatrices}

According to the previous section, the unfoldings of a matrix $A$ determine 
the determinantal deformations of the associated singularities. 
In practice, unfoldings of map germs are far easier to handle and to classify 
than deformations; for instance there is no obstruction 
theory\footnote{cf. \cite[Section 7.1.5]{Volume1}} to be 
taken into account. The paradigm of the following sections is 
to consider determinantal singularities as 
hybrid objects living in both the world of unfoldings of map germs 
and the world of deformations of space germs at the same time, 
and to consider unfoldings of matrices as a
handle for studying deformations of the associated determinantal singularities. 
The appropriate framework for this was developed by Damon in \cite{Damon84} where 
he defines the notion of a ``geometric subgroup'' of the contact group $\mathcal K$. 
The natural notion 
of equivalence for germs of matrices ---the so-called $\GL$-equivalence defined below--- 
leads to such a geometric subgroup as was observed in \cite{Haslinger01},
\cite{Bruce03}, \cite{BruceTari04}, \cite{GoryunovMond05} for square matrices 
and in \cite{Pereira10} for the general case.
In principal, all relevant theorems about finite determinacy, versality, etc. 
can be derived from that. 
For a general account on map germs, their unfoldings, finite determinacy 
and related techniques see for instance 
\cite{Wall81}, or \cite{BallesterosMond20}. 

In this note we will give the specific statements 
for $\GL$-equivalence of matrices in the complex analytic category with a 
focus on the explicit description of the infinitesimal theory. 
Many constructions and results can be carried over to the real analytic or 
even the real differentiable setup. Moreover, in 
\cite{BelitskiiKerner16} Belitskii and Kerner develop an analogous theory 
of unfoldings, equivalences, and finite determinacy in a purely algebraic 
fashion; particular applications to families of matrices can be found in 
\cite{BelitskiiKerner19} and \cite{Kerner19} and various further preprints are available.

\medskip
We start with the natural notion of equivalence of map germs in this context.

\begin{definition}
  Two matrices $A,B \in \CC\{x_1,\dots,x_p\}^{m\times n}$ are called $\GL$-equivalent, if there exist 
  matrices $P \in \GL(m;\CC\{x\})$ and $Q\in \GL(n;\CC\{x\})$ 
  and a biholomorphism $\Phi \colon (\CC^p,0) \to (\CC^p,0)$ 
  such that 
  \begin{equation}
    P(x) \cdot \left( A \circ \Phi^{-1}(x) \right) \cdot (Q(x))^{-1} = B(x).
    \label{eqn:GLEquivalence}
  \end{equation}
\end{definition}
Note that in (\ref{eqn:GLEquivalence}) we write $\Phi^{-1}(x)$ for the inverse 
of the map $\Phi$ applied to $x$ while $(Q(x))^{-1}$ denotes the inverse of the matrix $Q(x)$ 
at that point.

Similarly, we say that $A$ and $B$ are $\SL$-equivalent if $P$ and $Q$ only take 
values in the special linear groups. In the following, we will deliberately identify 
the space of matrices $\GL(m;\CC\{x\})$ with the space of map germs 
$P \colon (\CC^p,0) \to (\GL(m;\CC),P(0))$. Depending on the context, either one of the 
two notations has its advantages.

\begin{remark}
  \label{rem:GLEquivalenceOfMatricesGivesIsomorphicGerms}
  Note that in the above definition, neither $A$ nor $B$ are required to 
  define a determinantal singularity. If they do, however, then it follows 
  directly from Lemma \ref{lem:InvarianceOfDeterminantalIdealUnderLeftRightMultiplication}
  that any two associated determinantal singularities $(X_A^s,0)$ and $(X_{B}^s,0)$ 
  are isomorphic as germs.
\end{remark}

\begin{definition}
  \label{def:GLEquivalenceOfUnfoldings}
  Let $A \in \CC\{x_1,\dots,x_p\}^{m \times n}$ be a matrix. 
  Two unfoldings of $A$ 
  on $k$ parameters $t= t_1,\dots,t_k$ 
  given by 
  $\mathbf A(x,t), \mathbf B(x,t) \in \CC\{x,t\}^{m\times n}$ 
  are called $\GL$-equivalent if there exist unfoldings 
  \begin{eqnarray*}
    (\CC^p,0) \times (\CC^k,0) \to (\GL(m,\CC),\mathbf 1_m),&\quad&
    (x,t) \mapsto P_t(x)\\
    (\CC^p,0) \times (\CC^k,0) \to (\GL(n,\CC),\mathbf 1_n),&\quad&
    (x,t) \mapsto Q_t(x)\\
    (\CC^p,0)\times(\CC^k,0)\to (\CC^p,0),& \quad & 
    (x,t) \mapsto \Phi_t(x)
  \end{eqnarray*}
  of the identities $P_0 = \mathbf 1_m \in \GL(m;\CC\{x\})$, 
  $Q_0 = \mathbf 1_n \in \GL(n;\CC\{x\})$ 
  and $\Phi_0 = \Id_{\CC^p,0}$
  such that
  \[
    \mathbf B(x,t) = 
    P_t(x) \cdot \mathbf A( \Phi_t^{-1}(x),t) \cdot (Q_t(x))^{-1}
  \]
  An unfolding given by $\mathbf A(x,t)$ 
  is $\GL$-trivial (or just ``trivial'') if it is $\GL$-equivalent (as unfoldings) to 
  the map $\mathbf B(x,t) = (A(x),t)$.
\end{definition}
Again, we obtain the notion of $\SL$-equivalence of unfoldings by substituting 
$\SL$ for $\GL$ in the above definition. 

\begin{remark}
  \label{rem:GLEquivalentUnfoldingsGiveIsomorphicDeformations}
  Note that whenever $A$ defines a determinantal 
  singularity $(X_A^s,0) \subset (\CC^p,0)$ of type $s$ and $\mathbf A$ and $\mathbf B$ 
  are two $\GL$-equivalent unfoldings of $A$, the resulting flat families
  are isomorphic in the sense that there is a commutative diagram 
  \[
    \xymatrix{
      (\mathcal X_{\mathbf A}^s,0) \ar[rr]^\cong \ar[dr] & & 
      (\mathcal X_{\mathbf B}^s,0) \ar[dl] \\
      & (\CC^k,0). & 
    }
  \]
  In particular, a $\GL$-trivial unfolding $\mathbf A$ of $A$ gives rise to a 
  product $(\mathcal X_{\mathbf A}^s,0) \cong (X_A^s,0) \times (\CC^k,0)$.
\end{remark}

\medskip

We wish to classify all unfoldings of a given matrix $A \in \CC\{x\}^{m\times n}$
up to $\GL$-equivalence, thereby also capturing all possible determinantal deformations 
of the associated singularities $(X_A^s,0) \subset (\CC^p,0)$. While in the 
context of deformations of complex analytic germs this leads to the investigation 
of deformation functors over Artin rings developed by Schlessinger 
\cite{Schlessinger68} and in particular 
the \textit{first order deformations} (cf. \cite[Section 7.1.4]{Volume1}), the 
common viewpoint for unfoldings of map germs is to consider the action 
of an (infinite dimensional) algebraic group on them. 
The particular group in question for $\GL$-equivalence is the 
semi-direct product
\begin{equation}
  \mathcal G = \left(\GL(m;\CC\{x\}) \times \GL(n;\CC\{x\}) \right)\rtimes \mathrm{Diff}(\CC^p,0)
  \label{eqn:TheGroupG}
\end{equation}
with composition defined by 
\[
  (P,Q,\Phi) * (P',Q',\Psi) = 
  \left(P(x) \cdot P'(\Phi^{-1}(x)), 
  Q(x) \cdot Q'(\Phi^{-1}(x)), \Psi(\Phi(x))\right)
\]
so as to be compatible with the left action (\ref{eqn:GLEquivalence}) on the 
space of matrices $\CC\{x\}^{m\times n}$. 

As already mentioned in the introduction, 
the group $\mathcal G$ 
is a so called ``geometric subgroup'', 
a notion introduced by Damon in \cite{Damon84}, 
of the \textit{contact group} $\mathcal K$; see e.g. \cite[Proposition 2.5.1]{Pereira10}. 
This allows us to pursue a common path in the theory for unfoldings of map germs. 
The key object for further studies is the space $T^1_\GL(A)$ capturing 
the nontrivial \textit{infinitesimal} unfoldings of a given matrix $A$ up 
to $\GL$-equivalence. We shall introduce it now.

\medskip

For a given matrix $A \in \CC\{x_1,\dots,x_p\}^{m\times n}$, the trivial unfoldings 
of $A$ are those captured by the action of a $1$-parameter family 
$(P_t,Q_t,\Phi_t)$ in $\mathcal G$ with 
$(P_0,Q_0,\Phi_0) = (\mathbf 1_m, \mathbf 1_n,\Id_{\CC^p,0})$. Such unfoldings 
take the form 
\[
  \mathbf A(x,t) = P_t(x)\cdot A(\Phi^{-1}_t(x)) \cdot (Q_t(x))^{-1}.
\]
Differentiating with respect to $t$ at $t=0$ gives us the \textit{infinitesimally} trivial 
unfoldings
\[
  \frac{\D P}{\D t}(0) \cdot A - A \cdot \frac{\D Q}{\D t}(0)- 
  \sum_{i=1}^p \frac{\partial A}{\partial x_i} \cdot \frac{\D \Phi_i}{\D t}(0).
\]
These generate the so-called 
\textit{extended tangent space} (to the orbit) of $A$:
\begin{equation}
  T_e\mathcal G(A) = \fgl_m(\CC\{x\}) \cdot A - A \cdot \fgl_n(\CC\{x\}) + 
  \left\langle \frac{\partial A}{\partial x_1}, \dots, \frac{\partial A}{\partial x_p} \right\rangle
  \subset \CC\{x\}^{m\times n}, 
  \label{eqn:ExtendedTangentSpaceForGLEquivalence}
\end{equation}
cf. for instance \cite[Part I, Section 1]{Wall81}, 
or \cite[Proposition 2.5.1]{Pereira10}. We wrote $\fgl_m(\CC\{x\})$ for 
the space $\CC\{x\}^{m\times m}$ in which $\D P/\D t(0)$ lays, and vice versa for $\fgl_n(\CC\{x\})$.

Geometrically, the extended tangent space $T_e\mathcal G(A)$ can be described 
as follows. Consider the pullback of the tangent bundle $A^{*} (T\CC^{m\times n})$ 
of $\CC^{m\times n}$ along $A$. The sheaf of sections in this bundle 
is a free $\CC\{x\}$-module with stalk 
$\left( A^* T_{\CC^{m\times n}} \right)_0 \cong \CC\{x\}^{m\times n}$ at the origin. 
Then the extended tangent space is the submodule generated by 
the pullback of the specific 
logarithmic vector fields (\ref{eqn:LogarithmicVectorFieldsOfExponentialActions}),
\[\fgl_m(\CC\{x\}) \cdot A - A \cdot \fgl_n(\CC\{x\}),\]
and the image of the differential  
$\D A \colon T_{\CC^p,0} \to A^*T_{\CC^{m\times n},0}$ of $A$, 
\[
  \left\langle \frac{\partial A}{\partial x_1}, \dots, \frac{\partial A}{\partial x_p} \right\rangle.
\]
With $T_e\mathcal G(A) \subset \CC\{x\}^{m\times n}$ being the trivial unfoldinds, 
the quotient by this submodule 
\begin{equation}
  T^1_{\GL}(A) := \CC\{x\}^{m\times n}/T_e \mathcal G(A)
  \label{eqn:MatrixT1}
\end{equation}
captures the \textit{nontrivial} infinitesimal unfoldings of $A$ up to $\GL$-equivalence.
We will in the following refer to the dimension of $T^1\mathcal G(A)$ over $\CC$  
as the $\GL$-\textit{codimension} of $A$:
\begin{equation}
  \tau_{\GL}(A) = \dim_\CC T^1_\GL (A).
  \label{eqn:GLCodimension}
\end{equation}

\begin{remark}
  \label{rem:SLEquivalence}
  Besides $\GL$-equivalence some authors also consider other variants such 
  as $\SL$-equivalence, see for example \cite{GoryunovMond05}, or \cite{Goryunov21}.
  \begin{description}
    \item[\textbf{$\SL$-equivalence:}] For $\SL$-equivalence, the group $\mathcal S$ is 
      defined as in (\ref{eqn:TheGroupG}), only that 
      the matrices $P$ and $Q$ are restricted to take values in the subgroup 
      $\SL(m;\CC\{x\})$ and $\SL(n;\CC\{x\})$. Accordingly, for the extended tangent 
      space of a matrix $A$ one finds 
      \begin{equation}
	T_e \mathcal S(A) = \mathfrak{sl}_m(\CC\{x\}) \cdot A + A \cdot \mathfrak{sl}_n(\CC\{x\}) + 
	\left\langle \frac{\partial A}{\partial x_1}, \dots, \frac{\partial A}{\partial x_p} \right\rangle
	\subset \CC\{x\}^{m\times n}
	\label{eqn:ExtendedTangentSpaceForSLEquivalence}
      \end{equation}
      with $\mathfrak{sl}_m(\CC\{x\})$ the set of trace-free matrices in $\CC\{x\}^{m\times m}$ 
      and vice versa for the other term. In this case, we speak of the $\SL$-\textit{codimension} 
      of a matrix 
      \[
	\tau_{\SL}(A) = \dim_{\CC} T^1_\SL(A).
      \]

    \item[\textbf{Symmetric matrices:}] 
      Adaptations of $\GL$- and $\SL$-equivalence can also be made for symmetric matrices 
      $A \in \CC\{x_1,\dots,x_p\}^{m\times m}_\sym$ where one usually considers 
      the group 
      \begin{equation}
	\mathcal G_\sym = \GL(m;\CC\{x\}) \rtimes \mathrm{Diff}(\CC^p,0)
	\label{eqn:TheGroupGForSymmetricMatrices}
      \end{equation}
      with composition 
      \[
	(P,\Phi) * (P',\Psi) = (P \cdot (P'\circ \Phi^{-1}), \Phi\circ \Psi)
      \]
      and action on $\CC\{x\}^{m\times m}_\sym$ given by 
      \[
	((P,\Phi),A) \mapsto P \cdot (A \circ \Phi^{-1} ) \cdot P^T.
      \]
      In this case, the extended tangent space of a matrix $A$ is 
      \begin{equation}
	T_e \mathcal G_\sym (A) = 
	\left\langle M \cdot A + A \cdot M^T : M \in \fgl_m(\CC\{x\}) \right\rangle 
	+ \left\langle \frac{\partial A}{\partial x_1}, \dots, \frac{\partial A}{\partial x_p} \right\rangle.
	\label{eqn:ExtendedTangentSpaceForSymmetricGLEquivalence}
      \end{equation}

    \item[\textbf{Skew symmetric matrices:}]
      For skew-symmetric matrices $A \in \CC\{x\}^{2m \times 2m}_\sk$
      and their Pfaffian ideals it is customary to consider 
      the same equivalences as for symmetric matrices. Then the extended 
      tangent space for $\SL$-equivalence, for example, reads 
      \begin{equation}
	T_e \mathcal S_\sk (A) = 
	\left\langle M \cdot A + A \cdot M^T : M \in \fsl_m(\CC\{x\}) \right\rangle 
	+ \left\langle \frac{\partial A}{\partial x_1}, \dots, \frac{\partial A}{\partial x_p} \right\rangle
	\label{eqn:ExtendedTangentSpaceForSymmetricsLEquivalence}
      \end{equation}
      where again $\fsl_m(\CC\{x\})$ denotes the traceless matrices.
  \end{description}
  
\end{remark}

\subsection{Finite determinacy}
A natural question for map germs is whether or not they are 
finitely determined for a given notion of equivalence. 
Analogous to the case of holomorphic functions we define the $k$-jet of a matrix
$A \in \CC\{x\}^{m\times n}$ to be the Taylor expansion 
\[
  j^k(A) = A(0) + 
  \sum_{|\alpha|=1} \frac{x^\alpha}{\alpha!}\left(\frac{\partial }{\partial x}\right)^\alpha A|_{x=0}
  + \dots +
  \sum_{|\alpha|=k} \frac{x^\alpha}{\alpha!}\left(\frac{\partial }{\partial x}\right)^\alpha A|_{x=0} 
\]
of the entries of $A$ up to order $k$ modulo $\m^{k+1}\CC\{x\}^{m\times n}$.
As usual, $\alpha = (\alpha_1,\dots,\alpha_p) \in \NN_0^p$ denotes a multi-index with 
$\alpha! = \alpha_1! \cdots \alpha_p!$, $|\alpha| = \alpha_1 + \dots + \alpha_p$ and 
$(\partial/\partial x)^\alpha = (\partial/\partial x_1)^{\alpha_1} \cdots 
(\partial /\partial x_p)^{\alpha_p}$.

\begin{definition}
  \label{def:FiniteDeterminacy}
  A matrix $A \in \CC\{x_1,\dots,x_p\}^{m\times n}$ is $k$-determined (for $\GL$-equivalence)
  if for every other matrix $B$ an equality of jets $j^k(A) = j^k(B)$ implies that 
  $B$ is $\GL$-equivalent to $A$. 
\end{definition}
Straightforward adaptations can be given for the other groups discussed in 
Remark \ref{rem:SLEquivalence}.
We say that $A$ is finitely $\GL$-determined, if it is $k$-determined for some $k>0$. 
In particular, any finitely $\GL$-determined matrix is $\GL$-equivalent to a matrix with 
polynomial entries. 

The following explicit infinitesimal criterion for finite 
$\GL$-determinacy of matrices has been given by Pereira in \cite[Theorem 2.3.1]{Pereira10}:

\begin{theorem}
  Let $A \in \CC\{x_1,\dots,x_p\}^{m\times n}$ be a matrix and $k$ an integer such that 
  \[
    \m^{k+1} \CC\{x\}^{m\times n} \subset \m^2
    \left\langle 
    \frac{\partial A}{\partial x_1},\dots, \frac{\partial A}{\partial x_p} 
    \right\rangle
    + \m\cdot \left( \mathfrak{gl}_m(\CC\{x\}) \cdot A + A \cdot \mathfrak{gl}_n(\CC\{x\}) \right).
  \]
  Then $A$ is $k$-determined for $\GL$-equivalence.
  \label{thm:InfinitesimalCriterionForFiniteDeterminacy}
\end{theorem}

It should be pointed out, that Pereira also covered the 
real analytic case and the above theorem is the adapted citation for 
holomorphic matrices. While this criterion is useful for 
explicit computations, there is also another, more geometric criterion 
based on the transversality of maps. 
This was observed by Bruce \cite[Proposition 3.2]{Bruce03} for 
the specific case of symmetric matrices and later carried out explicitly 
by Pereira \cite[Theorem 2.4.1]{Pereira10} for $\GL$-equivalence 
for matrices of arbitrary size.

\medskip

We briefly recall the classical notion of 
transversality\footnote{Following Damon \cite{Damon01}, we will also refer to this as 
\textit{geometric transversality} in order to distinguish it from the 
\textit{algebraic transversality} which will be introduced in the next section.
}
for smooth maps of mani\-folds. 
Let $M$ and $U$ be smooth manifolds and $V \subset M$ a locally closed 
submanifold. We say that a map $f \colon U \to M$ is \textit{transverse} to $V$ 
at a point $p \in U$, if either $f(p) \notin V$, or $f(p) \in V$ and the 
tangent spaces $T_p U$ and $T_{f(p)} V$
\begin{equation}
  \D f(p) T_p U + T_{f(p)} V = T_{f(p)} M
  \label{eqn:GeometricTransversality}
\end{equation}
span the whole tangent space $T_{f(p)} M$ of the ambient manifold $M$; cf.  
\cite[Definition 4.2.11]{Volume1}.
We say that $f$ is transverse to $V$ on $U$ if this holds at every point in $U$.

\begin{figure}[h]
  \centering
  \includegraphics[scale=0.4,clip=true,trim=14cm 3cm 2cm 3cm]{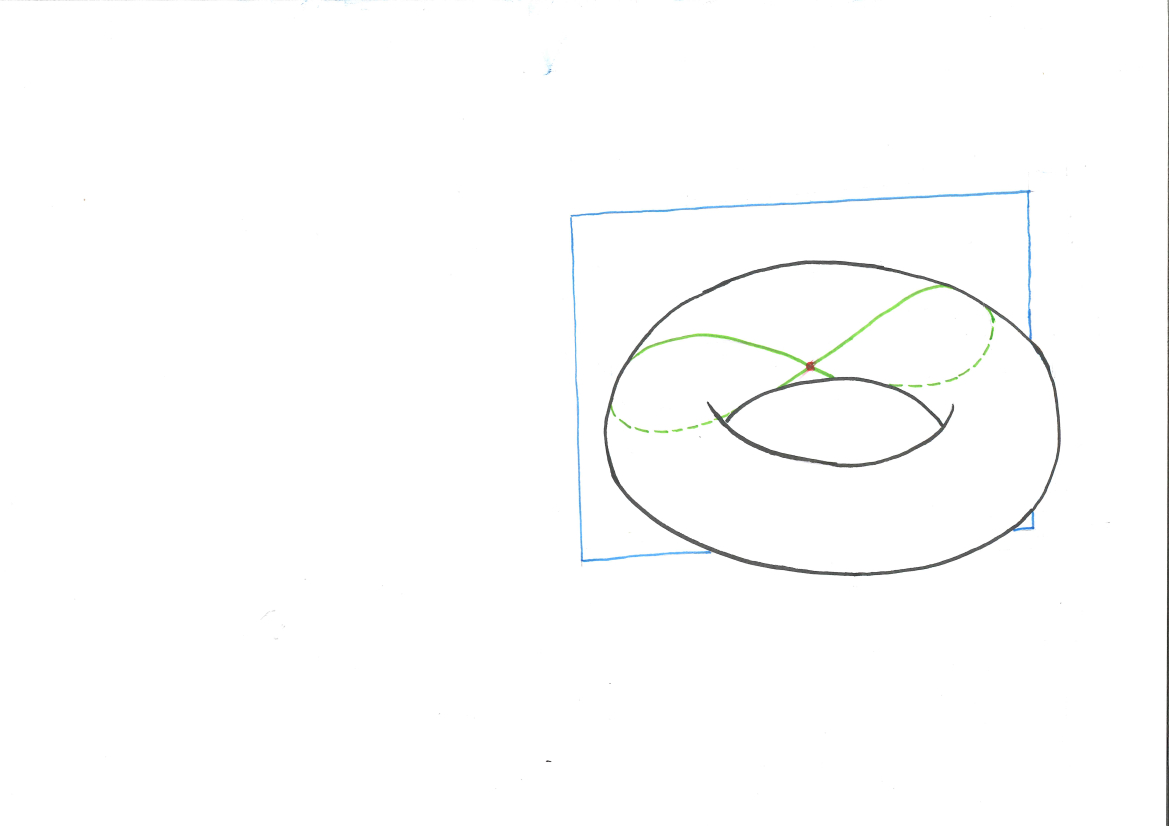}
  \includegraphics[scale=0.4,clip=true,trim=14cm 3cm 2cm 3cm]{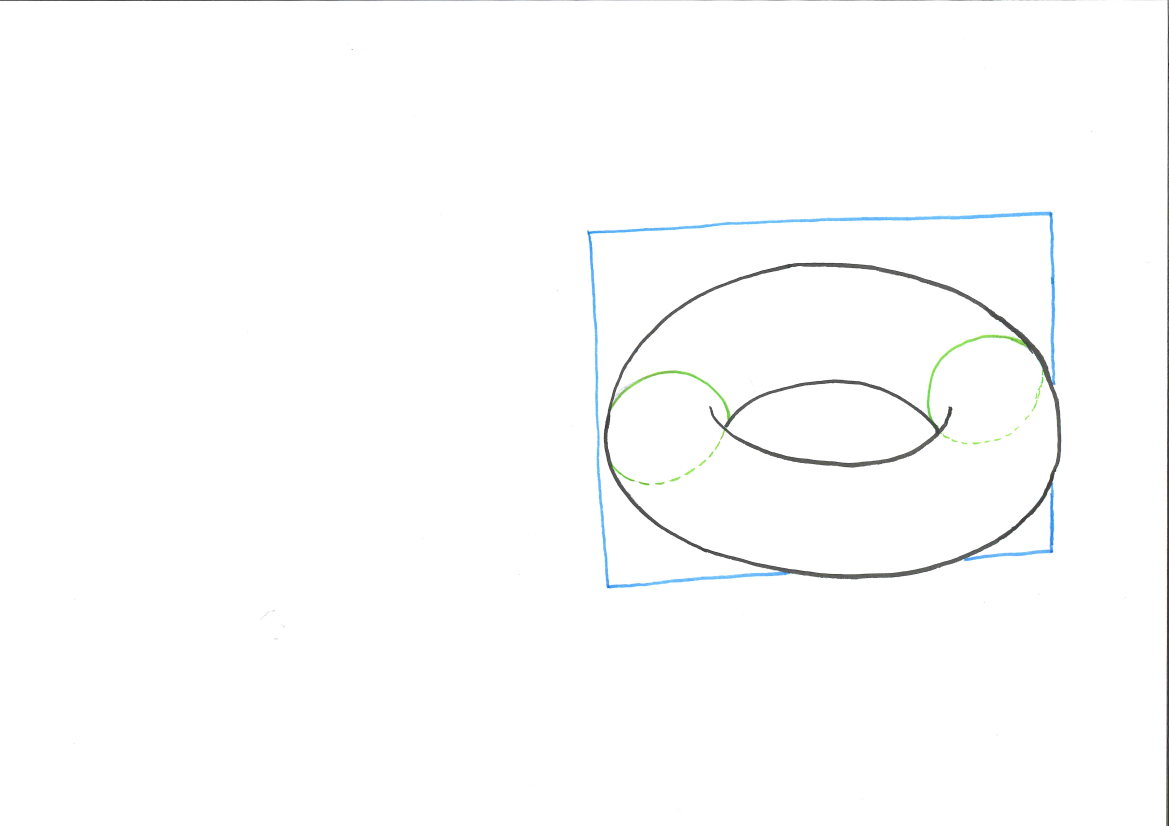}
  \caption{A transverse and a non-transverse intersection (green) of a torus 
    (black) in $\RR^3$ with the immersion of 
  a plane (blue).  }
  \label{fig:TransverseIntersectionExample}
\end{figure}

When $\bigcup_{i} V^i= M$ is a stratification of 
$M$ into smooth, locally closed submanifolds we say that a map $f$ as above 
is transverse to the stratification, if it is transverse to every stratum $V^i$. 

With these notions at hand we now have \cite[Theorem 2.4.1]{Pereira10} (cf. also 
\cite[Proposition 3.2]{Bruce03} for the symmetric and \cite[Proposition 3.2]{BruceTari04} 
for the arbitrary symmetric case):
\begin{theorem}
  Let $A \colon (\CC^p,0) \to (\CC^{m\times n},A(0))$ be a holomorphic map germ. 
  Then $A$ is finitely $\GL$-determined if and only if it is transverse 
  to the rank stratification on $\CC^{m\times n}$ in a punctured neighborhood of 
  the origin in $\CC^p$.
  \label{thm:GeometricCriterionForFiniteDeterminacy}
\end{theorem}

\subsection{Versal unfoldings}

An unfolding $F$ for a given map germ $f$ is called \textit{versal}, 
if every other unfolding $F'$ of $f$ can be written as a pullback from $F$ 
up to the underlying notion of equivalence of map germs. For an account 
on these notions, see e.g. \cite{Wall81} or \cite{BallesterosMond20}. 
In the explicit case of matrices and $\GL$-equivalence we can give the following 

\begin{definition}
  \label{def:GLVersality}
  Let $A \in \CC\{x_1,\dots,x_p\}^{m\times n}$ be a matrix. An unfolding 
  \[
    \mathbf A \colon (\CC^p,0) \times (\CC^k,0) \to (\CC^{m\times n},0)
  \]
  of $A$ on $k$ parameters $t$ 
  is called $\GL$-versal, if for every other unfolding 
  $\mathbf B$ of $A$ on $l$ parameters $s$ 
  there exists a holomorphic map germ 
  $h \colon (\CC^l,0) \to (\CC^k,0)$ such that $\mathbf B$ is $\GL$-equivalent 
  as an unfolding (Definition \ref{def:GLEquivalenceOfUnfoldings}) to the 
  unfolding of $A$ given by 
  \[
    A'(x,s) = A(x,h(s)).
  \]
  An unfolding $\mathbf A$ of $A$ is called 
  $\GL$-\textit{miniversal} (or semi-universal), if it is versal and 
  the number of parameters $k$ is minimal among all versal unfoldings of $A$.
\end{definition}

Versal unfoldings of matrices on a finite set of parameters as above 
do not necessarily exist. We saw earlier that the quotient $T^1_\GL(A)$ 
in (\ref{eqn:MatrixT1}) classifies all \textit{infinitesimal unfoldings} 
up to $\GL$-equivalence. If this space is of finite dimension, 
it can be used to construct miniversal unfoldings.

\begin{theorem}
  Let $A \in \CC\{x_1,\dots,x_p\}^{m\times n}$ be a matrix such that 
  $T^1_\GL(A)$
  has finite dimension $\tau$ as a 
  $\CC$-vector space. For any set of elements 
  $B_1,\dots,B_\tau \in \CC\{x\}^{m\times n}$ reducing to a 
  basis of $T^1_\GL(A)$, the unfolding on $\tau$ parameters given by 
  \[
    \mathbf A(x,t) = A(x) + t_1 \cdot B_1(x) + \dots + t_\tau \cdot B_\tau(x)
  \]
  is miniversal.
  \label{thm:ExistenceSemiUniversalUnfolding}
\end{theorem}
In particular, we see that the minimal number of parameters of a miniversal unfolding 
is always equal to $\dim_{\CC} T^1_\GL(A)$.

Theorem \ref{thm:ExistenceSemiUniversalUnfolding} is a special instance 
of the Unfolding Theorem of Damon \cite[Theorem 9.3]{Damon84} for 
geometric subgroups of $\mathcal K$.
As already remarked earlier, the theory developed by Damon allows one to also consider 
differentiable or real analytic setups.
For another, explicit proof of Theorem \ref{thm:ExistenceSemiUniversalUnfolding} 
which does not rely on Damon's work, see \cite[Theorem 1.4.10]{Zach17}. 

\medskip

In parallel to the previous section, we 
also give a geometric criterion for the condition $\dim_\CC T^1_\GL(A) < \infty$ 
to be satisfied,  analogous to Theorem \ref{thm:GeometricCriterionForFiniteDeterminacy}. 

\begin{definition}(Ebeling, Gusein-Zade \cite{EbelingGuseinZade09})
  \label{def:EssentiallyNonSingularPoint}
  Let $A \colon U \to \CC^{m\times n}$ be a holomorphic map on some open 
  subset $U \subset \CC^p$. A point $x\in U$ is called 
  \textit{essentially nonsingular}, if $A$ is transverse to the stratum 
  $V_{m,n}^r$ containing $A(x)$. 
\end{definition}
Here $r = \rank A(x)$ is the rank of the matrix $A$ at $x$. 
Note that the transversality of $A$ to the stratum $V_{m,n}^r$ 
at an essentially non-singular point $x$ 
and the Whitney-$(a)$-regularity of the rank stratification 
already imply that in a neighborhood of $x$ the 
map $A$ is stratified transversal to all 
strata $V_{m,n}^s$ for every $s\geq r$. In particular, 
the singularities $(A^{-1}(M_{m,n}^s),x) \subset (U,x)$
all have expected codimension. 

\begin{proposition}
  Let $A \colon (\CC^p,x) \to (\CC^{m\times n},A(x))$ be a holomorphic 
  map germ. Then $T^1_\GL(A) = 0$ if and only if $x$ is an essentially nonsingular 
  point of $A$.
  \label{prp:AlgebraicTransversality}
\end{proposition}
We include a brief proof of Proposition 
\ref{prp:AlgebraicTransversality} based on the notion of 
\textit{algebraic transversality} due to Damon in \cite{Damon96} 
(see also \cite{Damon01}). 

\medskip

Suppose $f \colon U \to M$ is a holomorphic map of complex manifolds and 
$V \subset M$ is a subvariety of $M$ with $I = I(V)$ the sheaf of ideals 
of functions vanishing on $V$. Then $f$ is said to be 
\textit{algebraically transverse} to $V$ at a point $p\in U$, if 
\begin{equation}
  \D f(p)T_pU + T^{\log}_{f(p)} V
  = T_{f(p)} M
  \label{eqn:AlgebraicTransversality}
\end{equation}
where $T^{\log}_q V$ is the logarithmic tangent space to $V$ at the 
point $q$, see Section \ref{sec:LogarithmicVectorFieldsRankStratification}.

\begin{proof}(of Proposition \ref{prp:AlgebraicTransversality})
  In the setup of Proposition \ref{prp:AlgebraicTransversality},
  and the generic determinantal varieties with $M_{m,n}^s \subset \CC^{m\times n}$ 
  in place of $V$ and $M$
  we find 
  \[
    T^{\log}_\varphi M_{m,n}^s = T_\varphi V_{m,n}^r
  \]
  for every $\varphi \in V_{m,n}^r \subset M_{m,n}^s$ according to 
  the holonomicity of the rank stratification, 
  Lemma \ref{lem:HolonomicityOfTheRankStratification}. Consequently, 
  a holomorphic map $A \colon U \to \CC^{m\times n}$ is algebraically 
  transversal to $M_{m,n}^s$ at a point $x\in U$ 
  if and only if it is geometrically transversal.
  Moreover, we saw in the proof of Lemma \ref{lem:HolonomicityOfTheRankStratification} 
  that the logarithmic tangent space $T_{\varphi}^{\log}M_{m,n}^s$ at 
  the point $\varphi = A(x)$ is spanned by the specific 
  vector fields $\xi = L \cdot \varphi - \varphi \cdot R$ introduced 
  in (\ref{eqn:LogarithmicVectorFieldsOfExponentialActions}).
  Using Nakayama's lemma, it is now easy to see 
  from the geometric description of the extended tangent space 
  $T_e \mathcal G(A)$ in (\ref{eqn:ExtendedTangentSpaceForGLEquivalence}) 
  that the module $T^1_\GL(A) = A^* T_{\CC^{m\times n}}/T_e \mathcal G(A)$ 
  is zero at $x$ if and only if $A$ is transverse to $M_{m,n}^s$ 
  (in either sense); the assertion follows.
\end{proof}

\begin{corollary}
  \label{cor:GeometricCriterionForSemiUniversality}
  Let $A \colon (\CC^p,0) \to (\CC^{m\times n},A(0))$ be a holomorphic map germ. 
  Then a miniversal unfolding of $A$ exists if and only if $A$ is transverse 
  to the rank stratification of $\CC^{m\times n}$ in a punctured neighborhood of 
  the origin.
\end{corollary}

\begin{proof}
  This follows directly from sheafification of the module $T^1_\GL(A)$, 
  cf. \cite[Section 4.5.1]{BallesterosMond20}. It 
  has finite length if and only if it is supported only at $0 \in \CC^p$ 
  and according to Proposition \ref{prp:AlgebraicTransversality} this is the case 
  whenever $A$ is transverse to the rank stratification off the origin.
\end{proof}

\begin{example}
  \label{exp:ThreeCoordinateAxisSemiUniversalUnfolding}
  Consider again the space curve singularity from Example 
  \ref{exp:DeformationOfThreeCoordinateAxisFlatness} given by the union of the three 
  coordinate axis in $(\CC^3,0)$. The defining matrix was 
  \[
    A=
    \begin{pmatrix} 
      x & 0 & z \cr
      0 & y & z 
    \end{pmatrix}
  \]
  and we will briefly indicate how to compute the space $T^1_\GL(A)$.

  To shorten notation, let $R = \CC\{x,y,z\}$. 
  We consider the extended tangent space of $A$ as a submodule 
  $T_e \mathcal G (A) \subset R^{2\times 3} \cong R^6$ and denote by 
  $E_{i,j}$ the generator of $R^{2\times 3}$ with $1$ at the $(i,j)$-th entry 
  and zeroes elsewhere. Then 
  the generators of $T_e \mathcal G(A)$ in the last summand 
  of (\ref{eqn:ExtendedTangentSpaceForGLEquivalence}) are 
  \[
    \frac{\partial A}{\partial x} = 
    \begin{pmatrix}
      1 & 0 & 0 \\
      0 & 0 & 0 
    \end{pmatrix}, \quad
    \frac{\partial A}{\partial y} = 
    \begin{pmatrix}
      0 & 0 & 0 \\
      0 & 1 & 0 
    \end{pmatrix}, \quad
    \frac{\partial A}{\partial z} = 
    \begin{pmatrix}
      0 & 0 & 1 \\
      0 & 0 & 1 
    \end{pmatrix}
  \]
  which allows us to reduce all multiples of $E_{1,1}$ and $E_{2,2}$ 
  to zero and write every multiple of $E_{1,3}$ as a multiple of $E_{2,3}$ 
  modulo $T_e\mathcal G(A)$. We may then proceed to show that for 
  every one of the remaining generators $E_{i,j}$ of $R^{2\times 3}$ we have 
  $\m \cdot E_{i,j} \subset T_e \mathcal G(A)$ where $\m = \langle x,y,z \rangle$.
  For instance for $E_{1,2}$ we find
  \[
    x\cdot E_{1,2} = A \cdot 
    \begin{pmatrix}
      0 & 1 & 0 \\
      0 & 0 & 0 \\
      0 & 0 & 0 \\
    \end{pmatrix},
    \, 
    y \cdot E_{1,2} = 
    \begin{pmatrix}
      -1 & 1 \\
      0 & 0 
    \end{pmatrix}
    \cdot 
    A
    + x\frac{\partial A}{\partial x}, \, 
    z \cdot E_{1,2} = A \cdot 
    \begin{pmatrix}
      0 & 0 & 0 \\
      0 & 0 & 0 \\
      0 & 1 & 0
    \end{pmatrix}
    - 
    z\frac{\partial A}{\partial y}
  \]
  and similarly for $E_{2,1}$ 
  \[
    x\cdot E_{2,1} = 
    \begin{pmatrix}
      0 & 0 \\
      1 & -1
    \end{pmatrix}
    \cdot A 
    + y \frac{\partial A}{\partial y}, \,
    y\cdot E_{2,1} = A \cdot
    \begin{pmatrix}
      0 & 0 & 0 \\
      1 & 0 & 0 \\
      0 & 0 & 0 
    \end{pmatrix},\,
    z\cdot E_{2,1} = A \cdot 
    \begin{pmatrix}
      0 & 0 & 0 \\
      0 & 0 & 0 \\
      1 & 0 & 0 \\
    \end{pmatrix} 
    - z \frac{\partial A}{\partial x}.
  \]
  The generators of $\m \cdot E_{2,3}$ are 
  \[
    x\cdot E_{2,3} = x\frac{\partial A}{\partial z} - A \cdot 
    \begin{pmatrix}
      0 & 0 & -1 \\
      0 & 0 & 0 \\
      0 & 0 & 0 \\
    \end{pmatrix},\,
    y\cdot E_{2,3} =
    A \cdot 
    \begin{pmatrix}
      0 & 0 & 0 \\
      0 & 0 & 1 \\
      0 & 0 & 0 \\
    \end{pmatrix},\,
    z\cdot E_{2,3} = 
    \begin{pmatrix}
      0 & 0 \\
      0 & 1
    \end{pmatrix}
    \cdot A - 
    y \frac{\partial A}{ \partial y}.
  \]
  Further calculations yield that a $\CC$-basis for 
  $T^1_\GL(A)$ is indeed given by the matrices 
  \[
    \begin{pmatrix}
      0 & 1 & 0 \\
      0 & 0 & 0 \\
    \end{pmatrix}, \quad
    \begin{pmatrix}
      0 & 0 & 0 \\
      1 & 0 & 0 \\
    \end{pmatrix}, \quad
    \begin{pmatrix}
      0 & 0 & 0 \\
      0 & 0 & 1 \\
    \end{pmatrix}
  \]
  so that we get a miniversal unfolding of $A$ by setting 
  \begin{equation}
    \mathbf A(x,y,z;t_1,t_2,t_3) = 
    \begin{pmatrix} 
      x & 0 & z \cr
      0 & y & z 
    \end{pmatrix} + 
    \begin{pmatrix}
      0 & t_1 & 0 \\
      t_2 & 0 & t_3
    \end{pmatrix}
    \label{eqn:MiniversalUnfoldingThreeCoordinateAxis}
  \end{equation}
  according to Theorem \ref{thm:ExistenceSemiUniversalUnfolding}.
\end{example}

\subsection{Discriminants of matrices}

Let $A \colon (\CC^p,0) \to (\CC^{m\times n},0)$ be a matrix and 
$\mathbf A \colon (\CC^p,0) \times (\CC^k,0) \to (\CC^{m\times n},0)$ an 
unfolding of $A$ on $k$ parameters $t_1,\dots,t_k$. For such an unfolding, one can 
define the \textit{relative}\footnote{cf. e.g. \cite[Definition 3.9]{BallesterosMond20}} space
$T^1_\GL(\mathbf A)$ of 
infinitesimal unfoldings of the fibers $A_t$ as the 
quotient of $\CC\{x,t\}^{m\times n}$ by the submodule 
\begin{equation}
  T_e\mathcal G(\mathbf A) = \CC\{x,t\}^{m\times m} \cdot \mathbf A 
  + \mathbf A \cdot \CC\{x,t\}^{n\times n} + 
  \left\langle \frac{\partial \mathbf A}{\partial x_1}, \dots, 
  \frac{\partial \mathbf A}{\partial x_p} \right\rangle
  \subset \CC\{x,t\}^{m\times n}. 
  \label{eqn:RelativeExtendedTangentSpaceForGLEquivalence}
\end{equation}
It is immediate from the definition, that 
$T^1_\GL(A) \cong T^1_\GL(\mathbf A)/\langle t_1,\dots,t_k\rangle T^1_\GL(\mathbf A)$.
Using the Weierstrass finiteness theorem\footnote{See e.g. \cite[Theorem 1.10]{GreuelLossenShustin07}.}
it is easy to see that $T^1_\GL(A)$ is a finite $\CC\{t\}$-module whenever 
$T^1_{\GL}(A)$ has finite dimension over $\CC$. 

\begin{definition}
  \label{def:DeterminantalDiscriminant}
  Let $A \in \CC\{x_1,\dots,x_t\}^{m\times n}$ be a matrix with $\dim_\CC T^1_\GL(A)< \infty$. 
  For an unfolding $\mathbf A$ of $A$ on $k$ parameters $t = (t_1,\dots,t_k)$, 
  the \textit{matrix discriminant} $(\Delta_{\mathbf A},0) \subset (\CC^k,0)$ is 
  the support of the $\CC\{t\}$-module $T^1_\GL(A)$.
  If $\mathbf A$ is a miniversal unfolding, then we also speak of the matrix 
  discriminant $(\Delta_A,0)$ of $A$ rather than the discriminant of the unfolding $\mathbf A$.
\end{definition}

Being the support of a finite analytic module, the matrix discriminant is a 
complex analytic set. In general we can decompose the discriminant into 
components 
\begin{equation}
  \Delta_{\mathbf A} = \bigcup_{0\leq r < \min\{m,n\}} \Delta_{\mathbf A}^r
  \label{eqn:DiscriminantDecomposition}
\end{equation}
as follows. 
Let $\mathbf A \colon U \times T \to \CC^{m\times n}$ be a 
representative of the unfolding and $t \in T \subset \CC^{k}$ a fixed parameter. 
Then $t \in \Delta_{\mathbf A}$ if and only if there are points $x \in U$ 
for which $A_t \colon (U,x) \to (\CC^{m\times n},A_t(x))$ is not transverse 
to the rank stratification. 
To any such point we can associate the rank $r$ of the stratum containing 
the critical value $A_t(x)$. Now $\Delta_{\mathbf A}^r$ is the component of 
$\Delta_{\mathbf A}$ whose \textit{generic} fiber $A_t$ has critical 
points of rank at most $r$. An example will be given below in \ref{exp:SurfaceDiscriminant}.

\begin{remark}
  \label{rem:DeterminantalDiscriminantStructure}
  It is easy to see that for either two miniversal unfoldings of $A$, the
  pullback maps in the parameter spaces take one matrix discriminant into the other. 
  However, since these pullbacks are not uniquely determined, the matrix discriminant 
  $(\Delta_A,0)$ of $A$ exists and is unique, but in general only up to non-unique isomorphism. 
\end{remark}

Sheafifying allows us to pass from $A = A_0$ to a nearby map 
$A_t \colon U \to \CC^{m\times n}$ defined on some suitable open subset $U \subset \CC^p$. 
Then $t\in \Delta_A$ lays in the matrix discriminant if and only if there 
are points $x \in U$ at which $A_t$ admits non-trivial unfoldings. 
Given Proposition \ref{prp:AlgebraicTransversality}, these are precisely 
the points at which $A_t$ is not transversal to the rank stratification of 
the target space $\CC^{m\times n}$. 
The Thom transversality theorem\footnote{See \cite{Thom54}} assures that any 
differentiable map $f \colon X \to Y$ can be deformed to a map $\tilde f$ 
which is transversal to a given submanifold $Z \subset Y$. In the holomorphic 
setting, this problem is more delicate due to the rigidity of 
holomorphic mappings. It has been addressed 
by Trivedi in \cite[Theorem 2.1]{Trivedi13} and \cite[Theorem 3.1]{Trivedi13}. 

For complex analytic germs such as $A \colon (\CC^p,0) \to (\CC^{m\times n},0)$ 
it also suffices to observe that the set of constant matrices 
$\varphi \in \CC^{m\times n}$ for which the map $A + \varphi$ is 
transverse to the rank stratification in a neighborhood of the origin, 
is dense in $\CC^{m\times n}$; cf. \cite[Lemma 2.2]{Trivedi13}.
As an immediate consequence we find:

\begin{corollary}
  \label{cor:DiscriminantsAreProperSubsets}
  Let $A$ be as in Definition 
  \ref{def:DeterminantalDiscriminant} and $\mathbf A$ a semi-universal 
  unfolding of $A$ (cf. Theorem \ref{thm:ExistenceSemiUniversalUnfolding}) on 
  $\tau = \dim_\CC T^1_\GL(A)$ parameters. Then the matrix 
  discriminant $(\Delta_A,0) \subset (\CC^\tau,0)$ is a \textit{proper} analytic subset of 
  $(\CC^\tau,0)$.
\end{corollary}
\noindent
For a more detailed discussion, 
the reader may also consult \cite[Section 2.2.1]{Zach17}. 

\medskip

In the following, we will refer to any map 
\begin{equation}
  A_t \colon U \to \CC^{m\times n}
  \label{eqn:Stabilization}
\end{equation}
which arises from an unfolding $\mathbf A$ of a finitely 
determined matrix $A$ and with $t \notin \Delta_{\mathbf A}$ as a 
(topological) \textit{stabilization} of $A$. 

\begin{remark}
  Using the results 
  by Trivedi \cite[Lemma 2.2]{Trivedi13}, it is easy to see that 
  stabilizations also exist for matrices $A$ which are not 
  finitely determined. In that case, however, a stabilization is 
  not uniquely determined by the original matrix $A$, similar 
  to the case of smoothings of non-isolated hypersurface singularities, 
  see Example \ref{exp:NonEquivalentSmoothingsForNonEIDS} below.
\end{remark}

\begin{example}
  \label{exp:DiscriminantThreeCoordinateAxis}
  Consider the miniversal unfolding (\ref{eqn:MiniversalUnfoldingThreeCoordinateAxis})
  of the matrix $A$ for the space 
  curve singularity given by the three coordinate axis in $(\CC^3,0)$ 
  from Examples \ref{exp:DeformationOfThreeCoordinateAxisFlatness} and 
  \ref{exp:ThreeCoordinateAxisSemiUniversalUnfolding}: 
  \[
    \mathbf A(x,y,z;t_1,t_2,t_3) = 
    \begin{pmatrix} 
      x & 0 & z \cr
      0 & y & z 
    \end{pmatrix} + 
    \begin{pmatrix}
      0 & t_1 & 0 \\
      t_2 & 0 & t_3
    \end{pmatrix}
  \]
  We claim that the $t_1$-axis in the parameter space $(\CC^3,0)$ is contained 
  in the discriminant $\Delta_A$. 

  Fix $t_2 = t_3 = 0$ and let $t_1 \neq 0$ be arbitrarily small. Then the 
  perturbed map 
  \[
    A_{(t_1,0,0)} \colon 0 \mapsto \varphi := A_{(t_1,0,0)}(0) = 
    \begin{pmatrix}
      0 & t_1, & 0 \\
      0 & 0 & 0
    \end{pmatrix}
  \]
  takes the origin to a matrix $\varphi$ of rank $1$. If we let 
  $y_{i,j}$ be the coordinates for the space of matrices $\CC^{2\times 3}$, then 
  locally at the point $\varphi$ the minor of the matrix $Y$ obtained 
  by deletion of the second column is 
  easily seen to be a superfluous generator of the ideal $\langle Y^{\wedge 2} \rangle$.
  The generic determinantal variety $M_{2,3}^2$ is a smooth local complete 
  intersection at $\varphi$ with tangent space 
  \[
    T_\varphi M_{2,3}^2 = \Span\left( 
    \begin{pmatrix}
      1 & 0 & 0 \\
      0 & 0 & 0 \\
    \end{pmatrix},\,
    \begin{pmatrix}
      0 & 1 & 0 \\
      0 & 0 & 0 \\
    \end{pmatrix},\,
    \begin{pmatrix}
      0 & 0 & 1 \\
      0 & 0 & 0 \\
    \end{pmatrix},\,
    \begin{pmatrix}
      0 & 0 & 0 \\
      0 & 1 & 0 \\
    \end{pmatrix}
    \right).
  \]
  But the image of the differential of $A_{(t_1,0,0)}$ at the origin is spanned by the 
  matrices 
  \[
    \D A_{(t_1,0,0)}(0)\,T_0\CC^3 = \Span\left( 
    \frac{\partial A}{\partial x},\,
    \frac{\partial A}{\partial y},\,
    \frac{\partial A}{\partial z}
    \right)
  \]
  which is clearly not transversal to $T_\varphi M_{2,3}^2$, since for $\lambda \neq 0$ none of the 
  matrices $\lambda \cdot E_{2,1}$ is contained in the sum of the two subspaces.

  This particular unfolding was already considered in 
  Example \ref{exp:DeformationOfThreeCoordinateAxisFlatness} and we saw that 
  the fibers $X_A^2(t_1,0,0)$ were not smooth. From the viewpoint taken in this example, 
  this is merely a consequence of the non-transversality of the map $A_{(t_1,0,0)}$ 
  at these points. This observation 
  will be generalized in Lemma \ref{lem:StratificationOfTheEssentialSmoothing}. 

  With the help of a computer algebra system it is easy to see that the full discriminant 
  $(\Delta_A,0) \subset (\CC^3,0)$ in the parameter space $(\CC^3,0)$ consists of the union 
  of the three coordinate \textit{hyperplanes} given by $t_1\cdot t_2 \cdot t_3 = 0$. 
  Hence, a simultaneous perturbation by $t_1 = t_2 = t_3 = t$ will lead to a smoothing 
  of $(X_A^2,0)$.
  \begin{figure}[h]
    \centering
    \includegraphics[scale=0.5,clip=true,trim=1cm 0.1cm 0.5cm 0.1cm]{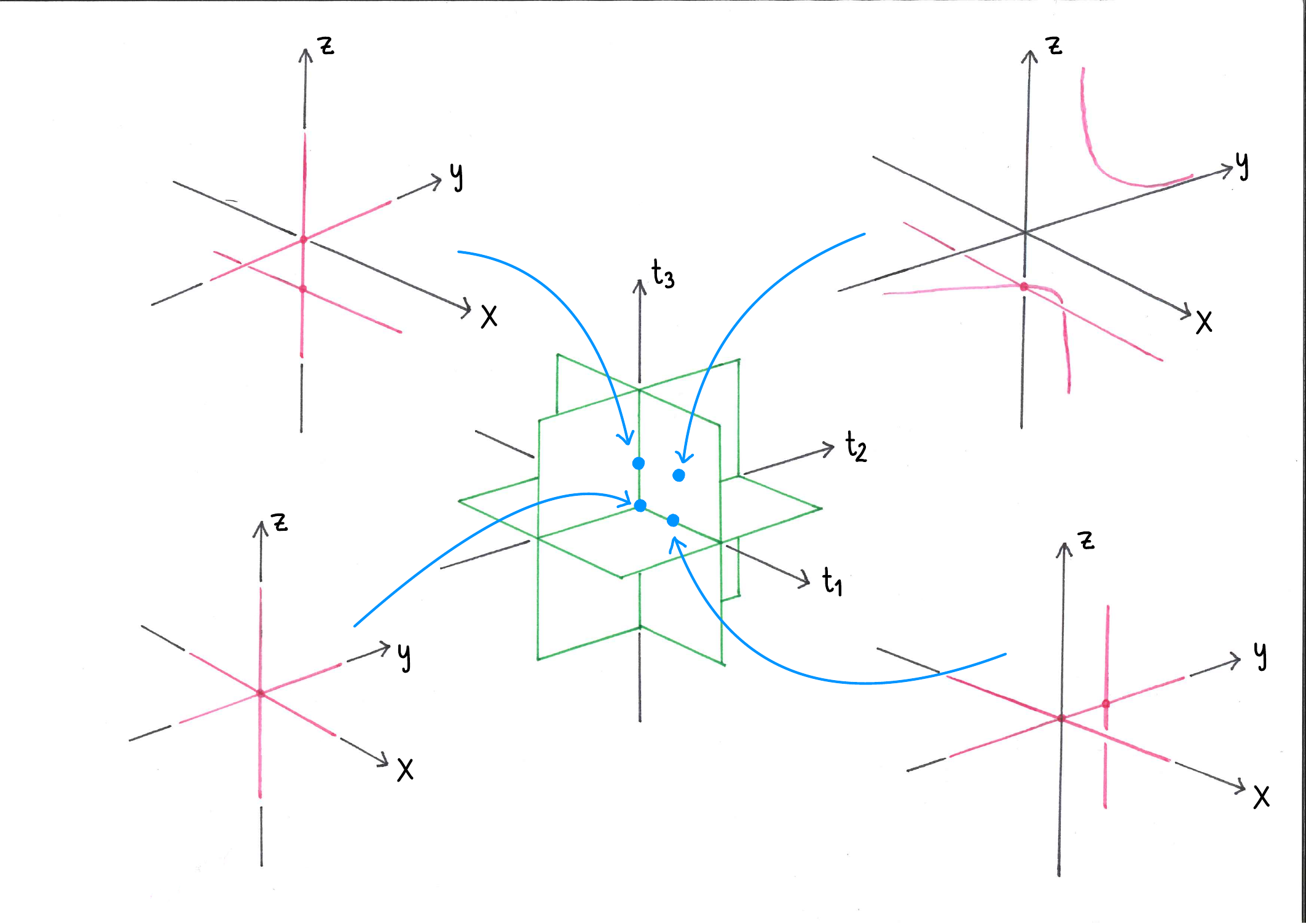}
    \caption{The discriminant of the space curve singularity in $\CC^3$ given by the union of 
    the three coordinate axis.}
    \label{fig:DiscriminantSpaceCurve}
  \end{figure}
\end{example}

\begin{example}
  \label{exp:SurfaceDiscriminant}
  The matrix discriminant is not always a hypersurface. For instance, the matrix
  \[
    A = 
    \begin{pmatrix}
      x & y & z \\
      y & z & w
    \end{pmatrix}
  \]
  has a miniversal unfolding given by 
  \[
    \mathbf A(x,y,z,w;t_0,t_1) = 
    \begin{pmatrix}
      x & y & z \\
      y & z & w
    \end{pmatrix} - 
    \begin{pmatrix}
      0 & 0 & 0 \\
      t_0 & t_1 & 0 \\
    \end{pmatrix}.
  \]
  Contrary to the previous example, the matrix discriminant $\Delta_A$ consists 
  only of the point $\{0\} \subset \CC^2$ in the parameter space. 

  This behaviour continues, as can be observed from the slightly more complicated 
  matrix 
  \[
    B = 
    \begin{pmatrix}
      x & y & z \\
      y^2 & z & w
    \end{pmatrix}
  \]
  This is a member of the second series from the list of simple isolated Cohen-Macaulay 
  codimension $2$ surface singularities in Table \ref{tab:FKNSurface}. 
  The miniversal unfolding can be realized on $3$ parameters as 
  \[
    \mathbf B(x,y,z,w;t_0,t_1,t_2) 
    = 
    \begin{pmatrix}
      x & y & z \\
      y^2 & z & w
    \end{pmatrix}
    + 
    \begin{pmatrix}
      0 & 0 & 0 \\
      t_0 + t_1 y & t_{2} & 0 
    \end{pmatrix}.
  \]
  The matrix discriminant decomposes as in (\ref{eqn:DiscriminantDecomposition}) into 
  components 
  $\Delta_B = \Delta_B^0 \cup \Delta_B^1$.
  The first one 
  \[
    \Delta_B^0 = \left\{ t_2 = t_0 = 0 \right\}
  \]
  is of codimension $2$. Geometrically it is characterized by the fact that 
  for a generic point $0 \neq t = (t_1,0,0) \in \Delta_B^0$ the map 
  $B_t \colon (\CC^4,0) \to (\CC^{2\times 3},0)$ $\GL$-equivalent to 
  the matrix $A$ above. In particular, it 
  is not transverse to the stratum $\{0\} = V_{2,3}^0$ of the 
  rank stratification so that the origin is a critical point of rank 
  $r=0$ for $B_t$. 
  
  The other component 
  \[
    \Delta_B^1 = \left\{ t_1^2 - 4t_0 \right\}
  \]
  is a divisor. For generic $t = \left( t_1^2/4, t_1, t_2 \right) \in \Delta_B^1$
  the map $B_t$ has a non-transverse point at $(x,y,z,w) = (0,-t_1/2,-t_2,0)$ 
  whose image 
  \[
    B_t(0,-t_1/2,t_2,0) = 
    \begin{pmatrix}
      0 & -t_1/2 & -t_2 \\
      0 & 0 & 0
    \end{pmatrix}
  \]
  is of rank $1$.
  Note that now the full discriminant $(\Delta_B,0)$ is of codimension $1$, but not a divisor.
\end{example}


\section{Essentially Isolated Determinantal Singularities and their deformations}
\label{sec:EIDSAndTheirDeformations}

In the preceding section we have discussed the singularities 
of matrices regarded as map germs. The geometric criteria for finite 
determinacy, Theorem \ref{thm:GeometricCriterionForFiniteDeterminacy},
and the existence of miniversal unfoldings, 
Proposition \ref{prp:AlgebraicTransversality},
motivated the definition of an essentially 
nonsingular point, Definition \ref{def:EssentiallyNonSingularPoint}, 
in a natural way. 
For the determinantal singularities 
associated to a matrix 
this leads to the following:

\begin{definition}(Ebeling, Gusein-Zade \cite{EbelingGuseinZade09})
  \label{def:EIDS}
  A determinantal singularity $(X_A^s,0) \subset (\CC^p,0)$ is called an 
  \textit{Essentially Isolated Determinantal Singularity} (EIDS) 
  of type $(m,n,s)$ 
  if the defining matrix $A \colon (\CC^p,0) \to (\CC^{m\times n},0)$
  has only essentially non-singular points in 
  a punctured neighborhood of the origin. 
\end{definition}

\begin{remark}
  \label{rem:EIDSAndFiniteDeterminacy}
  It follows directly from Proposition \ref{prp:AlgebraicTransversality} and 
  Theorem \ref{thm:GeometricCriterionForFiniteDeterminacy} that the 
  defining matrix of an EIDS is finitely $\GL$-determined. 
  Note that in general any such matrix gives rise to several EIDS:
  one for every 
  integer $s$ satisfying $0 \leq(m-s+1)(n-s+1) \leq p$ so that the intersection 
  of the image of $A$ with the variety $M_{m,n}^s$ 
  is generically nonempty.
\end{remark}

An EIDS is \textit{essentially isolated} in the sense that in general, 
the space $X_A^s$ has non-isolated singularities. However, due to 
the transversality condition in 
Definition \ref{def:EssentiallyNonSingularPoint} 
imposed on the map $A$, apart from the origin itself these singularities are 
locally products of the generic determinantal varieties with affine space: 

\begin{lemma}
  Let $A \in \CC\{x_1,\dots,x_p\}^{m\times n}$ be a matrix defining 
  an EIDS $(X_A^s,0)\subset (\CC^p,0)$. Then the preimages of the strata 
  $V^r_A = A^{-1}(V_{m,n}^r)$ with $r<s$
  form a Whitney stratification of $X_A^s\setminus\{0\}$. Moreover,
  at any point $x \in V_A^r\subset X_A^s\setminus\{0\}$ one has an isomorphism 
  \[
    (X_A^s,x) \cong (M_{m-r,n-r}^{s-r},0) \times (\CC^{p-(m-r)(n-r)},0).
  \]
  In particular, $X_A^s$ has isolated singularity 
  at the origin if and only if 
  \begin{equation}
    p\leq (m-s+2)(n-s+2) 
    \label{eqn:InequalityIsolatedSingularity}
  \end{equation}
  and it is smoothable by a determinantal deformation 
  if and only if this inequality is strict.
  \label{lem:WhitneyStratificationOfAnEIDSAndSmoothability}
\end{lemma}

\begin{proof}
  The first part of this lemma is a basic application of stratification theory, 
  see e.g. \cite{GoreskyMacPherson88}. For the question on smoothability 
  observe that in case $p=(m-s+2)(n-s+2)$, the space $(X_A^{s-1},0)$ is also 
  an EIDS, but of dimension zero. An unfolding of $A$ on a set of parameters $t$ 
  gives rise to a 
  (flat) deformation of this singularity and the principle of conservation of number 
  asserts that the total multiplicity of these points is preserved within 
  the family. It is now easy to see that for $t \neq 0$ any of these 
  points $x \in A_t^{-1}(M_{m,n}^{s-1}) \subset A_t^{-1}(M_{m,n}^s)$ 
  must be a singular point of the fiber $A_t^{-1}(M_{m,n}^s)$.
\end{proof}

\begin{figure}[h]
  \centering
  \includegraphics[page=1,scale=0.3,clip=true,trim=0cm 4cm 1cm 1cm]{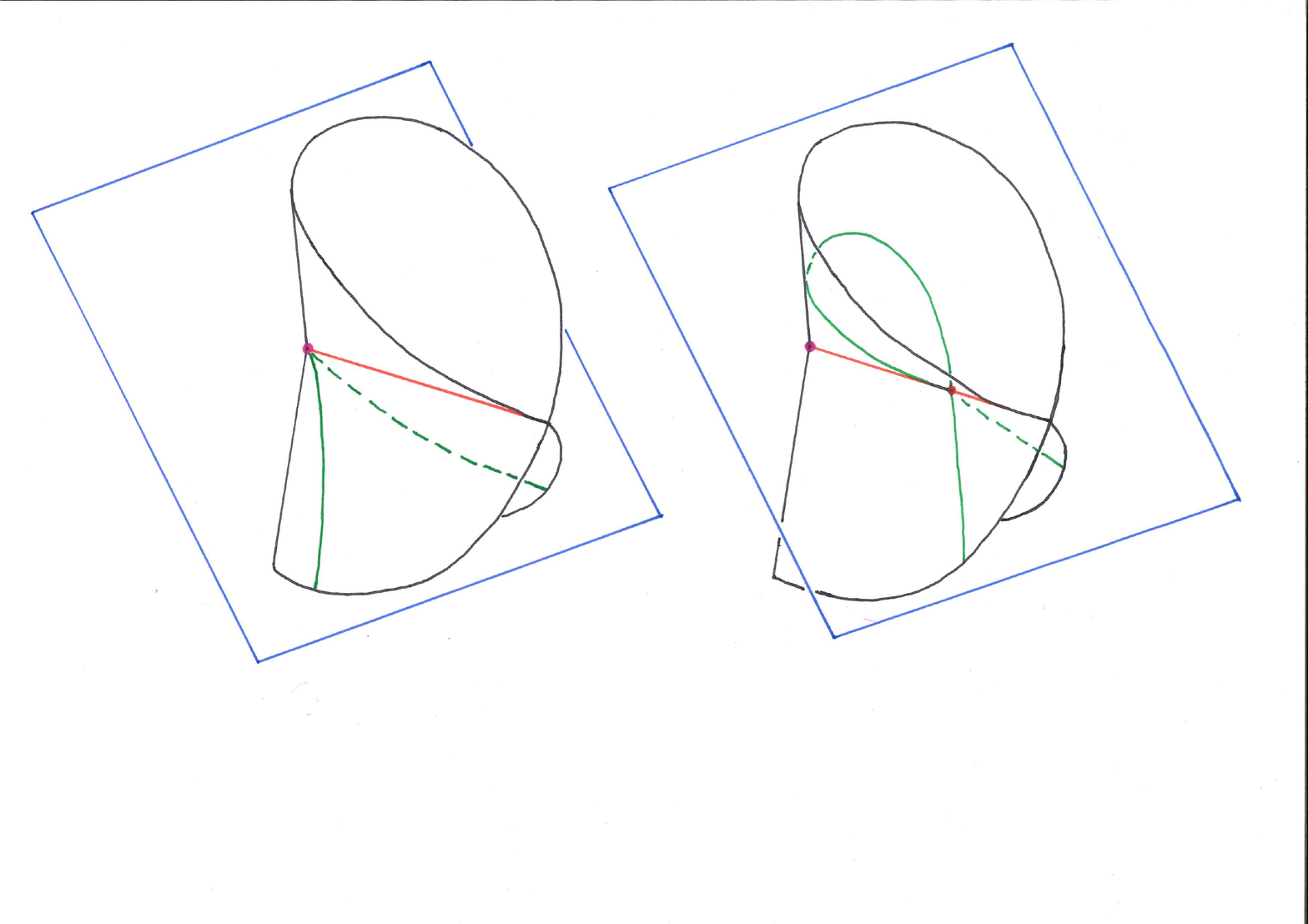}
  \caption{A stratified variety $M$ (black) with non-isolated singular locus (red) and its sections 
    (blue) via immersions of affine linear subspaces (green). On the left hand side, the immersion 
    is not transverse to $M$ in the stratified sense: The topology of the intersection changes when 
    perturbing the hyperplane slightly. On the right hand side, the map is transverse, but the 
    intersection of the plane with $M$ inherits a singular point in the nodal curve from the 
    singular locus of $M$.
  }
  \label{fig:WhitneyUmbrellaAndSections}
\end{figure}
Along the same lines, it is easy to prove:

\begin{lemma}
  \label{lem:StratificationOfTheEssentialSmoothing}
  In the setup of Lemma \ref{lem:WhitneyStratificationOfAnEIDSAndSmoothability} 
  let $A_t \colon U \to \CC^{m\times n}$ be a stabilization of $A$ defined on 
  some open neighborhood $U \subset \CC^p$ of the origin. Then the preimages 
  of the strata $V_A^r = A_t^{-1}(V_{m,n}^r)$ form a Whitney stratification of 
  the essential smoothing $M_A^s$.
\end{lemma}

\begin{remark}
  \label{rem:Laksov}
  In this note we confined ourselves mostly to the complex analytic setup. 
  However, the idea of essential smoothings works 
  equally well in the purely algebraic setting. This was discussed 
  by Laksov in \cite{Laksov75} where he develops the notions of transversality, 
  and determinantal deformations for affine determinantal schemes. 
\end{remark}

\subsection{The Tjurina transformation for EIDS}

We need to introduce one technical tool that has proved to be very useful 
in the study of determinantal singularities: The Tjurina transform and 
the Tjurina transform in family for arbitrary EIDS.
For the generic determinantal varieties $M_{m,n}^s \subset \CC^{m\times n}$ 
the Tjurina transform has already been introduced in 
(\ref{eqn:ResolutionsOfSingularitiesGenericDeterminantalVariety}). 
For an arbitrary EIDS we give the following definition.

\begin{definition}
  \label{def:TjurinaTransformForEIDS}
  Let $(X_A^s,0) \subset (\CC^p,0)$ be a determinantal singularity of type $s$ 
  given by a matrix $A \in \CC\{x\}^{m\times n}$. We define the \textit{Tjurina transform} 
  $\hat \nu \colon \hat X_A^s \to X_A^s$ of 
  $X_A^s$ to be the fiber product 
  \begin{equation}
    \xymatrix{ 
      X_A^s \times_{M_{m,n}^s} \hat M_{m,n}^s \ar[d]_{\hat \nu} \ar[r] &
      \hat M_{m,n}^s \ar[d]^{\hat \nu} \\
      X_A^s \ar[r]^A & 
      M_{m,n}^s
    }
    \label{eqn:DefinitionTjurinaTransform}
  \end{equation}
  where $\hat M_{m,n}^s$ is the Tjurina transform of the generic determinantal 
  variety (\ref{eqn:ResolutionsOfSingularitiesGenericDeterminantalVariety}).
  The \textit{strict Tjurina transform} $\overline{X}_A^s \subset \hat X_A^s$ is defined as 
  the closure 
  \begin{equation}
    \overline X_A^s = \overline{ (A \circ \hat \nu)^{-1}(V_{m,n}^{s-1})} \subset \hat X_A^s
    \label{eqn:DefinitionStrictTjurinaTransform}
  \end{equation}
  of the open set over the matrices of rank $s-1$.
\end{definition}

Using the properties already described for $\hat \nu \colon \hat M_{m,n}^s \to M_{m,n}^s$, 
it is easy to see that $\hat \nu \colon \hat X_A^s \to X_A^s$ is an isomorphism 
outside the singular locus $X_A^{s-1}$. 
In many practical cases with matrices of small size compared to the dimension of $(X_A^s,0)$, 
the Tjurina transform and the strict Tjurina transform coincide. For the general case, 
however, the Tjurina transform will have higher dimensional components in its 
exceptional set. While the strict Tjurina transform might be geometrically more 
intuitive, the definition in (\ref{eqn:DefinitionTjurinaTransform}) 
has the advantage that it provides explicit equations to work with. 

\begin{definition}
  \label{def:TjurinaTransformForEIDSInFamily}
  Let $(X_A^s,0) \hookrightarrow (\mathcal X_{\mathbf A}^s,0) \overset{\pi}{\longrightarrow} 
  (\CC^k,0)$ be a determinantal deformation of an EIDS induced from an unfolding 
  $\mathbf A(x,t)$ of the defining matrix $A$. 
  The \textit{Tjurina transformation in family} 
  \begin{equation}
    \xymatrix{
      \hat X^s_A \ar@{^{(}->}[r] \ar[d] & 
      \widehat{\mathcal X}^s_{\mathbf A} \ar[d]^{\hat \nu} \ar[r] &
      \hat M_{m,n}^s \ar[d]\\
      X_A^s \ar@{^{(}->}[r] \ar[d] & 
      \mathcal X_A^s \ar[d]^{\pi} \ar[r]^{\mathbf A} & 
      M_{m,n}^t\\
      \{0\} \ar@{^{(}->}[r] & 
      \CC^k& 
    }
    \label{eqn:DefinitionTjurinaTransformationInFamily}
  \end{equation}
  is obtained by applying the Tjurina transformation to the 
  total space $(\mathcal X_A^s,0)$ of the family.
\end{definition}

It is a priori not clear whether the family 
$\pi \circ \hat \nu \colon \widehat{\mathcal X}^s_{\mathbf A} \to \CC^k$ 
is well behaved (e.g. flat). In particular, the fibers of the family 
can not be expected to specialize to the strict Tjurina transform of $(X_A^s,0)$. 
However, using Lemma \ref{lem:StratificationOfTheEssentialSmoothing} it 
is not difficult to see the following:

\begin{proposition}
  If the family in Definition \ref{def:TjurinaTransformForEIDSInFamily} arises 
  from a stabilization of $A$, then for suitable representatives and 
  $t\notin \Delta_{\mathbf A}$ outside the discriminant of the deformation 
  the restriction to the fiber 
  \[
    \hat \nu \colon \hat X_{\mathbf A}^s(t) \to X_{\mathbf A}^s(t)
  \]
  is a resolution of singularities of $X_{\mathbf A}^s(t)$. In particular, this is an isomorphism 
  whenever $X_{\mathbf A}^s(t)$ is already smooth.
  \label{prp:PropertiesOfTheTjurinaTransformationInFamily}
\end{proposition}

\begin{remark}
  The same procedures can in principal also be applied to the dual 
  Tjurina transform and the Nash transform described in 
  (\ref{eqn:ResolutionsOfSingularitiesGenericDeterminantalVariety}). 
  This has, for instance, been done by Ebeling and Gusein-Zade in 
  \cite{EbelingGuseinZade09} in order to construct other resolutions
  of singularities for essential smoothings as in Proposition
  \ref{prp:PropertiesOfTheTjurinaTransformationInFamily}.
\end{remark}

\begin{example}
  \label{exp:TjurinaTransformationInFamilyThreeCoordinateAxis}
  Let again $(X_A^2,0) \subset(\CC^3,0)$ be the union of the three 
  coordinate axis and $\pi \colon (\mathcal X_{\mathbf A}^2,0) \to (\CC,0)$ 
  its smoothing induced from the unfolding 
  \[
    \mathbf A(x,y,z;t) = 
    \begin{pmatrix}
      x & 0 & z \\
      0 & y & z 
    \end{pmatrix}
    + 
    \begin{pmatrix}
      0 & t & 0 \\
      t & 0 & t 
    \end{pmatrix}
  \]
  as in the end of Example \ref{exp:DiscriminantThreeCoordinateAxis}.
  The strict Tjurina transform $\overline X_A^2 = L_x \dot \cup L_y \dot \cup L_z$ 
  consists of the three separated 
  coordinate axis only. The original Tjurina transform 
  $\hat X_A^2$, on the other hand, has as exceptional set $\hat E \cong \PP^1$ 
  a whole projective line as an additional component with the coordinate axis meeting 
  $\hat E$ at the points $0$, $\infty$, and $1$, respectively.

  For $t\neq 0$ the Tjurina transformation in family provides an identification 
  $\hat \nu \colon \hat X_{\mathbf A}^2(t) \overset{\cong}{\longrightarrow} X_{\mathbf A}^2(t)$ 
  of the smooth fibers. Observe that by construction, these fibers specialize 
  to the central fiber $\hat X_A^2 = \hat X_{\mathbf A}^2(0)$ but not to the strict 
  Tjurina transform of $(X_A^2,0)$. Moreover, one can easily verify from 
  explicit calculations, that the given family induces local, simultaneous smoothings 
  of the three $A_1$-singularities of $\hat X_A^2$ at the intersection points of 
  either $L_x$, $L_y$, and $L_z$ with $\hat E$. In particular, the family 
  $(\pi \circ \hat \nu ) \colon (\widehat{\mathcal X}_{\mathbf A}^2,0) \to (\CC,0)$ turns 
  out to be flat.

  For the \textit{dual} Tjurina transform we find that 
  \[
    \check X_A^2 = \PP^2 \cup L_x \cup L_y \cup L_z
  \]
  consists of an exceptional \textit{plane} $\check E \cong \PP^2$ 
  with the components of the strict transform meeting $\check E$ at the points 
  $(1:0:0)$, $(0:1:0)$, and $(0:0:1)$, respectively. Again, on the fibers 
  over $t \neq 0$, the projection $\check \nu \colon \check X_{\mathbf A}^2(t) \to X_{\mathbf A}^2(t)$ 
  is an isomorphism of curves with the fibers specializing to 
  $\check X_{\mathbf A}^2(0) = \check X_{A}^2$;
  only that this time the family can not be flat, given the jump in dimensions of the 
  fibers at $t=0$.
\end{example}

\subsection{Comparison of unfoldings and semi-universal deformations}
\label{sec:ComparisonOfUnfoldingsAndSemiUniversalDeformations}

Given a finitely $\GL$-determined matrix $A \in \CC\{x_1,\dots,x_p\}^{m\times n}$ 
together with its miniversal unfolding $\mathbf A$ on $k=\dim_\CC T^1_\GL(A)$ 
parameters as, for instance, in 
Theorem \ref{thm:ExistenceSemiUniversalUnfolding}, one can compare the 
unfoldings of $A$ with the deformations of the associated EIDS 
$(X_A^s,0)$ defined by $A$. To make this explicit, it is important that 
the germ $(X,0) = (X_A^s,0)$ given by the associated determinantal singularity has 
a \textit{semi-universal deformation}. An excellent overview for this topic 
with further references for details can be found in 
\cite[Chapter 7]{Volume1} and therefore we will restrict ourselves to 
briefly recalling the cornerstones of the theory.

Similar to the definition of versal unfoldings, 
a \textit{versal deformation} of a complex analytic germ $(X,0)$ 
is given by a flat family
\begin{equation}
  \xymatrix{
    (X,0) \ar@{^{(}->}[r] \ar[d] & 
    (\mathcal X,0) \ar[d]^{\pi} \\
    \{0\} \ar@{^{(}->}[r] &
    (S,0)
  }
  \label{eqn:SemiUniversalDeformationOfGrauert}
\end{equation}
such that any other deformation 
$(X,0) \hookrightarrow (\mathcal X',0) \to (S',0)$ can 
be written as a pullback from $(\mathcal X,0) \to (S,0)$ via some 
comparison map $\Phi \colon (S',0) \to (S,0)$.
A versal deformation is called \textit{semi-universal}, if 
the differential of the comparison map $\Phi$ is uniquely determined 
by the particular family $(\mathcal X',0) \to (S',0)$.
For a thourough discussion of these definitions we 
refer to \cite[Definition 7.1.13]{Volume1}.

This second condition on the differential of the comparison map is closely 
related to the space of \textit{first order deformations} $T^1_{X,0}$, 
see \cite[Section 7.1.4]{Volume1}, which 
is central to the construction of semi-universal deformations.
It can be defined as the set of isomorphism classes of deformations 
of $(X,0)$ over the formal algebra $\CC[t]/\langle t^2\rangle$, 
cf. \cite[Definition 7.1.26]{Volume1} and it is 
the Zariski tangent space of the base $(S,0)$ of a 
semi-universal deformation of $(X,0)$ --- provided such a semi-universal 
deformation exists. One can explicitly compute 
the space $T^1_{X,0}$ from the ideal $I \in \CC\{x\}$ for any given embedding  
$(X,0) \subset (\CC^p,0)$, cf. \cite[Proposition 7.1.33]{Volume1}: It 
appears in an exact sequence 
\begin{equation}
  \xymatrix{
    T_{\CC^p,0} \otimes_{\CC\{x\}} \CC\{x\}/I \ar[r]^\beta & 
    \Hom_{\CC\{x\}}(I,\CC\{x\}/I) \ar[r] & 
    T^1_{X,0} \ar[r] & 
    0
  }
  \label{eqn:DefiningSequenceT1}
\end{equation}
where $T_{\CC^p,0}$ is the module of germs of holomorphic vector fields on 
$(\CC^p,0)$, $\Hom_{\CC\{x\}}(I,\CC\{x\}/I)$ is the \textit{normal module} 
for the embedding of $(X,0)$, and $\beta$ is the map determined by 
\[
  \beta \colon 
  \frac{\partial}{\partial x_i} \mapsto \left( f \mapsto \frac{\partial f}{\partial x_i} \right)
\]
for $i = 1,\dots,p$. This endows $T^1_{X,0}$ with the structure of a 
$\CC\{x\}$-module and in particular, it is a $\CC$-vector space. 

The following theorem was proved by Grauert in \cite{Grauert72}, see 
also \cite[Theorem 7.1.14]{Volume1} and \cite[Remark 7.1.28]{Volume1}.

\begin{theorem}
  Let $(X,0)$ be a complex analytic singularity with $T^1_{X,0}$ of finite dimension. 
  Then a semi-universal deformation of $(X,0)$ exists.
  \label{thm:GrauertsTheorem}
\end{theorem}

The construction of such semi-universal deformations 
is quite technical and in general much more complicated than 
the construction of miniversal unfoldings for map germs as for instance 
in Theorem \ref{thm:ExistenceSemiUniversalUnfolding}. 
The general framework for this was developed by 
Schlessinger in his thesis \cite{Schlessinger68}. 
For details of the general case, we refer to 
\cite[Chapter 7]{Volume1}, ``Deformation and Smoothing of Singularities'' 
by Greuel.
For the specific purposes here, 
we only note the following, \cite[Definition 7.1.38]{Volume1}:

\begin{definition}
  \label{def:Unobstructed}
  A singularity $(X,0)$ with finite dimensional $T^1_{X,0}$ is called 
  \textit{unobstructed} if it has a semi-universal deformation with a 
  smooth base $(S,0)$.
\end{definition}

A criterion for a singularity $(X,0)$ to be unobstructed is that 
the module $T^2_{X,0}$ is zero. This is another coherent analytic 
module associated to the singularity that can be computed explicitly 
from the ideal defining the singularity. Again, we refer to 
\cite[Chapter 7]{Volume1} for details.

\medskip

Suppose that $(X,0) = (X_A^s,0)$ is an EIDS with 
isolated singularity (cf. Lemma \ref{lem:WhitneyStratificationOfAnEIDSAndSmoothability}).
Since every unfolding of $A$ gives rise to a deformation of 
$(X_A^s,0)$ due to Lemma \ref{lem:UnfoldingInducesDeformation}, we obtain a flat family 
over the parameter space of the unfolding 
and some comparison map
  $\Phi \colon (\CC^k,0) \to (S,0)$
to the base of the semi-universal deformation $\pi \colon (\mathcal X,0) \to (S,0)$ 
of $(X_A^s,0)$ in 
the sense of Grauert's theorem. Altogether this forms the commutative diagram 
\begin{equation}
  \xymatrix{
    (X_A^s,0) \ar@{^{(}->}[r] \ar[d] & 
    (\mathcal X_{\mathbf A}^s,0) \ar[d]^{\rho} \ar[r] & 
    (\mathcal X,0 ) \ar[d]^\pi\\
    \{0\} \ar@{^{(}->}[r] &
    (\CC^k,0)\ar[r]^\Phi & 
    (S,0).
  }
  \label{eqn:DeterminantalDeformationsVsSemiUniversalDeformations}
\end{equation}
Depending 
on the size of the matrix $A$ and the size of the minors $s$, the  
map $\Phi$ can take very different forms. We first give a list of 
particular examples and then discuss some cases where more structural 
results are known.

\begin{example}
  \label{exp:VersalDeterminantalDeformations}
  \begin{enumerate}
    \item Let $(X_A^2,0) \subset (\CC^3,0)$ be the determinantal hypersurface 
      singularity defined by the matrix 
      \[
	A = 
	\begin{pmatrix}
	  x & y \\ z & x
	\end{pmatrix}.
      \]
      The ideal $I$ of $(X_A^2,0)$ is thus generated by the equation $f = x^2-yz$ and 
      we recognize the well-known $A_1$-surface singularity.
      A basis of $T^1_\GL(A)$ is given by 
      \[
	\begin{pmatrix}
	  1 & 0 \\ 0 & -1
	\end{pmatrix}
      \]
      and hence if we let $t$ be the deformation parameter in the semi-universal 
      unfolding of $A$, then 
      the induced deformation of the space germ $(X_A^2,0)$ comes from 
      a perturbation of $f$ by $-t^2$. 

      The semi-universal deformation of $(X_A^2,0)$ as a space germ on the other hand 
      is given by the perturbation of $f$ by a constant $u$. It follows that 
      the comparison map (\ref{eqn:DeterminantalDeformationsVsSemiUniversalDeformations}) takes 
      the form 
      \[
	\Phi : (\CC,0) \to (\CC,0), \quad t \mapsto u = t^2.
      \]
      In other words: The base of the miniversal unfolding of $A$ is a 
      $2:1$ cover of the base of the semi-universal deformation of $(X_0,0)$. 

    \item (Pinkham, \cite{Pinkham73}) 
      Recall from Example \ref{exp:VariousDeterminantalStructures} 
      that the ideal $I \subset \CC\{x_0,\dots,x_4\}$ in Pinkham's example \cite{Pinkham73}
      was given as $I = \langle A^{\wedge 2} \rangle = \langle B^{\wedge 2}\rangle$ 
      for the two matrices 
      \[
	A = 
	\begin{pmatrix}
	  x_0 & x_1 & x_2 & x_3 \\
	  x_1 & x_2 & x_3 & x_4
	\end{pmatrix}
	\quad\textnormal{and}\quad
	B = 
	\begin{pmatrix}
	  x_0 & x_1 & x_2 \\
	  x_1 & x_2 & x_3 \\
	  x_2 & x_3 & x_4
	\end{pmatrix}.
      \]
      The germ $(X,0) \subset (\CC^5,0)$ defined by $I$ is the cone over the 
      rational normal curve of degree $4$ in $\PP^4$. 
      A direct computation of $T^1_\GL(A)$ and application of Theorem 
      \ref{thm:ExistenceSemiUniversalUnfolding} 
      yields that the unfolding of $A$ given by 
      \begin{equation}
	\mathbf A(x,t) = A(x) - 
	\begin{pmatrix}
	  0 & 0 & 0 & 0 \\
	  t_1 & t_2 & t_3 & 0
	\end{pmatrix}
	\label{eqn:PinkhamsExampleFirstMiniversalUnfolding}
      \end{equation}
      on the parameters $t = (t_1,t_2,t_3)$ is miniversal. 
      Similarly, one has a miniversal unfolding on a single parameter $u$ 
      for the symmetric matrix $B$ given by 
      \begin{equation}
	\mathbf B(x,u) = B(x) - 
	\begin{pmatrix}
	  0 & 0 & u \\
	  0 & 0 & 0 \\
	  u & 0 & 0
	\end{pmatrix}.
	\label{eqn:PinkhamsExampleSecondMiniversalUnfolding}
      \end{equation}
      In \cite{Pinkham73} Pinkham shows explicitly that the 
      semi-universal deformation of $(X,0)$ as a complex analytic germ has a base 
      $(S,0) \subset (\CC^4,0)$ of the following form. Let $t_1,t_2,t_3,u$ be 
      the coordinates of $\CC^4$. Then $(S,0)$ consists of two components: 
      The plane $H = \{u = 0 \}$ and the line $L = \{t_1= t_2=t_3=0\}$. 
      An explicit computation of $(S,0)$ can also be found in \cite[Example 7.1.41]{Volume1}.

      Indeed, the comparison map for the miniversal unfolding of $A$ identifies 
      the parameter space 
      \[
	\Phi_A \colon (\CC^3,0) \overset{\cong}{\longrightarrow} (H,0) \subset(S,0), 
	\quad t \mapsto (t,0)
      \]
      of $\mathbf A$ with the hyperplane $(H,0)$ in $(S,0)$ and similarly for $B$ one has an 
      isomorphism 
      \[
	\Phi_B \colon (\CC,0) \overset{\cong}{\longrightarrow} (L,0) \subset (S,0), 
	\quad u \mapsto (0,u).
      \]
    \item Consider the $A_1$ threefold singularity in $(\CC^4,0)$ as a determinantal 
      singularity of type $(2,2,2)$ via the matrix 
      \[
	A = 
	\begin{pmatrix}
	  x & y \\
	  z & w 
	\end{pmatrix}.
      \]
      There are no nontrivial unfoldings of this matrix. 
      For the space germ on the other hand we find that the perturbation of 
      $f = \det A$ by a constant is semi-universal so that we have $T^1_{X_0,0} \cong \CC$. 
      Therefore, the comparison map takes the form 
      \[
	\Phi : \{ \textnormal{pt} \} \to (\CC,0).
      \]
      
    \item This example is taken from Schaps \cite{Schaps83} and it also appears 
      in \cite{Buchweitz81}. Explicit computations can be found in \cite{FruehbisKrueger18}. 

      Let $(X_A^2,0) \subset (\CC^4,0)$ be the union of the four coordinate axis. 
      This is a determinantal singularity via any matrix 
      \[
	\begin{pmatrix}
	  x_1 &
	  \alpha \cdot x_2 &
	  \beta \cdot x_3 &
	  \gamma \cdot x_4 \\
	  0 &
	  x_2 &
	  x_3 &
	  x_4 &
	\end{pmatrix}
      \]
      for general values $\alpha, \beta, \gamma \in \CC$. 
      Using row and column operations and local coordinate changes, one can always 
      bring this matrix to the form 
      \[
	A:=
	\begin{pmatrix}
	  x_1 & 0 & x_3 & \gamma' \cdot x_4 \\
	  0 & x_2 & x_3 & x_4
	\end{pmatrix}
      \]
      with $\gamma' \notin \{0,1\}$.
      One can show that the following matrices give a $\CC$-basis of $T^1_\GL(A)$: 
      \begin{eqnarray*}
	&
	\begin{pmatrix}
	  0 & 0 & 0 & 0 \\
	  1 & 0 & 0 & 0
	\end{pmatrix}, \quad
	\begin{pmatrix}
	  0 & 1 & 0 & 0 \\
	  0 & 0 & 0 & 0
	\end{pmatrix}, \quad
	\begin{pmatrix}
	  0 & 0 & 1 & 0 \\
	  0 & 0 & 0 & 0
	\end{pmatrix}, &
	\\
	&
	\begin{pmatrix}
	  0 & 0 & 0 & 1 \\
	  0 & 0 & 0 & 0
	\end{pmatrix}, \quad
	\begin{pmatrix}
	  0 & 0 & 0 & x_4 \\
	  0 & 0 & 0 & 0
	\end{pmatrix}.& 
      \end{eqnarray*}
      Hence, the base of the miniversal unfolding of $A$ is $(\CC^5,0)$. 
      Let $t_1,\dots,t_5$ be the parameters associated to these matrices 
      as in Theorem \ref{thm:ExistenceSemiUniversalUnfolding}.

      Computations of Rim\footnote{The computations are attributed to Rim in 
	\cite{Schaps83} without further reference} 
      and independently of Buchweitz \cite{Buchweitz81}
      have shown that the base $(S,0)$ 
      of the semi-universal deformation 
      of $(X_A^2,0)$ is isomorphic to the cone of the Segre embedding of 
      $\PP^1 \times \PP^3$ into $\PP^7$ and thus also of dimension $5$. 
      Consider the comparison map 
      \[
	\Phi : (\CC^5,0) \to (S,0).
      \]
      It is easy to see that the perturbation by $t_5$ alone 
      does not change the ideal generated by the $2$-minors of $A$ in $\CC\{x\}$: 
      This is a non-trivial deformation of the map germ $A$ which induces 
      a trivial deformation of the underlying space germ!
      Accordingly, as the computations by Fr\"uhbis-Kr\"uger show,
      $\Phi$ is a contraction of the $t_5$-axis
      but a local diffeomorphism away from it; just as if $\Phi$ was 
      a local chart in a resolution of singularities for $(S,0)$.
  \end{enumerate}
\end{example}

\subsection{Complete intersections}
\label{sec:DeformationsOfICIS}

It has been pointed out 
earlier that any (isolated) complete intersection singularity 
$(X,0) \subset (\CC^p,0)$ of codimension $c$ is determinantal 
of type $1$ for some $1\times c$-matrix $F = (f_1,\dots,f_c)$. 
For this class of singularities, the deformations of $(X,0)$ 
coincide with the unfoldings of $F$ up to $\GL$-equivalence. 

This starts with an explicit identification of the space of 
first order deformations $T^1_{X,0}$ of $(X,0)$ as in 
(\ref{eqn:DefiningSequenceT1}) and the infinitesimal unfoldings 
$T^1_\GL(F)$ of $F$ from (\ref{eqn:MatrixT1}). 

\begin{lemma}
  Let $(X,0) \subset (\CC^p,0)$ be a complete intersection singularity 
  defined by a regular sequence $F = (f_1,\dots,f_c)$ in $\CC\{x\}$. 
  Then there is an explicit isomorphism $T^1_\GL(F) \cong T^1_{X,0}$ 
  of the infinitesimal unfoldings of $F$ considered as a $1\times c$-matrix
  and the first order deformations of $(X,0)$.
  \label{lem:IsomorphismOfFirstOrderDeformationsForICIS}
\end{lemma}

The construction of this isomorphism builts on the description 
of $T^1_{X,0}$ for complete intersections as in \cite[Remark 7.1.35]{Volume1}.

\begin{proof}
  For a complete intersection ideal $I = \langle f_1,\dots,f_c\rangle$ 
  in $\CC\{x_1,\dots,x_p\}$ the resolution of $\CC\{x\}/I$ by the 
  Koszul complex (\ref{eqn:SpecializedKoszulComplex}) 
  can be used to show that the normal module 
  $\Hom_{\CC\{x\}}(I,\CC\{x\}/I)$ is a free $\CC\{x\}/I$-module 
  in generators $e_1,\dots,e_c$ which are dual to the $f_i$'s.
  An element $g = g_1 \cdot e_1 + \dots + g_c \cdot e_c 
  \in \Hom_{\CC\{x\}}(I,\CC\{x\}/I)$ corresponds to a formal 
  deformation of $(X,0)$ over $\CC[t]/\langle t^2\rangle$ given 
  by the $c$ equations 
  \[
    \mathbf F(x,t) = F(x) + t \cdot 
    \begin{pmatrix}
      g_1(x) & \dots & g_c(x)
    \end{pmatrix} = 0
  \]
  in the ring $\CC\{x\}[t]/\langle t^2\rangle$.
  The obvious translation to 
  unfoldings yields an isomorphism 
  \begin{eqnarray*}
    \Hom_{\CC\{x\}}(I,\CC\{x\}/I) &\overset{\cong}{\longrightarrow}&
    \CC\{x\}^{1\times c}/ \langle \CC\{x\}^{1\times 1} \cdot F + F \cdot \CC\{x\}^{c\times c} \rangle, \\
    \sum_{i=1}^c g_i \cdot e_i 
    & \mapsto & 
    \begin{pmatrix}
      g_1 & \dots & g_c
    \end{pmatrix} 
  \end{eqnarray*}
  where the right hand side is a quotient of $\CC\{x\}^{1\times c}$ by 
  a submodule of the extended tangent space 
  (\ref{eqn:ExtendedTangentSpaceForGLEquivalence}) of $F$ up to $\GL$-equivalence. 
  It is now easy to see that the remaining relations to be added, namely 
  \[
    \left\langle
    \frac{\partial F}{\partial x_1},\dots,\frac{\partial F}{\partial x_p} 
    \right\rangle
  \]
  coincide with the image of the map $\beta$ in (\ref{eqn:DefiningSequenceT1}). 
\end{proof}

Once the space of first order deformations has been computed, 
the construction of a semi-universal deformation of an isolated 
complete intersection singularity $(X,0) \subset (\CC^p,0)$ 
is straightforward, see \cite{Tjurina69}, \cite{KasSchlessinger72}; 
cf. also \cite[Theorem 7.1.22]{Volume1}. First 
note that for a complete intersection singularity $T^1_{X,0}$ 
is finite dimensional if and only if $(X,0)$ has isolated 
singularity.
Now according to the above cited theorems, 
any set of elements $G_1,\dots,G_\tau \in \CC\{x\}^{1\times c}$ reducing to a 
$\CC$-basis of $T^1_{X,0}$ in the above description gives rise to a semi-universal 
deformation $(X,0) \hookrightarrow (\mathcal X,0) \overset{\pi}{\longrightarrow} 
(\CC^\tau,0)$
where the germ $(\mathcal X,0)\subset (\CC^p,0) \times (\CC^\tau,0)$ is 
defined by the equations 
\[
  \mathbf F(x,t) = F + t_1\cdot G_1 + \dots + t_\tau \cdot G_\tau = 0. 
\]
This turns out to be the same procedure as for the construction of the 
miniversal unfolding of $F \in \CC\{x\}^{1\times c}$ described in 
Theorem \ref{thm:ExistenceSemiUniversalUnfolding}. 

\begin{corollary}
  \label{cor:SemiUniversalDeformationsAndUnfoldingsCoincide}
  For an isolated complete intersection singularity $(X,0) \subset (\CC^p,0)$ 
  defined by a regular sequence $F = (f_1,\dots,f_c)$ in $\CC\{x\}$ 
  the flat family (\ref{eqn:DeterminantalDeformationsVsSemiUniversalDeformations}) induced from 
  a miniversal unfolding of $F \in \CC\{x\}^{1\times c}$ is a 
  semi-universal deformation of $(X,0)$.
\end{corollary}
In particular, the comparison map $\Phi$ in 
(\ref{eqn:DeterminantalDeformationsVsSemiUniversalDeformations}) is an isomorphism.

\subsection{Cohen-Macaulay codimension $2$ singularities}

\label{sec:DeformationsOfICMC2Singularities}

The particular interest in Cohen-Macaulay singularities of 
codimension $2$ stems from the celebrated Hilbert-Burch theorem, 
see \cite{Hilbert90} and \cite{Burch68}, or \cite{Eisenbud95} for 
a modern textbook account:

\begin{theorem}[Hilbert-Burch]
  Let $I \subset \CC\{x_1,\dots,x_p\}$ be an ideal of codimension $2$ such that $\CC\{x\}/I$ is 
  Cohen-Macaulay. Then the minimal resolution of $\CC\{x\}/I$ as a $\CC\{x\}$-module takes the 
  form 
  \begin{equation}
    \xymatrix{
      0 \ar[r] &
      \CC\{x\}^{m} \ar[r]^A & 
      \CC\{x\}^{m+1} \ar[r]^f & 
      \CC\{x\} \ar[r] & 
      \CC\{x\}/I \ar[r] & 
      0 
    }
    \label{eqn:ResolutionCMC2}
  \end{equation}
  for some matrix $A \in \CC\{x\}^{m \times (m+1)}$ and 
  $I = \langle A^{\wedge m} \rangle$
  as ideals in $\CC\{x\}$. \\
  Conversely, suppose $A\in \CC\{x\}^{m\times (m+1)}$ is 
  any matrix such that the ideal 
  $I := \langle A^{\wedge m} \rangle$ has codimension $2$.  
  If we let 
  \[
    f = 
    \begin{pmatrix}
      \delta_1 & \cdots & \delta_{m+1}
    \end{pmatrix},
  \]
  where $\delta_i$ is $(-1)^i$ times the determinant of $A$ after deleting 
  the $i$-th row, then
  (\ref{eqn:ResolutionCMC2}) gives a minimal free resolution of $\CC\{x\}/I$.
  \label{thm:HilbertBurch}
\end{theorem}

In other words: Any Cohen-Macaulay singularity $(X,0) \subset (\CC^p,0)$ of 
codimension $2$ is determinantal for some $m\times (m+1)$-matrix in a canonical 
way. This differs drastically from the general case where one needs to specify 
the matrix in order to turn $(X,0)$ into a determinantal singularity; 
cf. Pinkham's example in Example \ref{exp:VersalDeterminantalDeformations}.

\medskip

Schaps has observed that this theorem can also be very well applied in the 
context of deformations. In \cite{Schaps77} and \cite{Schaps83} she pursues Schlessinger's approach 
to deformation theory \cite{Schlessinger68} for affine algebraic determinantal schemes, in particular 
those which are Cohen-Macaulay of codimension $2$. Her results can easily be adapted 
to the case of complex analytic singularities and, rephrasing them accordingly, 
she establishes the following, cf. \cite[Corollary 1]{Schaps77}: 

\begin{proposition}
  Let $(X,0) \subset (\CC^p,0)$ be Cohen-Macaulay of codimension $2$ endowed 
  with its canonical determinantal structure for some matrix $A \in \CC\{x\}^{m\times (m+1)}$. 
  Then a family $(X,0) \hookrightarrow (\mathcal X,0) \overset{\pi}{\longrightarrow} (\Spec B,0)$ 
  over some Artinian ring $B$ with special fiber $(X,0)$ is flat, if and only if 
  there exists a matrix $\mathbf A \in (\CC\{x\}\otimes_\CC B)^{m\times (m+1)}$ 
  such that the germ $(\mathcal X,0)$ is determinantal with matrix $\mathbf A$.
  \label{prp:SchapsResultForCMC2Singularities}
\end{proposition}

In other words: Any formal, infinitesimal deformation of a Cohen-Macaulay 
singularity of codimension $2$ is determinantal for its canonical determinantal 
structure. In particular, this holds for the first order deformations which 
leads to an explicit description of the space $T^1_{X,0}$ in (\ref{eqn:DefiningSequenceT1}) 
in its ``matrix form'':

\begin{corollary}
  For $(X,0) \subset (\CC^p,0)$ as in Proposition \ref{prp:SchapsResultForCMC2Singularities}
  one has a canonical isomorphism 
  \[
    T^1_\GL(A) \cong T^1_{X,0}.
  \]
  \label{cor:FruehbisKruegerMatrixT1}
\end{corollary}
 
\begin{proof}
  This was explicitly carried out by the first named author in \cite[Lemma 2.6]{FruehbisKrueger99} 
  and \cite[Lemma 2.7]{FruehbisKrueger99}.
\end{proof}

The deformation theory of Cohen-Macaulay codimension $2$ singularities is 
unobstructed, cf. Definition \ref{def:Unobstructed}, 
so there exists a semi-universal deformation 
for every such $(X,0)\subset (\CC^p,0)$ 
with $\dim T^1_{X,0} = \tau <\infty$ over a smooth base $(\CC^\tau,0)$, 
cf. \cite[Proposition 7.1.37]{Volume1}. Again, 
this semi-universal deformation can be derived from any 
$\CC$-basis of $T^1_{X,0}$ in the same way as for isolated complete 
intersection singularities in the previous section. 
Using the explicit identification from Corollary \ref{cor:FruehbisKruegerMatrixT1} 
we find:

\begin{corollary}
  Let $(X,0) \subset (\CC^p,0)$ be an isolated Cohen-Macaulay codimension $2$ singularity 
  with its canonical determinantal structure for a matrix $A \in \CC\{x\}^{m\times (m+1)}$. 
  Then the flat family (\ref{eqn:DeterminantalDeformationsVsSemiUniversalDeformations}) induced from 
  a miniversal unfolding of $A$ is a semi-universal deformation of $(X,0)$.
  \label{cor:SemiUniversalDeformationOfICMC2}
\end{corollary}
Again, the comparison map $\Phi$ in 
(\ref{eqn:DeterminantalDeformationsVsSemiUniversalDeformations}) 
can be chosen to be an isomorphism.

\subsection{Gorenstein singularities in codimension $3$} 
\label{sec:GorensteinSingularitiesInCodimension3}

Similar to the case of Cohen-Macaulay codimension $2$ singularities, 
Gorenstein singularities of codimension $3$ are equipped with a canonical 
Pfaffian structure. This was established by Buchsbaum and Eisenbud 
in \cite[Theorem 2.1]{BuchsbaumEisenbud77}: 

\begin{theorem}
  \label{thm:BuchsbaumEisenbudGorensteinInCodimension3}
  Let $(R,\m)$ be a Noetherian local ring.
  Let $n>0$ be an integer and $A \in R^{(2n+1)\times (2n+1)}$ a matrix 
  with entries in $\m$. Suppose that the ideal $I = \langle A^{\wedge n}_\sk \rangle$
  has grade $3$. Then $R/I$ is Gorenstein and $I$ is minimally generated by 
  $2n+1$ elements. \\
  Conversely, every ideal $I$ of grade $3$ with $R/I$ Gorenstein arises this way.
  %
\end{theorem}
They furthermore show that the free resolution of an ideal $I$ as 
in Theorem \ref{thm:BuchsbaumEisenbudGorensteinInCodimension3}
is given by the complex 
\begin{equation}
  \xymatrix{
    0 \ar[r] &
    R \ar[r]^{f^T} & 
    R^{2n+1} \ar[r]^A & 
    R^{2n+1} \ar[r]^f & 
    R \ar[r] & 
    R/I \ar[r] & 
    0
  }
  \label{eqn:ResolutionOfGorensteinIdealsOfCodimension3}
\end{equation}
where $f = (f_1,\dots,f_{2n+1})$ is the $1\times (2n+1)$-matrix 
with the generators of $\langle A^{\wedge n}_\sk\rangle$ as entries.
One may therefore follow the same arguments as in the Cohen-Macaulay
codimension $2$ case (cf. also the proof of Lemma 
\ref{lem:UnfoldingInducesDeformation}) to show that every deformation 
of a Gorenstein singularity of codimension $3$ is determinantal. 


Waldi has used the results by Buchsbaum and Eisenbud to show that 
the deformation theory of Gorenstein ideals of codimension $3$ 
is unobstructed, see \cite[Satz 1]{Waldi79}. 
Whenever such an ideal defines an isolated singularity, 
Waldi therefore finds \cite[Satz 2]{Waldi79}:
\begin{theorem}
  \label{thm:WaldiSmoothabilityForGorensteinCase}
  Let $(X,0) \subset (k^p,0)$ be an algebraic Gorenstein singularity 
  over an algebraically closed field $k$ which is analytically irreducible 
  and of codimension $3$ with $p\leq 9$. Then the semi-universal deformation of 
  $(X,0)$ has a smooth base and the generic fiber of this deformation is 
  also smooth. 
\end{theorem}
However, we should note that in this setting -- at least to the author's knowledge -- 
the interplay of $T^1_\GL(A)$ and $T^1_{X,0}$ has not yet been investigated. 
In particular, no classification of simple Gorenstein singularities in 
codimension $3$ has been done as of this writing.
\subsection{Rational surface singularities}
\label{sec:RationalSurfaceSingularities}

Determinantal deformations also appear in the study of rational surface singularities. 
Recall that a singularity $(X,0) \subset (\CC^p,0)$ of dimension $d\geq 2$ 
is called \textit{rational} if there exists a \textit{resolution of singularities} 
\[
  \rho \colon (Z,E) \to (X,0)
\]
such that one has 
\begin{equation}
  R^i\rho_* \OO_Z = 
  \begin{cases}
    \OO_X & \textnormal{ if } i = 0 \\
    0 & \textnormal{ otherwise }
  \end{cases}
  \label{eqn:DefinitionRationality}
\end{equation}
for the higher direct images of the structure sheaf $\OO_Z$. Here, $Z$ is smooth, 
$E = \rho^{-1}(X_{\mathrm{sing}})$ is the preimage of the singular locus of $X$ and 
$\rho$ is an isomorphism outside $E$. We refer to \cite[Chapter 3]{Volume1} 
for a discussion of resolutions of singularities. 

\medskip

For surface singularities, it is customary to require 
$\rho \colon (Z,E) \to (X,0)$ to be a \textit{good resolution}.
This means that the set $E = \bigcup_i E_i$ 
is a \textit{simple normal crossing divisor} with a decomposition into smooth, irreducible 
components $E_i$. To each component, one can assign two integers, the 
\textit{genus} of $E_i$ and its \textit{self intersection} in $Z$. 
A good resolution is called \textit{minimal} if there are no 
components $E_i \cong \PP^1$ with self intersection $-1$ which  
meet only one or two other components of the exceptional set.
Whenever such components occur in an arbitrary, 
good resolution, they can 
be blown down to a minimal good resolution.
For a discussion of the existence and uniqueness of such 
minimal good resolutions for surfaces we refer to \cite[Chapter 2]{Volume1} 
``The topology of surface singularities'' by Michel.

The genera and intersection multiplicities of the components $E_i$ 
of a good resolution can be summarized in its \textit{dual graph}.
Surface singularities are often described in terms of this resolution 
graph rather than by explicit equations; for instance, 
such dual graphs have already appeared 
in \cite[Chapter 10]{Volume1} on ``Finite dimensional Lie algebras in singularities''.
Note, however, that in general such a 
dual graph does not determine the singularity up to analytic 
isomorphism\footnote{Those normal surface singularities for which this is the 
  case are called ``taut'', see \cite{Laufer73}
}.

\medskip
Let $(X,0)$ be a rational surface singularity,  
$(X,0) \hookrightarrow (\mathcal X,0) \overset{\pi}{\longrightarrow} (S,0)$ 
its semi-universal deformation, and suppose $\rho \colon Z \to X$ is a 
minimal good resolution of singularities for $(X,0)$ as above.
Artin has shown in \cite[Theorem 3]{Artin74} that there exists a smooth space 
$(R,0)$ parametrizing those deformations of $Z$ that \textit{blow down} to 
deformations of $(X,0)$, i.e. the resolution $\rho$ extends to a projection 
of the total space $\mathcal Z$ of the deformation of $Z$ giving rise 
to another total space $\rho(\mathcal Z)$ of a deformation of $X$. 
This provides a commutative diagram 
\begin{equation}
  \xymatrix{
    Z \ar@{^{(}->}[r] \ar[d]^{\rho} & 
    \mathcal Z \ar[d]^{\rho} & 
    \\
    X \ar@{^{(}->}[r] \ar[d] & 
    \rho(\mathcal Z) \ar[r] \ar[d] & 
    \mathcal X \ar[d]^\pi \\
    \{ 0 \}  \ar@{^{(}->}[r] &
    R \ar[r]^{\Phi} &
    S.
  }
  \label{eqn:SimultaneousResolutionRationalSingularities}
\end{equation}
where $\Phi \colon (R,0) \to (S,0)$ is the comparison map to the base 
of the semi-universal deformation.
Artin has shown that $\Phi$ is finite and maps 
surjectively onto an irreducible component $(S',0)$ of $(S,0)$ 
which is now called the \textit{Artin component} of the base $(S,0)$.

Since the central fiber $Z$ in (\ref{eqn:SimultaneousResolutionRationalSingularities}) 
is smooth and proper over $X$, 
the family $\mathcal Z \to R$ is topologically trivial due to Ehresmann's 
Lemma (cf. Lemma \ref{lem:DiffeomorphismForResolutionInFamily} below). 
In particular, this entails that $\rho \colon Z_t \to X_t$ is a resolution of 
singularities for every fiber over $t\in R$ in a neighborhood of $0$. 
For this reason, a diagram like (\ref{eqn:SimultaneousResolutionRationalSingularities}) 
is also called a \textit{resolution in family} or a \textit{simultaneous resolution}.

\begin{example}
  \label{exp:SimultaneousResolutionDoublePoint}
  One instance of diagram (\ref{eqn:SimultaneousResolutionRationalSingularities}) 
  can be constructed for the $A_1$-singularity $(X,0) \subset (\CC^3,0)$ 
  which has already appeared in Example \ref{exp:VersalDeterminantalDeformations}, 1.
  It is defined by the equation $f = x^2-yz$, but can also be regarded as a determinantal 
  singularity of type $2$ for the matrix 
  \[
    A = 
    \begin{pmatrix}
      x & y \\ z & x
    \end{pmatrix}.
  \]
  One can check that a resolution of singularities for $(X,0) = (X_A^2,0)$ is given by 
  the Tjurina transform $\hat X_A^2 \subset \CC^3 \times \PP^1$. Furthermore, the 
  Tjurina transformation in family for the deformation induced by the unfolding 
  \[
    \mathbf A(x,y,z;t) = A(x,y,z) - t \cdot 
    \begin{pmatrix}
      1 & 0 \\
      0 & -1
    \end{pmatrix}
  \]
  on a parameter $t$ 
  blows down to the deformation of $(X,0)$ given by the perturbation $f - t^2$. This  
  furnishes the left hand side of (\ref{eqn:SimultaneousResolutionRationalSingularities}). 
  As predicted by Artin's theorem, 
  the comparison map $\Phi \colon (\CC,0) \to (\CC,0)$ taking $t$ to $t^2$ is a 
  $2:1$-cover of the base of the semi-universal deformation of $(X,0)$.
\end{example}

\begin{figure}[h]
  \centering
  \includegraphics[origin=c,page=1,scale=0.4,clip=true,trim=0cm 2cm 2cm 1cm]{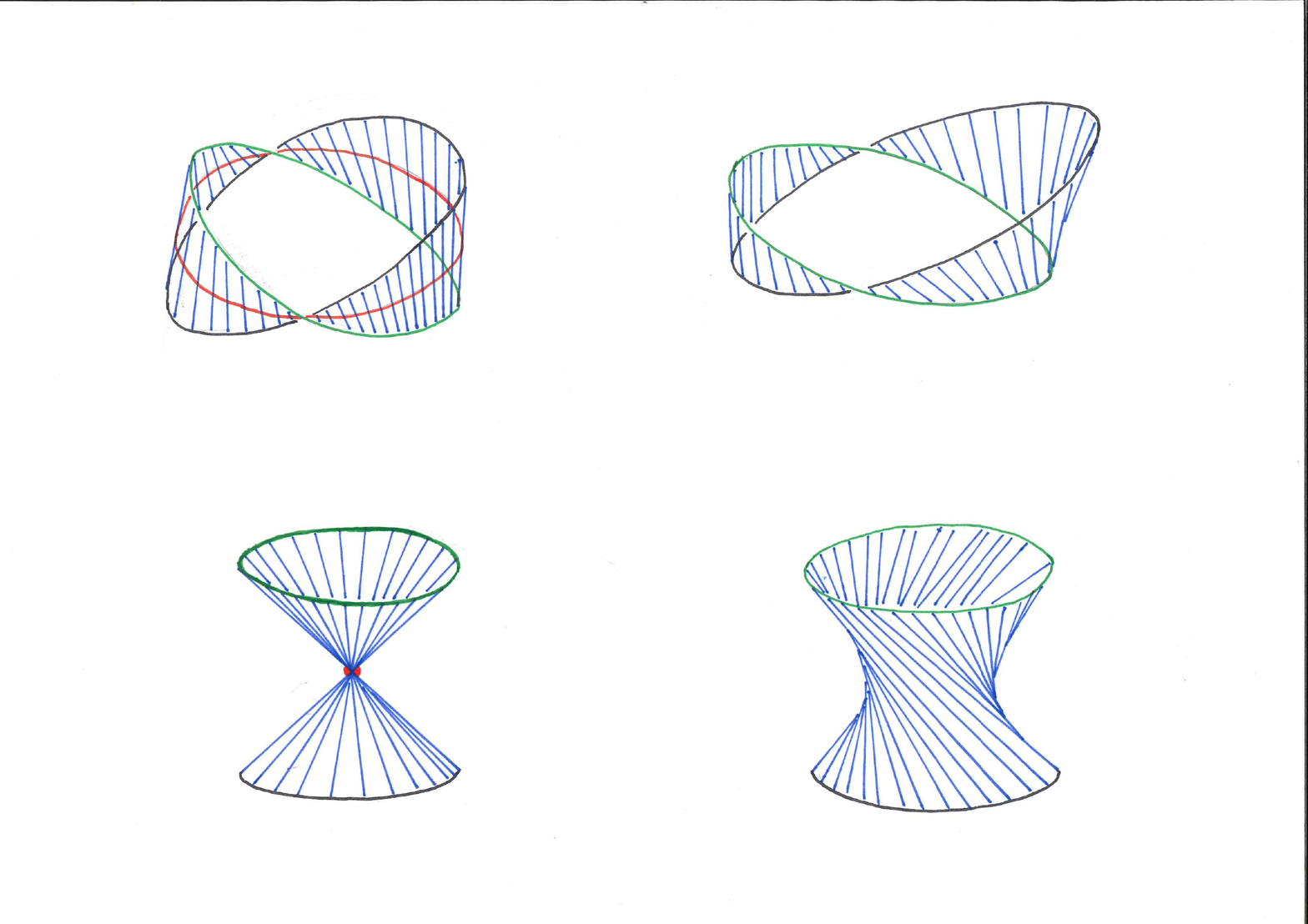}
  \caption{A resolution in family for the $A_1$-surface singularity in the positions of 
    the upper left square in (\ref{eqn:SimultaneousResolutionRationalSingularities}). The deformation 
    of the Tjurina transform does not change the topology, but forgets about the zero section (red) 
    of the disc bundle. The projection of the smooth fibers is an intrinsic isomorphism despite 
    the different embeddings.
  }
  \label{fig:ResolutionInFamilyDoublePoint}
\end{figure}

The previous example is also an instance of a general result on simultaneous 
resolutions for 
rational double points due to Brieskorn \cite{Brieskorn71}, which was independently 
proved by Tjurina \cite{Tjurina70}. 
We cite the summarized form from \cite{Wahl77}, 
cf. also \cite[Theorem 2]{Artin74}.

\begin{theorem}
  The versal deformation of a rational double point resolves simultaneously 
  after a Galois base change.
  \label{thm:BrieskornRationalDoublePoints}
\end{theorem}
Note that Brieskorn's and Tjurina's constructions do not 
necessarily involve the choice of a matrix structure and Tjurina transformation 
in family as in Example \ref{exp:SimultaneousResolutionDoublePoint} above. 
However, a resolution of singularities can at times be constructed using 
only Tjurina modifications; see \cite{Pedersen21} for further examples.


In \cite{Wahl77} Wahl has addressed the 
question how to obtain equations and even free resolutions of the 
defining ideals for rational singularities, given their dual graphs.
He shows in \cite[Proposition 3.2]{Wahl77}:

\begin{proposition}
  \label{prp:WahlDeterminantalRationalSingularities}
  Suppose $(X,0)\subset(\CC^p,0)$ is a rational surface singularity of embedding 
  dimension $p\geq 4$. If $(X,0)$ is determinantal, then it is determinantal of 
  type $(2,p-2,2)$.
\end{proposition}
Using this result, he then obtains in \cite[Theorem 3]{Wahl77}:

\begin{theorem}
  Let $(X,0)$ be a determinantal rational surface singularity 
  of embedding dimension $p$. Then 
  the dual graph of $(X,0)$ consists of one $-(p-1)$ curve and (possibly) 
  some $-2$ curves. 
  \label{thm:WahlDeterminantalRationalSingularities}
\end{theorem}
The converse of this theorem was later established by R\"ohr \cite{Roehr92} by de Jong 
\cite{deJong98} (who also produces explicit matrices)
so that indeed every rational surface singularity 
with this particular configuration in its dual graph is 
in fact determinantal.

\begin{example}
  \label{exp:TjurinaTransformRationalSurfaceExample}
  Consider the normal surface singularity $(X_A^2,0) \subset (\CC^4,0)$ given by the 
  matrix 
  \[
    A = 
    \begin{pmatrix}
      x & y & z \\
      y^2 & z & w
    \end{pmatrix}
  \]
  from Example \ref{exp:SurfaceDiscriminant}. Let $(u:v)$ be the homogeneous coordinates 
  of $\PP^1$ and $\hat X_A^2 \subset \CC^4\times \PP^1$ the Tjurina transform of 
  $(X_A^2,0)$ defined by the equations
  \[
    \begin{pmatrix}
      u & v
    \end{pmatrix} \cdot
    \begin{pmatrix}
      x & y & z \\
      y^2 & z & w
    \end{pmatrix} = 
    \begin{pmatrix}
      0 & 0 & 0
    \end{pmatrix}.
  \]
  Using the explicit equations it is easy to see that on the chart $\{ u \neq 0\}$ 
  the variables $x,y$ and $z$ can be eliminated so that 
  $\hat X_A^2 \cap \{ u \neq 0 \}$ is smooth. In the other chart $\{ v \neq 0 \}$ we 
  find a single $A_1$-hypersurface singularity at the origin, which is given by the 
  equation 
  \[
    y^2 + x \cdot \frac{u}{v} = 0
  \]
  after elimination of $z$ and $w$. 
  \begin{figure}[h]
    \centering
    \includegraphics[scale=0.25,clip=true,trim=2cm 7cm 2cm 2cm,angle=90,origin=c]{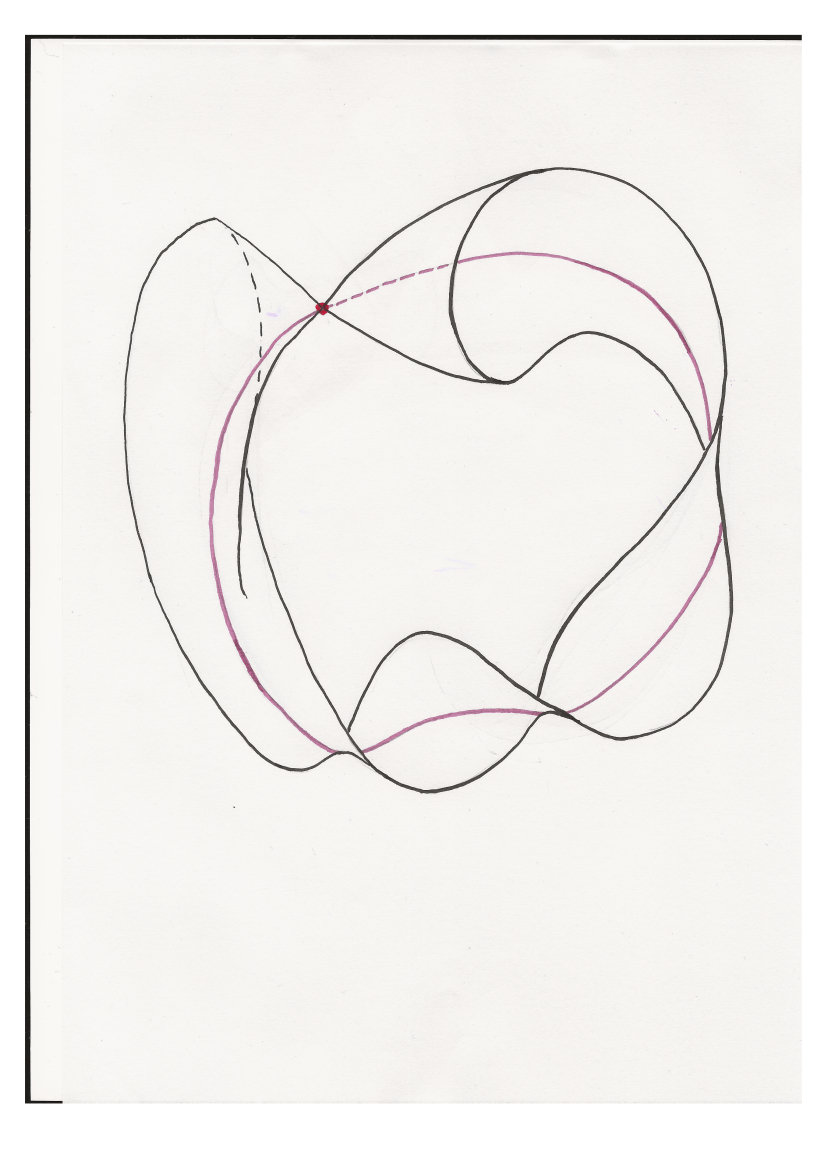}
    \includegraphics[scale=0.25,clip=true,trim=3cm 0.1cm 0cm 7cm,angle=-90,origin=c]{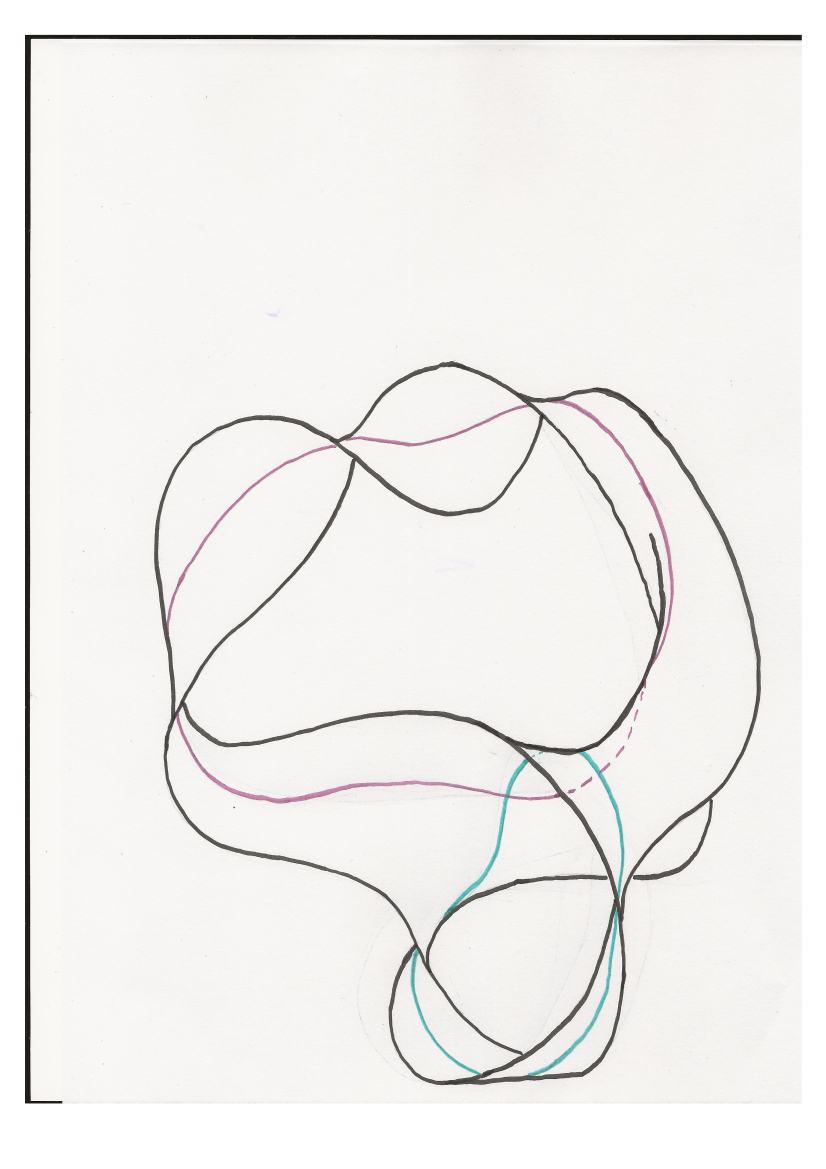}
    \caption{The real picture of the Tjurina transform (left) of the rational 
      double point 
      in Example \ref{exp:TjurinaTransformRationalSurfaceExample} and its full resolution (right), 
      together with their exceptional sets.
    }
    \label{fig:RDPExampleTjurinaTransformAndFullResolution}
  \end{figure}
  This singularity can be resolved 
  by a single classical blowup so that the overall exceptional set is a 
  normal crossing divisor. The dual graph of the full resolution is pictured 
  in Figure \ref{fig:DualGraphSimpleSurfaceExample}.
  \begin{figure}[h]
    \centering
    \includegraphics[scale=0.8,clip=true,trim=0cm 19cm 9cm 7cm]{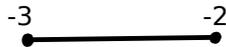}
    \caption{The dual graph of the rational double point in Example 
      \ref{exp:TjurinaTransformRationalSurfaceExample}; the vertices 
      stand for the exceptional curves (purple and green in Figure 
      \ref{fig:RDPExampleTjurinaTransformAndFullResolution}), the 
      numbers indicate their self intersections (the number of half 
      twists of their normal bundles in Figure 
      \ref{fig:RDPExampleTjurinaTransformAndFullResolution}), and 
      the edge corresponds to the existence of the intersection point of 
      the two curves. 
    }
    \label{fig:DualGraphSimpleSurfaceExample}
  \end{figure}
\end{example}

In a subsequent article \cite{Wahl79}, Wahl then investigates the deformation theoretic 
behaviour of rational surface singularities and finds in \cite[Theorem 3.2]{Wahl79}: 

\begin{theorem}
  For a determinantal rational surface singularity $(X,0)$ of 
  multiplicity $e\geq 3$ the Artin component 
  $(S',0)$ of deformations of $(X,0)$ that admit a resolution in family
  consists precisely of the determinantal deformations. 
  \label{thm:WahlArtinComponent}
\end{theorem}
Note that for any rational surface singularity, the multiplicity $e$ is 
equal to $p-1$ where $p$ is its embedding dimension, cf. \cite[Corollary 6]{Artin66}.
Wahl furthermore observes that for the determinantal deformations, the 
resolution in family always factors through the Tjurina transformation 
in family. In this context, the Tjurina transform $\hat X_A^2$ of $(X,0)$ 
is obtained from the resolution $(Z,E)$ by blowing down all the $-2$ configurations, 
cf. Theorem \ref{thm:WahlDeterminantalRationalSingularities}. 

\begin{example}
  We continue with Example \ref{exp:TjurinaTransformRationalSurfaceExample}. 
  Given the $A_1$-hypersurface singularity in the chart $\{v \neq 0\}$ 
  we already saw in Example \ref{exp:SimultaneousResolutionDoublePoint} how 
  to construct a resolution in family: Write the local 
  equation for the singularity as the determinant of a $2\times 2$-matrix
  \[
    y^2 + \frac{u}{v} \cdot x = \det 
    \begin{pmatrix}
      y & x \\
      -\frac{v}{u} & y
    \end{pmatrix}
  \]
  and resolve it by Tjurina modification rather than the standard blowup.
  One can check that the local resolution in family is compatible with the Tjurina 
  transformation in family so that we indeed obtain a full resolution in family for 
  any stabilization of the defining matrix $A$. 
\end{example}

\subsection{Further references and techniques} 
\label{sec:FurtherReferencesAndTechniques}

\subsubsection{Construction of ``Versal determinantal deformations''}

In view of the three particular cases of complete intersection, 
Cohen-Macaulay codimension $2$, and Gorenstein singularities of 
codimension $3$ discussed above and given the fact that the 
theory of unfoldings of map germs is 
in many respects much simpler than the theory of deformations of 
germs $(X,0) \subset(\CC^p,0)$ in general, one might be tempted 
to develop a theory of ``semi-universal determinantal deformations'' 
for arbitrary determinantal singularities $(X_A^s,0)\subset(\CC^p,0)$ 
by fixing the defining matrix $A \in \CC\{x\}^{m\times n}$ and restricting 
to those deformations coming from unfoldings of $A$, i.e. the 
image of the comparison map $\Phi$ in 
(\ref{eqn:DeterminantalDeformationsVsSemiUniversalDeformations}). 
We already saw in Pinkham's example and, more generally, the discussion of the Artin component 
in Section \ref{sec:RationalSurfaceSingularities} that 
$\Phi$ is not surjective in general. But for determinantal singularities 
it is nevertheless a natural question, whether the image of $\Phi$ 
can be reconstructed directly from the germ $(X_A^s,0)$ and its 
fixed determinantal structure itself rather than from 
the unfoldings of the defining matrix $A$. 


An attempt along these lines has been made by Schaps in 
\cite{Schaps83} where she follows the classical approach to deformation 
theory due to Schlessinger \cite{Schlessinger68} in order 
to define a functor of ``$A$-determinantal deformations''
for determinantal schemes with a fixed matrix $A$.
The goal is then to prove that this functor has a pro-representable hull. 
However, one runs into technical difficulties caused by the fact in the 
general case, an infinitesimal determinantal deformation of $(X_A^s,0)$ does not lift 
uniquely to an infinitesimal unfolding of $A$. One instance of this 
phenomenon can be found in Example \ref{exp:VersalDeterminantalDeformations}, 4 
where we have a continuous modulus $\gamma'$ of the defining matrix, 
the variation of which results in a trivial deformation 
of the associated determinantal singularity. This observation leads Schaps to define a 
``unique lifiting proporty'': 

\begin{definition}
  Let $B' \to B$ be a surjection of Artin rings and 
  $A \in B[x]$ a matrix with $x = x_1,\dots,x_p$ a fixed set of variables. 
  Then $A$ is said to satisfy the unique lifting property, if 
  for any two liftings $A_1,A_2 \in B'[x]$ of $A$, whose minors generate the same ideal, 
  one has that $A_1$ is $\GL$-equivalent to $A_2$. 
\end{definition}
Schaps then proceeds to show, \cite[Proposition 1]{Schaps83}:
\begin{proposition}
  \label{eqn:SchapsResultOnProrepresentabilityOfTheDeformationFunctor}
  Let $(X,0) \subset (\CC^p,0)$ be an isolated determinantal singularity with 
  a defining matrix $A \in \CC\{x\}^{m\times n}$ satisfying the unique lifting 
  property. Then the functor of $A$-determinantal deformations of $(X,0)$ has 
  a prorepresentable hull.
\end{proposition}

\subsubsection{Buchweitz' criterion for deformations to be determinantal}

The question as to when or under which conditions all deformations 
of a given singularity $(X_A^s,0)$ are determinantal for the defining 
matrix $A$, 
has been addressed more generally 
by Buchweitz in his thesis in \cite[Section 4.3]{Buchweitz81} 
headed ``Deformations d'un type donn\'e et d\'eploiements''.
Determinantal singularities appear as a special case in 
\cite[Exemples 4.3.2 b)]{Buchweitz81}. Before we can state the main 
theorem \cite[Theorem 4.7.1]{Buchweitz81} in its adapted version for 
determinantal singularities, we need to introduce some further 
mathematical notions.

Recall, that
given a morphism of analytic local $k$-algebras $\varphi: R \longrightarrow S$
and an $R$-module $M$ of finite type, $\varphi$ is said to be transversal to 
$M$, if $\operatorname{Tor}_i^R(M,S)$ vanishes for all $i>0$. For an 
$R$-module of finite type, the Auslander module $D(M)$ of $M$ is defined as the
cokernel of the dual of a free presentation of $M$. More precisely, let
$$F_1 \stackrel{\psi}{\longrightarrow} F_0  
        \longrightarrow M  \longrightarrow  0$$
be a free presentation of $M$, then dualizing 
$$0 \longrightarrow M^\vee \longrightarrow F_0^\vee 
    \stackrel{\psi^\vee}{\longrightarrow} F_1 ^\vee \longrightarrow D(M)
    \longrightarrow 0$$
provides a presentation of $D(M)$.
For an analytic space germ $(X,0) \subset ({\mathbb C}^p,0)$ with ideal
$I_{X,0} \subseteq {\mathbb C}\{x\}$ the \textit{Auslander module of $(X,0)$} 
is then the module $D(I_{X,0}/I_{X,0}^2)$ whose homological properties are
independent of the choice of embedding and free resolution, 
cf. \cite[Section 4.5]{Buchweitz81}. With this notation we have the 
following adaptation of \cite[Theorem 4.7.1]{Buchweitz81}:

\begin{theorem}
  \label{thm:Buchweitz}
  Let $(X,0) = (X_A^s,0) \subseteq ({\mathbb C}^p,0)$ be a determinantal 
  singularity of type $s$ 
  defined by a matrix $A \in \CC\{x\}^{m\times n}$.
  Then the following properties are equivalent:
  \begin{enumerate}
    \item Each first order deformation of $(X,0)$ is determinantal for $A$.
    \item Each analytic deformation of $(X,0)$ is determinantal for $A$.
    \item The generic determinantal singularity $(M_{m,n}^s,0) \subset (\CC^{m\times n},0)$
      is rigid and $(X,0)$ is unobstructed.
    \item The Auslander module of the generic determinantal singularity $(M_{m,n}^s,0)$
      is transversal to the map 
      $A: (X,0) \longrightarrow (M_{m,n}^s,0)$ and $(M_{m,n}^s,0)$ is rigid.
  \end{enumerate}
  If additionally the vector space dimension of $T^1_{X,0}$ is finite, then this is
  also equivalent to the property that there is an unfolding as in 
  Lemma \ref{lem:UnfoldingInducesDeformation} with smooth analytic base which induces a 
  versal deformation of $(X,0)$.
\end{theorem}

\begin{remark}
  \label{rem:SvanesResultOnRigidity}
  The generic determinantal singularities $(M_{m,n}^s,0) \subset (\CC^{m\times n},0)$ are 
  all rigid except for the generic determinantal hypersurface singularities 
  $(M_{m,m}^m,0)$ defined by $f = \det = 0$ (or $\Pf = 0$ for skew symmetric matrices).
  This was proved independently by Svanes \cite{Svanes72} and J\"ahner \cite{Jaehner74}; 
  see \cite[Chapter 15.C]{BrunsVetter88} for a textbook account.
  Therefore, given characterization 3. in Theorem \ref{thm:Buchweitz}, 
  the only condition on a given 
  determinantal singularity $(X_A^s,0)$ that really needs to be checked is the unobstructedness. 
\end{remark}

\begin{remark}
  Theorem \ref{thm:Buchweitz} suggests to review 
  Corollary \ref{cor:SemiUniversalDeformationsAndUnfoldingsCoincide} for 
  complete intersections 
  and Corollary \ref{cor:SemiUniversalDeformationOfICMC2} for isolated Cohen-Macaulay 
  codimension $2$ singularities. Given the 
  explicit identifications $T^1_\GL(A) \cong T^1_{X,0}$ from 
  Lemma \ref{lem:IsomorphismOfFirstOrderDeformationsForICIS} and respectively 
  Corollary \ref{cor:SemiUniversalDeformationsAndUnfoldingsCoincide}, 
  the coincidence of a semi-universal deformation of the determinantal 
  singularity with the deformation induced from a miniversal unfolding of the defining 
  matrix now follows directly from the implications 1. $\Rightarrow$ 3. $\Rightarrow$ 2. 
  in Theorem \ref{thm:Buchweitz}. 

  It is expected that the construction of a ``matrix-$T^1$'' as in 
  Corollary \ref{cor:FruehbisKruegerMatrixT1} for ICMC2 singularities is possible 
  in many more cases. For non-maximal minors, however, one can in general only 
  expect surjective maps $T^1_\GL(A) \to T^1_{X_A^s,0}$ which have a non-trivial kernel. 
  This can, for instance, be observed from Table \ref{tab:FKN4Fold} for the 
  simple ICMC2 fourfold singularities given by matrices 
  $A \colon (\CC^6,0) \to (\CC^{2\times 3},0)$. In this case not only $(X_A^2,0)$, but 
  also the zero-dimensional complete intersections $(X_A^1,0)$ are determinantal singularities. 
  The $\GL$-codimensions 
  $\tau_\GL(A)= \dim_\CC T^1_\GL(A)$ of the defining matrices and the Tjurina numbers
  $\tau(X_A^1,0) = \dim_\CC T^1_{X^1_A,0}$ of the singularities 
  $(X_A^1,0)$ are listed in the fourth and the fifth column, respectively. 
  We find series for which these numbers coincide and others where 
  $\tau(X^1_A,0) \leq \tau_\GL(A)$ with a strict inequality in general. 

  A geometric reason for this discrepancy for non-maximal minors is the following. The 
  $\GL$-equivalence of matrices is sensitive to the position of the image of the defining 
  matrix $A$ relative to all strata $V_{m,n}^r$ of the rank stratification. For 
  non-maximal minors of size $s< \min\{m,n\}$ only the strata $V_{m,n}^r$ 
  of $M_{m,n}^s$ with $r<s$ are relevant for the deformations of $(X_A^s,0)$. Those 
  inifitesimal unfoldings varying only the position of the image of $A$ relative to higher dimensional 
  strata will therefore be discarded by the comparison map $\Phi \colon T^1_\GL(A) \to T^1_{X_A^s,0}$.
\end{remark}

\section{Classification of simple singularities}

Recall that a singularity is called \textit{simple} if only a finite number of
non-equivalent singularities appear in its versal family. 
Arnold gives 
a complete list of simple isolated hypersurface singularities of arbitrary 
dimension in his article \cite{Arnold72} from 1972 and the 
classification of all simple isolated complete intersection singularities 
was completed  by Giusti \cite{Giu} 
in the mid 1980s. The question which singularities are simple 
is still not fully answered for 
determinantal singularities, but there exist classifications 
for symmetric square matrices by Bruce in \cite{Bruce03},
for square matrices by Bruce and Tari in \cite{BruceTari04},
and for skew-symmetric
ones by Haslinger in \cite{Haslinger01} (incomplete), as well as for isolated Cohen-Macaulay 
codimensions $2$ singularities 
by the first named author and Neumer in \cite{FruehbisKrueger99} and \cite{FruehbisKruegerNeumer10}. 
Note that,
in contrast to simple hypersurfaces or complete intersections, a simple 
determinantal singularity does not need to be smoothable or isolated, as there
are rigid non-isolated determinantal singularities and any
rigid singularity is simple for trivial reasons. We give a brief overview on 
known classification results, leaving the explicit tables to the appendix of
this article. 

\subsection{Singularities of square matrices} 
\label{sec:SimpleSquareMatrices}

In \cite{BruceTari04} Bruce and 
Tari classify all simple singularities for square matrices 
\[
  A \colon (\mathbb K^p,0) \to (\mathbb K^{m\times m},0)
\] 
up to $\GL$-equivalence, where either $\mathbb K = \RR$ and $A$ smooth, 
or $\mathbb K = \CC$ and $A$ holomorphic, the treatement of the real case made possible 
by Damon's theory, \cite[Remark 2.8]{BruceTari04}.
In the Tables 
\ref{tab:SquareSymm3p2}, \ref{tab:NormalFormsSquareMatricesSize2In2Space}, 
\ref{tab:NormalFormsSquareMatricesSize2In3Space}, and 
\ref{tab:NormalFormsSquareMatricesSize3In2Space} 
the consideration of the real case occasionally leads to 
a $\pm$-sign with different singularities over $\RR$. 
Over $\CC$ these signs can be omitted.

Associated to the matrices $A$ as above, Bruce and Tari also consider 
\textit{determinantal hypersurface singularities}\footnote{In 
  \cite{BruceTari04} these are called the 
\textit{discriminants} of the matrix.} (abbreviated by DHS) in the following)
\[
  (X_A^m,0) = (\{\det A = 0 \},0) \subset (\CC^p,0)
\]
which are defined by the equation $f = \det A$.
Moreover, they also consider the Tjurina 
transformation\footnote{Tjurina transforms are called \textit{criminants} in \cite{BruceTari04}.} 
and, as it turns out, the simple singularities of $A$ are closely related to simple 
singularities of both, the determinantal hypersurface $(X_A^m,0)$ and of its 
Tjurina transforms $\hat X_A^m$.
Besides the simple isolated hypersurface singularities, the well known 
A-D-E-singularities that were classified 
by Arnold up to right- or $\mathcal R$-equivalence in \cite{Arnold72}, there also appear the 
simple singularities from the classification for functions on manifolds 
with boundary up to $\mathcal R_\delta$-equivalence from \cite{Arnold78}. 

As usual, Bruce and Tari assume all matrices $A$ to have entries in the 
maximal ideal $\m$ of $\CC\{x_1,\dots,x_p\}$ so that the \textit{corank} 
of the matrix $A(0)$, i.e. the codimension of its image, is equal to the 
size of the matrix $m$.
A key handle to the classification is to also consider the corank of the differential 
\[
  \D A(0)\colon T_0 \mathbb K^p \to T_0 \mathbb K^{m\times m}
\]
of $A$ at the origin.

With all this notation at hand, we can now reproduce their main results:

\begin{theorem}
  \label{thm:ClassificationOfSimpleSquareMatrices}
  \begin{enumerate}
    \item When $p=1$, all finitely $\GL$-determined germs are simple and $\GL$-equivalent 
      to a germ of the form $\mathrm{diag}(x^{\alpha_1}, x^{\alpha_2},\dots,x^{\alpha_m})$ 
      where $\alpha_1\leq \alpha_2 \leq \dots \leq \alpha_m$. This germ has 
      $\GL$-codimension $\sum_{i=1}^m (2(m-i)+1)\alpha_i -1$. Its Tjurina transform is 
      empty.
    \item When the corank of the differential $\D A(0)$ is zero there is a normal form 
      \[
	A \colon (\mathbb K^{m^2},0) \times (\mathbb K^l,0) \to (\mathbb K^{m\times m},0),
	\quad 
	(x,z) \mapsto A(x), \quad a_{i,j} = x_{i,j}.
      \]
      where the variables $z_1,\dots,z_l$ are redundant. 
    \item When the corank of the differential is one, there are two cases. 
      \begin{enumerate}
	\item A normal form $A \colon (\mathbb K^{m^2-1},0) \times (\mathbb K^l,0) \to 
	  (\mathbb K^{m\times m},0)$ given by
	  \[
	    \left( \sum_{i=2}^m x_{i,i} + f(z) \right) \cdot E_{1,1} + 
	    \sum_{(i,j) \neq (1,1)} x_{i,j} \cdot E_{i,j}
	  \]
	  where $E_{i,j}$ is the matrix with a $1$ in the $(i,j)$-th place and zeroes 
	  elsewhere, $\{x_{i,j}\}_{(i,j) \neq (1,1)}$ are the coordinates of the first factor
	  $\mathbb K^{m^2-1}$, and $f \colon (\mathbb K^l,0) \to (\mathbb K,0)$ is 
	  one of Arnold's $\mathcal R$-simple germs (see Table \ref{tab:ArnoldR} in the appendix). 
	  The $\GL$-codimension of $A$ 
	  coincides with the Tjurina number of $f$. 
	\item A normal form $A$ as above given by 
	  \[
	    \left( \sum_{i=2}^{m-1} x_{i,i} + f(x_{m,m},z) \right)\cdot E_{1,1} 
	    + \sum_{(i,j) \neq (1,1)} x_{i,j} \cdot E_{i,j}
	  \]
	  where $f \colon (\mathbb K,0) \times (\mathbb K^{l},0) \to (\mathbb K,0)$ 
	  is one of Arnold's $\mathcal R_\delta$-simple germs of singularities of 
	  functions on manifolds with boundary (see Table \ref{tab:ArnoldRDelta}). 
	  The $\GL$-codimension of $A$ coincides 
	  with the $\mathcal R_\delta$-codimension of $f$.
      \end{enumerate}
      In both cases, the Tjurina transforms of the singularities are smooth.
    \item When the corank of the differential is two, then $m=3$ and $p=7$. 
      The simple matrices in this class can be derived from the symmetric 
      ones\footnote{
        Goryunov has informed us that there were mistakes in this part of the  
	original classification in \cite[Paragraph C, p 757]{BruceTari04}. 
	The corrected description given here is taken from \cite[Theorem 3.7]{Goryunov21}.
      } in Table \ref{tab:SquareSymm3p4} by addition of the matrix
      \[
	U = 
	\begin{pmatrix}
	  0 & u_{12} & u_{13} \\
	  -u_{12} & 0 & u_{23} \\
	  -u_{13} & -u_{23} & 0 
	\end{pmatrix}.
      \]
      The Tjurina transforms are all smooth. 
    \item When $m = 2$, the simple germs that are not covered by the preceding items 
      are given in Table \ref{tab:NormalFormsSquareMatricesSize2In2Space} and 
      Table \ref{tab:NormalFormsSquareMatricesSize2In3Space}. 
    \item When $m =3$ the simple germs that are not covered by the preceding items 
      are given in Table \ref{tab:SquareSymm3p2}\footnote{
	Again, the normal forms given in Table \ref{tab:SquareSymm3p2} are not the original ones 
	found by Bruce and Tari. We were informed by Goryunov about a mistake in the original 
	classification and the given table has the correct matrices up to $\GL$-equivalence.
      }.
  \end{enumerate}

\end{theorem}

The following 
classification of simple \textit{symmetric} matrices by Bruce \cite[Theorem 1.1]{Bruce03} precedes 
the above classification by one year. The setting is slightly different, 
since only \textit{holomorphic} matrices 
\[
  A \colon (\mathbb C^p,0) \to (\CC^{m\times m}_\sym,0)
\] 
are considered up to symmetric $\GL$-equivalence, see Remark \ref{rem:SLEquivalence}.

\begin{theorem}
  \label{thm:ClassificationOfSimpleSymmetricSquareMatrices}
  \begin{enumerate}
    \item When $p=1$, all finitely $\GL$-determined germs are simple and $\GL$-equivalent 
      to a germ of the form $\mathrm{diag}(x^{\alpha_1}, x^{\alpha_2},\dots,x^{\alpha_m})$ 
      where $\alpha_1\leq \alpha_2 \leq \dots \leq \alpha_m$. This germ has 
      $\GL$-codimension $\sum_{i=1}^m ((m-i)+1)\alpha_i-1$. 
    \item When the corank of the differential $\D A(0)$ is zero there is a normal form 
      \[
	A \colon (\mathbb C^{N},0) \times (\mathbb C^l,0) \to (\CC^{m\times m}_\sym,0),
	\quad 
	(x,z) \mapsto A(x), \quad a_{i,j} = x_{i,j}.
      \]
      where $N=\frac{m(m+1)}{2}$ and the variables $z_1,\dots,z_l$ are 
      redundant. 
    \item When the corank of the differential is one, there are two cases. 
      \begin{enumerate}
	\item A normal form $A \colon (\mathbb C^{N-1},0) \times (\mathbb C^l,0) \to 
	  (Sym_m(\mathbb C),0)$ given by
	  \[
	    \left( \sum_{i=2}^m x_{i,i} + f(z) \right) \cdot E_{1,1} + 
	    \sum_{(i,j) \neq (1,1)} x_{i,j} \cdot E_{i,j}
	  \]
	  where $E_{i,j}$ is the matrix with a $1$ in the $(i,j)$-th place and zeroes 
	  elsewhere, $\{x_{i,j}\}_{(i,j) \neq (1,1)}$ are the coordinates of the first factor
	  $\mathbb C^{N-1}$, and $f \colon (\mathbb C^l,0) \to (\mathbb C,0)$ is 
	  one of Arnold's $\mathcal R$-simple germs. The $\GL$-codimension of $A$ 
	  coincides with the Tjurina number of $f$. 
	\item A normal form $A$ as above given by 
	  \[
	    \left( \sum_{i=2}^{m-1} x_{i,i} + f(x_{m,m},z) \right)\cdot E_{1,1} 
	    + \sum_{(i,j) \neq (1,1)} x_{i,j} \cdot E_{i,j}
	  \]
	  where $f \colon (\mathbb C,0) \times (\mathbb C^l,0) \to (\mathbb C,0)$ 
	  is one of Arnold's $\mathcal R_\delta$-simple germs of singularities of 
	  functions on manifolds with boundary. 
	  The $\GL$-codimension of $A$ coincides 
	  with the $\mathcal R_\delta$-codimension of $f$.
      \end{enumerate}
    \item When $m=p=2$ the $\GL$-simple germs are given in Table \ref{tab:SquareSymm2p2}.
    \item If $m=3$, $p=2$ the $\GL$-simple germs are listed in Table \ref{tab:SquareSymm3p2}.
    \item If $m=3$, $p=4$ the $\GL$-simple germs are given in Table \ref{tab:SquareSymm3p4}.
\end{enumerate}
\end{theorem}

For skew-symmetric square matrices, only a partial classification of simple
singularities is known from the thesis of Haslinger \cite{Haslinger01}. 
\begin{theorem}
  Let $A \colon (\CC^p,0) \to (\CC^{m\times m}_\sk,0)$ be a germ of skew-symmetric 
  matrices. 
  \begin{enumerate}
    \item (\cite[Proposition 4.6.1]{Haslinger01}) 
      When $p=1$, then $A$ is $\mathcal G_\sk$-equivalent to a matrix of the form 
      \[
	\bigoplus_{j=1}^{m'} x^{k_j} \cdot I^\sk_{s_j}
      \]
      for some $m'\leq m$ and integers $0<k_1<k_2<\dots<k_{m'}$
      where $I^\sk_{s} = \sum_{i=1}^s \left(E_{2i-1,2i} - E_{2i,2i-1} \right)$ is 
      the skew symmetric analogue of the $2s\times 2s$-identity matrix 
      and the direct sum stands for taking successive diagonal blocks.
    \item (\cite[Lemma 4.6.2 (i)]{Haslinger01}) When $m = 2$ any two matrices 
      \[
	A = 
	\begin{pmatrix}
	  0 & a \\ -a & 0
	\end{pmatrix}, \quad
	B = 
	\begin{pmatrix}
	  0 & b \\ -b & 0
	\end{pmatrix}
      \]
      are $\mathcal G_\sk$-equivalent if and only if $a$ and $b$ are $\mathcal K$-equivalent. 
      In particular, the $\mathcal G_\sk$-simple matrices are given by the classification 
      of $\mathcal R$-simple functions $(\CC^p,0) \to (\CC,0)$. 
    \item 
      (\cite[Lemma 4.6.2 (ii)]{Haslinger01})
      When $m=3$ then $A$ is finitely $\mathcal G_\sk$-determined if and only if 
      the corresponding map germ $(\CC^p,0) \to (\CC^3,0)$ is finitely $\mathcal K$-determined 
      and two matrices are $\mathcal G_\sk$-equivalent if and only if this holds 
      for the corresponding map germ. In particular, the $\mathcal G_\sk$-simple 
      matrices are given by the $\mathcal K$-classification of simple 
      germs $(\CC^p,0) \to (\CC^3,0)$. 
    \item (\cite[Theorem 6.1.1]{Haslinger01}) 
      The $\mathcal G_\sk$-simple germs for $p=2$ and $m=4$ are listed 
      in Table \ref{tab:SquareSkewm4}. 
  \end{enumerate}
  \label{thm:ClassificationOfHaslinger}
\end{theorem}

Arnold's simple hypersurface singularities from the A-D-E-classification in 
Table \ref{tab:ArnoldR} and the boundary singularities in 
Table \ref{tab:ArnoldRDelta} make numerous appearances among the simple matrices 
in Theorems \ref{thm:ClassificationOfSimpleSquareMatrices}, 
\ref{thm:ClassificationOfSimpleSymmetricSquareMatrices}, and 
\ref{thm:ClassificationOfHaslinger}. While this is not really a surprise, 
the general theme of the interplay between the different lists
is still far from being completely understood. 

\begin{remark}
  \label{rem:SimpleHypersurfacesInTheMatrixSingularities}

  In \cite{GoryunovZakalyukin03} Goryunov and Zakalyukin 
  relate the simple symmetric matrix singularities in two variables 
  from Theorem \ref{thm:ClassificationOfSimpleSymmetricSquareMatrices}  
  (Tables \ref{tab:SquareSymm2p2} and \ref{tab:SquareSymm3p2})
  to Arnold's simple hypersurfaces as follows. 
  
  As their name suggests, the A-D-E-singularities correspond to 
  \textit{Dynkin diagrams} which have already appeared in
  \cite[Chapter 10]{Volume1} 
  and \cite[Chapter 8]{Volume1}. 
  These diagrams encode
  the dual resolution graphs of the simple hypersurface singularities 
  in dimension $2$, see \cite[Section 10.2.1]{Volume1}, 
  and, as suggested by the resolution in family 
  for these singularities, Theorem \ref{thm:BrieskornRationalDoublePoints}, 
  they can also be used to describe the vanishing homology and monodromy 
  operations, see e.g. \cite[Chapter 8]{Volume1}.
  
  Originally, the Dynkin diagrams classify the Weyl groups of the semi-simple Lie algebras. 
  These are subgroups of the linear automorphisms of Euclidean space 
  generated by reflections on hyperplanes. 
  In a case-by-case analysis, Goryunov and Zakalyukin found that 
  the simple symmetric matrices from Theorem \ref{thm:ClassificationOfSimpleSymmetricSquareMatrices} 
  correspond uniquely to pairs $(X;Y)$ of a Weyl group $X$ and a subgroup $Y$ 
  obtained by omitting certain reflections. This can be stated more 
  precisely in terms of Dynkin 
  diagrams: One obtains the \textit{affine} Dynkin diagram of $Y$ 
  by deletion of either two $1$-vertices for the corank $2$ families 
  (Table \ref{tab:SquareSymm2p2}) 
  or one $2$-vertex for the corank $3$ families (Table \ref{tab:SquareSymm3p2} from the affine 
  Dynkin diagram of $X$. These pairs of Weyl groups are listed in 
  the last columns of Table \ref{tab:SquareSymm2p2} and 
  Table \ref{tab:SquareSymm3p2}, respectively.

  Interestingly, the correspondence of simple symmetric matrices $A$ with pairs of 
  Weyl groups $(X;Y)$ is established via a study of the miniversal 
  unfoldings of the matrix $A$ and the associated hypersurface $f = \det A$.
  It turns out that 
  for all the singularities
  in question one has $\mu = \mu_f = \tau_{\GL}(A)$, 
  cf. Theorem \ref{thm:GoryunovMondSmoothableHypersurfaces} below.
  The space $\CC^\mu$ can be regarded as the configuration 
  space for the reflections in $X$ and $Y$ and 
  then the bases of the respective miniversal unfoldings can be 
  identified with $\CC^{\tau_{\GL}(A)} \cong\CC^\mu/X$ and 
  $\CC^{\mu_f} \cong \CC^\mu/Y$, respectively so that 
  the comparison map $\Phi$ from \ref{eqn:DeterminantalDeformationsVsSemiUniversalDeformations} 
  completes the quotient maps to a diagram 
  \[
    \xymatrix{
      (\CC^\mu, \mathcal A_X) \ar[rr]^{/Y} \ar[dr]_{/X} & & 
      (\CC^{\mu},\Sigma) \ar[dl]^\Phi \\
      & (\CC^{\mu}, \Delta)
    }
  \]
  where $\mathcal A_X$ is the configuration of mirrors of $X$ in $\CC^\mu$, 
  $\Sigma \subset \CC^\mu$ is the discriminant of $f$, 
  and $\Delta \subset \CC^\mu$ the matrix discriminant of $A$, 
  see \cite[Corollary 3.10]{GoryunovZakalyukin03}.
  Then $\Phi$ is a finite covering of order $|X:Y|$, 
  branched over the discriminant 
  $\Delta$ of the function $f$. 
\end{remark}

\begin{remark}
  \label{rem:SimpleHypersurfacesInSimpleMatrices}
  Recently in \cite{Goryunov21}, the investigation of the discriminants, 
  bifurcation diagrams, and the monodromy of simple matrix 
  singularities has been pushed further by Goryunov to also comprise those matrices $A$ for 
  which the associated hypersurface $f = \det(A)$ is not 
  smoothable via determinantal deformations or has non-isolated singular locus 
  (cf. Lemma \ref{lem:WhitneyStratificationOfAnEIDSAndSmoothability}).
  This has revealed a further correspondence of the simple symmetric/square matrices 
  in Table \ref{tab:SquareSymm3p4} with $\mathcal R_{\mathrm{odd}}$-simple 
  singularities of 
  odd functions\footnote{A function $f \colon (\CC^p,0)\to (\CC,0)$ is \textit{odd} if it 
    changes sign under the central symmetry of $(\CC^p,0)$. These are classified 
    up to right equivalence by diffeomorphisms of $(\CC^p,0)$ which preserve this 
    symmetry.
  } on $(\CC^2,0)$ and $\mathcal K_{\mathrm{odd}}$-simple ICIS\footnote{
    These are defined as centrally symmetric isolated complete 
    intersection curves $(C,0) \subset(\CC^3,0)$ and classified by the 
    subgroup of $\mathcal K$ in which the group of diffeomorphisms 
    $\mathcal R$ in the source is replaced by those preserving the symmetry. 
  } in $(\CC^3,0)$. These associated odd functions and symmetric ICIS 
  are listed\footnote{
    The notation used extends Giusti's list in Table \ref{tab:GiustiCurve} 
    in a natural way (note that $U_{11}$ and $U_{13}$ are not simple if the 
    symmetry condition is removed).
  }
  in the last two columns of Table \ref{tab:SquareSymm3p4}. 
  In this case, the correspondence is established by a natural identification 
  of the \textit{bifurcation diagrams} of the matrix with those of 
  the respective odd functions and symmetric space curves, see 
  \cite[Proposition 6.7]{Goryunov21}. 
\end{remark}

\subsection{Cohen-Macaulay codimension $2$ singularities}

In \cite{FruehbisKrueger99} the first named author classified the simple isolated space curve 
singularities up to isomorphism of space germs. Later she extended this
together with Neumer to all isolated Cohen-Macaulay codimension $2$ 
singularities in \cite{FruehbisKruegerNeumer10}.
By the Hilbert-Burch theorem (Theorem \ref{thm:HilbertBurch}), 
this amounts to classifying singularities of matrices 
$$A:({\mathbb C}^p,0) \longrightarrow ({\mathbb C}^{m\times (m+1)},0)$$
up to $\GL$-equivalence. 

Chronologically preceding the criteria for finite determinacy 
based on Damon's work, these 
articles use a weighted determinacy criterion followed by a matrix variant 
of the Arnold's rotating ruler method. 

The main results can be summarized as follows

\begin{theorem}
Isolated Cohen-Macaulay codimension $2$ singularities with defining matrix 
$A \colon (\CC^p,0) \to (\CC^{m\times (m+1)},0)$ exist only in the following cases:
\begin{enumerate}
\item Fat points in $({\mathbb C}^2,0)$: $p=2$, $m\in \{1,2\}$; 
Giusti's list of complete intersections ($m=1$) is included in Table 
\ref{tab:GiustiFatPoints}, the one for $m=2$ in Table \ref{tab:FKNFatPoints}.
\item Space curves in $({\mathbb C}^3,0)$: $p=3$, $m\in \{1,2\}$;
  the complete list can be found in Tables \ref{tab:GiustiCurve} and \ref{tab:FKCurve}.
\item Normal surfaces in $({\mathbb C}^4,0)$: $p=4$, $m=2$;
the list is in Table \ref{tab:FKNSurface}.
\item $3$-folds in $({\mathbb C}^5,0)$: $p=5$, $m=2$;
see Table \ref{tab:FKN3Fold}.
\item $4$-folds in $({\mathbb C}^6,0)$: $p=6$, $m=2$;
see Table \ref{tab:FKN4Fold}.
\end{enumerate}
\end{theorem}

\noindent Concerning these lists, there are a few noteworthy observations.
All singularities in the lists, except the $4$-folds, are smoothable, but only
fat points and surfaces exhibit complete intersection singularities in their 
versal family, which are not of hypersurface type. 
For normal surfaces, the list reproduces the list of rational triple point
singularities found by Tjurina in \cite{Tjurina68}, where she proceeds by a 
completely different approach, cf. Sections \ref{sec:RationalSurfaceSingularities}
and \ref{sec:ICMC2Topology}.
For $3$-fold singularities there is one infinite series, which holds 
as one matrix entry precisely the equations of hypersurface singularities 
in $2$ variables and also mimics their deformation behaviour (i.e. their
adjacencies).
For $4$-folds, we find the generic determinantal singularity of this type 
in the list and each of the other simple singularities is adjacent to it.

\section{Stabilizations and the topology of Essential Smoothings} 

In this section we will be concerned with \textit{essential smoothings}
of (essentially isolated) determinantal singularities. This is the generic object 
that a given EIDS can deform to using only 
determinantal deformations, see Definition \ref{def:EssentialSmoothing} below. 
They can be regarded as a generalization of the 
Milnor fiber for isolated hypersurface and complete intersection singularities and 
they coincide with the latter whenever a complete intersection singularity
defined by a regular sequence $F = (f_1,\dots,f_c)$ 
is considered as a determinantal singularity for the matrix $F \in \CC\{x\}^{1\times c}$.
For the construction and properties of Milnor fibers the reader may consult 
\cite[Chapter 6]{Volume1}.

Besides the topology of the essential smoothing itself, we will in the following 
also be interested in the interplay of topological invariants with analytic invariants 
arising from the deformation theory of the singularity. This is motivated from the 
classical results relating the Milnor with the Tjurina number for isolated 
hypersurface and complete intersection singularities. A treatment of this subject 
can be found in \cite[Section 7.2.4]{Volume1}.

To motivate this discussion, recall that the Milnor fiber 
\[
  M_f := B_\varepsilon \cap f^{-1}(\{\delta\}), \quad 
  1\gg \varepsilon \gg |\delta| >0,
\]
of an isolated hypersurface 
singularity $f \colon (\CC^p,0) \to (\CC,0)$ 
is homotopy-equivalent to a bouquet of spheres 
\begin{equation}
  M_f \cong_{\mathrm{ht}} \bigvee_{i=1}^{\mu(f)} S^{p-1}
  \label{eqn:BouquetDecompositionDHS}
\end{equation}
of real dimension $p-1$. The number of these spheres is the 
\textit{Milnor number}\footnote{See e.g. \cite[Definition 6.5.2]{Volume1}.} 
$\mu(f)$ of $f$. This was first described 
by Milnor in \cite{Milnor68} where he also shows that this number 
can be computed as the length 
\begin{equation}
  \mu = \mu(f)= \dim_{\CC} \left(\CC\{x_1,\dots,x_p\}/\left\langle \frac{\partial f}{ \partial x_1},\dots,
  \frac{\partial f}{\partial x_p} \right\rangle \right)
  \label{eqn:MilnorsFormulaForIHS}
\end{equation}
of the so-called \textit{Milnor algebra}, i.e. the 
quotient of $\CC\{x\}$ by the \textit{Jacobian ideal} 
$\Jac(f) = \langle \partial f /\partial x_1,\dots, \partial f /\partial x_p\rangle$,
see \cite{Milnor68} or \cite[Theorem 6.5.3]{Volume1}. The latter space naturally 
coincides with the space $T^1_{\mathcal R}(f)$ of non-trivial unfoldings of 
$f$ up to right- or $\mathcal R$-equivalence, see for instance 
\cite[Definition 3.3]{BallesterosMond20}, and therefore, the Milnor number of $f$ 
is equal to its $\mathcal R$-codimension. This fundamental result already suggests 
a close connection between topological invariants of the smoothing and 
analytic invariants related to unfoldings and deformations. 

If instead of the function $f$ and its unfoldings, one considers the germ of the 
hypersurface $(X,0) = (f^{-1}(\{0\}),0) \subset (\CC^p,0)$ and 
its deformations, one is naturally lead to the 
\textit{Tjurina module} $T^1_{X,0}$ which was already discussed 
earlier and which classifies the 
non-trivial first order deformations of the germ $(X,0)$. 
For isolated hypersurface singularities, this module becomes 
\[
  T^1_{X,0} = \CC\{x_1,\dots,x_p\}/\left\langle 
  \frac{\partial f}{\partial x_1},\dots,\frac{\partial f}{\partial x_p}, f
  \right\rangle.
\]
From this explicit form it is immediately clear that $T^1_{X,0}$ has finite 
dimension over $\CC$ if the Milnor number $\mu(f)<\infty$ is finite. The converse is 
also true but not so obvious, see for instance \cite[Lemma 2.3]{GreuelLossenShustin07}. 
The length of the Tjurina module is called the 
\textit{Tjurina number} of $(X,0)$:
\begin{equation}
  \tau = \dim_\CC T^1_{X,0}
  \label{eqn:DefinitionTjurinaNumber}
\end{equation}
and one has the famous inequality
\begin{equation}
  \mu \geq \tau
  \label{eqn:MuVsTauForIHS}
\end{equation}
for isolated hypersurface singularities.

Suppose that $f$ is weighted homogeneous of some degree $e = \deg_w f$ 
for weights $w_i = \deg_w x_i >0$, i.e. $f$ is a polynomial 
that can be written as a linear combination of monomials 
$x^\alpha = x_1^{\alpha_1} x_2^{\alpha_2} \cdots x_p^{\alpha_p}$ 
which all satisfy 
\[
  \deg_w x^\alpha = 
  \sum_{i=1}^p \deg_w x_i^{\alpha_i} = 
  w_1 \cdot \alpha_1 + w_2 \cdot \alpha_2 + \dots + w_p \cdot \alpha_p = e.
\]
Using the so-called \textit{Euler vector field} 
\[
  \theta_w = \sum_{i=1}^p w_i \cdot x_i \frac{\partial}{\partial x_i}
\]
on $\CC^p$ associated to these weights, it is easy to see that $e \cdot f = \theta_w(f)$. 
So in this case, the function $f$ is already contained in the 
Jacobian ideal of $f$, and the Milnor algebra and the Tjurina 
module are isomorphic. It follows that $\mu = \tau$ if $f$ is quasi-homogeneous. 

The converse implication was proven by Saito in \cite{Saito71} in '71: 
Starting from $f$ with $\mu = \tau$ he establishes the existence of a 
change of coordinates of $(\CC^p,0)$ such that $f$ is quasi-homogeneous 
for some weights in this new coordinate system. 

The colloquial form of this result
\begin{equation}
  \textnormal{``$\mu \geq \tau$ with equality iff $f$ is quasi-homogeneous''}
  \label{eqn:MuEqualsTauColloquial}
\end{equation}
has been abundant in complex analytic singularity theory ever since and much effort has been 
invested in order to generalize it to isolated complete intersection singularities 
defined by maps $f \colon (\CC^p,0) \to (\CC^c,0)$. 

The generalization of 
the bouquet decomposition (\ref{eqn:BouquetDecompositionDHS}) for ICIS 
is due to Hamm \cite{Hamm72}, which settles the definition of a Milnor number in this 
context. The Tjurina number is given by the length of $T^1_{X,0}$ as described 
in Section \ref{sec:DeformationsOfICIS}. With these definitions at hand, the journey went 
on for more than thirty years:
\begin{itemize}
  \item The equality $\mu = \tau$ has been established in '80 for quasi-homogeneous ICIS of 
    positive dimension\footnote{Note that $\mu=\tau$ does not hold for quasi-homogeneous 
      ICIS of dimension zero,
      see for instance Table \ref{tab:GiustiFatPoints}.} by 
    Greuel in \cite{Greuel80}.
  \item The general inequality $\mu \geq \tau$ for ICIS of positive dimension 
    was shown in '85 by Looijenga and 
    Steenbrink in \cite{LooijengaSteenbrink85}.
  \item In the same year, the full statement (\ref{eqn:MuEqualsTauColloquial}) has been proved for 
    Gorenstein curve singularities by Greuel, Martin, and Pfister in 
    \cite{GreuelMartinPfister85},
  \item while Wahl proved  (\ref{eqn:MuEqualsTauColloquial}) for ICIS of 
    codimension $2$ in \cite{Wahl85}.
  \item Finally, the part ``$\mu = \tau$ implies quasi-homogeneity'' was 
    established in '02 by Vosegaard \cite{Vosegaard02}.
\end{itemize}
After the question as to whether or not (\ref{eqn:MuEqualsTauColloquial}) holds 
has been established for isolated complete intersection singularities, it 
seems natural to ask to which extent this result can be generalized to arbitrary 
EIDS and we will report on what is known in this regard. For now, let us mention 
that in general, this question is wide open. For instance, one has the following conjecture 
by Wahl in \cite{Wahl16}:

\begin{conjecture}
  \label{con:WahlsConjecture}
  Let $(X,0) \subset (\CC^4,0)$ be a normal surface singularity which is not a complete 
  intersection. Then $\mu \geq \tau -1$ with equality if and only if $(X,0)$ is 
  quasi-homogeneous.
\end{conjecture}

These singularities fall into the category of isolated Cohen-Macaulay codimension 
$2$ singularities discussed earlier and 
the precise definitions of the Milnor and Tjurina number will 
be given below. That quasi-homogeneity of $(X,0)$ implies equality was already 
shown by Wahl in \cite{Wahl16},
but the converse implication is not settled as 
of this writing\footnote{A counterexample given by the first named author in \cite{FruehbisKrueger18} 
turned out to be wrong.}.

\subsection{Construction of Essential Smoothings}

We briefly describe the construction of the essential smoothing of 
an EIDS, parallel to the well known Milnor-L\^e-fibration for isolated 
hypersurface and complete intersection singularities, cf. 
\cite{Volume1}. Another description based on the transformation into a 
complete intersection on a singular ambient space, 
cf. Remark \ref{rem:GraphTransformationAndRelativeCompleteIntersections},
will be given later in Corollary \ref{cor:EssentialSmoothingAsICISMilnorFiber}. 

\medskip

Let $(X_A^s,0) \subset (\CC^p,0)$ be an EIDS defined by a finitely 
$\GL$-determined matrix $A \in \CC\{x\}^{m\times n}$. Choose 
a miniversal unfolding $\mathbf A(x,t)$ on 
$\tau = \dim_\CC T^1_\GL(A)$ parameters $t_1,\dots,t_{\tau}$ and 
let 
\[
  \mathbf A \colon U \times T \to \CC^{m\times n}
\]
be a representative thereof defined on some product of open neighborhoods 
of the origin $U \subset \CC^p$ and $T \subset \CC^\tau$. As usual, let 
$\mathcal X_{\mathbf A}^s = \mathbf A^{-1}(M_{m,n}^s) \subset U\times T$ 
be the total space of the induced deformation 
$(X_A^s,0) \hookrightarrow (\mathcal X_{\mathbf A}^{s},0) \overset{\pi}{\longrightarrow} (T,0)$
of the determinantal singularity. 

Recall that, according to Lemma \ref{lem:WhitneyStratificationOfAnEIDSAndSmoothability}, 
any EIDS $(X_A^s,0)$ has a canonical Whitney stratification by the strata 
$V_{A}^r = A^{-1}(V_{m,n}^r)$. This allows us to choose a \textit{Milnor sphere}: 
For $\varepsilon_0>0$ sufficiently small, the intersection of the sphere 
$S_\varepsilon = \partial B_\varepsilon \subset U$ with the $X_A^s$ is transverse 
for every $\varepsilon_0 \geq \varepsilon >0$, see e.g. \cite[Theorem 6.10.1]{Volume1}. 
The \textit{real link} of $(X_A^s,0)$ can then be defined as the transverse intersection 
\begin{equation}
  \mathcal K_A^s := S_\varepsilon \cap X_A^s.
  \label{eqn:DefinitionRealLinkForEIDS}
\end{equation}
Note that, since $(X_A^s,0)$ has non-isolated singularities in general, the 
real link will also be singular, but endowed with a canonical Whitney regular 
stratification. Using Thom's first isotopy 
lemma\footnote{See for instance \cite{GoreskyMacPherson88}.} it is easy to 
see that due to the various transversalities, small determinantal deformations 
of $X_A^s$ do not change the compact stratified space $\mathcal K_A^s$ 
up to homeomorphism. Hence, after shrinking $T$ to some small 
disk $D_\delta \subset \CC^\tau$ if necessary, 
the restriction of the projection 
\begin{equation}
  \pi \colon (S_\varepsilon \times T) \cap \mathcal X_A^s \cong_{\mathrm{homeo}} 
  \mathcal K_A^s \times T \to T
  \label{eqn:TrivialFibrationOnTheBoundary}
\end{equation}
to the parameter space of the miniversal unfolding of $A$ 
is a trivial topological fiber bundle. Here one makes 
essential use of the fact that $\mathcal K_A^s$ is compact. 

Due to the characterization of transversality 
given in Proposition \ref{prp:AlgebraicTransversality}, 
the matrix discriminant $\Delta_A \subset T$ 
(Definition \ref{def:DeterminantalDiscriminant}) 
consists of those parameters $t$ for which 
$A_t \colon U \to \CC^{m\times n}$ is not transversal to the 
rank stratification. It is easy to see that for $t \notin \Delta_A$, the projection 
\begin{equation}
  \pi \colon (B_\varepsilon \times (T\setminus\Delta_A)) \cap \mathcal X_A^s 
  \to T \setminus \Delta_A
  \label{eqn:MilnorFibrationForEssentialSmoothings}
\end{equation}
is a stratified submersion along every fiber 
$\pi^{-1}(\{t\}) = A_t^{-1}(M_{m,n}^s) \subset U$. 
By choice of the ball $B_\varepsilon$, these fibers 
are again compact and canonically stratified, so that 
$\pi$ is proper. 
Another application of Thom's first isotopy lemma can now be made 
in order to show that (\ref{eqn:MilnorFibrationForEssentialSmoothings}) 
is a topological fiber bundle extending the fibration on the boundary 
(\ref{eqn:TrivialFibrationOnTheBoundary}) to the interior over $T\setminus \Delta_A$. 

\begin{definition}
  \label{def:EssentialSmoothing}
  The essential smoothing of an EIDS $(X_A^s,0)\subset (\CC^p,0)$ defined by a matrix 
  $A \in \CC\{x_1,\dots,x_p\}^{m\times n}$ is the generic fiber 
  \[
    M_A^s:= A_t^{-1}(M_{m,n}^s)
  \]
  of (\ref{eqn:MilnorFibrationForEssentialSmoothings}) over the base 
  of the miniversal unfolding of $A$ for $t \notin \Delta_A$ 
  outside the matrix discriminant.
\end{definition}


\begin{example}[Smoothing of a space curve]
  \label{exp:SmoothingOfSpaceCurve}
  We return to the study of the semi-universal deformation 
  of the three coordinate axis in $(\CC^3,0)$ from Example 
  \ref{exp:ThreeCoordinateAxisSemiUniversalUnfolding}:
  \[
    \mathbf A(x,y,z;t_1,t_2,t_3) = 
    \begin{pmatrix} 
      x & 0 & z \cr
      0 & y & z 
    \end{pmatrix} + 
    \begin{pmatrix}
      0 & t_1 & 0 \\
      t_2 & 0 & t_3
    \end{pmatrix}.
  \]
  It had already been discussed in Example \ref{exp:DiscriminantThreeCoordinateAxis} 
  that the simultaneous perturbation by $t_1=t_2=t_3 = t$ leads to a smoothing 
  $M_A^2$ of $(X_A^2,0)$. Note that due to the homogeneity of the singularity, we may choose the 
  Milnor ball $B_\varepsilon$ arbitrarily large so that we can consider the whole affine 
  varieties as suitable representatives of the singularity and its fibers in a deformation.
  
  In order to determine the topological type of $M_A^2$ 
  we can exploit the adjacencies of $(X_A^2,0)$: First we deform only along the 
  $t_1$-axis as in Example \ref{exp:DeformationOfThreeCoordinateAxisFlatness} to 
  obtain a configuration of lines consisting of $L_x'$, $L_y$, and $L_z$, with 
  $L'_x$ and $L_z$ meeting transversally in the point $(0,0,t_1)$ and 
  $L_y$ and $L_z$ at the origin, respectively.
  Since all these fibers for arbitrary $t_1$ sit over the discriminant $\Delta_A$ in a 
  semi-universal deformation of $(X_A^2,0)$, the (global) smoothing of this configuration 
  must be diffeomorphic to $M_A^2$. 
  \begin{figure}[h]
    \centering
    \includegraphics[scale=0.3,clip=true,trim=2cm 3cm 1cm 2cm]{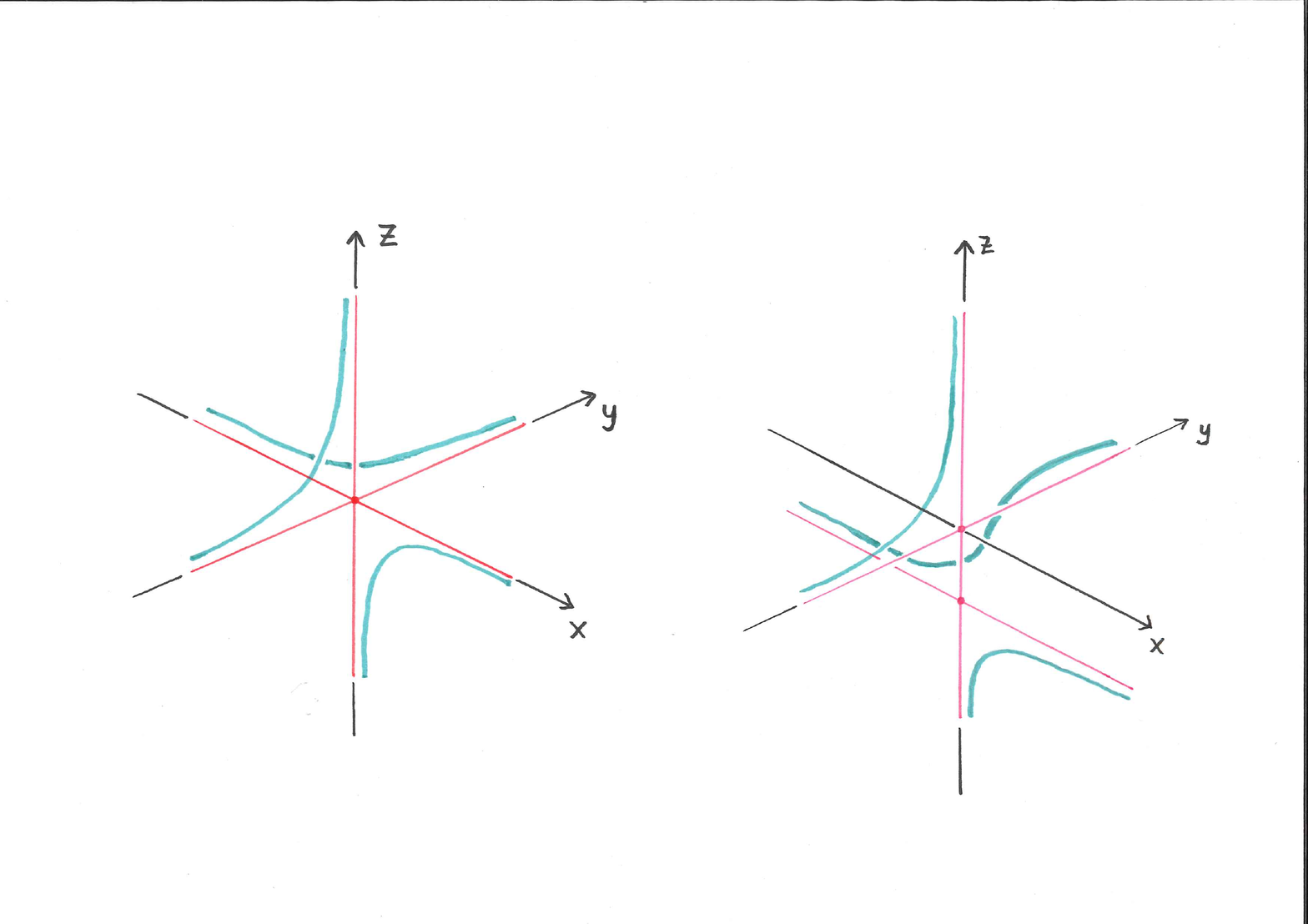}
    \caption{Real picture of a deformation of the three coordinate axis in $(\CC^3,0)$ 
      and its smoothing, passing 
    through the adjacency with two $A_1$-hypersurface singularities.}
    \label{fig:SmoothingThreeCoordinateAxisWithAdjacencies}
  \end{figure}

  The topology of the local smoothings of the $A_1$-singularities 
  at the intersection points is known: Over the complex numbers, 
  a double cone is replaced by a tube bounding the two circles in its boundary. 
  Now it is easy to see that in fact 
  \[
    M_A^2 \cong_{\mathrm{ht}} S^1 \vee S^1
  \]
  is homotopy equivalent to a bouquet of two circles. In parallel to 
  the definition of Milnor numbers for IHS and ICIS, one would say that 
  the Milnor number for $(X_A^2,0)$ is two in this case.
  \begin{figure}[h]
    \centering
    \includegraphics[scale=0.2,clip=true,trim=0 3cm 12cm 2cm]{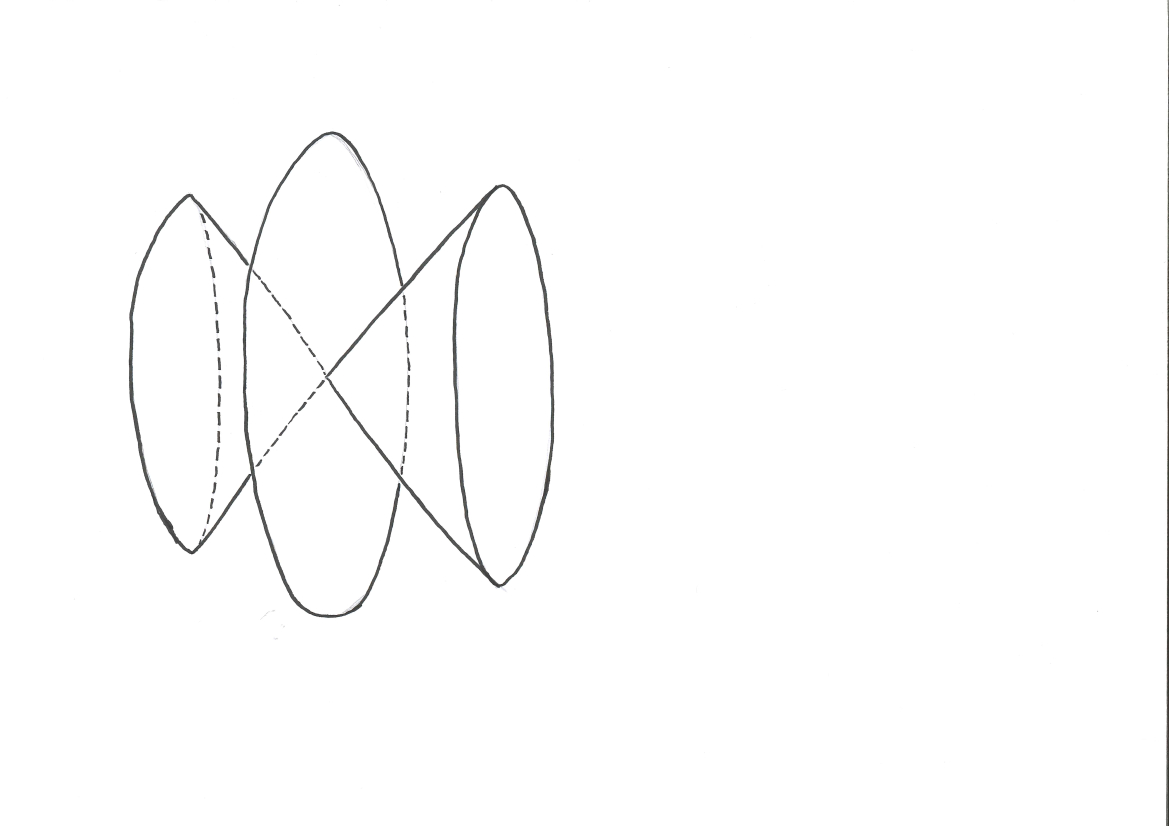}
    \includegraphics[scale=0.2,clip=true,trim=0 3cm 12cm 2cm]{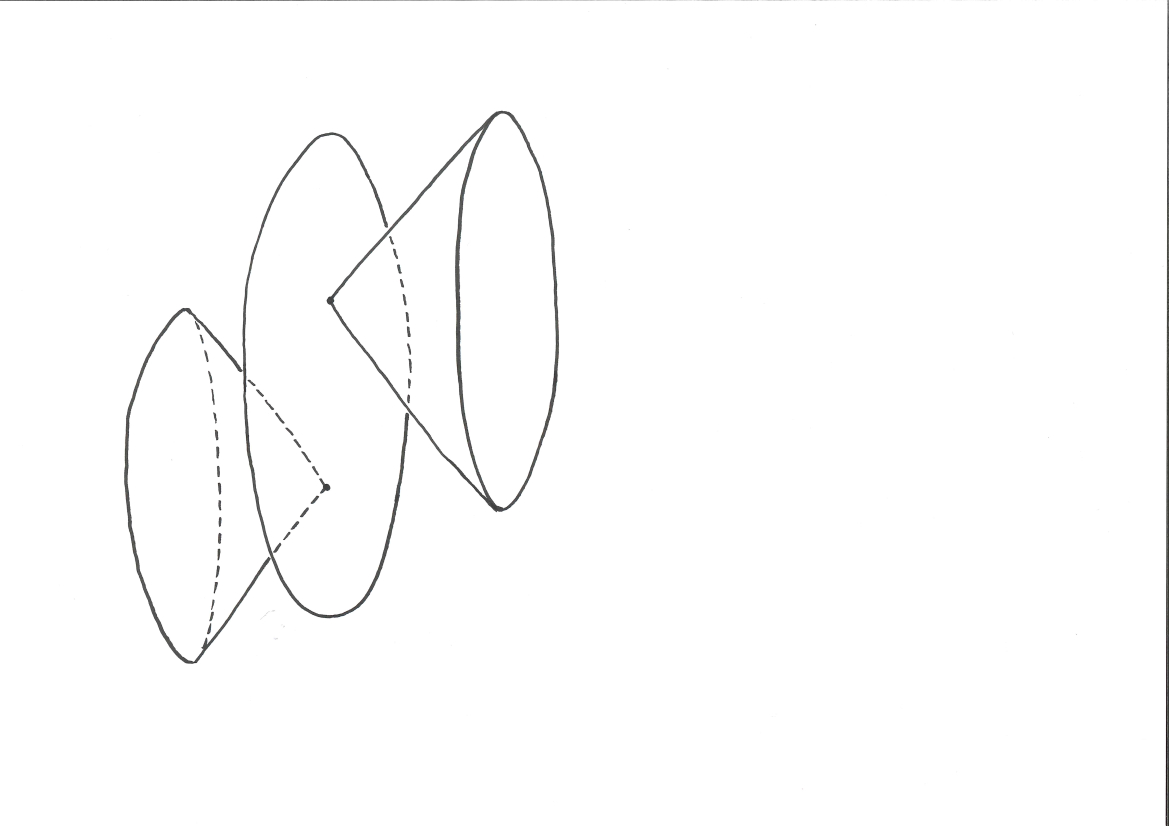}
    \includegraphics[scale=0.2,clip=true,trim=0 3cm 12cm 2cm]{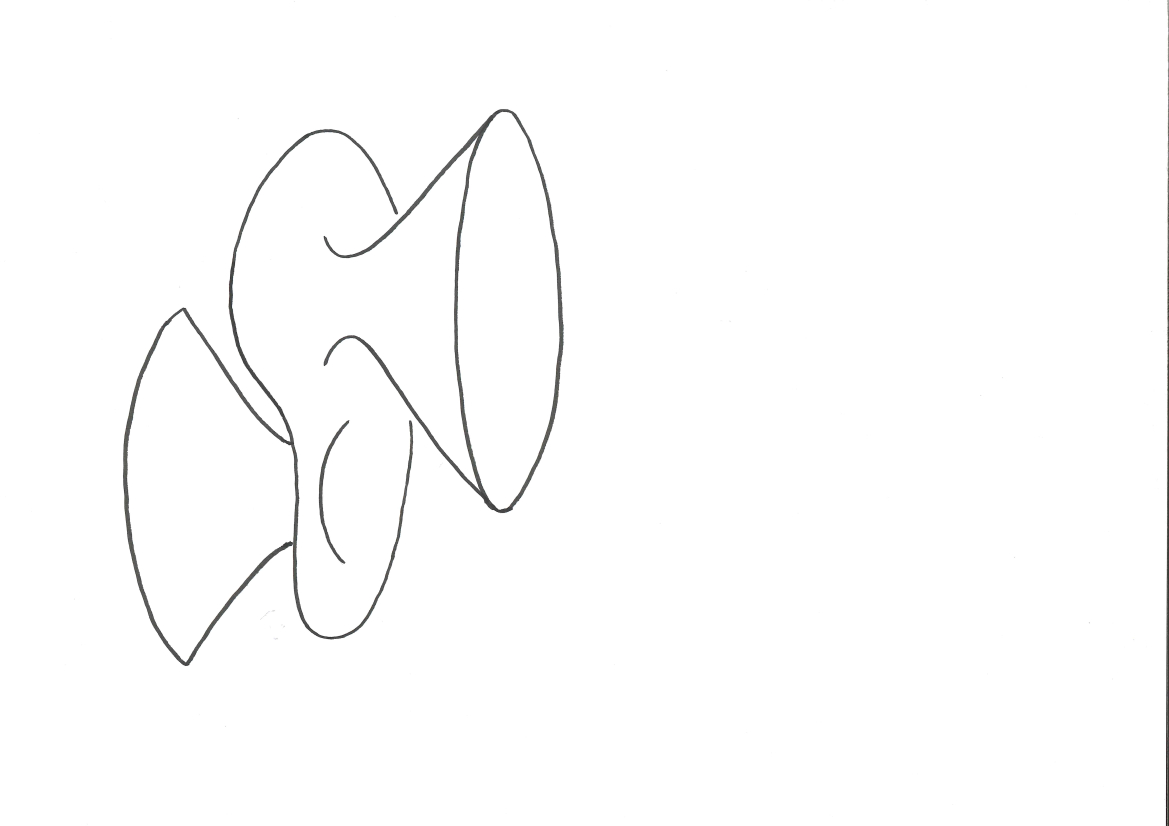}
    \caption{The change in topology for a smoothing of the 
      three coordinate axis in $(\CC^3,0)$ passing 
    through the adjacency with two $A_1$-hypersurface singularities.}
    \label{fig:SmoothingThreeCoordinateAxisWithAdjacencies}
  \end{figure}
\end{example}

\begin{remark}
  \label{rem:EssentialSmoothingsOfNonIsolatedSings}
  Essential smoothings can also be defined for determinantal singularities 
  $(X_A^s,0) \subset (\CC^p,0)$ which are not EIDS; i.e. those 
  for which the defining matrix $A\in \CC\{x\}^{m\times n}$ is not 
  finitely $\GL$-determined. In the parallelism with hypersurface 
  and complete intersection singularities, these correspond to non-isolated 
  singularities and, as in the classical case, the essential smoothings 
  are not uniquely determined by the singularity $(X_A^s,0)$ itself anymore, 
  but they can differ depending on the underlying unfolding of $A$. 

  To give a definition, let $\mathbf A \colon (\CC^p,0) \times (\CC,0) \to 
  (\CC^{m\times n},0)$ be a $1$-parameter unfolding of $A$. We say that 
  $\mathbf A$ is a \textit{stabilization} of $A$ if for some representative 
  \[
    \mathbf A \colon U \times T \to \CC^{m\times n}
  \]
  one has that $A_t \colon U \to \CC^{m\times n}$ is transversal to the 
  rank stratification for every $t\in T\setminus\{0\}$. The existence 
  of such stabilizations can for example be derived from  
  \cite[Theorem 2.2 and Theorem 3.1]{Trivedi13}. 
  In particular, the set of constant matrices $C \in \CC^{m\times n}$ 
  for which $\mathbf A(x,t) = A(x) + t \cdot C$ is a stabilization of 
  $A$, is dense in $\CC^{m\times n}$. 

  Given any suitable representative of a stabilization $\mathbf A$ of $A$ 
  as above,
  we can consider the total space of the induced deformation 
  $\mathcal X_{\mathbf A}^s = \mathbf A^{-1}(M_{m,n}^s) \subset U \times T$ together 
  with its  projection 
  $\pi \colon \mathcal X_{\mathbf A}^s \to T$ 
  to the parameter space.
  General fibration theorems\footnote{See \cite{Le77}, cf. also \cite[Theorem 6.10.3]{Volume1}.}
  can be used to show that for a sufficiently 
  small ball $B_\varepsilon \subset U$ of radius $\varepsilon>0$ around 
  the origin $0 \in U$ and some subsequently chosen, sufficiently 
  small disc $D_\delta \subset T$, 
  the projection 
  \begin{equation}
    \pi \colon B_\varepsilon \times (D_\delta\setminus \{0\}) \cap \mathcal X_{\mathbf A}^s \to 
    D_\delta\setminus \{0\}
    \label{eqn:FibrationForEssentialSmoothingForNonIsolatedSings}
  \end{equation}
  is a topological fiber bundle with fiber 
  \begin{equation}
    M_{\mathbf A}^s := B_\varepsilon \cap A_t^{-1}(M_{m,n}^s).
    \label{eqn:EssentialSmoothingForNonIsolatedSings}
  \end{equation}
  This is the \textit{essential smoothing} of $(X_A^s,0)$ defined 
  by the stabilization $\mathbf A$. 
  
  It is easy to see using the 
  properties of miniversal unfoldings that this notion of essential smoothing 
  coincides with the previous one given for EIDS in case $A$ is finitely 
  $\GL$-determined.

  Note that the fibration (\ref{eqn:FibrationForEssentialSmoothingForNonIsolatedSings}) 
  always exists for $1$-parameter unfoldings, regardless of whether or not $A_t$ 
  is a stabilization for $t \neq 0$. However, we shall refer to the fiber 
  of this fibration as an essential smoothing of $(X_A^s,0)$ only if this is the case.
\end{remark}

\begin{example}
  \label{exp:NonEquivalentSmoothingsForNonEIDS}
  Consider the space curve $(X_A^2,0) \subset (\CC^3,0)$ given by the matrix 
  \[
    A = 
    \begin{pmatrix}
      x & 0 & z \\
      0 & y & z^2
    \end{pmatrix}.
  \]
  Set theoretically this coincides with the union of the three coordinate 
  axis discussed in several previous examples. But a primary decomposition of the 
  ideal $I = \langle A^{\wedge 2} \rangle$ reveals that 
  \[
    I = \langle y, z^2 \rangle \cap 
    \langle x,z \rangle \cap \langle x,y \rangle
  \]
  so that $X_A^2$ consists of the double $x$-axis in the 
  $x$-$z$-plane and 
  the $y$- and the $z$-axis.
  Consequently, the module $T^1_\GL(A)$ is indeed supported 
  along the whole $x$-axis. 

  We describe two distinct unfoldings leading to topologically different smoothings 
  of $(X_A^2,0)$. 
  The first one is given by 
  \[
    \mathbf A(x,y,z;u) = 
    \begin{pmatrix}
      x & 0 & z \\
      0 & y & z^2-u^2
    \end{pmatrix}.
  \]
  For $u\neq 0$ the variety defined by the ideal $\langle A_u^{\wedge 2}\rangle$ 
  in $\CC^3$ consists of \textit{four} lines with two of them being the now split 
  parallels to the $x$-axis passing through the points $(0,0,\pm u)$.
  A global smoothing of the three singular points of this variety is homotopy 
  equivalent to a bouquet of three real spheres:
  \[
    M_{\mathbf A}^2 \cong_{\mathrm{ht}} S^1\vee S^1 \vee S^1.
  \]

  \begin{figure}[h]
    \centering
    \includegraphics[scale=0.35,clip=true,trim=0.5cm 1cm 13cm 1cm]{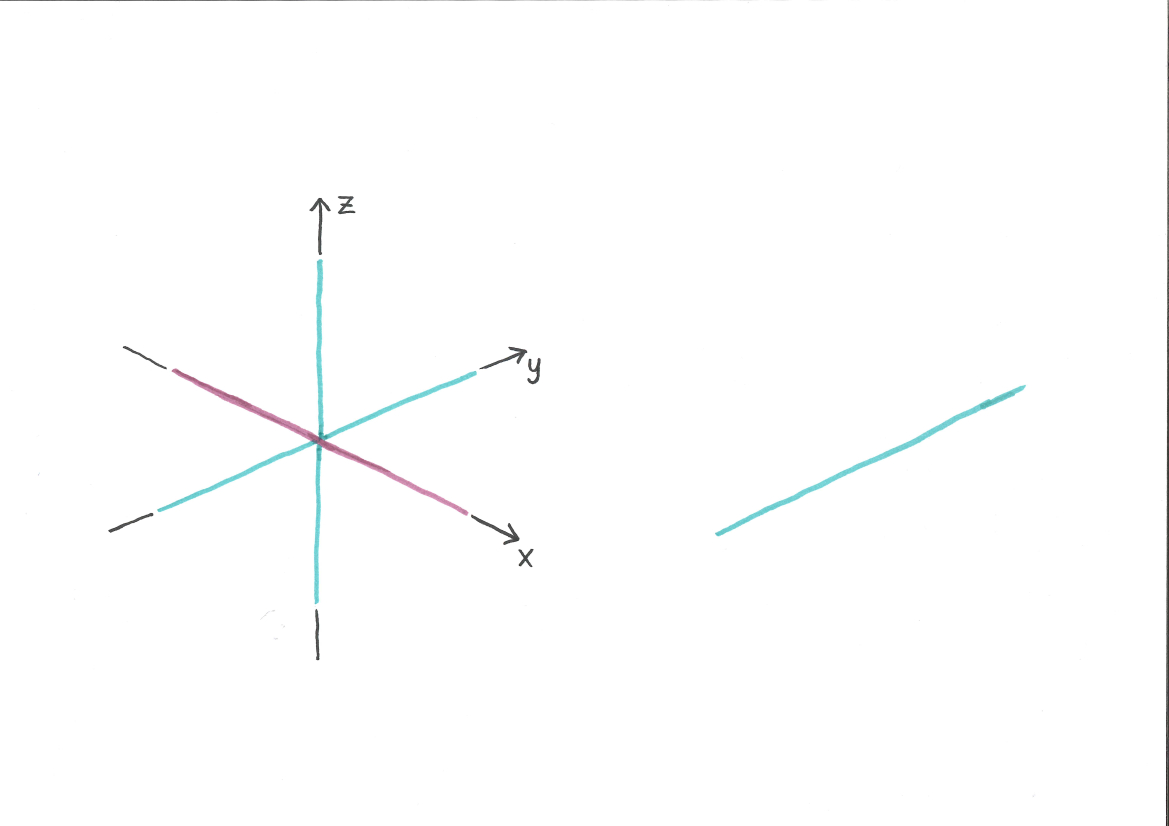}
    \includegraphics[scale=0.35,clip=true,trim=0.5cm 1cm 13cm 1cm]{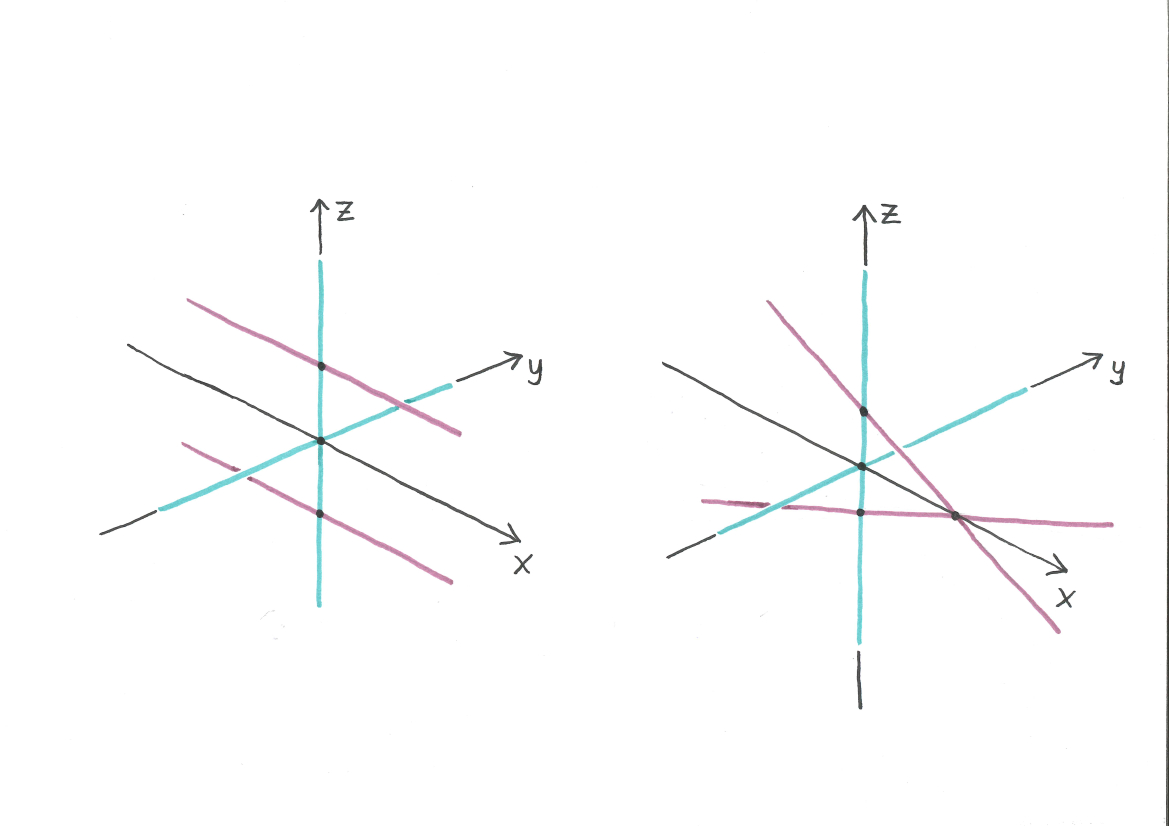}
    \caption{The first deformation of a non-reduced space curve singularity (purple and green) 
    with non-isolated singular locus (purple)}
    \label{fig:NonIsolatedSpaceCurveFirstDeformation}
  \end{figure}

  The other deformation that we wish to consider is induced by the unfolding
  \[
    \mathbf B(x,y,z;v) = 
    \begin{pmatrix}
      x & 0 & z \\
      0 & y & z^2-v^2(x^2-2z-v^2) 
    \end{pmatrix}.
  \]
  For this deformation, the variety defined by $\langle B_v^{\wedge 2} \rangle$ for 
  $v\neq 0$ consists again of four lines, only that this time the double $x$-axis does 
  not split into parallels, but opens up like a scissor which is pulled backwards at 
  the same time:
  \[
    \langle B_v^{\wedge 2} \rangle = 
    \langle y, xv-(z-v^2) \rangle \cap
    \langle y, xv+(z-v^2) \rangle \cap
    \langle x,z \rangle \cap \langle x,y \rangle.
  \]
  \begin{figure}[h]
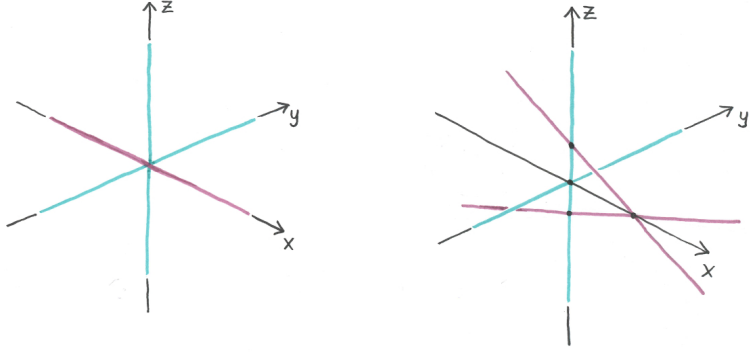

    \centering
    \includegraphics[scale=0.35,clip=true,trim=0.5cm 1cm 13cm 1cm]{Pictures/Non-isolated_space_curve_singularity_page_2.jpg}
    \includegraphics[scale=0.35,clip=true,trim=15cm 1cm 1cm 1cm]{Pictures/Non-isolated_space_curve_singularity_page_1.jpg}
    \caption{The second deformation of a the same non-reduced space curve singularity.
    }
    \label{fig:NonIsolatedSpaceCurveFirstDeformation}
  \end{figure}
  Thus for $v \neq 0$ the four lines meet transversally in pairs at four 
  points in total. The global smoothing of this variety is then 
  homotopy equivalent to a bouquet of \textit{five} real spheres:
  \[
    M_{\mathbf B}^2 \cong_{\mathrm{ht}} S^1 \vee S^1 \vee S^1 \vee S^1 \vee S^1.
  \]
  Choosing appropriate combinations of either one of these two deformations 
  with their respective global smoothings, it is easy to see that also the original 
  singularity $(X_A^2,0)$ can be deformed directly into $M_{\mathbf A}^2$, 
  but also into $M_{\mathbf B}^2$. Therefore, a smoothing is not unique 
  for this non-essentially isolated determinantal singularity.
\end{example}

\subsection{Determinantal hypersurfaces}
\label{sec:DeterminantalHypersurfaces}

Determinantal hypersurfaces have already appeared in the context of the classification 
of simple square matrices by Bruce and Tari \ref{sec:SimpleSquareMatrices}. 
For this section we will need this notion to also comprise the symmetric and 
the skew symmetric cases.

\begin{definition}
  \label{def:DeterminantalHypersurface}
  Let $A \in \CC\{x_1,\dots,x_p\}^{m\times m}$ be a non-constant matrix, 
  either arbitrary square, symmetric, 
  or skew symmetric. Then the singularity $(X_A^m,0) \subset (\CC^p,0)$ defined by 
  the equation $f = \det A$ in the first two, or by $f = \Pf A$ in the third case, 
  is called a \textit{determinantal hypersurface singularity} with defining equation 
  $f$.
\end{definition}
Depending on the case, we will consider the matrix $A$ up to either $\mathcal G$-, 
$\mathcal G_{\sym}$-, or $\mathcal G_\sk$-equivalence. Since the bound on the rank 
$m$ is always equal to the size of the defining matrix $A$ for determinantal hypersurface 
singularities, we will usually not mention $m$ explicitly throughout this section 
and omit it from our notation. Note that, moreover, $m$ is necessarily even in 
the skew-symmetric case since $\Pf(A) = 0$ for matrices of odd size.

\subsubsection{The singular Milnor fiber}
For determinantal hypersurfaces we always have two different deformation theories at 
hand: The deformations arising from unfoldings of $f$ and the determinantal 
deformations induced from unfoldings of the matrix $A$.
Correspondingly, there are two notions of ``Milnor fiber'' in this setup which 
differ in general. The first one is the classical Milnor fiber 
\begin{equation}
  M_f := B_\varepsilon(0) \cap f^{-1}(\{\delta\}),
  \label{eqn:MilnorFiberForHypersurfaces}
\end{equation}
$1\gg \varepsilon \gg |\delta| >0$ (see \cite{Milnor68}, cf. also \cite[Chapter 6]{Volume1}),
of the hypersurface singularity determined by $f$. 
The other one is the \textit{essential smoothing} of $(X_A,0)$ 
given by 
\begin{equation}
  M_A := B_\varepsilon(0) \cap A_t^{-1}(M_{m,m}^m) = B_\varepsilon(0) \cap \{ \det A_t = 0 \},
  \label{eqn:DeterminantalMilnorFiberHypersurfaceCase}
\end{equation}
$1\gg \varepsilon \gg |t| >0$ with the appropriate substitution for the 
skew-symmetric case.
Depending on how large the dimension $p$ is compared to the codimension of 
the singular locus of the set of degenerate matrices, 
these spaces either coincide or differ:

\begin{lemma}
  Let $(X_A,0) \subset (\CC^p,0)$ be a determinantal hypersurface defined by 
  a finitely determined square matrix $A \in \CC\{x_1,\dots,x_p\}^{m\times m}$ 
  which is either square, symmetric, or skew-symmetric. 
  Depending on the case, let $c'=4, 3$, or $6$ be the 
  codimension of the singular locus of the respective set 
  of degenerate matrices.
  
  \medskip

  \noindent
  When $p<c'$, the essential smoothing is in fact smooth and the 
  manifolds $M_A \cong_{\mathrm{diff}} M_f$ are diffeomorphic.

  \medskip 

  \noindent
  When $p\geq c'$, the essential smoothing 
  $M_A$ is singular with singular locus of dimension $m^2-4$ and $M_f$ is diffeomorphic 
  to a global smoothing of $M_A$. 
  \label{lem:MilnorFibersOfDeterminantalHypersurfaces}
\end{lemma}
\begin{proof}
  In the first case when $p<c'$ this is clear since the smoothing 
  $M_f$ is unique. 
  The second case must be split into two further subcases, namely 
  $p=c'$ so that $f = \det A$ has \textit{isolated} singularity, 
  and $p>c'$ in which case $f$ has \textit{non-isolated singularities}.

  When $p=c'$, there exists again a semi-universal deformation of the 
  hypersurface germ $(X_A,0) = (\{f=0\},0)$ over some base 
  $(\CC^{\tau},0)$ with $\tau = \tau(X_A,0)$ the Tjurina number 
  of $(X_A,0)$, and a comparison map 
  $\Phi \colon (\CC^{\tau_\GL(A)},0) \to (\CC^{\tau(f)},0)$ 
  as in the diagram (\ref{eqn:DeterminantalDeformationsVsSemiUniversalDeformations}) 
  where $(\CC^{\tau_\GL(A)},0)$ is the parameter space of a 
  $\GL$-miniversal unfolding $\mathbf A$ of $A$. 
  But the image of $\Phi$ must be contained in the discriminant 
  $(\Delta_f,0) \subset (\CC^{\tau(f)},0)$ of $f$ since $(X_A,0)$ 
  does not admit any determinantal smoothing.
  Choosing an appropriate representative
  $\mathbf A \colon U \times T \to \CC^{m\times m}$ of the 
  $\GL$-miniversal unfolding of $A$, it is easy to 
  see that the generic fiber $M_A \subset U$ has only 
  isolated singularities and hence posesses a unique smoothing. 
  Since we can think of the fiber $M_A$ as a fiber in the 
  semi-universal deformation of $(X_A,0)$ by virtue of 
  Diagram (\ref{eqn:DeterminantalDeformationsVsSemiUniversalDeformations}), 
  this smoothing of $M_A$ must coincide with $M_f$ as this is the 
  only smooth nearby fiber of $M_A$ in the semi-universal 
  deformation $\pi \colon \mathcal X \to T'$ of $(X_A,0)$.

  In the case $p>c'$ there does not exist as semi-universal 
  deformation of $(X_A,0)$ anymore but one has the function 
  $F = \det \mathbf A$ for a $\GL$-miniversal unfolding 
  $\mathbf A$ of $A$ which assigns a canonical smoothing 
  \[
    B_\varepsilon \cap \left( \det A_t \right)^{-1}(\{\delta\}), \quad 
    \varepsilon \gg |t| \gg |\delta|>0
  \]
  to every fiber $X_{\mathbf A}(t) = \{ \det A_t = 0\}$ for 
  $t \in \CC^{\tau_{\GL}(A)}$ sufficiently small.
  Details for the treatment of this setup can for example be found in
  \cite{Siersma91_2}, see also \cite[Proposition 4.5]{DamonMond91}.
\end{proof}

\subsubsection{The smoothable case}

The first case $p<c'$ has been studied by Goryunov and Mond in \cite{GoryunovMond05}. 
In this case the function $f \colon (\CC^p,0) \to (\CC,0)$ defines an 
isolated hypersurface singularity in the classical sense. 
Hence, on the topological side, the theory is in principal covered by Milnor's 
results \cite{Milnor68}: 
The Milnor fiber 
\[
  M_f \cong_{\mathrm{ht}} \bigvee_{i=1}^{\mu(f)} S^{p-1}
\]
is homotopy equivalent to a bouquet of $\mu(f)$ spheres of real dimension $p-1$ 
where $\mu(f)$ is the Milnor number of the singularity defined by $f$. 

For an arbitrary isolated hypersurface singularity, one would compare the 
Milnor with the Tjurina number of $f$. However, the classification of 
simple singularities of the square matrices 
in Table \ref{tab:NormalFormsSquareMatricesSize2In3Space} and of the
symmetric matrices in Table \ref{tab:SquareSymm3p2} 
shows an equality $\mu(f) = \tau_\GL(A)$ (resp. $\mu(f) = \tau_\GL^\sym(A)$) 
of the Milnor numbers with the $\GL$-Tjurina numbers of 
the defining matrices. This peculiarity motivated Goryunov and Mond to closer 
inspections and they revealed in \cite[Corollary 4.4]{GoryunovMond05} 
the following, more general phenomenon:

\begin{theorem}
  Let $(X_A,0) \subset (\CC^p,0)$ be a smoothable determinantal hypersurface 
  singularity for some matrix $A \in \CC\{x\}^{m\times m}$, either 
  symmetric, arbitray square, or skew symmetric with $m$ even, with defining equation  
  $f$. Then
  \begin{eqnarray*}
    \tau_{\SL}^{\sym}(A) = \mu( f ) & \quad & \textnormal{ if $A$ is symmetric and $p = 2$;}\\
    \tau_{\SL}^{\sq}(A) = \mu( f ) & \quad & \textnormal{ if $A$ is arbitrary square and $p=3$;}\\
    \tau_{\SL}^{\sk}(A) = \mu( f ) & \quad & \textnormal{ if $A$ is skew-symmetric and $p=5$.}
  \end{eqnarray*}
  \label{thm:GoryunovMondSmoothableHypersurfaces}
\end{theorem}

For determinantal hypersurface singularities defined by quasi-homogeneous matrices
(such as all the simple ones listed above),
the $\GL$- and the $\SL$-Tjurina numbers coincide, see e.g.
\cite[Proposition 1.1]{Goryunov21}. In the general case, however, the 
$\SL$-Tjurina numbers seem to be the prefered choice for comparison 
with topological invariants. 

\begin{remark}
  \label{rem:GoryunovMondTorsionCorrections}

  It is remarkable that the proof of Theorem \ref{thm:GoryunovMondSmoothableHypersurfaces} 
  in \cite{GoryunovMond05} 
  relies on an argument which is similar to the use of generic perfection introduced 
  in Section \ref{sec:DeterminantalIdeals}, but where the bound on the 
  grade in Theorem \ref{thm:InheritanceOfPerfection} is \textit{not} attained. 

  More precisely, Goryunov and Mond
  consider the complex $K(m,m,m-1)$ from (\ref{eqn:GenericDeterminantalComplex}), 
  (respectively $K^\sym(m,m-1)$, or $K^\sk\left(m,\frac{m}{2}-1\right)$) over the ring $\CC\{y\}$ 
  of convergent power series at the origin $0 \in \CC^{m\times m}$ in the target 
  of $A$. These provide a resolution of the Jacobian ideal of the function $\det$ 
  (resp. $\Pf$) 
  which is generated by
  \[
    \frac{\partial \det}{\partial y_{i,j}} = (-1)^{i+j} Y^{\wedge m-1}_{\hat i, \hat j}
  \]
  where $\hat i$ denotes the multiindex $(1,2,\dots,\hat i,\dots, m-1,m)$ obtained by deleting 
  the $i$-th entry. The pullback $A^* K(m,m,m-1)$ of these complexes 
  admit a morphism 
  \[
    \Phi \colon 
    \Kosz\left( \frac{\partial f}{\partial x_1},\dots,\frac{\partial f}{\partial x_p}; \CC\{x\} \right)
    \to 
    A^* K(m,m,m-1)
  \]
  lifting the natural map on the generators 
  \[
    \frac{\partial f}{\partial x_k} = \frac{\partial \det \circ A}{\partial x_k} \mapsto 
    \sum_{(i,j)} \frac{\partial \det}{\partial y_{i,j}} \frac{\partial A_{i,j}}{\partial x_k}
    = \sum_{(i,j)} (-1)^{i+j} \frac{\partial A_{i,j}}{\partial x_k} A^{\wedge m-1}_{\hat i, \hat j}
  \]
  induced by the chain rule.
  Completing this to a short exact sequence of complexes via the mapping cone of $\Phi$, 
  the associated long exact sequence reads 
  \begin{eqnarray*}
    \dots \to
    H_1\left( \Kosz\left(\frac{\partial f}{\partial x}; \CC\{x\} \right)\right) &\to&
    H_1\left( A^* K(m,m,m-1) \right) \to 
    T^1_\SL(A) \to \\
    &\to & \CC\{x\}/\left\langle \frac{\partial f}{\partial x} \right\rangle \to 
    H_0( A^* K(m,m,m-1)) \to 
    0,
  \end{eqnarray*}
  cf. \cite[Theorem 1.2]{GoryunovMond05} and \cite[Corollary 1.3]{GoryunovMond05}. 
  In case $f$ has isolated singularity at the origin, the Koszul complex 
  $\Kosz\left( \frac{\partial f}{ \partial x}; \CC\{x\} \right)$ is exact 
  in degrees $>0$ so that this restricts to an exact four-term complex. 
  Writing $\beta_i$ for the length of the module
  \[
    H_i(A^*K(m,m,m-1)) = 
    \Tor_i^{\CC\{y\}}\left( \CC\{x\}, \CC\{y\}/\langle Y^{\wedge m-1} \rangle \right)
  \]
  Goryunov and Mond obtain a more general formula 
  \begin{equation}
    \tau_{\SL}(A) = \mu(f) - \beta_0 + \beta_1
    \label{eqn:GoryunovMondTorsionCorrectionMuVsTau}
  \end{equation}
  for arbitrary determinantal hypersurface singularities with 
  \textit{isolated} singularity. 

  When the dimension of the source $p$ is equal to the expected codimension 
  $c' = \grade \langle Y^{\wedge m-1} \rangle$ of the singular locus of the set of 
  degenerate matrices in $\CC^{m\times m}$, Theorem \ref{thm:InheritanceOfPerfection} 
  applies: $\beta_0 = \dim_\CC \CC\{x\}/\langle A^{\wedge m-1}\rangle$ and 
  $\beta_1 = 0$. This leads to Theorem \ref{thm:GoryunovMondNonSmoothableIsolatedHypersurfaces} 
  below. 

  For Theorem \ref{thm:GoryunovMondSmoothableHypersurfaces} observe that for a 
  $1$-parameter stabilization $\mathbf A(x,t)$ of $A$ the complex
  $\mathbf A^* K(m,m,m-1)$ is exact. From the long exact sequence induced 
  by multiplication with the unfolding parameter $t$ on $\mathbf A^*K(m,m,m-1)$ one 
  may then infer $\beta_1 = \beta_0$ which yields the desired result.
\end{remark}

\subsubsection{Isolated hypersurface singularities without determinantal smoothing}

Most of the simple singularities of square matrices do not admit determinantal 
smoothings, i.e. they are subsumed under the second case of Lemma 
\ref{lem:MilnorFibersOfDeterminantalHypersurfaces}; the essential smoothing 
$M_A$ will be singular along the subspace $M_A^{m-1}$ and in particular 
different from the Milnor fiber $M_f$. This raises the question for the 
correct notion of Milnor number for the determinantal hypersurface singularity 
in this setting. The answer is given by the following proposition which is 
based on a modified version of a theorem by L\^e \cite{Le86}:

\begin{proposition}
  \label{prp:BouquetDecompositionSingularMilnorFiberHypersurface}
  Let $(X_A,0) \subset (\CC^p,0)$ be an essentially isolated 
  determinantal hypersurface singularity defined by a matrix 
  $A \in \CC\{x_1,\dots,x_p\}$. Then the essential smoothing 
  \begin{equation}
    M_A \cong_{\mathrm{ht}} \bigvee_{i=1}^{\mu(A)} S^{p-1}
    \label{eqn:BouquetDecompositionSingularMilnorFiberHypersurface}
  \end{equation}
  of $(X_A,0)$ is homotopy equivalent to a bouquet 
  of spheres of real dimension $p-1$. 
\end{proposition}

\begin{proof}
  We give a sketch of the proof; full details can be found in \cite{Le86} and 
  \cite{Siersma91_2}, cf. also \cite{DamonMond91}. 
  For a suitable representative of a generic $1$-parameter unfolding 
  $\mathbf A \colon U \times T \to \CC^{m\times m}$ of $A$, the 
  \textit{polar locus} 
  \[
    \Gamma = \{ (x,t) \in U \times T \colon x \textnormal{ is a critical point of } \det A_t \}
  \]
  is a finite, branched covering over the parameter space $T$. 
  Denote the fiber over $t\in T$ by $\Gamma_t$.
  By genericity, we may assume that for $t\neq 0$, the function 
  $\det A_t$ has a complex Morse singularity at every point 
  $x \in \Gamma_t$. Now we can use the real valued function 
  $|\det A_t|$ as a (stratified) Morse function on the complement 
  $B_\varepsilon \setminus M_A$ of $M_A$ in a 
  Milnor ball $B_\varepsilon$, to find 
  that $B_\varepsilon$ is obtained from the essential smoothing 
  by attaching handles of Morse index $p$ at the points of 
  $\Gamma_t$. Since $B_\varepsilon$ itself is contractible, 
  the claim follows.
\end{proof}

\begin{definition}
  The \textit{singular Milnor number} of $(X_A,0)$ 
  is the number of spheres in the bouquet decomposition 
  (\ref{eqn:BouquetDecompositionSingularMilnorFiberHypersurface})
  of the essential smoothing.
\end{definition}

\medskip

When the number of variables 
$p$ is equal to the codimension $c'$ of the singular locus of the respective 
set of degenerate matrices, the determinantal hypersurface singularity $(X_A,0)$
is still isolated, but its essential smoothing $M_A$ will retain isolated 
singularities at some points $x_1,\dots,x_s \in M_A$. We are therefore in a
boundary setting where we have \textit{two} Milnor numbers at hand: The classical 
one, $\mu(f)$, for the isolated hypersurface singularity defined by $f = \det A$ 
(resp. $f = \Pf A$), and the singular Milnor number $\mu(A)$. 
In order to compare the two, let $\mathbf A(x,t)$ be an unfolding of the defining 
matrix $A$ of $(X_A,0)$ over some parameter space $T \subset \CC^k$.
As was already pointed out in the proof of Lemma 
\ref{lem:WhitneyStratificationOfAnEIDSAndSmoothability}, the (relative) 
singular locus of the fibers $X_{\mathbf A}(t)$ is again determinantal 
for the $(m-1)$-minors (resp. Pfaffians) of the defining matrices.
The principle of conservation of number assures that the total 
multiplicity of the singular points of the fibers $X_{\mathbf A}(t)$ 
is preserved in 
any family and therefore equal to the multiplicity 
of the singular locus of the central fiber 
$e:=\dim_\CC \CC\{x\}/\langle A^{\wedge m-1} \rangle$. 

When the unfolding $\mathbf A(x,t)$ is sufficiently generic so that 
$A_t$ is a stabilization for $t\neq 0$, we may assume that all these 
singular points $x_1,\dots,x_e \in M_A = X_{\mathbf A}(t)$ are 
Morse critical points of $f$ on $\CC^p$. Given that, according 
to Lemma \ref{lem:MilnorFibersOfDeterminantalHypersurfaces}, 
$M_f$ is a global smoothing of $M_A$, it is then 
not difficult to see that on the topological side, $M_A$ is 
homotopy equivalent to a suspension of exactly $e$ spheres 
in the Milnor fiber $M_f$. Therefore, the number $e$ measures 
the difference 
\begin{equation}
  e = \dim_\CC \CC\{x\}/\langle A^{\wedge m-1} \rangle = \mu(f) - \mu(A)
  \label{eqn:DifferenceOfMilnorNumbers}
\end{equation}
of the classical and the singular Milnor number when 
$p=3$ in the symmetric, $p=4$ in the arbitrary square, or $p=6$ in the 
skew-symmetric case. 

In these terms, Goryunov and Mond found that 
the $\SL$-Tjurina number is equal to the \textit{singular} Milnor number, 
cf. \cite[Theorem 4.6]{GoryunovMond05} and \cite[Corollary 4.2]{GoryunovMond05}:

\begin{theorem}
  In the same setting as Theorem \ref{thm:GoryunovMondSmoothableHypersurfaces} one has 
  \begin{eqnarray*}
    \tau_{\SL}^{\sym}(A) = \mu( f ) - \dim_\CC \CC\{x\}/\langle A^{\wedge m-1}\rangle & 
    & \textnormal{ for $A$ symmetric and $p = 3$;}\\
    \tau_{\SL}^{\sq}(A) = \mu( f ) - \dim_\CC \CC\{x\}/\langle A^{\wedge m-1}\rangle & 
    & \textnormal{ for $A$ arbitrary square and $p=4$;}\\
    \tau_{\SL}^{\sk}(A) = \mu( f ) - \dim_\CC \CC\{x\}/\langle A^{\wedge \frac{m}{2}-1}_\sk\rangle &
    & \textnormal{ for $A$ skew-symmetric and $p=6$.}
  \end{eqnarray*}
  \label{thm:GoryunovMondNonSmoothableIsolatedHypersurfaces}
\end{theorem}
Note that in the skew-symmetric case we assume $m$ to be even. 

\subsubsection{Determinantal hypersurfaces with non-isolated singularities}

For values of $p>c'$ in Lemma \ref{lem:MilnorFibersOfDeterminantalHypersurfaces}, 
the hypersurface singularities defined by $f$ are non-isolated, so 
the classical Milnor number is no longer defined. In fact, Kato and 
Matsumoto have established a lower bound on the connectivity of 
the Milnor fiber in \cite{KatoMatsumoto75}: 
When $s$ denotes the dimension of the critical locus of 
$f$, then $M_f$ is only $(p-s-2)$-connected and this bound is sharp in 
general. In this section we will survey recent results due to Goryunov 
\cite{Goryunov21} and Damon \cite{Damon16}, \cite{Damon18}, on both the 
singular Milnor fiber and its smoothing, respectively.

\medskip

Starting with the singular Milnor fiber of a square matrix 
$A \in \CC\{x_1,\dots,x_p\}^{m\times m}_{(*)}$ with $(*)$ either $\sym$, $\sq$, or $\sk$,
Proposition 
\ref{prp:BouquetDecompositionSingularMilnorFiberHypersurface}
assures that despite the presumably low connectivity of the classical Milnor fiber 
due to Kato and Matsumoto, the essential smoothing $M_A$ of $(X_A,0)\subset(\CC^p,0)$ is 
again homotopy equivalent to a bouquet of real spheres of equal dimension $p-1$.
In particular, the singular Milnor fiber is defined. 

The equality 
\begin{equation}
  \mu(A) = \tau_{\SL}(A)
  \label{eqn:MuEqualsTauDHS}
\end{equation}
(with $\tau_\SL(A)$ the appropriate Tjurina number in the symmetric, square, or skew-symmetric case)
was established by Goryunov in 
\cite{Goryunov21} 
for a wider class of singularities comprising all the known 
simple singularities of square matrices with non-isolated singularities, 
\cite[Corollary 6.3]{Goryunov21}, i.e. those listed in 
Theorem \ref{thm:ClassificationOfSimpleSquareMatrices}, 3 and 4 
(Table \ref{tab:SquareSymm3p4}), 
Theorem \ref{thm:ClassificationOfSimpleSymmetricSquareMatrices}, 3 and 
6 (also Table \ref{tab:SquareSymm3p4}), and Table \ref{tab:SquareSkewm4}.
In particular, he considers singularities of square matrices $A$ whose 
differential $\D A(0)$ is of corank $1$. Such singularities have already 
appeared in the classification by Bruce and Tari, 
Theorem \ref{thm:ClassificationOfSimpleSquareMatrices}, 3. 
The same reasons leading to the distinction of the two cases $a$ and $b$ listed there 
yield the classes of matrices of the form (\ref{eqn:GoryunovsSpecialSingularitiesI}) 
and (\ref{eqn:GoryunovsSpecialSingularitiesII}) below.

Let $x_{i,j}$ be independent variables for $0<i,j\leq m$ with $(i,j) \neq (1,1)$ 
and $z_1,\dots,z_q$ an additional set of variables. Then for $m\times m$-matrices of the form 
\begin{equation}
  A = \left( g(z) - \sum_{i=2}^m x_{i,i} \right) \cdot E_{1,1} 
  + \sum_{(i,j)\neq (1,1)} x_{i,j} \cdot E_{i,j}
  \label{eqn:GoryunovsSpecialSingularitiesI}
\end{equation}
where $g\in \CC\{z_1,\dots,z_q\}$ defines an isolated hypersurface singularity, 
Goryunov establishes in \cite[Theorem 3.5]{Goryunov21} that 
$\mu(A) = \tau_\SL(A) = \mu(g)$ whenever the $\SL$-Tjurina 
number is finite. 

Similarly, for matrices of the form 
\begin{equation}
  A = \left( h(x_{2,2},z) - \sum_{i=3}^m x_{i,i} \right) \cdot E_{1,1} 
  + \sum_{(i,j)\neq (1,1)} x_{i,j} \cdot E_{i,j}
  \label{eqn:GoryunovsSpecialSingularitiesII}
\end{equation}
where $h \in \CC\{x_{2,2},z_1,\dots,z_q\}$ defines a boundary singularity 
in the sense of Arnold \cite{Arnold78}, Goryunov establishes 
(\ref{eqn:MuEqualsTauDHS}) in \cite[Theorem 3.8]{Goryunov21}, together 
with a further equality in the quasi-homogeneous case (cf. \cite[Remark 3.9]{Goryunov21}) 
to the ``boundary Milnor number'' $\mu_\partial(h)$ 
for the boundary singularities.
The analogous results also hold, with the appropriate adaptations, for 
symmetric and skew-symmetric matrices. 

Given that no counterexamples to (\ref{eqn:MuEqualsTauDHS}) have been encountered 
among the known simple singularities, Goryunov has conjectured:

\begin{conjecture}[\cite{Goryunov21}]
  \label{con:GoryunovConjecture}
  Let $A \colon (\CC^p,0) \to (\CC^{m\times m}_{(*)},0)$ be a holomorphic 
  map germ with $(*)$ either $\sq$, $\sym$, or $\sk$ with finite $\SL$-codimension 
  and the dimension $p$ sufficiently large such that the 
  associated hypersurface $(X_A^m,0)$ has non-isolated singular locus. Then 
  \[
    \mu(A) = \tau_{\SL}(A).
  \]
\end{conjecture}

\medskip

We now turn to the study of smoothings of essentially isolated 
determinantal hypersurface singularities with non-isolated singularities.
These have occupied central stage in Damon's study of 
prehomogeneous vector spaces and exceptional orbit varieties, \cite{Damon16} 
and \cite{Damon18}. 

First, he studies the Milnor fibration of the functions 
$f = \det$ in the arbitrary square and symmetric, and $f=\Pf$ in the skew-symmetric 
case on the space of matrices $\CC^{m\times m}_\sq$, $\CC^{m\times m}_\sym$, 
resp. $\CC^{m\times m}_{\sk}$ with $m=2n$ even. We shall refer to these as the 
\textit{generic determinantal hypersurface singularities}. 
The homogeneity of all these singularities allows one to identify the local 
Milnor fibration at the origin with a fibration of affine manifolds
\begin{equation}
  \xymatrix{
    F_m \ar@{^{(}->}[r] \ar[d] &
    E_m \ar[d]^f \\
    \{1\}\ar@{^{(}->}[r] &
    S^1
  }
  \label{eqn:GlobalMilnorFibrationDamon}
\end{equation}
where $S^1\subset \CC^*$ is the unitary subgroup. Damon observed 
in \cite[Theorem 3.1]{Damon16} that the total 
spaces $E_m$ and the fibers $F_m$ are homotopy equivalent to 
``symmetric spaces'' in the sense of Cartan. In particular, 
this allows for an explicit computation of their cohomology 
rings with coefficients in a field $k$ of characteristic zero. 
The results are listed in Table \ref{tab:Damon16Table1}.

\begin{table}
  \centering
  \begin{tabular}{|l|c|c|l|}
    \hline
    Case & $F_m$ & Symmetric space & Cohomology of $F_m$ \\
    \hhline{|=|=|=|=|}
    $\sym$, $m$ odd & 
    $\SL(m;\CC)/\mathrm{SO}(m;\CC)$ & 
    $\mathrm{SU}(m;\CC)/\mathrm{SO}(m;\CC)$ & 
    $\bigwedge k\left\langle e_5, e_9, \dots, e_{2m-1}\right\rangle$ \\
    \hline
    $\sym$, $m$ even & 
    $\SL(m;\CC)/\mathrm{SO}(m;\CC)$ & 
    $\mathrm{SU}(m;\CC)/\mathrm{SO}(m;\CC)$ & 
    $\{1,e_m\}\cdot \bigwedge k\left\langle e_5, e_9, \dots, e_{2m-3}\right\rangle$ \\
    \hline
    $\sq$ & 
    $\SL(m;\CC)$ & 
    $\mathrm{SU}(m;\CC)$ & 
    $\bigwedge k\left\langle e_3, e_5, \dots, e_{2m-1}\right\rangle$ \\
    \hline
    $\sk$, $m=2n$ even & 
    $\SL(m;\CC)/\mathrm{Sp}(n;\CC)$ & 
    $\mathrm{SU}(m;\CC)/\mathrm{Sp}(n;\CC)$ & 
    $\bigwedge k\left\langle e_5, e_9, \dots, e_{2m-3}\right\rangle$ \\
    \hline

  \end{tabular}
  \caption{Milnor fibers of the generic determinantal hypersurfaces, their 
    associated symmetric spaces, and cohomology rings according to Damon \cite[Table 1]{Damon16}.}
  \label{tab:Damon16Table1}
\end{table}

Following Damon, we denote by $k\langle e_{i_1},e_{i_2},\dots\rangle$ 
the graded vector space generated by elements $e_{i_l}$ which 
are homogeneous of degree $i_l$, respectively. By $\bigwedge M$ we denote 
the full exterior algebra of $M$ and $\{1,e_m\}\cdot \bigwedge M$ denotes 
the free module generated by the elements $1$ and $e_m$ over $\bigwedge M$.

The results in Table \ref{tab:Damon16Table1} are also valid when replacing 
the coefficients $k$ by $\ZZ$, except in the symmetric case where the 
cohomology has $2$-torsion. In that case, the cohomology of $F_m$ 
with coefficients in 
$\ZZ_2 = \ZZ/2\ZZ$ is 
\begin{equation}
  H^\bullet(F_m; \ZZ_2) \cong H^\bullet(\mathrm{SU}(m;\CC)/\mathrm{SO}(m;\RR);\ZZ_2) 
  \cong \bigwedge \ZZ_2\left\langle s_2,s_3,\dots,s_m\right\rangle
  \label{eqn:DamonMod2Cohomology}
\end{equation}
with $s_j$ of degree $j$.

A ``Schubert decomposition'' of the global Milnor fibers $F_m$ was presented 
by Damon in \cite{Damon18}. He shows that the associated \textit{Schubert cycles} 
in the homology of $F_m$ are dual to the generators in cohomology from 
Table \ref{tab:Damon16Table1} in many cases and it is conjectured that this 
always holds. The decomposition itself is too complicated to be reproduced here 
and we refer to \cite{Damon18} for details.

It is already clear from Table \ref{tab:Damon16Table1} that the smoothings 
of the generic determinantal hypersurface singularities are hardly ever homotopy equivalent 
to a bouquet of spheres of the same dimension. Damon also provides insight to 
the homotopy groups of the $F_m$ in certain ranges in \cite[Theorem 3.5]{Damon16}
by considering the associated symmetric spaces as subspaces of the respective 
infinite dimensional symmetric spaces 
\begin{eqnarray*}
  \mathbb S \mathbb U &=&  \bigcup_{m=1}^\infty \mathrm{SU}(m;\CC), \\
  \mathbb S \mathbb U/ \mathbb S \mathbb O &=& 
  \bigcup_{m=1}^\infty \mathrm{SU}(m;\CC)/\mathrm{SO}(m;\CC), \\
  \mathbb S \mathbb U/ \mathbb S p &=&  
  \bigcup_{n=1}^\infty \mathrm{SU}(2n;\CC)/\mathrm{Sp}(n;\CC).
\end{eqnarray*}
The homotopy groups of the spaces on the left hand side are known, 
\begin{table}
  \centering
  \begin{tabular}{|c|p{0.3cm}|p{0.3cm}|p{0.3cm}|p{0.3cm}|p{0.3cm}|p{0.3cm}|p{0.3cm}|p{0.3cm}|p{0.3cm}|p{0.5cm}|}
    \hline
    $j=$ & 
    $0$ & $1$ & $2$ & $3$ & $4$ & $5$ & $6$ & $7$ & $8$ & $9$ \\
    \hhline{|=|=|=|=|=|=|=|=|=|=|=|}
    $\pi_j(\mathbb S \mathbb U)$ & 
    $0$	& $0$	& $0$ 	& $\ZZ$ & $0$ 	& $\ZZ$ & $0$ 	& $\ZZ$ & $0$	& $\ZZ$ \\
    $\pi_j(\mathbb S \mathbb U/ \mathbb S \mathbb O)$ & 
    $0$	& $0$	& $\ZZ_2$ & $\ZZ_2$ & $0$ 	& $\ZZ$ & $0$ 	& $0$ & $0$	& $\ZZ$ \\
    $\pi_j(\mathbb S \mathbb U/ \mathbb S p)$ & 
    $0$	& $0$	& $0$ & $0$ & $0$ 	& $\ZZ$ & $\ZZ_2$ & $\ZZ_2$ & $0$	& $\ZZ$ \\
    \hline
  \end{tabular}
  \caption{The stable homotopy groups of the infinite dimensional symmetric spaces. 
  They are periodic of period $8$ starting at $j=2$.}
  \label{tab:Damon16Table3}
\end{table}
see Table \ref{tab:Damon16Table3}, and there are certain \textit{stable ranges} in which 
these homotopy groups coincide with those of their respective subspaces:
\begin{theorem}[\cite{Damon16}]
  The homotopy groups of the Milnor fibers $F_m$ in the different cases are: \\
  \begin{tabular}[h]{llrcll}
    $\sq$,& $m$ arbitrary: &
      $\pi_j(F_m)$ &$\cong$& $\pi_j(\mathrm{SU}(m;\CC)) \cong \mathbb S \mathbb U$ & for $j<2m$;\\
    $\sym$, & $m$ arbitrary: & 
      $\pi_j(F_m)$ &$\cong$& $\pi_j(\mathrm{SU(m;\CC)}/\mathrm{SO}(m;\CC))\cong \mathbb S \mathbb U/
      \mathbb S \mathbb O$ &
      for $j<m-1$; \\
    $\sk$, & $m=2n$ even: &
      $\pi_j(F_m)$ &$\cong$ & $\pi_j(\mathrm{SU}(2n;\CC)/\mathrm{Sp}(n;\CC)) \cong \mathbb S \mathbb U/
      \mathbb S p$&
      for $j<2m-2$.\\
  \end{tabular}
  \label{thm:Damon16Theorem3.5}
\end{theorem}

\medskip 

After studying the smoothings of the generic determinantal hypersurface 
singularities, Damon turns to the study of essentially isolated hypersurfaces 
given by map germs 
\[
  A \colon (\CC^p,0) \to (\CC^{m\times m}_{(*)},0)
\]
and their smoothings $M_f = B_\varepsilon \cap (\det \circ A)^{-1}(\{\delta\})$, 
$1 \gg \varepsilon \gg |\delta| >0$, where again $(*)$ denotes either $\sq$, $\sym$, 
or $\sk$, and $\det$ is replaced by $\Pf$ in the latter case. 

Suppose $A\colon U \to \CC^{m\times m}$ is a suitable representative. Then 
by construction the map $A$ restricts to  
\[
  M_f \overset{A}{\longrightarrow} F_m
\]
where $F_m$ is the Milnor fiber of the generic determinantal hypersurface singularity.
In this setup, Damon gives the following definition:

\begin{definition}
  \label{def:CharacteristicCohomology}
  The image $\mathcal A_{M_{m,n}^s}(A)$ of the pullback in cohomology 
  $A^* \colon H^\bullet(F_m) \to H^\bullet(M_f)$
  is the \textit{characteristic cohomology} of $M_f$.
\end{definition}

According to Damon, it is an open question to determine the structure of 
$H^\bullet(M_f)$ as a module over $H^\bullet(F_m)$. He notes that in the two 
extremal cases where either i) the singularity $(X_A,0)$ admits a determinantal 
smoothing, or ii) when $A$ is the germ of a submersion, one has 
\begin{equation}
  H^\bullet(M_f) \cong \mathcal A_{M_{m,n}^s}(A) \oplus k^\mu[p-1]
  \label{eqn:DamonsCohomologyDecomposition}
\end{equation}
where $k$ denotes the chosen ring of coefficients for cohomology (cf. Table \ref{tab:Damon16Table3} 
and (\ref{eqn:DamonMod2Cohomology})) and $[p-1]$ the shift in cohomological degree by $p-1$. 
He remarks that in case i), the characteristic cohomology 
$\mathcal A_{M_{m,n}^s}(A)$ consists of the degree-zero-part only, so that $\mu$ is in fact the classical 
Milnor number of the associated isolated hypersurface singularity and in case ii) 
the second summand is trivial. For all other cases he asks:
\begin{center}
  \textit{``How generally valid is (\ref{eqn:DamonsCohomologyDecomposition}) for matrix 
singularities of the three types?''}
\end{center}
We will see in Theorem \ref{thm:BouquetDecompositionForEIDS} 
below that similar decompositions for the cohomology can be observed for 
essential smoothings of EIDS defined by non-square matrices.

\subsection{Isolated Cohen-Macaulay codimension $2$ singularities} 
\label{sec:ICMC2Topology}

As was already discussed earlier, isolated Cohen-Macaulay codimension $2$ singularities, 
which are not complete intersections, 
arise in a range of dimensions from $0$ up to dimension $4$.
In either case, the miniversal unfolding of the defining 
matrix\footnote{Since the determinantal structure of a 
  given ICMC2 singularity is unique, we will in what follows usually 
  supress the matrix $A$ from the notation and simply write $(X,0)$ 
  for the singularity, $M$ for its (essential) smoothing, $\hat X$ 
  for its Tjurina transform, etc. 
}
$A \in \CC\{x\}^{m\times (m+1)}$ induces a semi-universal deformation 
of the singularity $(X,0) = (A^{-1}(M_{m,m+1}^m),0)$ and there is 
a unique essential smoothing $M=M_A^m$ which is in fact smooth in 
dimensions $d =\dim (X,0)\leq 3$.

\medskip

The zero-dimensional case is rather trivial: Since determinantal singularities 
are Cohen-Macaulay, the multiplicity of a fat point is preserved under 
deformations. A stabilization of the defining matrix will therefore split 
any zero-dimensional determinantal singularity into a collection of finitely 
many simple points, the number of which is equal to the 
multiplicity of the singularity. 

\medskip

The $1$-dimensional case is known as ``space curve singularities''. 
For curves, a lot of theory has been developed beyond the complete intersection 
case already, in particular for the study of their topology. We refer 
to \cite[Section 7.2.6]{Volume1}
for an account 
on smoothings of curves and deliberately exclude them from our 
further discussion here. 

\subsubsection{From simple ICMC2 surfaces towards higher dimensions}

In dimension two, the isolated Cohen-Macaulay codimension $2$ singularities 
are the normal surface singularities in $(\CC^4,0)$. 
Just like the curve case, normal surface singularities have been extensively studied.
But even in the particular case of 
codimension $2$ it is, for example, still unknown whether or not 
a smoothing of the singularity is always simply connected. 

Similar to the case of determinantal hypersurfaces, 
the known results that we wish to present here for surfaces and threefolds 
are again mostly motivated by observations made for the lists of simple singularities 
in Table \ref{tab:FKNSurface} and Table \ref{tab:FKN3Fold}.
Starting with surfaces, we already noted that the list of simple isolated 
Cohen-Macaulay codimension $2$ singularities coincides with the list of 
\textit{rational triple points} that have already appeared 
in Section \ref{sec:RationalSurfaceSingularities}. The rational 
triple points were classified by Artin 
\cite{Artin74} in terms of the dual graphs of their resolution and 
explicit equations for their embeddings in $(\CC^4,0)$ have been given 
by Tjurina \cite{Tjurina68}. 

A resolution of these singularities can be constructed
as follows.
Tjurina has shown in \cite{Tjurina68} 
that any rational triple point $(X,0)$ is determinantal of type 
$(2,3,2)$ in $(\CC^4,0)$ for some matrix $A \in \CC\{x_1,\dots,x_4\}^{2\times 3}$ 
(cf. Proposition \ref{prp:WahlDeterminantalRationalSingularities} and 
Theorem \ref{thm:WahlDeterminantalRationalSingularities}) 
and the matrix $A$ has rank $1$ along the smooth locus $X_{\mathrm{reg}}$ 
of $(X,0)$. Blowing up the associated rational map
\[
  (X,0) \dashrightarrow \PP^1,\quad x \mapsto \Span A(x)
\]
then provides her with the \textit{Tjurina transform} 
$\hat X \overset{\hat \nu}{\longrightarrow} X$ which 
has at most A-D-E-singularities along its exceptional set $\hat E = \hat\nu^{-1}(\{0\})$.
The A-D-E-singularities are known to be rational with explicitly given resolutions. 
It can then be shown that the resolution of singularities $Z \overset{\rho}{\longrightarrow} X$
obtained from resolving the singular points of $\hat X$ in fact satisfies the 
requirements (\ref{eqn:DefinitionRationality}) for rationality.

\medskip

This construction is of course not the original approach pursued by Tjurina, 
since she started with a configuration of exceptional divisors 
coming from a resolution of singularities in the first place. 
Just as in Theorem \ref{thm:WahlArtinComponent} by Wahl, 
she obtained the space $\hat X$ by \textit{blowing down} the 
A-D-E configurations on $(Z,E)$.
The reason for the presentation given here is that the above construction is 
what generalizes to simple ICMC2 threefold singularities, as was shown 
by the authors in \cite{FruehbisKruegerZach20}. Applying the Tjurina transformation 
to simple ICMC2 singularities of higher dimension they obtain:

\begin{theorem}[\cite{FruehbisKruegerZach20}]
  All \textit{simple}
  ICMC2 singularities $(X,0)$ of dimension $d\geq 2$ 
  have at most A-D-E-singularities in their Tjurina transform $(\hat X,\hat E)$. 
  Moreover, they admit a resolution of singularities $\rho \colon Z \to X$ 
  factoring through $\hat \nu \colon \hat X\to X$ such that 
  \[
    R^k \rho_* \OO_Z =
    \begin{cases}
      \OO_{X} & \text{ if } k = 0, \\
      0 & \text{ otherwise, }
    \end{cases}
  \]
  i.e. they are rational.
  \label{thm:SimpleICMC2sAreRational}
\end{theorem}
Note that the exceptional set $\hat E = \hat \nu^{-1}(\{0\})$ 
is a $\PP^1$ and therefore not a divisor in $\hat X$ in dimensions $d>2$.

\medskip

Interestingly, the proof of Theorem \ref{thm:SimpleICMC2sAreRational} 
in dimensions $3$ and $4$ 
is based on the compatibility of the Tjurina 
transformation with deformations for the simple singularities. 
For the rational triple points of surfaces this has already been noted 
(cf. Theorem \ref{thm:WahlArtinComponent}) and exploited in order to construct a 
\textit{resolution in family} factoring through the Tjurina transformation 
in family, i.e. a commutative diagram
\begin{equation}
  \xymatrix{ 
    Z \ar@{^{(}->}[r] \ar[d] \ar@/_2pc/[dd]_(.3)\rho&
    \mathcal Z \ar[d] \ar@/_2pc/[dd]_(.3)\rho &
    \\
    \hat X \ar@{^{(}->}[r] \ar[d]^{\hat \nu} &
    \widehat{\mathcal X}'\ar[d]^{\hat \nu} &
    \\
    X \ar@{^{(}->}[r] \ar[d] &
    \mathcal X' \ar[d]^{\pi'} \ar[r] &
    \mathcal X \ar[d]^\pi \\
    \{0\} \ar@{^{(}->}[r] &
    R \ar[r]^{\Phi} &
    S
  }
  \label{eqn:ResolutionInFamily}
\end{equation}
where, as in (\ref{eqn:SimultaneousResolutionRationalSingularities}),
$Z\to X$ is a minimal good resolution of singularities, 
$\mathcal Z \to R$ a flat family deforming 
$Z$, $\mathcal X' = \rho(\mathcal Z)$ the blowdown of the deformation of $Z$, 
$\widehat{ \mathcal X}'$ the Tjurina transformation in family 
for the induced deformation of $X$, and $\mathcal X \to S$ a representative 
of the semi-universal deformation of $(X,0)$.
Given the results of Artin, Wahl, and de Jong summarized in 
Section \ref{sec:RationalSurfaceSingularities}, 
such resolutions in family can be constructed for all 
rational surface singularities which are determinantal and 
the family over the Artin component $\Phi(R) = S' \subset S$ 
is always a smoothing of $(X,0)$.

\begin{lemma}
  Suppose $(X,0)$ is a determinantal rational surface singularity and 
  (\ref{eqn:ResolutionInFamily}) a minimal good resolution in family. 
  Then the generic fiber over the Artin component 
  $X_s = \pi^{-1}(\{s\})$ is diffeomorphic to the central fiber $Z$ 
  of the resolution.
  \label{lem:DiffeomorphismForResolutionInFamily}
\end{lemma}

\begin{proof}
  We may assume $X$ to be embedded in some open domain $U \subset \CC^p$ 
  containing the origin.
  Let $B_\varepsilon$ be a Milnor ball for the singularity of $X$ at $0$ 
  and $\rho^{-1}(B_\varepsilon)$ its preimage in $Z$. Since $(R,0)$ 
  is smooth, we may replace $R$ by a small disc $D_\delta\subset \CC^k$ for 
  some $k$ and assume $\mathcal X'$ to be embedded in the product $U \times \CC^k$. 
  We then extend the Milnor ball to a tube $B_\varepsilon \times D_\delta$ and its preimage 
  $\rho^{-1}(B_\varepsilon \times D_\delta)$ in the total space $\mathcal Z$.

  Using Ehresmann's lemma and the transversality of the intersection 
  $\partial B_\varepsilon \cap X$ we see that the deformation of 
  $X$ by $t \in D_\delta$ does not alter its boundary up to diffeomorphism 
  for $1\gg \delta >0$ sufficiently small.
  Since the projection $\rho \colon Z \to X$ is an isomorphism off $0$, the 
  same holds for $Z$ and its boundary $\rho^{-1}(\partial B_\varepsilon) \cap Z$ 
  in the family $\mathcal Z \to D_\delta$. But by construction, the space 
  $Z = (\pi' \circ \rho)^{-1}(\{0\})$ is already smooth and hence 
  $Z \cap \rho^{-1}(B_\varepsilon)$ a smooth compact manifold with boundary. 
  Another application of Ehresmann's lemma shows that, again for $\delta>0$ 
  sufficiently small, the family $\mathcal Z \to D_\delta$ is a trivial 
  fiber bundle. 

  Since the deformation of $(X,0)$ over the Artin component $(S',0)$ is a smoothing 
  and $\Phi$ is finite and surjective onto $(S',0)$, the fiber 
  $X_s \cong \pi'^{-1}(\{t\})$ at a point $s = \Phi(t)$ is smooth for 
  generic $t \in D_\delta$.
  But then the restriction of 
  $\rho \colon Z_t \cap \rho^{-1}(B_\varepsilon) \to \pi'^{-1}(\{t\}) \cap B_\varepsilon$ 
  must be an isomorphism for $\rho$ was a resolution in family and the resolution 
  in the central fiber was minimal. The assertion now follows from the identifications
  \[
    X_s \cong_{\mathrm{diff}} \pi'^{-1}(\{t\}) \cong_{\mathrm{diff}} 
    (\pi'\circ \rho)^{-1}(\{t\}) \cong_{\mathrm{diff}} Z
  \]
  up to diffeomorphism.
\end{proof}

\begin{corollary}
  \label{cor:TopologyOfSimpleICMC2Surfaces}
  The smoothing of a determinantal rational surface singularity over 
  the Artin component is homotopy equivalent to a bouquet of $2$-dimensional 
  spheres.
\end{corollary}
\begin{proof}
  Due to Laufer's criterion 
  on rationality (\cite[Theorem 4.2]{Laufer72}), all components $E_i$ 
  of the exceptional set $E \subset Z$ in a good resolution of singularities 
  must be smooth rational curves $E_i \cong \PP^1 \cong S^2$ 
  intersecting transversally. Since the dual graph of a rational singularity 
  can have no cycles, the statement follows from the fact that $E \hookrightarrow Z$ 
  is a deformation retract of an appropriately chosen representative of $Z$, 
  cf. \cite{Lojasiewicz64}.
\end{proof}

\begin{remark}
  \label{rem:SmoothingsOfICMC2SurfacesAndDeterminantalRationalSingularities}
  For an \textit{arbitrary} ICMC2 surface singularity $(X,0)$ not 
  everything is known about the homotopy type of its smoothing $M$ or even 
  the homology groups with integer coefficients. Due to a 
  result by Greuel and Steenbrink (\cite{GreuelSteenbrink83}, 
  cf. \cite[Theorem 7.2.16]{Volume1}) the smoothing is connected
  and furthermore its first Betti number vanishes 
  (\cite{GreuelSteenbrink83}, \cite[Theorem 7.2.18]{Volume1}). But 
  for instance it is not known whether or not $\pi_1(M) = 0$ in general.  
\end{remark}

\subsubsection{Vanishing homology for ICMC2 threefolds}

For simple ICMC2 \textit{threefold} singularities there is no full resolution in family. 
However, the Tjurina transformation in family applied to a 
smoothing of $(X,0)$ provides us with a truncated 
version of the diagram (\ref{eqn:ResolutionInFamily}): Just delete the top 
row $Z \hookrightarrow \mathcal Z$ and observe that due to the 
results by Schaps, Corollary \ref{cor:SemiUniversalDeformationOfICMC2}, 
the map $\Phi$ and consequently also $\mathcal X' \to \mathcal X$ are isomorphisms. 

\begin{theorem}[\cite{FruehbisKruegerZach20}]
  \label{thm:MainTheoremFromFKZ}
    Let $M$ be the smoothing of an ICMC2 
    threefold singularity $(X,0) \subset (\CC^5,0)$, 
    which has only isolated 
    singularities in the Tjurina transform $\hat X$. Then the singularities 
    of $\hat X$ are isolated complete intersection singularities 
    and the homology groups 
    with integer coefficients 
    of $M$ are 
    \[
        H_0(M) = \ZZ, \quad 
	H_1(M) = 0, \quad 
	H_2(M) = \ZZ, \quad 
	H_3(M) = \ZZ^r,
    \]
    where $r \in \mathbb N_0$ is the sum of the Milnor numbers 
      $r = \sum_{p\in \Sigma(\hat X)} \mu(\hat X,p)$
    of the singularities in $\hat X$ if the defining matrix is of size $2\times 3$, or 
    $r = \sum_{p\in \Sigma(\hat X)} \mu(\hat X,p) - 1$
    otherwise.
\end{theorem}

\begin{proof}
  We sketch a proof for the case of simple singularities. The full proof 
  can be found in \cite{FruehbisKruegerZach20}.
  
  According to Theorem \ref{thm:SimpleICMC2sAreRational} the singularities 
  in the Tjurina transform $\hat X$ of $(X,0)$ are at most A-D-E-singularities. 
  Choose a Milnor tube $B_\varepsilon \times D_\delta$ 
  as in the proof of Corollary \ref{cor:TopologyOfSimpleICMC2Surfaces}. 
  One can show that a smoothing of $(X,0)$ induces a smoothing of all singularities 
  of $\hat X$ and that the projection $\rho \colon \hat X(t) \to X(t)$ is an isomorphism 
  over smooth fibers $X(t) \cong_{\mathrm{diff}} M$. 
  Now the central fiber $\hat X=\hat X(0)$ of the Tjurina transform retracts onto its exceptional 
  set $\hat E \cong \PP^1 \cong S^2$ and the class of that sphere freely generates 
  $H^2(\hat X)$. Passing to the smooth fiber $\hat X(t) \cong_{\mathrm{diff}} M$, 
  one can show that this cycle survives 
  while the vanishing cycles from the smoothings of the A-D-E-singularities in $\hat X$ 
  freely generate $H_3(M)$.
\end{proof}

\begin{remark}
  \label{rem:HomologyForArbitraryICMC2Threefolds}
  In \cite{Zach18} 
  the second named author shows that the homology groups of the smoothing take 
  the same form for all threefolds defined by matrices of size $2\times 3$, i.e. 
  also those with 
  \textit{non-isolated} singularities in the Tjurina transform. In that case, however, 
  there is no formula for the computation of the rank $r$ of the middle homology 
  group.
  Recently, the second named author also announced that the same decomposition of the homology 
  holds for smoothings of arbitrary ICMC2 threefold singularities, i.e. independent 
  of the size of the defining matrix. 
\end{remark}

\begin{remark}
  \label{rem:MilnorNumberForThreefolds}
  What is interesting in Theorem \ref{thm:MainTheoremFromFKZ}
  is the non-vanishing of the second Betti number. 
  It had already been observed by Damon and Pike in \cite{DamonPike14} 
  that for some simple ICMC2 threefolds the reduced Euler characteristic 
  of the smoothing was negative. Since the first Betti number of $M_A^m$ 
  was known to be zero due to results by Greuel and 
  Steenbrink\footnote{See \cite{GreuelSteenbrink83}, cf. also \cite[Theorem 7.2.17]{Volume1}}, 
  this implied 
  that at least $H_2(M_A^m)$ had to be nontrivial. The fact that there is more 
  than one degree with nontrivial vanishing homology raises the question 
  as to what is the correct notion of Milnor number for ICMC2 threefolds. 
  As of this writing, this question has not been settled in a satisfactory way. 
  However, the independence of the second Betti number from the defining matrix 
  $A$ and the more recent Bouquet decomposition for essential smoothings, 
  see Theorem \ref{thm:BouquetDecompositionForEIDS} suggest that 
  $b_3 = \rank H_3(M_A^m)$ is a good candidate.
\end{remark}

\subsubsection{Comparison of Milnor and Tjurina numbers}

Coming back to the case of surface singularities again, 
there is Wahl's conjecture \ref{con:WahlsConjecture} which was already stated 
in the beginning of this section. For normal surface singularities, the first Betti 
number of any smoothing is zero due to the previously cited result of Greuel and Steenbrink, 
see \cite{GreuelSteenbrink83}, or also \cite[Theorem 2.7.18]{Volume1}. 
Consequently, the Milnor number $\mu$ of the singularity is defined to be the second Betti number. 

Since the ``if''-part of Wahl's conjecture has been established already and 
all \textit{simple} ICMC2 surface singularities are quasi-homogeneous, we see that 
in this case 
\begin{equation}
  \mu = b_2(M_A^2) = \tau(X_A^2,0) -1 = \tau_\GL(A) -1
  \label{eqn:MuVsTauForSimpleICMC2Surfaces}
\end{equation}
where as usual $A \in \CC\{x,y,z,w\}^{2\times 3}$ is the defining matrix of the singularity.

For the simple ICMC2 \textit{threefolds} on the other hand, 
the authors observed a different behaviour in \cite{FruehbisKruegerZach20}: The Tjurina number seems to be rather 
unrelated to the topology of the smoothing. For instance, singularities defined 
by the matrices 
\[
  \begin{pmatrix}
    w & y & x \\
    z & w & y + v^k
  \end{pmatrix},
\]
the so-called $\Pi_k$-family, have Tjurina number $\tau = 2k -1$. Their smoothings, 
however, are all homotopy equivalent to the sphere $S^2$, independent of $k$. 
More generally, they establish an equality 
\begin{equation}
  \tau(X,0) = h^1(\hat X, T^0_{\hat X}) + \sum_{p \in \Sigma(\hat X)} \tau(\hat X,p)
  \label{eqn:TauForSimpleICMC2Threefolds}
\end{equation}
for all ICMC2 threefolds defined by matrices of size $2\times 3$ with only 
isolated singularities in the Tjurina transform. Here $T^0_{\hat X}$ denotes the 
tangent sheaf of $\hat X$ and $\tau(\hat X,p)$ the local Tjurina numbers at the 
singularities $p\in \Sigma(\hat X)$. These local Tjurina numbers can then be related 
to the local Milnor numbers and the topology of the smoothing. But the example given by 
the $\Pi_k$-family shows that the ``correction term'' given by $h^1(\hat X, T^0_{\hat X})$ 
can become arbitrary big. 

\subsection{Arbitrary EIDS and some further, particular cases}

\subsubsection{Bouquet decomposition for essential smoothings}

For arbitrary EIDS, the second named author has established a bouquet decomposition 
of the essential smoothing as an application of a result due to Tib{\u a}r 
\cite{Tibar95} (see also \cite[Theorem 6.10.6]{Volume1}) 
to determinantal singularities. Tib{\u a}r's bouquet theorem 
itself is based on the carousel construction exhibited 
in \cite[Chapter 6]{Volume1}.

\begin{theorem}{\cite{Zach18BouquetDecomp}}
  Let $(X_A^s,0)\subset (\CC^p,0)$ be an EIDS of positive dimension and 
  of type $(m,n,s)$ defined by a matrix $A \in \CC\{x_1,\dots,x_p\}^{m\times n}$. 
  Then the essential smoothing $M_A^s$ of $(X_A^s,0)$ is homotopy equivalent to 
  \begin{equation}
    M_A^s \cong_{\mathrm{ht}} L_{m,n}^{s,p} \vee 
    \bigvee_{0\leq r < s} \bigvee_{i=1}^{\lambda(r,A)} S^{p-(m-r)(n-r)+1}( 
      L_{m-r,n-r}^{s-r-1,(m-r)(n-r)-1}
      )
    \label{eqn:BouquetDecomposition}
  \end{equation}
  where $L_{m,n}^{s,k}$ is the intersection 
    $L_{m,n}^{s,k} = H^k \cap M_{m,n}^s$
  of a $k$-dimensional hyperplane $H^k$ in general position off the origin 
  with the generic determinantal variety.
  \label{thm:BouquetDecompositionForEIDS}
\end{theorem}
In this formula, $S^k(\cdot)$ denotes the $k$-fold repeated suspension 
of a topological space with the convention that $S(\emptyset) = S^0$  
is the $0$-dimensional sphere, $S^0(X) = X$ for every $X$, 
and $S^k(X) = \emptyset$ for negative $k$. Note that 
$L_{m,n}^{s,mn-1}$ is nothing but the \textit{complex link} of the 
generic determinantal variety $M_{m,n}^s$ at the origin, see 
e.g. \cite[Section 5.9.3]{Volume1}. For more general values of 
$p$, the space $L_{m,n}^{s,p}$ is the essential smoothing of 
a linear EIDS determined by a generic linear map $\CC^p \to \CC^{m\times n}$.

\begin{proof}
  We give a rough outline of the proof in order to indicate 
  its synopsis with the other methods. Recall from Remark 
  \ref{rem:GraphTransformationAndRelativeCompleteIntersections} 
  that $(X_A^s,0)$ can be realized as a complete intersection 
  on $(\CC^p,0)\times(M_{m,n}^s,0)$ by considering the graph 
  $\Gamma_A = \{ (x,\varphi) \in \CC^p \times \CC^{m\times n} : \varphi = A(x) \}$ 
  which is given by the $m\cdot n$ equations $h_{i,j} = y_{i,j} - a_{i,j}(x)$. 
  Then 
  \[
    X_A^s \cong \Gamma_A \cap \CC^p \times M_{m,n}^s \subset \CC^p\times \CC^{m\times n}
  \]
  and these equations form a regular sequence on the coordinate ring 
  $\CC\{x\}[y]/\langle Y^{\wedge s}\rangle$ of the 
  variety $(\CC^p,0) \times M_{m,n}^s$. Note that the latter variety 
  is canonically Whitney stratified by the product of the rank stratification 
  with $\CC^p$. It is now not too difficult to 
  see that the transversality conditions imposed on $A$ off the origin 
  translate to an appropriate notion of transversality of the functions $h_{i,j}$ 
  that allows one to say that these equations define an isolated complete intersection 
  singularity on $(\CC^p,0)\times (M_{m,n}^s,0)$ in the \textit{stratified sense}.
  
  It is now easy to see that by construction the essential smoothing $M_A^s$ 
  coincides with the Milnor fiber\footnote{See e.g. \cite[Theorem 6.10.3]{Volume1}}
  of this isolated complete intersection on 
  $(\CC^p,0)\times(M_{m,n}^s,0)$. 
  The remainder of the proof now follows from Tib{\u a}r's decomposition theorem 
  \cite[Bouquet Theorem]{Tibar95} and a modification of \cite[Corollary 4.2]{Tibar95}.
\end{proof}

\begin{corollary}
  \label{cor:EssentialSmoothingAsICISMilnorFiber}
  The essential smoothing of an EIDS $(X_A^s,0) \subset (\CC^p,0)$ is homeomorphic 
  to the Milnor fiber of the complete intersection defined by 
  $y_{i,j} - a_{i,j}$ for $0<i\leq m$, $0<j\leq n$ on $(\CC^p,0) \times (M_{m,n}^s,0)$.
\end{corollary}

In certain particular cases (see below), Theorem \ref{thm:BouquetDecompositionForEIDS} 
leads to a full understanding of the essential smoothing 
$M_A^s$ up to homotopy: The numbers $\lambda(r,A)$ can, in principal, 
be computed from the so-called Cerf-diagrams involved in the 
carousel construction, see e.g. \cite[Theorem 6.6.6]{Volume1} and 
\cite[Section 6.7]{Volume1}. What is more difficult, is the 
investigation of the complex links and, more generally, 
the spaces $L_{m,n}^{s,p}$ for the various 
values of $m,n,s$ and $p$. These are the atomic building blocks for 
the topology of essential smoothings that can not be turned into 
singularities of complete intersections anymore. 

\medskip

The Euler characteristic of the complex links of generic determinantal 
varieties has been computed by Ebeling and Gusein-Zade in \cite{EbelingGuseinZade09}. 
Without loss of generality, one may assume that $m\leq n$. Then 
the reduced Euler characteristic of the complex link of $M_{m,n}^s$ at 
the origin is
\begin{equation}
  \overline \chi( L_{m,n}^{s,mn-1} ) = (-1)^s { m-1 \choose s-1 }.
  \label{eqn:EulerCharacteristicComplexLink}
\end{equation}
This formula has been used by Gaffney, Grulha, and Ruas to compute the 
\textit{local Euler obstructions}\footnote{The local Euler obstruction was introduced 
  by MacPherson in \cite{MacPherson74}.} of the generic determinantal varieties in 
\cite{GaffneyGrulhaRuas19}. Again, for $m \leq n$ they find
\begin{equation}
  \Eu(M_{m,n}^s,0) = { m \choose s-1 }.
  \label{eqn:EulerObstructionGenericDeterminantalVarieties}
\end{equation}
This formula has then been generalized for the generic determinantal 
varieties over arbitrary fields by Zhang \cite{Zhang17}.

Both these invariants are of fundamental importance in the study 
of stratified Morse theory on determinantal varieties, cf. 
\cite[Chapter 5]{Volume1}, and \cite{GoreskyMacPherson88}. 
In particular, they allow for the computation of the Euler characteristic 
of the essential smoothing $M_A^s$ of an arbitrary 
EIDS $(X_A^s,0)$ as in Theorem \ref{thm:BouquetDecompositionForEIDS} 
(again assuming $m\leq n$)
\begin{equation}
  \overline \chi(M_A^s) = \overline \chi(L_{m,n}^{s,p}) + 
  \sum_{0\leq r < s} 
  (-1)^{p+s-r-(m-r)(n-r)} \cdot { m-r-1 \choose s-r-2 }\cdot 
  \lambda(r,A) 
  \label{eqn:EulerCharacteristicFromBouquetTheorem}
\end{equation}
up to the term $\overline \chi(L_{m,n}^{s,p})$, which is constant and 
independent of the specific matrix $A$.

\begin{remark}
  In the spirit of Definition \ref{def:CharacteristicCohomology} for the 
  ``characteristic cohomology'' for smoothings of 
  determinantal hypersurface singularities due to Damon 
  we can make the following definition for 
  smoothings of EIDS defined 
  by non-square matrices $A \in \CC\{x_1,\dots,x_p\}^{m\times n}$.
  If $(X_A^s,0) \subset (\CC^p,0)$ is 
  smoothable, then a stabilization 
  \[
    A_t \colon B_\varepsilon \to \CC^{m\times n} \supset M_{m,n}^s
  \]
  will not meet the singular locus $M_{m,n}^{s-1}$ of the generic 
  determinantal variety $M_{m,n}^s$ so that the intersection 
  $A_t(B_\varepsilon) \cap M_{m,n}^s$ will be completely contained 
  in the stratum $V_{m,n}^{s-1}$. In analogy to Damon's definition 
  we may set 
  \[
    \mathcal A_{M_{m,n}^s}(A) = A_t^*( H^\bullet(V_{m,n}^{s-1})) \subset H^\bullet(M_A^s)
  \]
  to be the image 
  of the pullback in cohomology of $A_t \colon M_A^s \to V_{m,n}^{s-1}$ 
  to the essential smoothing $M_A^s = A_t^{-1}(M_{m,n}^s)$.
  In the smoothable case the wedge sum in (\ref{eqn:BouquetDecomposition}) 
  simplifies to $r = s-1$ so that 
  \[
    M_A^s \cong_{\mathrm{ht}} L_{m,n}^{s,p} \vee \bigvee_{i=1}^\lambda S^{p-(m-s+1)(n-s+1)}.
  \]
  We conjecture that $\mathcal A_{M_{m,n}^s}(A)$ is precisely the 
  contribution of $H^\bullet(L_{m,n}^{s,p})$ so that, given the 
  analogy of the definitions, the number $\lambda$ plays the 
  r\^ole of $\mu$ in (\ref{eqn:DamonsCohomologyDecomposition}).
\end{remark}

\subsubsection{Formulae for the vanishing Euler characteristic using polar varieties}

Various other methods have been developed in order to compute the 
Euler characteristics of essential smoothings, including the space 
$L_{m,n}^{s,p}$ itself. One of these methods is to study successive 
hyperplane sections of a given determinantal singularity $(X_A^s,0)$ 
and its deformation, 
together with the associated \textit{polar varieties}\footnote{Polar 
  varieties for complex analytic germs were introduced in \cite{LeTeissier81}. 
  See also \cite{Teissier82}.} and their multiplicities 
$m_i(X_A^s,0)$ for $0\leq i < \dim(X_A^s,0) = d$. In addition to that, 
one 
needs the so-called $m_d$-multiplicity $m_d(X_A^s,0)$ introduced by Gaffney
in \cite{Gaffney93} to also capture the behaviour of the singularity in families. 
For a given singularity, this multiplicity measures the number 
of critical points of a generic linear form on the smooth locus of 
the nearby stable object\footnote{Note that in order for $m_d(X_A^s,0)$ 
to be an invariant of the singularity $(X_A^s,0)$ itself rather than the 
given family, the nearby stable object needs to be uniquely determined 
by $(X_A^s,0)$. 
This is guaranteed for EIDS by the existence of the miniversal unfolding 
and the uniqueness of the essential smoothing.}. To give a concise formula, 
we furthermore use the convention that $m_i(X,0) = 0$ for negative $i$.

\begin{theorem}
  Let $M_A^s$ be the essential smoothing of an EIDS $(X_A^s,0) \subset (\CC^p,0)$ 
  defined by a matrix $A \in \CC\{x_1,\dots,x_p\}^{m\times n}$ with $m \leq n$.
  Then 
  \begin{equation}
    \chi(M_{A}^s) = 
    \sum_{0\leq r < s} \left( 
    \sum_{j=0}^{d(r)} (-1)^{d(r)-j} m_{d(r)-j}(X_A^{r+1},0)
    \right)\cdot (-1)^{s-r-1} { m-r \choose s-r-1 }
    \label{eqn:EulerCharacteristicEssentialSmoothingPolarMethods}
  \end{equation}
  where $d(r) = \dim(X_A^{r+1},0) = p - (m-r)(n-r)$.
  \label{thm:EulerCharacteristicViaPolarMethods}
\end{theorem}
This theorem and its variants have appeared in several places such as 
\cite{EbelingGuseinZade09}, \cite{BraChaRua15}, 
\cite{Pereira10}, \cite{BallesterosOreficeTomazella13}, 
\cite{GaffneyGrulhaRuas19}, and \cite{Zach17}.

\begin{proof}
  Again, we only sketch the proof to illustrate the related ideas. 
  Details can be found in the various sources cited above. 

  Let $\mathbf A \colon U \times T \to \CC^{m\times n}$ be a suitable 
  representative of a $1$-parameter unfolding of the defining matrix 
  $A$ on a parameter $t$ such that $A_t$ is a stabilization for $t\neq 0$ 
  and choose a sequence of linear forms $l_1,\dots,l_{d(s)}$ on 
  the ambient space $\CC^p$ of the singularity. We denote by $D_j$ 
  the hyperplane of codimension $j$ defined by $l_1= \dots = l_j = 0$. 
  If the linear forms have been chosen sufficiently general, 
  then either one of the singularities $(X_A^s \cap D_j,0) \subset (D_j,0)$ is 
  again an EIDS and the restriction of the unfolding $\mathbf A$ to 
  $D_j$ induces an essential smoothing 
  thereof. Furthermore, for every $t\neq 0$ sufficiently small, 
  every $r<s$, and every $j>0$ the function $l_{j+1}$ has 
  only complex Morse singularities on the interior of 
  the fiber 
  \[
    D_j \cap A_t^{-1}(M_{m,n}^{r+1}) \cong 
    D'_j \cap A_t^{-1}(M_{m,n}^{r+1}) \cong 
    D'_j \cap A_0^{-1}(M_{m,n}^{r+1}).
  \]
  where $D'_j$ denotes a hyperplane parallel to $D_j$ off 
  the origin. The above isomorphisms come from wiggeling either 
  $D_j$ or $A_t$ and parallel transport of the associated fibers 
  by virtue of Thom's isotopy lemma. 
  The number of the Morse critical points is precisely 
  $m_{d(r)-j}(X_A^{r+1},0)$, see \cite{LeTeissier81}. 
  Finally, for $j=0$, the multiplicity $m_{d(r)}(X_A^{r+1},0)$ 
  counts the number of critical points of $l_1$ on the regular 
  locus of the essential smoothing $M_A^{r+1}$, see \cite{Gaffney93}.
  
  By downwards induction on the codimension $j$ of the hyperplanes
  one can now rebuild the original essential 
  smoothing $M_A^s = A_t^{-1}(M_{m,n}^s)$ from its hyperplane sections 
  starting with $j=d(s-1)$.
  The Morse critical points on the regular loci of the various hyperplane 
  sections of $A_t^{-1}(M_{m,n}^{r+1})$ above are then \textit{stratified} 
  Morse critical points of the fiber $D_j \cap A_t^{-1}(M_{m,n}^s)$ located 
  on the respective strata $D_j \cap A_t^{-1}(V_{m,n}^r)$. At any such point,
  the factor 
  \[
    (-1)^{s-r-1}{m-r\choose s-r-1 } = 
    1 - \chi( L_{m-r,n-r}^{s-r,(m-r)(n-r)-1} )
  \]
  accounts for the \textit{normal Morse datum}\footnote{For the discussion 
    of normal Morse data in complex stratified Morse theory, see \cite{GoreskyMacPherson88}.}
  associated to that point which is determined by the complex links of the 
  generic determinantal varieties according to 
  Lemma \ref{lem:WhitneyStratificationOfAnEIDSAndSmoothability} and 
  (\ref{eqn:TrivializationOfTheRankStratification}).
\end{proof}

\begin{remark}
  \label{rem:MultiplicitiesVsTjurinaNumbers}
  The multiplicities $m_i(X_A^{r+1},0)$ appearing in 
  (\ref{eqn:EulerCharacteristicEssentialSmoothingPolarMethods}) 
  can in principal be computed by various different methods. 
  For $i< d = \dim(X_A^{r+1},0)$ this is possible directly from their definitions 
  in \cite{LeTeissier81}. For the $m_d$-multiplicity, Gaffney and 
  Ruas have given a formula in \cite[Proposition 4.6]{GaffneyRuas21} 
  as the sum of the multiplicity of the pair of modules given 
  by the Jacobian module and the normal module of the singularity 
  and the intersection number of the image of the defining matrix 
  with a certain polar variety of the generic determinantal variety.
\end{remark}

\begin{remark}
  \label{rem:TopologyVsDeformations}
  In general, the hard task is to relate the multiplicities in 
  (\ref{eqn:EulerCharacteristicEssentialSmoothingPolarMethods}) 
  or the numbers $\lambda(r;A)$ in 
  (\ref{eqn:BouquetDecomposition}) to deformation 
  theoretic invariants such as the $\GL$-Tjurina number of the defining matrix. 
  In the special cases described earlier, as for example in 
  Theorem \ref{thm:MainTheoremFromFKZ}, one can, of course, also 
  apply Theorem \ref{thm:BouquetDecompositionForEIDS} or
  Theorem \ref{thm:EulerCharacteristicViaPolarMethods} and then 
  compare $m_i(X_A^s,0)$ and $\lambda(r;A)$ to the Milnor and Tjurina numbers 
  of the singularities. 

  Yet another approach to the computation of the 
  vanishing Euler characteristic has been pursued by Damon and Pike in 
  \cite{DamonPike15} and \cite{DamonPike14}. For certain determinantal 
  varieties (including the determinantal hypersurfaces and type $(2,3,2)$), 
  they succeed to embed the associated generic determinantal variety in a 
  special arrangement of so-called ``$H$-holonomic divisors'' $W_i$:
  \[
    M_{m,n}^s \subset \bigcup_{i} W_i.
  \]
  These \textit{free completions} of $M_{m,n}^s$
  are extracted from certain group representations on 
  the space of matrices which are closely related to 
  Cholesky decomposition. 

  A nonlinear section 
  $A \colon (\CC^p,0) \to (\CC^{m\times n},0) \supset W$ of an 
  $H$-holonomic divisor $W$ which is algebraically transverse to $W$ 
  off the origin, gives rise to an ``almost free divisor'' 
  $(W_A,0) = (A^{-1}(W),0) \subset (\CC^p,0)$. 
  Similar to the case of determinantal varieties, one considers 
  the deformations of $(W_A,0)$ induced from unfoldings of $A$. 
  Stabilizations $A_t$ of $A$ lead to the analogue of essential smoothings 
  $A_t^{-1}(W)$ for 
  $(W_A,0)$ and these spaces are homotopy equivalent to a bouquet 
  of spheres for the same reasons as outlined in Lemma 
  \ref{lem:MilnorFibersOfDeterminantalHypersurfaces} for 
  determinantal hypersurfaces --- 
  only that in this case the number $\mu$ of these spheres is equal 
  to the $\mathcal K_{H,e}$-codimension of $A$ 
  (\cite[Theorem 3.1]{DamonPike14}). This $\mathcal K_{H,e}$-equivalence 
  for non-linear sections of $H$-holonomic free divisors 
  is similar to the $\SL$-equivalence of matrices 
  discussed in Remark 
  \ref{rem:SLEquivalence}, see e.g. \cite{GoryunovMond05}. 

  Using the additivity of the topological Euler characteristic on complex 
  analytic sets, Damon and Pike can then infer the Euler characteristic of 
  the essential smoothing of the determinantal singularity from 
  its relative position in the stabilizations 
  \[
    M_A^s = A_t^{-1}(M_{m,n}^s)  
    \subset \bigcup_{i} A_t^{-1}(W_i).
  \]
  This machinery is constructed for the purpose of relating the vanishing 
  Euler characteristic of essential smoothings to deformation theoretic 
  invariants. For various specific configurations
  closed formulas are given for 
  $\overline \chi(M_{A}^s)$ in terms of the various $\mathcal K_{H,e}$-codimensions, 
  cf. e.g. \cite[Theorem 8.1]{DamonPike14}.
  However, it is not clear how to determine such formulae in general 
  and how the various 
  $\mathcal K_{H,e}$-codimensions for the $W_i$ relate to the 
  $\GL$- or $\SL$-Tjurina numbers of the original matrix singularity. 
\end{remark}

\subsubsection{Some explicitly known complex links}

To conclude this section, we will summarize a few special cases in 
which the topology of the essential smoothings of generic linear 
EIDS, i.e. the basic building blocks in (\ref{eqn:BouquetDecomposition}), 
is known. 

The complex links of the degenerate square matrices have been described by 
Goryunov \cite[Theorem 2.1]{Goryunov21}:

\begin{theorem}
  Let $L_{(*)}$ be the complex link of the generic determinantal hypersurface in 
  $\CC^{m\times m}_{(*)}$ where $(*)$ denotes either 
  $\sq$, $\sym$, or $\sk$. Then in every case $L_{(*)}$ is homotopy equivalent 
  to a single sphere $S^{N-2}$ where $N$ equals either $m^2$, $\frac{1}{2}m(m+1)$, 
  or $\frac{1}{2}m(m-1)$ with $m = 2n$ even, respectively. Depending on the case, these spheres 
  can be chosen to be 
  \begin{description}
    \item[$\sq:$] all degenerate Hermitian matrices in $\CC^{m\times m}$ 
      with trace $1$ and all eigenvalues non-negative.
    \item[$\sym:$] all degenerate real matrices in $\CC^{m\times m}_\sym$ 
      with trace $1$ and all eigenvalues non-negative.
    \item[$\sk:$] all degenerate quaternionic matrices in $\CC^{m\times m}_\sk$
      with skew trace $1$ and all skew eigenvalues non-negative.
  \end{description}
  \label{thm:GoryunovComplexLink}
\end{theorem}
The homogeneity of these particular singularities allow for scaling of these spaces 
to arbitrarily small sizes close to the origin. For the skew-symmetric matrices $A$ 
of even size $m = 2n$ Goryunov has set the \textit{skew trace} to be $\sum_{i=1}^n a_{2i-1,2i}$ 
and the skew eigenvalues $\lambda_i$ the solutions of the equation 
\[
  \Pf\left(  A-\lambda \cdot \sum_{i=1}^n \left( E_{2i-1,2i} - E_{2i,2i-1} \right)\right) = 0.
\]

\medskip

For non-square matrices it is possible to determine 
the spaces $L_{m,n}^{s,p}$ up to homotopy for all values of $n$ and $p$ when $m = s = 2$, 
see \cite[Section 4.1]{Zach18BouquetDecomp}: 
\begin{equation}
  L_{2,n}^{2,p} \cong_{\mathrm{ht}} 
  \begin{cases}
    \{\mathrm{pt.}\} & \textnormal{ if } p \geq 2n,\\
    S^2 & \textnormal{ if } n < p < 2n,\\
    \bigvee_{i=1}^{n-1} S^1 & \textnormal{ if } p = n,\\
    \{ n \,\mathrm{ points } \} & \textnormal{ if } p= n-1.
  \end{cases}
  \label{eqn:DeterminantalLinks2xn}
\end{equation}

\begin{remark}
  Using (\ref{eqn:DeterminantalLinks2xn}) in combination with 
  Theorem \ref{thm:BouquetDecompositionForEIDS}, it is possible to sharpen 
  Theorem \ref{thm:MainTheoremFromFKZ} for the case of $2\times 3$ matrices 
  so that the smoothing is in fact homotopy equivalent to a bouquet 
  of one $2$-sphere and $r$ spheres of dimension $3$.
\end{remark}


\section{Appendix: Lists of simple singularities}

\subsection{Arnold's lists}

For reader's convenience and to complement Theorems
\ref{thm:ClassificationOfSimpleSquareMatrices} and
\ref{thm:ClassificationOfSimpleSymmetricSquareMatrices}, we give the
lists due to Arnold (see \cite{ArnoldGuseinZadeVarchenkoVolI}) which
are mentioned there. Note that the $\mathcal R$-simple germs in Table 
\ref{tab:ArnoldR} are stated as plane curve singularities in cases $D_k$ and 
$E_k$ and as fat point singularities in the case $A_k$. This is the smallest 
dimension in which they occur, but they also exist as simple singularities 
in any higher dimension by stable equivalence (i.e. using the generalized Morse lemma).

\begin{longtable}{|c|c|c|c|c|}
\hline Type & Normal form  & $\mu$ & $\tau$ & \\ \hline
$A_k$ & $x^{k+1}$  & $k$ & $k$ & $k \geq 1$ \\
\hline
$D_k$ & $x^2y+y^{k-1}$ & $k$ & $k$ & $k \geq 4$\\
\hline
$E_6$ & $x^3+y^4$ & 6 & 6 &\\
\hline
$E_7$ & $x^3+xy^2$ & 7 & 7 &\\
\hline
$E_8$ & $x^3+y^5$ & 8 & 8 &\\
\hline
\caption{The $\mathcal R$-simple germs from \cite{Arnold72}}
\label{tab:ArnoldR}
\end{longtable}

\noindent
A \textit{boundary singularity} is given by a germ $f \colon (\CC^n,0) \to (\CC,0)$ 
together with a ``boundary'' specified by first 
coordinate $(\CC^{n-1},0) = (\{x_1=0\},0) \subset (\CC^n,0)$. 
Associated to $f$ we now have two Milnor fibers $M_f$ and $M_{f|\{x_1=0\}}$ 
with a natural inclusion 
\[
  M_{f|\{x_1=0\}} \subset M_f.
\]
When both $f$ and $f|\{x_1=0\}$ have isolated singularity, then the 
Milnor fibers are homotopy equivalent to bouquets of spheres of dimension 
$n-1$ and $n-2$, respectively. Arnold has shown in 
\cite[Theorem 3]{Arnold78} that also the factor space $M_f/M_{f|\{x_1=0\}}$ 
has the homotopy type of a bouquet of $\mu_\partial(f)$ spheres of dimension $n-1$ 
where $\mu_\partial(f)$ is the \textit{boundary Milnor number}
\[
  \mu_\partial(f) = 
  \dim_{\CC} \CC\{x_1,\dots,x_n\}/
  \left\langle 
  x_1 \cdot \frac{\partial f}{\partial x_1},
  \frac{\partial f}{\partial x_2}, 
  \dots,
  \frac{\partial f}{\partial x_n}
  \right\rangle.
\]
All these Milnor numbers satisfy $\mu(f) + \mu(f|\{x_1=0\}) = \mu_\partial(f)$. 

\begin{longtable}{|c|c|c|c|c|c|}
  \hline
  Name	 & Normal form 	& $\mu(f)$ & $\mu(f|\{x=0\})$ & $\mu_\partial(f)$ &\\
  \hline
  $B_k$ &
  $\pm x^k \pm y^2$ &
  $k-1$ & $1$ & $k$ & $k \geq 2$\\
  \hline

  $C_k$ & 
  $xy \pm y^k $ & 
  $1$ & $k-1$ & $k$ & $k \geq 2$\\
  \hline

  $F_4$ & 
  $\pm x^2 + y^3$ & 
  $2$ & $2$ & $4$ &\\
  \hline
  \caption{The $\mathcal R_\delta$-simple germs from \cite{Arnold78}. 
    The boundary is given by $x=0$.
  }
\label{tab:ArnoldRDelta}
\end{longtable}

\subsection{Complete interesections}

Giusti proved in \cite{Giu} that apart from hypersurfaces simple complete 
intersections can only occur in two settings: fat points in the plane and
curves in $3$-space. He gave exhaustive lists of the simple singularities 
in these cases.  For completeness of the ICMC2 case below, we also include 
these tables here: 

\begin{longtable}{|c|c|c|c|c|}
\hline Type & Normal form  & $\mu$ & $\tau$ &  \\ \hline
$A_k$ & $\left<y,x^{k+1}\right>$  & $k$ & $k$ & $k\geq 1$ \\
\hline
$F_{q+r-1}^{q,r}$ & $\left<xy,x^q+y^r\right>$ & $q+r-1$ & $q+r$   &  $q,r\geq 2$ \\
\hline
$G_5$ & $\left<x^2,y^3\right>$ & 5 & 7 & \\
$G_7$ & $\left<x^2,y^4\right>$ & 7 & 10 & \\
\hline 
$H_{q+3}$ & $\left<x^2+y^{q},xy^2\right>$ & $q+3$ & $q+5$ &
$q\geq 3$\\
\hline $I_{2q-1}$ & $\left<x^2+y^3,y^q\right>$ & $2q-1$ & $2q+1$ &
$q\geq 4$\\
$I_{2r+2}$ & $\left<x^2+y^3,xy^r\right>$ & $2r+2$ & $2r+4$ & $r\geq 3$\\
\hline
\caption{Simple ICIS fat point singularities\label{tab:GiustiFatPoints}}
\end{longtable}

\begin{longtable}{|c|c|c|c|c|}
\hline Type & Normal form & $\mu$ & $\tau$ & 
\\ \hline
$S_{n+3}$ & $(x^2+y^2+z^{n},yz)$ & $n+3$ & $n+3$ &  $n\geq 2$ \\
\hline
$T_7$ & $(x^2+y^3+z^3,yz)$ & 7 & 7 &\\
\hline
$T_8$ & $(x^2+y^3+z^4,yz)$ & 8 & 8 &\\
\hline
$T_9$ & $(x^2+y^3+z^5,yz)$ & 9 & 9 & \\
\hline
$U_7$ & $(x^2+yz,xy+z^3)$ & 7 & 7 & \\
\hline
$U_8$ & $(x^2+yz+z^3,xy)$ & 8 & 8 & \\
\hline
$U_9$ & $(x^2+yz,xy+z^4)$ & 9 & 9 & \\
\hline
$W_8$ & $(x^2+z^3,y^2+xz)$ & 8 & 8 & \\
\hline
$W_9$ & $(x^2+yz^2,y^2+xz)$ & 9 & 9 & \\
\hline
$Z_9$ & $(x^2+z^3,y^2+z^3)$ & 9 & 9 & \\
\hline
$Z_{10}$ & $(x^2+yz^2,y^2+z^3)$ & 10 & 10 & \\
\hline
\caption{Simple ICIS space curve singularities\label{tab:GiustiCurve}}
\end{longtable}

\subsection{Simple square and symmetric matrices}

%
%

We list 
the tables mentioned in the results of Bruce in \cite{Bruce03} and Bruce and Tari \cite{BruceTari04}
on simple determinantal singularities defined by symmetric and arbitrary square matrices.
The original tables have been extended by several auxiliary results of related publications 
such as \cite{GoryunovZakalyukin03} and \cite{Goryunov21}. Also, we were informed 
by V. Goryunov about certain mistakes in the original classifications and we adopt 
his unified and corrected exposition from \cite{Goryunov21}.
The notations were at times streamlined with those from the tables in the other sections. \\

All simple square matrices of size $m=3$ in two variables turn out to be 
symmetric, whence Table \ref{tab:SquareSymm3p2} simultaneously states both
cases. The list of simple square matrices of size $m=3$ in seven variables
is also very closely related to a list of simple symmetric matrices: It can
be obtained from the simple symmetric matrices in four variables by adding
the skew symmetric matrix $U$ in three new variables as in Theorem 
\ref{thm:ClassificationOfSimpleSquareMatrices} 4. For square matrices the 
entries in Table \ref{tab:SquareSymm3p4} have to be understood in this sense.

\begin{longtable}{|c|c|c|c|c|c|}
  \hline 
  Normal form 	& associated 	& Tjurina 	& $\GL$-codim.	& $\GL$-codim.	& pair of \\
  		& hypersurface	& transform	& as symm. mat.	& as square mat.& Weyl groups \\
  \hline

  $
  \begin{pmatrix}
    y^k & x & 0 \\
    x & \pm y^l & 0 \\
    0 & 0 & y
  \end{pmatrix} 
  $ &
  $D_{k+l+2}$ & 
  $
  \begin{matrix}
    \textnormal{$A_0+ A_0 + A_0$ for $k=l$}\\
    \textnormal{$A_0 + A_{|l-k|-1}$ for $k\neq l$} 
  \end{matrix}
  $ & 
  $k+l+2$ &
  $2k+l+4$ &
  $(D_{k+l+2};D_{k+1} \oplus D_{l+1})$ \\
  \hline

  $
  \begin{pmatrix}
    0 & x & y \\
    x & y & 0 \\
    y & 0 & x^2
  \end{pmatrix}
  $ & 
  $E_6$ & 
  $A_0$ & 
  $6$ & 
  $9$ & 
  $(E_6;A_5 \oplus A_1)$ \\
  \hline

  $
  \begin{pmatrix}
    0 & x & y \\
    x & y & 0 \\
    y & 0 & xy
  \end{pmatrix}
  $ &
  $E_7$ & 
  $A_1$ & 
  $7$ & 
  $10$ & 
  $(E_7;D_6 \oplus A_1)$ \\
  \hline

  $
  \begin{pmatrix}
    0 & x & y \\
    x & y & 0 \\
    y & 0 & x^3
  \end{pmatrix}
  $ & 
  $E_8$ & 
  $A_2$ & 
  $8$ &
  $11$ &
  $(E_8;E_7 \oplus A_1)$ \\
  \hline

  $
  \begin{pmatrix}
    x & 0 & 0 \\
    0 & y & x \\
    0 & x & y^2
  \end{pmatrix}
  $ & 
  $E_7$ & 
  $A_0 + A_0$ & 
  $7$ &
  $11$ & 
  $(E_7;A_7)$ \\
  \hline

  $
  \begin{pmatrix}
    x & 0 & y^2 \\
    0 & y & x \\
    y^2 & x & 0
  \end{pmatrix}
  $ & 
  $E_8$ & 
  $A_0$ & 
  $8$ & 
  $12$ &
  $(E_8;D_8)$\\
  \hline

  \caption{Simple singularities of symmetric and square matrices of size $m=3$ in 
    two variables}
  \label{tab:SquareSymm3p2}
\end{longtable}

\begin{longtable}{|c|c|c|c|c|c|}
  \hline
  Name 	& Normal form &	$\GL$-codim.	& $\GL$-codim. of & odd   	& ass. symm. \\
       	&             &	as symm. mat.	& ass. square mat.& function	& ICIS \\
  \hline
  $
  \begin{matrix}
    I_{k+1}, \\
    k\geq 1
  \end{matrix}
  $ & 
  $
  \begin{pmatrix}
    x & 0 & z \\
    0 & y + x^k& w \\
    z & w & y
  \end{pmatrix}
  $ & 
  $k+1$ & 
  $k+1$ & 
  $
  \begin{matrix}
    D_{2k+2}/\ZZ_2\colon \\
    ac^2 + a^{2k+1}
  \end{matrix}
  $ & 
  $
  \begin{matrix}
    S_{2k+3}\colon \\
    c^2 + 2bc+a^{2k} \\
    ab
  \end{matrix}
  $\\
  \hline

  $II_{4}$ & 
  $
  \begin{pmatrix}
    x & w^2 & y \\
    w^2 & y & z \\
    y & z & w
  \end{pmatrix}
  $ & 
  $4$ & 
  $4$ &
  $
  \begin{matrix}
    E_8/\ZZ_2 \colon \\
    b^3+c^5
  \end{matrix}
  $ & 
  $
  \begin{matrix}
    U_9\colon \\
    b^2 - ac + c^4 \\
    ab -c^4
  \end{matrix}
  $\\
  \hline

  $II_{5}$ & 
  $
  \begin{pmatrix}
    x & 0 & y+w^2 \\
    0 & y & z \\
    y+w^2 & z & w
  \end{pmatrix}
  $ & 
  $5$ & 
  $5$ &
  $
  \begin{matrix}
    J_{10}/\ZZ_2 \colon \\
    b^3 - bc^4
  \end{matrix}
  $ & 
  $
  \begin{matrix}
    U_{11}\colon \\
    b^2 - ac + c^4 \\
    ab
  \end{matrix}
  $ \\
  \hline

  $II_{6}$ & 
  $
  \begin{pmatrix}
    x & w^3 & y \\
    w^3 & y & z \\
    y & z & w
  \end{pmatrix}
  $ & 
  $6$ &
  $6$ &
  $
  \begin{matrix}
    E_{12}/\ZZ_2 \colon \\
    b^3 + c^7
  \end{matrix}
  $ &
  $
  \begin{matrix}
    U_{13} \colon \\
    b^2 - ac \\
    ab -c^6
  \end{matrix}
  $ \\
  \hline
\caption{Simple singularities of symmetric matrices of size $m=3$ in four 
  variables (also providing simple singularities of square matrices of the
  same size in seven variables)}
\label{tab:SquareSymm3p4}
\end{longtable}

Table \ref{tab:NormalFormsSquareMatricesSize2In2Space} is dealing with curve
singularities and hence the Tjurina transform needs to be considered with 
special care, as the (strict) Tjurina transform is of the same dimension as the
exceptional locus 
(see Example \ref{exp:TjurinaTransformationInFamilyThreeCoordinateAxis}). The strict
Tjurina transforms are listed in the third column of this table. If the relative
position of the exceptional locus w.r.t. the strict Tjurina transform is of interest,
the singularity type of the reduced structure of the total Tjurina transform
is stated in brackets. The naming of all singularities is according to the
tables of simple plane curve and space curve singularities, i.e. Tables 
\ref{tab:ArnoldR}, \ref{tab:GiustiCurve} and  \ref{tab:FKCurve}. The only
exception is the notation $A_0$ for a smooth branch meeting the exceptional 
locus transversally, which we adopted from the original table of Bruce and 
Tari. 

For instance the entry ``$A_0,A_{k-3}$'' in the third column of the second row
indicates that there is an $A_{k-3}$ singularity in the strict Tjurina 
transform and additionally a smooth branch, both of which meet the exceptional
locus, but not in the same point. The entry ``$E_6(1) \;\;\;\; (U_7)$'' in 
the fifth row indicates that the strict Tjurina transform is a space curve
of type $E_6(1)$ from Table \ref{tab:FKCurve} (with parametrization 
$(t^3,t^4,t^5)$) and that it meets the exceptional locus to form an $U_7$
singularity from Giusti's list \ref{tab:GiustiCurve}.

  \begin{longtable}{|l|c|c|c|}
    \hline
    Normal form 	& hypersurface	& Tj. transf.	& $\GL$-codim. \\
    \hline
    $
    \begin{pmatrix}
      x & y^k \\ \pm y^l & x
    \end{pmatrix}
    ,\,
    1\leq k \leq l
    $ &
    $A_{k+l-1}$ &
    $
    \begin{matrix}
      \textnormal{$A_{l-k-1}$ for $k \neq l$ }\\
      \textnormal{2 $A_0$ for $k=l$,}
    \end{matrix} 
    $ &
    $2k+l-1$ \\
    \hline

    $
    \begin{pmatrix}
      x & y \\ x^2 \pm y^k & 0
    \end{pmatrix},\,
    2\leq k
    $ & 
    $D_{k+2}$ & 
    $A_0, A_{k-3}$ & 
    $k+3$ \\
    \hline

    $
    \begin{pmatrix}
      x & x^2\pm y^k \\ y & 0
    \end{pmatrix},\, 
    2\leq k
    $ & 
    $D_{k+2}$ & 
    $A_{k-1} \vee L \;\;\;\; (S_{k+3})$ & 
    $k+3$ \\
    \hline

    $
    \begin{pmatrix}
      x & y \\ y^3 & x^2
    \end{pmatrix}
    $ & 
    $E_6$ & 
    $A_0$ & 
    $7$ \\
    \hline

    $
    \begin{pmatrix}
      x & y^3 \\ y & x^2
    \end{pmatrix} 
    $ & 
    $E_6$ & 
    $E_6(1)  \;\;\;\;(U_7)$ & 
    $7$ \\
    \hline

    $
    \begin{pmatrix}
      x & y \\ xy^2 & x^2
    \end{pmatrix}
    $ & 
    $E_7$ & 
    $A_1$ & 
    $8$ \\
    \hline

    $
    \begin{pmatrix}
      x & xy^2 \\ y & x^2 
    \end{pmatrix}
    $ & 
    $E_7$ & 
    $E_7(1) \;\;\;\;(U_8)$ & 
    $8$ \\
    \hline

    $
    \begin{pmatrix}
      x & y \\ y^4 & x^2
    \end{pmatrix}
    $ & 
    $E_8$ & 
    $A_2$ & 
    $9$ \\
    \hline

    $
    \begin{pmatrix}
      x & y^4 \\ y & x^2 
    \end{pmatrix}
    $& 
    $E_8$ & 
    $E_8(1)\;\;\;\;(U_9)$ & 
    $9$ \\
    \hline

    $
    \begin{pmatrix}
      x & 0 \\ 0 & y^2 \pm x^k
    \end{pmatrix},\,
    2\leq k
    $ &
    $D_{k+2}$ & 
    $A_0,A_{k-1}$ & 
    $k+4$ \\
    \hline

    $
    \begin{pmatrix}
      x & 0 \\ 0 & xy + y^k
    \end{pmatrix},\,
    3\leq k
    $ & 
    $D_{2k}$ & 
    $A_0,A_1$ & 
    $3k$ \\
    \hline

    $
    \begin{pmatrix}
      x & y^k \\ \pm y^l & xy
    \end{pmatrix} ,\,
    3 \leq k \leq l
    $ & 
    $D_{k+l+1}$ & 
    $
    \begin{matrix}
      \textnormal{ $A_1$ for $k=l$, } \\
      \textnormal{3 $A_0$ for $k+1 = l$ } \\
      \textnormal{ $A_0+A_{l-k-2}$ for $k+1 < l$ } 
    \end{matrix}
    $ & 
    $2k+l+1$ \\
    \hline

    $
    \begin{pmatrix}
      x & \pm y^l \\ y^k & xy
    \end{pmatrix},\,
    3 \leq k < l$ & 
    $D_{k+l+1}$ & 
    $A_{l-k} \vee L \;\;\;\;(D_{l+k+3} \vee L)$ & 
    $2k+l+1$ \\
    \hline

    $
    \begin{pmatrix}
      x & y^2 \\ y^2 & x^2 
    \end{pmatrix}
    $ & 
    $E_6$ & 
    $A_2$ & 
    $8$ \\
    \hline 

    $
    \begin{pmatrix}
      x & y^2 \\ 0 & x^2 +y^3
    \end{pmatrix}
    $ &
    $E_7$ & 
    2 $A_0$ & 
    $9$ \\
    \hline 

    $
    \begin{pmatrix}
      x & 0 \\ y^2 & x^2 +y^3
    \end{pmatrix}
    $ & 
    $E_7$ & 
    $A_2 \vee L \;\;\;\;(S_6)$ & 
    $9$ \\
    \hline 

    $
    \begin{pmatrix}
      x & 0 \\ 0 & x^2 +y^3
    \end{pmatrix}
    $ & 
    $E_7$ & 
    $A_0, A_2$ & 
    $10$ \\
    \hline

    $
    \begin{pmatrix}
      x & y^2 \\ y^3 & x^2
    \end{pmatrix}
    $ & 
    $E_8$ & 
    $A_0$ & 
    $10$ \\
    \hline

    $
    \begin{pmatrix}
      x & y^3 \\ y^2 & x^2 
    \end{pmatrix}
    $& 
    $E_8$ &
    $E_6(1) \;\;\;\;(W_9)$ & 
    $10$ \\
    \hline
    \caption{Simple singularities of square matrices in two variables of size $m=2$, 
      \cite[Table 2]{BruceTari04}. }
      \label{tab:NormalFormsSquareMatricesSize2In2Space}
    \end{longtable}

    \begin{longtable}{|l|c|c|c|}
      \hline
      Normal form 	& hypersurface 	& Tj. transf.	& $\GL$-codimension \\
      \hline 
      $
      \begin{pmatrix}
	x & z^k \\ \pm z^l & y
      \end{pmatrix}, \,
      1 \leq k \leq l$ & 
      $A_{k+l-1}$ & 
      $A_{k-1},A_{l-1}$ & 
      $k+l-1$ \\
      \hline 

      $
      \begin{pmatrix}
	x & -y \\ y+z^k & x
      \end{pmatrix}
      $ & 
      $A_{2k-1}$ & 
      $2A_0$ & 
      $2k-1$ \\
      \hline 

      $
      \begin{pmatrix}
	x & y \\ z^2 \pm y^k & x
      \end{pmatrix},\,
      2\leq k$ & 
      $D_{k+2}$ &
      $A_0, D_{k+1}$ & 
      $k+2$ \\
      \hline

      $
      \begin{pmatrix}
	x & y \\ y^2 & x + z^2 
      \end{pmatrix} 
      $ & 
      $E_6$ &
      $A_0, D_5$ & 
      $6$ \\
      \hline

      $
      \begin{pmatrix}
	x & y \\ y^2 + z^3 & x
      \end{pmatrix} 
      $ & 
      $E_7$ & 
      $A_0, E_6$ & 
      $7$ \\
      \hline

      $ 
      \begin{pmatrix}
	x & y \\ yz & x+z^k
      \end{pmatrix},\,
      2\leq k$ & 
      $D_{2k+1}$ & 
      $A_0, A_{k-1}$ & 
      $2k+1$ \\
      \hline

      $
      \begin{pmatrix}
	x & y \\ yz + z^k & x
      \end{pmatrix}
      ,\, 3 \leq k$ & 
      $D_{2k}$ & 
      $A_0, A_{k-1}$ & 
      $2k$ \\
      \hline
      \caption{Simple singularities of square matrices in three variables of size $m=2$,
	\cite[Table 3]{BruceTari04}}
	\label{tab:NormalFormsSquareMatricesSize2In3Space}
      \end{longtable}

      \begin{longtable}{|l|c|c|c|}
	\hline
	Normal form 	& hypersurface 	& Tj. transf. 	& $\GL$-codimension \\
	\hline 
	$
	\begin{pmatrix}
	  x & y^k & 0 \\
	  \pm y^l & x & 0 \\
	  0 & 0 & y
	\end{pmatrix},\,
	1\leq k \leq l $ &
	$D_{k+l+2}$ & 
	$ 
	\begin{matrix}
	  \textnormal{2 $A_0$ for $l=k$ }\\
	  \textnormal{$A_0 + A_{l-k-1}$ for $l \neq k$} 
	\end{matrix} 
	$ & 
	$2k+l+4$ \\
	\hline

	$
	\begin{pmatrix}
	  x & y & 0 \\
	  0 & x & y \\
	  y^2 & 0 & x
	\end{pmatrix}
	$ & 
	$E_6$ & 
	$A_0$ & 
	$9$ \\
	\hline

	$
	\begin{pmatrix}
	  x & y & 0 \\
	  0 & x & y \\
	  xy & 0 & x
	\end{pmatrix}
	$ & 
	$E_7$ & 
	$A_1$ & 
	$10$ \\
	\hline 

	$
	\begin{pmatrix}
	  x & y & 0 \\
	  y^2 & x & 0 \\
	  0 & 0 & x
	\end{pmatrix}
	$ & 
	$E_7$ & 
	2 $A_0$ & 
	$11$ \\
	\hline 

	$ 
	\begin{pmatrix}
	  x & y & 0 \\
	  0 & x & y^2 \\
	  y^2 & 0 & x
	\end{pmatrix}
  $ & 
  $E_8$ & 
  $A_0$ & 
  $12$ \\
  \hline
\caption{Simple singularities of square matrices in two variables of size $m=3$,
  \cite[Table 4]{BruceTari04}}
\label{tab:NormalFormsSquareMatricesSize3In2Space}
\end{longtable}

%
%



\begin{longtable}{|c|c|c|c|c|c|}
  \hline 
  Normal form	& associated	& Tjurina	& $\mathcal K$-type	& $\GL$-codim.	& Pair of 	\\
             	& hypersurface	& transform	&                  	&             	& Weyl groups 	\\
  \hline

  $
  \begin{pmatrix}
    y^k & x \\ x & y^l
  \end{pmatrix},\,
  k\geq 1, l \geq 2$ &
  $A_{k+l-1}$ &
  $
  \begin{matrix}
    \textnormal{ $A_0 + A_0$ for $k=l$ }\\
    \textnormal{ $A_{|l-k|-1}$ for $k \neq l$} 
  \end{matrix}
  $ & 
  $A_{p-1}$ & 
  $k+l-1$ & 
  $(A_{k+l-1}; A_{k-1} \oplus A_{l-1} )$ \\
  \hline

  $
  \begin{pmatrix}
    x & 0 \\ 0 & y^2 + x^k
  \end{pmatrix},\,
  k\geq 2$ &
  $D_{k+2}$ & 
  $A_0 + A_{k-1}$ & 
  $A_1$ & 
  $k+2$ & 
  $(D_{k-2};D_{k-3})$ \\
  \hline

  $
  \begin{pmatrix}
    x & 0 \\ 0 & xy +y^k
  \end{pmatrix},\,
  k\geq 2$ & 
  $D_{2k}$ & 
  $A_0 + A_1$ & 
  $A_{k-1}$ & 
  $2k$ & 
  $(D_{2k};A_{2k-1})$ \\
  \hline

  $
  \begin{pmatrix}
    x & y^k \\ y^k & xy
  \end{pmatrix},\,
  k\geq 2$ &
  $D_{2k+1}$ & 
  $A_1$ & 
  $A_{k-1}$ & 
  $2k+1$ & 
  $(D_{2k+1}; A_{2k})$ \\
  \hline

  $
  \begin{pmatrix}
    x & y^2 \\ y^2 & x^2
  \end{pmatrix}
  $ & 
  $E_6$ & 
  $A_2$ & 
  $A_1$ & 
  $6$ & 
  $(E_6;D_5)$ \\
  \hline

  $
  \begin{pmatrix}
    x & 0 \\
    0 & x^2+y^3
  \end{pmatrix}
  $ & 
  $E_7$ & 
  $A_0 + A_2$ & 
  $A_2$ & 
  $7$ & 
  $(E_7;E_6)$ \\
  \hline

\caption{Simple singularities of symmetric matrices of size $m=2$ in two variables,
  \cite[Theorem 1.1]{Bruce03}, extended according to \cite[Table 1]{GoryunovZakalyukin03}}
\label{tab:SquareSymm2p2}
\end{longtable}

%
%

\subsection{Skew-symmetric matrices}

For Haslinger's result \cite{Haslinger01}, we only reproduce his table of 
simple singularities given by skew-symmetric $4 \times 4$ matrices in $2$
variables. 
All matrices in the table are of the block form 
\[
  A = 
  \left( 
  \begin{array}[h]{cc}
    0 & B \\
    -B^T & 0 
  \end{array}
  \right)
\]
for some $2\times2$-matrix $B$ so that the associated 
hypersurface singularity is given by $f = \Pf A = \det B$. 

\begin{longtable}{|c|c|c|c|c|}
\hline
Name & matrix $B$ & & hypersurface & $\mathcal G_\sk$-codimension \\
\hline
$B_{kl}$ & 
$
\begin{pmatrix}
  x & y^k \\ y^l & x
\end{pmatrix}
$ & 
$ 1 \leq k \leq l $ & 
$A_{k+l-1}$ & 
$4k+l-1$ \\

$S_k$ & 
$
\begin{pmatrix}
  x & xy \\ y & x^k 
\end{pmatrix}
$ & 
$k \geq 2$ & 
$D_{k+2}$ & 
$k+5$ \\

$M_9$ & 
$
\begin{pmatrix}
  x & y^3 \\ y & x^2 
\end{pmatrix}
$ & &
$E_6$ & 
$9$ \\

$M_{10}$ & 
$
\begin{pmatrix}
  x & xy^2 \\ y & x^2 
\end{pmatrix} 
$ & &
$E_7$ & 
$10$ \\

$M_{11}$ &  
$
\begin{pmatrix}
  x & y^4 \\ 
  y & x^2
\end{pmatrix}
$ & & 
$E_8$ & 
$11$ \\

$F_k$ & 
$
\begin{pmatrix}
  x & 0 \\ 
  0 & y^2 + x^k 
\end{pmatrix}
$ & 
$k \geq 2$ & 
$D_{k+2}$ & 
$k+8$ \\

$G_k$ & 
$
\begin{pmatrix}
  x & 0 \\
  0 & xy + y^k
\end{pmatrix}
$ & 
$k\geq3$ & 
$D_{2k}$ & 
$5k$ \\

$H_{kl}$ & 
$
\begin{pmatrix}
  x & y^k \\ y^l & xy
\end{pmatrix}
$ & 
$2 \leq k \leq l$ & 
$D_{k+l+1}$ & 
$4k+l+1$ \\

$T_{12}$ & 
$
\begin{pmatrix}
  x & y^2 \\ y^2 & x^2
\end{pmatrix} 
$ & &
$E_6$ & 
$12$ \\

$T_{13}$ & 
$
\begin{pmatrix}
  x & y^2 \\ 0 & x^2 + y^3
\end{pmatrix}
$ & &
$E_8$ & 
$14$ \\

$T_{16}$ & 
$
\begin{pmatrix}
  x & 0 \\ 0 & x^2 + y^3
\end{pmatrix}
$ & & 
$E_7$ & 
$16$\\

\hline
\caption{Simple singularities, skew-symmetric matrices, case m=4}
\label{tab:SquareSkewm4}
\end{longtable}

Note that the classification of determinantal singularities of
skew-symmetric matrices is incomplete and hence the list is only exhaustive
for the given case, but does not exclude the existence of simple singularities
for other values of $m$ and $p$.

\subsection{Cohen-Macaulay codimension $2$ singularities}

In this case the lists reproduced here are extracted from the articles of 
Fr\"uhbis-Kr\"uger \cite{FruehbisKrueger99} and of Fr\"uhbis-Kr\"uger and Neumer 
\cite{FruehbisKruegerNeumer10}.
Together with Giusti's lists of simple ICIS, these lists are complete for the 
simple isolated Cohen-Macaulay codimension $2$ singularities:

\begin{longtable}{|c|c|c|c|}
\hline
$\Xi_k$ & $\begin{pmatrix} x & y & 0 \\ 0 & x^k & y \end{pmatrix}$ & $k+3$ & $k\geq 1$  \\
\hline
\caption{Simple non-ICIS fat point singularities in the plane\label{tab:FKNFatPoints} }
\end{longtable}

The list of simple space curve singularities from \cite{FruehbisKrueger99} reads: 
\begin{longtable}{|c|c|c|c|}
\hline
Type & Normal form  & $\mu$ & $\tau$  
\\ \hline
\parbox{1.5cm}{\begin{center}$A_{k-3}\vee L$\\ $k\geq 4$\end{center}} &
$\begin{pmatrix} z & y & x^{k-3} \\ 0 & x & y \end{pmatrix}$ & $k-2$ & $k-1$  \\
\hline
$E_6(1)$ &  $\begin{pmatrix} z & y & x^2 \\ x & z & y \end{pmatrix}$  & 4 & 5  \\
\hline
$E_7(1)$ & $\begin{pmatrix} z+x^2 & y & x \\ 0 & z & y \end{pmatrix}$ & 5 & 6   \\
\hline
$E_8(1)$ & $\begin{pmatrix} z & y & x^3 \\ x & z & y \end{pmatrix}$ & 6 & 7  \\
\hline
$J_{2,0}(2)$ &
$\begin{pmatrix} z+x^2 & y & x^{2} \\ 0 & z & y \end{pmatrix}$ & 6 & 7 \\
\hline
$J_{2,1}(2)$ &
$\begin{pmatrix} z+x^2 & y & x^{3} \\ 0 & z & y \end{pmatrix}$ & 7 & 8 \\
\hline
$E_{12}(2)$ & $\begin{pmatrix} z & y & x^3 \\ x^2 & z & y \end{pmatrix}$ & 8 & 9 \\
\hline
\parbox{1.5cm}{\begin{center}$D_{k+4}\vee L$ \\ $k\geq 0$ \end{center} } & $\begin{pmatrix} z & 0
& x^{k+2}-y^2 \\ 0 & x & y
\end{pmatrix}$
 &  $k+5$ & $k+6$ \\
\hline
$E_{6}\vee L$ & $\begin{pmatrix} z & -y^2 & -x^3 \\ 0 & x & y \end{pmatrix}$ &  7 & 8   \\
\hline
$E_{7}\vee L$ & $\begin{pmatrix} z & x^3-y^2 & 0 \\ 0 & x & y \end{pmatrix}$ &  8 & 9 \\
\hline
$E_{8}\vee L$ & $\begin{pmatrix} z & -y^2 & -x^4 \\ 0 & x & y \end{pmatrix}$ &  9 & 10 \\
\hline
$S_6^*$ & $\begin{pmatrix} z & x & y \\ 0 & y & x^2-z^2 \end{pmatrix}$ & 6 & 7 \\
\hline
$T_7^*$ & $\begin{pmatrix} z & x & y \\ 0 & y & x^2-z^3 \end{pmatrix}$ &  7 & 8  \\
\hline
$U_7^*$ & $\begin{pmatrix} z & xy & x^2 \\ x & z & y \end{pmatrix}$ &  7 & 8  \\
\hline
$W_8^*$ & $\begin{pmatrix} z & y^2 & x^2 \\ x & z & y \end{pmatrix}$ &  8 & 9 \\
\hline
\caption{Simple non-ICIS space curve singularities \label{tab:FKCurve}}
\end{longtable}

The list of simple ICMC2 surface singularities from \cite{FruehbisKruegerNeumer10} is 
given in Table \ref{tab:FKNSurface} below. The Tjurina number $\tau$ is equal to both the $\GL$-Tjurina number 
of the defining matrix $A \in \CC\{x,y,z,w\}^{2\times 3}$, as well as the Tjurina number 
of the associated germ $(X_A^2,0) \subset (\CC^4,0)$ according to Corollary 
\ref{cor:FruehbisKruegerMatrixT1}. This list recovers the rational triple points 
from \cite{Tjurina68} and we give the corresponding name in the last column. 
The last three entries of the list are the nameless sporadic members  
from \cite{Artin66} in the same order as there and in \cite{Tjurina68}.
The Milnor numbers of the smoothings (i.e. the second Betti number of $M_A^2$) 
can be computed as $\mu = \tau -1$ by virtue 
of the ``if''-part of Wahl's conjecture \ref{con:WahlsConjecture}.
\begin{longtable}{|c|l|l|c|l|}
\hline
 Type
    & Normal form
    &
    & $\tau$ 
    & Name of Triple
\\  
    &
    &
    &
    & Point in \cite{Tjurina68}
\\ \hline
$\Lambda_{1,1}$
    & $\begin{pmatrix}w & y & x \\ z & w & y \end{pmatrix}$
    &
    & 2
    & $A_{0,0,0}$
\\ \hline
$\Lambda_{k,1}$
    & $\begin{pmatrix}w & y & x  \\ z & w & y^k \end{pmatrix}$
    & $k\geq 2$
    & $k+1$
    & $A_{0,0,k-1}$
\\ \hline
$\Lambda_{k,l}$
    & $\begin{pmatrix}w^l & y & x \\ z & w & y^k \end{pmatrix}$
    & $k\geq l\geq 2$
    & $k+l$
    & $A_{0,l-1,k-1}$
\\ \hline

    & $\begin{pmatrix}z & y & x \\ x & w & y^2+z^k \end{pmatrix}$
    & $k\geq 2$
    & $k+3$
    & $C_{k+1,0}$
\\ \hline
    & $\begin{pmatrix}z & y & x \\ x & w & yz+y^kw \end{pmatrix}$
    & $k\geq 1$
    & $2k+4$
    & $B_{2k+2,0}$
\\
    & $\begin{pmatrix}z & y & x \\ x & w & yz+y^k \end{pmatrix}$
    & $k\geq 3$
    & $2k+1$
    & $B_{2k-1,0}$
\\ \hline
    & $\begin{pmatrix}z & y & x \\ x & w & z^2+yw \end{pmatrix}$
    &
    & 7
    & $D_0$
\\ \hline
    & $\begin{pmatrix}z & y & x \\ x & w & z^2+y^3 \end{pmatrix}$
    &
    & 8
    & $F_0$
\\ \hline
    & $\begin{pmatrix} z & y+w^l & w^m \\ w^k & y & x \end{pmatrix}$
    & $k,l,m \geq 2$
    & $k+l+m-1$
    & $A_{k-1,l-1,m-1}$
\\ \hline
    & $\begin{pmatrix}z & y & x^l+w^2 \\ w^k & x & y\end{pmatrix}$
    & $k,l\geq 2$
    & $k+l+2$
    & $C_{l+1,k-1}$
\\ \hline
    & $\begin{pmatrix}z & y+w^l & xw \\ w^k & x & y \end{pmatrix}$
    & $k,l\geq 2$
    & $k+2l+1$
    & $B_{2l,k-1}$
\\
    & $\begin{pmatrix}z & y & xw+w^l \\ w^k & x & y \end{pmatrix}$
    & $k\geq 2, l\geq 3$
    & $k+2l$
    & $B_{2l+1,k-1}$
\\ \hline
    & $\begin{pmatrix}z & y+w^2 & x^2 \\ w^k & x & y \end{pmatrix}$
    & $k\geq 2$
    & $k+6$
    & $D_{k-1}$
\\ \hline
    & $\begin{pmatrix}z & y & x^2+w^3 \\ w^k & x & y \end{pmatrix}$
    & $k\geq 2$
    & $k+7$
    & $F_{k-1}$
\\ \hline
    & $\begin{pmatrix}z & y & xw+w^k  \\ y & x & z \end{pmatrix}$
    &
    & $3k+1$
    & $H_{3k}$
\\
    & $\begin{pmatrix}z & y & xw \\ y & x & z+w^k \end{pmatrix}$
    &
    & $3k+2$
    & $H_{3k+1}$
\\ 
    & $\begin{pmatrix}z & y & xw \\ y+w^k & x & z \end{pmatrix}$
    &
    & $3k+3$
    & $H_{3k+2}$
\\ \hline
    & $\begin{pmatrix}z & y & w^2  \\ y & x & z+x^2 \end{pmatrix}$
    &
    & 8
    & 
\\ \hline
    & $\begin{pmatrix}z & y & x^2  \\ y & x & z+w^2 \end{pmatrix}$
    &
    & 9
    & 
\\ \hline
    & $\begin{pmatrix}z & y & x^3+w^2  \\ y & x & z \end{pmatrix}$
    &
    & 9
    &
\\ \hline
\caption{Simple normal surface singularities in $({\mathbb C}^4,0)$\label{tab:FKNSurface}}
\end{longtable}
\noindent

\medskip

For the simple ICMC2 threefold singularities we extend the list from 
\cite{FruehbisKruegerNeumer10} by the middle Betti number of the smoothing 
as in \cite{FruehbisKruegerZach20}. Recall that the second Betti number is always 
equal to $1$, cf. Theorem \ref{thm:MainTheoremFromFKZ}.
\begin{center}
\begin{longtable}[h]{|c|c|p{2cm}|c|}
    \hline
    Transpose of the & Tjurina & Singularities &  \\
    presentation matrix $A$ & number $\tau$ & in Tj.-transf. &  
    $b_3(M_A^2)$  \\
    \hline
    $
    \begin{pmatrix}
        x & y & z \\
        v & w & x
    \end{pmatrix}
    $ &
    $1$ &
    - &
    $0$ \\
    \hline

    $
    \begin{pmatrix}
        x & y & z \\
        v & w & x^{k+1} + y^2
    \end{pmatrix}
    $ & 
    $k+2$ & 
    $A_k$ & 
    $k$  \\
    $
    \begin{pmatrix}
        x & y & z \\
        v & w & x y^2 + x^{k-1} 
    \end{pmatrix}
    $ & 
    $k+2$ & 
    $D_{k}$ &
    $k$\\

    $
    \begin{pmatrix}
        x & y & z \\
        v & w & x^3 + y^4
    \end{pmatrix}
    $ & 
    $8$ & 
    $E_6$ & $6$ \\

    $
    \begin{pmatrix}
        x & y & z \\
        v & w & x^3+xy^3
    \end{pmatrix}
    $ & 
    $9$ & 
    $E_7$ &
    $7$ \\

    $
    \begin{pmatrix}
        x & y & z \\
        v & w & x^3 + y^5
    \end{pmatrix}
    $ & 
    $10$ & 
    $E_8$ &
    $8$ \\
    \hline

    $
    \begin{pmatrix}
        w & y & x \\
        z & w & y+v^k
    \end{pmatrix}
    $ &
    $2k-1$ & 
    - & 
    $0$ \\

    $
    \begin{pmatrix}
        w & y & x \\
        z & w & y^k + v^2
    \end{pmatrix}
    $ &
    $k+2$ & 
    $A_{k-1}$ &
    $k-1$ \\

    $
    \begin{pmatrix}
        w & y & x \\
        z & w & yv + v^k 
    \end{pmatrix}
    $ &
    $2k$ & 
    $A_1$ &
    $1$ \\

    $
    \begin{pmatrix}
        w+v^k & y & x \\
        z & w & yv
    \end{pmatrix} 
    $ &
    $2k+1$ &
    $A_1$ &
    $1$ \\

    $
    \begin{pmatrix}
        w+v^2 & y & x \\
        z & w & y^2 + v^k
    \end{pmatrix}
    $ & 
    $k+3$ & 
    $A_{k-1}$ &
    $k-1$ \\

    $
    \begin{pmatrix}
        w & y & x \\
        z & w & y^2 + v^3
    \end{pmatrix}
    $ &
    $7$ &
    $A_2$ & 
    $2$\\

    \hline

    $
    \begin{pmatrix}
        v^2 + w^k & y & x \\
        z & w & v^2 + y^l
    \end{pmatrix}
    $
    & $k+l+1$ &
    $A_{k-1}$, $A_{l-1}$ & 
    $k+l-2$ \\

    $
    \begin{pmatrix}
        v^2 + w^k & y & x \\
        z & w & yv
    \end{pmatrix}
    $
    & $k+4$ &
    $A_{k-1}$, $A_1$ &
    $k$ \\

    $
    \begin{pmatrix}
        v^2 + w^k & y & x \\
        z & w & y^2+v^l
    \end{pmatrix}
    $
    & $k+l+2$ &
    $A_{k-1}$, $A_{l-1}$ &
    $k+l-2$ \\ 

    \hline
    $
    \begin{pmatrix}
        wv+ v^k & y & x \\
        z & w & yv + v^k
    \end{pmatrix}
    $
    & $2k+1$ &
    $A_1$, $A_1$ & 
    $2$ \\

    $
    \begin{pmatrix}
        wv + v^k & y & x \\
        z & w & yv
    \end{pmatrix}
    $
    & $2k+2$ &
    $A_1$, $A_1$ & 
    $2$ \\
    \hline

    $
    \begin{pmatrix}
        wv + v^3 & y & x \\
        z & w & y^2+v^3
    \end{pmatrix}
    $
    & $8$ &
    $A_1$, $A_2$ & 
    $3$ \\

    $
    \begin{pmatrix}
        wv & y & x \\
        z & w & y^2 + v^3
    \end{pmatrix}
    $
    & $9$ &
    $A_1$, $A_2$ &
    $3$ \\

    $
    \begin{pmatrix}
        w^2 + v^3 & y & x \\
        z & w & y^2+v^3
    \end{pmatrix}
    $
    & $9$ &
    $A_2$, $A_2$ &
    $4$ \\
    \hline

    $
    \begin{pmatrix}
        z & y & x \\
        x & w & v^2 + y^2 + z^k
    \end{pmatrix}
    $ 
    & 
    $k+4$ &
    $D_{k+1}$ &
    $k+1$ \\

    $
    \begin{pmatrix}
        z & y & x \\
        x & w & v^2 + yz + y^k w
    \end{pmatrix}
    $ & 
    $2k+5$ & 
    $A_{2k+2}$ &
    $2k+2$ \\

    $
    \begin{pmatrix}
        z & y & x \\
        x & w & v^2 + yz + y^{k+1}
    \end{pmatrix}
    $ & 
    $2k+4$ & 
    $A_{2k+1}$  &
    $2k+1$ \\

    $
    \begin{pmatrix}
        z & y & x \\
        x & w & v^2 + yw + z^2
    \end{pmatrix}
    $ & 
    $8$ &
    $D_5$ & 
    $5$ \\

    $
    \begin{pmatrix}
        z & y & x \\
        x & w & v^2 + y^3 + z^2
    \end{pmatrix}
    $ & 
    $9$ &
    $E_6$ &
    $6$ \\

    $
    \begin{pmatrix}
        z & y & x + v^2 \\
        x & w & vy + z^2
    \end{pmatrix}
    $ & 
    $7$ &
    $D_3$ &
    $3$ \\

    $
    \begin{pmatrix}
        z & y & x + v^2 \\
        x & w & vz + y^2
    \end{pmatrix}
    $ & 
    $8$ & 
    $A_4$ &
    $4$ \\

    $
    \begin{pmatrix}
        z & y & x + v^2 \\
        x & w & z^2 + y^2
    \end{pmatrix}
    $ & 
    $9$ &
    $D_5$ &
    $5$ \\
    \hline
\caption{Simple 3-fold singularities in $({\mathbb C}^5,0)$}
  \label{tab:FKN3Fold}
\end{longtable}
\end{center}

\begin{longtable}{|c|l|c|c|c|}
\hline
     Type
    & Normal form $A$
    &
    & $\tau_\GL(A)$
    & $\tau(X_{A}^1,0)$
\\ \hline
$\Omega_1$
    & $\begin{pmatrix}x & y & v \\ z & w & u \end{pmatrix}$
    &
    & 0
    & 0
\\ \hline
$\Omega_k$
    & $\begin{pmatrix}x & y & v \\ z & w & x+u^k \end{pmatrix}$
    & $k\geq 2$
    & $k-1$
    & $k-1$ 
\\ \hline
$A_k^\sharp$
    & $\begin{pmatrix}x & y & z \\ w & v & u^2+x^{k+1}+y^2 \end{pmatrix}$
    & $k\geq 1$
    & $k+2$
    & $1$ 
\\ \hline
   $D_k^\sharp$
    & $\begin{pmatrix}x & y & z \\ w & v & u^2+xy^2+x^{k-1} \end{pmatrix}$
    & $k\geq 4$
    & $k+2$
    & $1$ 
\\ \hline
   $E_6^\sharp$
    & $\begin{pmatrix}x & y & z \\ w & v & u^2+x^3+y^4 \end{pmatrix}$
    &
    & 8
    & $1$ 
\\ \hline
   $E_7^\sharp$
    & $\begin{pmatrix}x & y & z \\ w & v & u^2+x^3+xy^3 \end{pmatrix}$
    &
    & 9
    & $1$ 
\\ \hline
   $E_8^\sharp$
    & $\begin{pmatrix}x & y & z \\ w & v & u^2+x^3+y^5 \end{pmatrix}$
    &
    & 10
    & $1$ 
\\ \hline
    & $\begin{pmatrix}x & y & z \\ w & v & ux+y^k+u^l \end{pmatrix}$
    & $k\geq 2, l\geq 3$
    & $k+l-1$
    & $l-1$
\\ \hline
    & $\begin{pmatrix}x & y & z \\ w & v & x^2+y^2+u^3 \end{pmatrix}$
    &
    & 6
    & $2$ 
\\ \hline
   $F_{q,r}^\sharp$
    & $\begin{pmatrix}w & y & x \\ z & w+vu & y+v^q+u^r \end{pmatrix}$
    & $q,r\geq 2$
    & $q+r$
    & $q+r$ 
\\ \hline
   $G_5^\sharp$
    & $\begin{pmatrix}w & y & x \\ z & w+v^2 & y+u^3 \end{pmatrix}$
    &
    & 7
    & $7$ 
\\ \hline
   $G_7^\sharp$
    & $\begin{pmatrix}w & y & x \\ z & w+v^2 & y+u^4 \end{pmatrix}$
    &
    & 10
    & $10$ 
\\ \hline
   $H_{q+3}^\sharp$
    & $\begin{pmatrix}w & y & x \\ z & w+v^2+u^q & y+vu^2 \end{pmatrix}$
    & $q\geq 3$
    & $q+5$
    & $q+5$
\\ \hline
   $I_{2q-1}^\sharp$
    & $\begin{pmatrix}w & y & x \\ z & w+v^2+u^3 & y+u^q \end{pmatrix}$
    & $q\geq 4$
    & $2q+1$
    & $2q+1$
\\ \hline
   $I_{2r+2}^\sharp$
    & $\begin{pmatrix}w & y & x \\ z & w+v^2+u^3 & y+vu^r \end{pmatrix}$
    & $r\geq 3$
    & $2r+4$
    & $2r+4$
\\ \hline
    & $\begin{pmatrix}w & y & x \\ z & w+v^{k_1}+u^{k_2} & y^l+uv \end{pmatrix}$
    & $k_1,k_2,l\geq 2$
    & $k_1+k_2+l-1$
    & $k_1 + k_2$
\\ \hline
    & $\begin{pmatrix}w & y & x \\ z & w+v^2 & u^2+yv \end{pmatrix}$
    &
    & 6
    & $4$
\\ \hline
    & $\begin{pmatrix}w & y & x \\ z & w+uv & u^2+yv+v^k \end{pmatrix}$
    & $k\geq 3$
    & $k+4$
    & $k+2$
\\ \hline
    & $\begin{pmatrix}w & y & x \\ z & w+v^k & u^2+yv+v^3 \end{pmatrix}$
    & $k\geq 3$
    & $2k+2$
    & $2k+1$
\\
    & $\begin{pmatrix}w & y & x \\ z & w+uv^k & u^2+yv+v^3 \end{pmatrix}$
    & $k\geq 2$
    & $2k+5$
    & $2k+4$
\\ \hline
    & $\begin{pmatrix}w & y & x \\ z & w+v^3 & u^2+yv \end{pmatrix}$
    &
    & 9
    & $7$
\\ \hline
    & $\begin{pmatrix}w & y & x \\ z & w+v^k & u^2+y^2+v^3 \end{pmatrix}$
    & $k\geq 3$
    & $2k+3$
    & $2k+1$
\\
    & $\begin{pmatrix}w & y & x \\ z & w+uv^k & u^2+y^2+v^3 \end{pmatrix}$
    & $k\geq 2$
    & $2k+6$
    & $2k+4$
\\ \hline
\caption{Simple 4-fold singularities in $({\mathbb C}^6,0)$\label{tab:FKN4Fold}}
\end{longtable}



\section*{Acknowledgements}

The authors wish to thank V. Goryunov for careful reading, several hints and 
useful discussion on parts of an earlier version of the manuscript, and also 
G. Pe\~nafort Sanchis, W. Ebeling, H. M\o{}ller Pedersen and D. Kerner for 
fruitful discussions.

\bibliographystyle{plain}
\bibliography{./Bibtex/sources} 

\begin{thebibliography}{100}

\bibitem{AkinBuchsbaumWeyman81}
K.~Akin, D.A. Buchsbaum, and J.~Weyman.
\newblock Resolutions of determinantal ideals: the submaximal minors.
\newblock {\em Adv. in Math.}, 39(1):1--30, 1981.

\bibitem{Arnold72}
V.~Arnold.
\newblock Normal forms for functions near degenerate critical point.
\newblock {\em FAA}, 6:254--272, 1972.

\bibitem{ArnoldGuseinZadeVarchenkoVolI}
V.~Arnold, S.M. Gusein-Zade, and A.~Varchenko.
\newblock {\em Singularities of Differentiable Maps, Volume I}.
\newblock Birkh\"auser, 1985.

\bibitem{Arnold78}
V.~I. Arnold.
\newblock Critical points of functions on a manifold with boundary, the simple
  {L}ie groups {$B_{k}$}, {$C_{k}$}, {$F_{4}$} and singularities of evolutes.
\newblock {\em Uspekhi Mat. Nauk}, 33(5(203)):91--105, 237, 1978.

\bibitem{Artin66}
M.~Artin.
\newblock On isolated rational singularities of surfaces.
\newblock {\em Amer. J. Math.}, 88:129--136, 1966.

\bibitem{Artin74}
M.~Artin.
\newblock Algebraic construction of {B}rieskorn's resolutions.
\newblock {\em J. Algebra}, 29:330--348, 1974.

\bibitem{BelitskiiKerner16}
G.~Belitskii and D.~Kerner.
\newblock Group actions on filtered modules and finite determinacy. {F}inding
  large submodules in the orbit by linearization.
\newblock {\em C. R. Math. Acad. Sci. Soc. R. Can.}, 38(4):113--153, 2016.

\bibitem{BelitskiiKerner19}
G.~Belitskii and D.~Kerner.
\newblock Finite determinacy of matrices over local rings. {T}angent modules to
  the miniversal deformation for {$R$}-linear group actions.
\newblock {\em J. Pure Appl. Algebra}, 223(3):1288--1321, 2019.

\bibitem{BraChaRua15}
J.P. Brasselet, N.~Chachapoyas, and M.A.S Ruas.
\newblock Generic sections of essentially isolated determinantal singularities.
\newblock 2015.

\bibitem{Brieskorn71}
E.~Brieskorn.
\newblock Singular elements of semi-simple algebraic groups.
\newblock In {\em Actes du {C}ongr\`es {I}nternational des {M}ath\'{e}maticiens
  ({N}ice, 1970), {T}ome 2}, pages 279--284. 1971.

\bibitem{Bruce03}
J.~W. Bruce.
\newblock On families of symmetric matrices.
\newblock volume~3, pages 335--360, 741. 2003.
\newblock Dedicated to Vladimir I. Arnold on the occasion of his 65th birthday.

\bibitem{BruceTari04}
J.~W. Bruce and F.~Tari.
\newblock On families of square matrices.
\newblock {\em Proc. London Math. Soc. (3)}, 89(3):738--762, 2004.

\bibitem{BrunsHerzog93}
W.~Bruns and J.~Herzog.
\newblock {\em Cohen-Macaulay rings}.
\newblock Cambridge University Press, 1993.

\bibitem{BrunsVetter88}
W.~Bruns and U.~Vetter.
\newblock {\em Determinantal Rings}, volume 1327 of {\em Lecture Notes in
  Mathematics}.
\newblock Springer, 1988.

\bibitem{Buchsbaum79}
D.A. Buchsbaum.
\newblock {A} {N}ew {C}onstruction of the {E}agon-{N}orthcott {C}omplex.
\newblock 34:58--79, 1979.

\bibitem{BuchsbaumEisenbud73}
D.A. Buchsbaum and D.~Eisenbud.
\newblock Remarks on ideals and resolutions.
\newblock In {\em Symposia {M}athematica, {V}ol. {XI} ({C}onvegno di {A}lgebra
  {C}ommutativa, {INDAM}, {R}ome, 1971)}, pages 193--204. 1973.

\bibitem{BuchsbaumEisenbud77}
D.A. Buchsbaum and D.~Eisenbud.
\newblock Algebra structures for finite free resolutions, and some structure
  theorems for ideals of codimension {$3$}.
\newblock {\em Amer. J. Math.}, 99(3):447--485, 1977.

\bibitem{Buchweitz81}
R.~Buchweitz.
\newblock {\em Contributions {\`a} la th{\'e}orie des singularit{\'e}s}.
\newblock Th{\`e}se, 1981.

\bibitem{Burch68}
L.~Burch.
\newblock On ideals of finite homological dimension in local rings.
\newblock {\em Proc. Cambridge Philos. Soc.}, 64:941--948, 1968.

\bibitem{Volume1}
J.L. Cisneros~Molina, L\^e~D{\~u}ng Tr{\'a}ng, and J.~Seade, editors.
\newblock {\em {H}andbook of {G}eometry and {T}opology of {S}ingularities {I}}.
\newblock {S}pringer {V}erlag, 2020.

\bibitem{Pereira10}
M.~da~Silva~Pereira.
\newblock {\em Variedades determinantais e singularidades de matrizes}.
\newblock Phd-thesis, 2010.

\bibitem{Damon84}
J.~Damon.
\newblock {\em The unfolding and determinacy theorems for subgroups of A and
  K}, volume~50 of {\em Memoirs of the American Mathematical Society}.
\newblock American Math. Soc., 1984.

\bibitem{Damon87}
J.~Damon.
\newblock Deformations of sections of singularities and {G}orenstein surface
  singularities.
\newblock {\em Amer. J. Math.}, 109(4):695--721, 1987.

\bibitem{Damon96}
J.~Damon.
\newblock Higher multiplicities and almost free divisors and complete
  intersections.
\newblock {\em Mem. Amer. Math. Soc.}, 123(589):x+113, 1996.

\bibitem{Damon01}
J.~Damon.
\newblock Nonlinear sections of nonisolated complete intersections.
\newblock In {\em New Developments in Singularity Theory}, volume~21 of {\em
  NATO Conf. Series}, pages 405--445. Kluwer, 2001.

\bibitem{Damon16}
J.~Damon.
\newblock {T}opology of {E}xceptional {O}rbit {H}ypersurfaces of
  {P}rehomogeneous {S}paces.
\newblock {\em {J}ournal of {T}opology}, 9:797--825, 2016.

\bibitem{Damon18}
J.~Damon.
\newblock Schubert decomposition for {M}ilnor fibers of the varieties of
  singular matrices.
\newblock {\em J. Singul.}, 18:358--396, 2018.

\bibitem{DamonMond91}
J.~Damon and D.~Mond.
\newblock A-codimension and the vanishing topology of discriminants.
\newblock {\em {I}nvent. {M}ath.}, (106):217--242, 1991.

\bibitem{DamonPike14}
J.~Damon and B.~Pike.
\newblock Solvable groups, free divisors and nonisolated matrix singularities
  {II}:vanishing topology.
\newblock {\em Geom. Topol.}, 18:911--962, 2014.

\bibitem{DamonPike15}
J.~Damon and B.~Pike.
\newblock Solvable groups, free divisors and nonisolated matrix singularities
  {I}: {T}owers of free divisors.
\newblock {\em Ann. Inst. Fourier (Grenoble)}, 65(3):1251--1300, 2015.

\bibitem{deJong98}
T.~de~Jong.
\newblock Determinantal rational surface singularities.
\newblock {\em Compositio Math.}, 113(1):67--90, 1998.

\bibitem{EagonNorthcott62}
J.~Eagon and D.G. Northcott.
\newblock Ideals defined by matrices, and a certain complex associated to them.
\newblock {\em Proc. Royal Soc.}, 269:188--204, 1962.

\bibitem{EagonNorthcott67}
J.A. Eagon and D.G. Northcott.
\newblock Generically acyclic complexes and generically perfect ideals.
\newblock 299 A:147--172, 1967.

\bibitem{EbelingGuseinZade09}
W.~Ebeling and S.~Gusein-Zade.
\newblock On indices of 1-forms on determinantal singularities.
\newblock 267:113--124, 2009.

\bibitem{Eisenbud95}
D.~Eisenbud.
\newblock {\em Commutative Algebra with a View Toward Algebraic Geometry}.
\newblock Springer-Verlag New York, 1995.

\bibitem{FruehbisKrueger99}
A.~Fr{\"u}hbis-Kr{\"u}ger.
\newblock Classification of simple space curve singularities.
\newblock {\em Comm. Algebra}, 27(8):3993--4013, 1999.

\bibitem{FruehbisKrueger18}
A.~Fr\"{u}hbis-Kr\"{u}ger.
\newblock On discriminants, {T}jurina modifications and the geometry of
  determinantal singularities.
\newblock {\em Topology Appl.}, 234:375--396, 2018.

\bibitem{FruehbisKruegerNeumer10}
A.~{Fr\"uhbis-Kr\"uger} and A.~{Neumer}.
\newblock {Simple Cohen-Macaulay Codimension 2 Singularities}.
\newblock {\em Comm. in Alg.}, 38(2):454--495, 2010.

\bibitem{FruehbisKruegerZach20}
A.~Fr\"uhbis-Kr\"uger and M.~Zach.
\newblock {O}n the {V}anishing {T}opology of {I}solated {C}ohen-{M}acaulay
  codimension 2 {S}ingularities.
\newblock 2020.
\newblock To appear in {G}eometry \& {T}opology.

\bibitem{Gaffney93}
T.~Gaffney.
\newblock Polar multiplicities and equisingularity of map germs.
\newblock {\em Topology}, 32(1):185--223, 1993.

\bibitem{GaffneyGrulhaRuas19}
T.~Gaffney, N.~G. Grulha, Jr., and M.~A.~S. Ruas.
\newblock The local {E}uler obstruction and topology of the stabilization of
  associated determinantal varieties.
\newblock {\em Math. Z.}, 291(3-4):905--930, 2019.

\bibitem{GaffneyRuas21}
T.~Gaffney and M.A.S. Ruas.
\newblock Equisingularity and {EIDS}.
\newblock {\em Proc. Amer. Math. Soc.}, 149(4):1593--1608, 2021.

\bibitem{Giu}
M.~Giusti.
\newblock Classification des singularit\'es isol\'ees simples d'intersections
  compl\`etes.
\newblock {\em Proc. Symp Pure Math.}, 40:457--494, 1983.

\bibitem{GoreskyMacPherson88}
M.~Goresky and R.~MacPherson.
\newblock {\em Stratified Morse Theory}.
\newblock Springer Verlag, 1988.

\bibitem{Goryunov21}
V.~Goryunov.
\newblock Vanishing cycles of matrix singularities.
\newblock {\em {J}. {L}ondon {M}ath. {S}oc.}, 103:991 -- 1015, 2021.

\bibitem{GoryunovMond05}
V.~Goryunov and D.~Mond.
\newblock Tjurina and milnor numbers of matrix singularities.
\newblock {\em J. London Math. Soc.}, 72((2)):205 -- 224, 2005.

\bibitem{GoryunovZakalyukin03}
V.~V. Goryunov and V.~M. Zakalyukin.
\newblock Simple symmetric matrix singularities and the subgroups of {W}eyl
  groups {$A_\mu$}, {$D_\mu$}, {$E_\mu$}.
\newblock volume~3, pages 507--530, 743--744. 2003.
\newblock Dedicated to Vladimir I. Arnold on the occasion of his 65th birthday.

\bibitem{Grauert72}
H.~Grauert.
\newblock \"uber die deformation isolierter singularit\"aten analytischer
  mengen.
\newblock 15:171 -- 198, 1972.

\bibitem{GreuelMartinPfister85}
G.-M. Greuel, B.~Martin, and G.~Pfister.
\newblock Numerische {C}harakterisierung quasihomogener
  {G}orenstein-{K}urvensingularit\"{a}ten.
\newblock {\em Math. Nachr.}, 124:123--131, 1985.

\bibitem{GreuelSteenbrink83}
G.-M. Greuel and J.~Steenbrink.
\newblock On the topology of smoothable singularities.
\newblock In {\em Singularities, {P}art 1 ({A}rcata, {C}alif., 1981)},
  volume~40 of {\em Proc. Sympos. Pure Math.}, pages 535--545. Amer. Math.
  Soc., Providence, R.I., 1983.

\bibitem{Greuel80}
G.M. Greuel.
\newblock Dualit\"at in der lokalen {K}ohomologie isolierter
  {S}ingularit\"aten.
\newblock {\em Math. Ann.}, 250(2):157--173, 1980.

\bibitem{GreuelLossenShustin07}
G.M. Greuel, C.~Lossen, and E.~Shustin.
\newblock {\em Introduction to Singularities and Deformations}.
\newblock Springer Monographs in Mathematics. Springer, 2007.

\bibitem{GulliksenNegard72}
T.H. Gulliksen and O.G. Neg\.{a}rd.
\newblock Un complexe r\'{e}solvant pour certains id\'{e}aux
  d\'{e}terminantiels.
\newblock {\em C. R. Acad. Sci. Paris S\'{e}r. A-B}, 274:A16--A18, 1972.

\bibitem{Hamm72}
H.~Hamm.
\newblock Lokale topologische eigenschaften komplexer r\"aume.
\newblock {\em Math. Ann.}, 191:235--252, 1972.

\bibitem{Haslinger01}
G.~Haslinger.
\newblock {\em Families of skew-symmetric matrices}.
\newblock Phd-thesis, 2001.

\bibitem{Hilbert90}
D.~Hilbert.
\newblock Ueber die {T}heorie der algebraischen {F}ormen.
\newblock {\em Math. Ann.}, 36(4):473--534, 1890.

\bibitem{EagonHochster71}
M.~Hochster and J.A. Eagon.
\newblock Cohen-macaulay rings, invariant theory, and the generic perfection of
  determinantal loci.
\newblock 93:1020--1058, 1971.

\bibitem{Jaehner74}
U.~J{\"a}hner.
\newblock {\em Beispiele starrer analytischer Algebren}.
\newblock Dissertation, 1974.

\bibitem{Jozefiak78}
T.~J\'{o}zefiak.
\newblock Ideals generated by minors of a symmetric matrix.
\newblock {\em Comment. Math. Helv.}, 53(4):595--607, 1978.

\bibitem{JozefiakPragacz79}
T.~J\'{o}zefiak and P.~Pragacz.
\newblock Ideals generated by {P}faffians.
\newblock {\em J. Algebra}, 61(1):189--198, 1979.

\bibitem{KasSchlessinger72}
A.~Kas and M.~Schlessinger.
\newblock On the versal deformation of a complex space with an isolated
  singularity.
\newblock {\em Math. Ann.}, 196:23--29, 1972.

\bibitem{KatoMatsumoto75}
M.~Kato and Y.~Matsumoto.
\newblock On the connectivity of the milnor fiber of a holomorphic function at
  a critical point.
\newblock {\em Proc. Internat. Conf., Tokyo}, 1975.

\bibitem{Kerner19}
D.~Kerner.
\newblock Group actions on matrices over local rings.
\newblock 2019.

\bibitem{Kirby74}
D.~Kirby.
\newblock A sequence of complexes associated with a matrix.
\newblock {\em J. London Math. Soc. (2)}, 7:523--530, 1974.

\bibitem{KleppeLaksov80}
H.~Kleppe and D.~Laksov.
\newblock The algebraic structure and deformation of {P}faffian schemes.
\newblock {\em J. Algebra}, 64(1):167--189, 1980.

\bibitem{Kurano89}
Kazuhiko Kurano.
\newblock The first syzygies of determinantal ideals.
\newblock {\em J. Algebra}, 124(2):414--436, 1989.

\bibitem{Kutz74}
R.E. Kutz.
\newblock Cohen-{M}acaulay rings and ideal theory in rings of invariants of
  algebraic groups.
\newblock {\em Trans. Amer. Math. Soc.}, 194:115--129, 1974.

\bibitem{Laksov75}
D.~Laksov.
\newblock Deformation of determinantal schemes.
\newblock {\em Compositio Mathematica}, 30:273--292, 1975.

\bibitem{Lascoux78}
A.~Lascoux.
\newblock Syzygies des vari\'et\'es d\'eterminantales.
\newblock 30:202--237, 1978.

\bibitem{Laufer72}
H.B. Laufer.
\newblock On rational singularities.
\newblock {\em Amer. J. Math.}, 94:597--608, 1972.

\bibitem{Laufer73}
H.B. Laufer.
\newblock Taut two-dimensional singularities.
\newblock {\em Math. Ann.}, 205:131--164, 1973.

\bibitem{Le77}
{L\^e D\~ung Tr\'ang}.
\newblock Some remarks on relative monodromy.
\newblock {\em Proc. Ninth Nordic Summer School/NAVF Sympos. Math., Oslo 1976},
  pages 397--403, 1977.

\bibitem{Le86}
{L\^e D\~ung Tr\'ang}.
\newblock Le concept du singularit\'e isol\'ee de fonction analytique.
\newblock {\em Adv. Stud. Pure Math.}, 8:215 -- 227, 1986.

\bibitem{Lojasiewicz64}
S.~Lojasiewicz.
\newblock Triangularitions of semi analytic sets.
\newblock {\em Ann. Scoula Norm. Sup. Pisa}, 18(3):449--474, 1964.

\bibitem{LooijengaSteenbrink85}
E.~Looijenga and J.~Steenbrink.
\newblock Milnor number and {T}jurina number of complete intersections.
\newblock {\em Math. Ann.}, 271(1):121--124, 1985.

\bibitem{Ma93}
Yonghao Ma.
\newblock The first syzygies of determinantal ideals.
\newblock {\em J. Pure Appl. Algebra}, 85(1):73--103, 1993.

\bibitem{Macaulay16}
F.~S. Macaulay.
\newblock {\em The algebraic theory of modular systems}.
\newblock Cambridge Mathematical Library. Cambridge University Press,
  Cambridge, 1994.
\newblock Revised reprint of the 1916 original, With an introduction by Paul
  Roberts.

\bibitem{MacPherson74}
R.~D. MacPherson.
\newblock Chern classes for singular algebraic varieties.
\newblock {\em Ann. of Math. (2)}, 100:423--432, 1974.

\bibitem{Milnor68}
J.W. Milnor.
\newblock {\em Singular points of complex hypersurfaces}.
\newblock Princenton University Press, 1968.

\bibitem{Pedersen21}
H.~M{\o}ller~Pedersen.
\newblock On {T}jurina {T}ransform and {R}esulotion of {D}eterminantal
  {S}ingularities.
\newblock In J.F de~Bobadilla, T.~L\'aszl\'o, and A.~Stipsicz, editors, {\em
  Singularities and {T}heir {i}nteraction with {G}eometry and {L}ow
  {D}imensional {T}opology}, pages 289 -- 312. Birkh\"auser Verlag, 2021.

\bibitem{BallesterosMond20}
J.J. Nu\~no Ballesteros and D.~Mond.
\newblock {\em Singularities of {M}appings}, volume 357 of {\em {G}rundlehren
  der mathematischen {W}issenschaften}.
\newblock Springer Verlag, 2020.

\bibitem{BallesterosOreficeTomazella13}
J.J. Nu\~no Ballesteros, B.~Or\'efice-Okamoto, and J.N. Tomazella.
\newblock The vanishing {E}uler characteristic of an isolated determinantal
  singularity.
\newblock {\em Israel J. Math.}, 197:475--495, 2013.

\bibitem{Pinkham73}
H.C. Pinkham.
\newblock Deformations of cones with negative grading.
\newblock {\em Journal of Algebra}, 30:92--102, 1974.

\bibitem{PragaczWeyman86}
P.~Pragacz and J.~Weyman.
\newblock On the construction of resolutions of determinantal ideals: a survey.
\newblock In {\em S\'{e}minaire d'alg\`ebre {P}aul {D}ubreil et {M}arie-{P}aule
  {M}alliavin, 37\`eme ann\'{e}e ({P}aris, 1985)}, volume 1220 of {\em Lecture
  Notes in Math.}, pages 73--92. Springer, Berlin, 1986.

\bibitem{Roehr92}
A.~R{\"o}hr.
\newblock {\em Formate rationaler {F}l{\"a}chensingularit{\"a}ten}.
\newblock Dissertation, 1992.

\bibitem{Saito71}
Kyoji Saito.
\newblock Quasihomogene isolierte {S}ingularit\"{a}ten von {H}yperfl\"{a}chen.
\newblock {\em Invent. Math.}, 14:123--142, 1971.

\bibitem{Saito80}
Kyoji Saito.
\newblock Theory of logarithmic differential forms and logarithmic vector
  fields.
\newblock {\em J. Fac. Sci. Univ. Tokyo Sect. IA Math.}, 27(2):265--291, 1980.

\bibitem{Sard42}
A.~Sard.
\newblock The measure of the critical values of differentiable maps.
\newblock {\em Bulletin of the American Mathematical Society}, 48(12):883--890,
  1942.

\bibitem{Schaps77}
M.~Schaps.
\newblock Deformations of {C}ohen-{M}acaulay schemes of codimension 2 and
  nonsingular deformations of space curves.
\newblock {\em Amer. J. Math.}, 99:669 --684, 1977.

\bibitem{Schaps83}
M.~Schaps.
\newblock Versal determinantal deformations.
\newblock {\em Pacific Journal of Mathematics}, 107(1), 1983.

\bibitem{Schlessinger68}
M.~Schlessinger.
\newblock Functors of {A}rtin rings.
\newblock {\em Trans. AMS}, 130:536--577, 1968.

\bibitem{Siersma91_2}
D.~Siersma.
\newblock Vanishing cycles and special fibres.
\newblock In {\em Singularity theory and its applications, {P}art {I}
  ({C}oventry, 1988/1989)}, volume 1462 of {\em Lecture Notes in Math.}, pages
  292--301. Springer, Berlin, 1991.

\bibitem{Svanes72}
T.~Svanes.
\newblock {\em C{OHERENT} {COHOMOLOGY} {ON} {FLAG} {MANIFOLDS} {AND}
  {RIGIDITY}}.
\newblock ProQuest LLC, Ann Arbor, MI, 1972.
\newblock Thesis (Ph.D.)--Massachusetts Institute of Technology.

\bibitem{Teissier82}
Bernard Teissier.
\newblock Vari\'{e}t\'{e}s polaires. {II}. {M}ultiplicit\'{e}s polaires,
  sections planes, et conditions de {W}hitney.
\newblock In {\em Algebraic geometry ({L}a {R}\'{a}bida, 1981)}, volume 961 of
  {\em Lecture Notes in Math.}, pages 314--491. Springer, Berlin, 1982.

\bibitem{Thom54}
R.~Thom.
\newblock Quelques propri\'{e}t\'{e}s globales des vari\'{e}t\'{e}s
  diff\'{e}rentiables.
\newblock {\em Comment. Math. Helv.}, 28:17--86, 1954.

\bibitem{Tibar95}
M.~Tib{\u a}r.
\newblock Bouquet decomposition of the {M}ilnor fibre.
\newblock {\em Topology}, 35:227--241, 1995.

\bibitem{Tjurina68}
G.~N. Tjurina.
\newblock Absolute isolatedness of rational singularities and triple rational
  points.
\newblock {\em Func. Anal. Appl.}, 2:324--333, 1968.

\bibitem{Tjurina69}
G.~N. Tjurina.
\newblock Locally semi-universal flat deformations of isolated singularities of
  complex spaces.
\newblock {\em Izv. Akad. Nauk SSSR Ser. Mat.}, 33:1026--1058, 1969.

\bibitem{Tjurina70}
G.~N. Tjurina.
\newblock Resolution of singularities of flat deformations of double rational
  points.
\newblock {\em Funkcional. Anal. i Prilo\v{z}en.}, 4(1):77--83, 1970.

\bibitem{LeTeissier81}
L\^e~Dung Tr\`ang and B.~Teissier.
\newblock Vari\'{e}t\'{e}s polaires locales et classes de {C}hern des
  vari\'{e}t\'{e}s singuli\`eres.
\newblock {\em Ann. of Math. (2)}, 114(3):457--491, 1981.

\bibitem{Trivedi13}
S.~Trivedi.
\newblock Stratified transversality of holomorphic maps.
\newblock {\em International Journal of Mathematics}, 24(13), 2013.

\bibitem{Vosegaard02}
H.~Vosegaard.
\newblock A characterization of quasi-homogeneous complete intersection
  singularities.
\newblock {\em J. Algebraic Geom.}, 11(3):581--597, 2002.

\bibitem{Wahl16}
J.~Wahl.
\newblock Milnor and {T}jurina numbers for smoothings of surface singularities.
\newblock 2016.

\bibitem{Wahl77}
J.M. Wahl.
\newblock Equations defining rational singularities.
\newblock {\em Ann. Sci. \'{E}cole Norm. Sup. (4)}, 10(2):231--263, 1977.

\bibitem{Wahl79}
J.M. Wahl.
\newblock Simultaneous resolution and discriminantal loci.
\newblock {\em Duke Math. J.}, 46(2):341--375, 1979.

\bibitem{Wahl85}
Jonathan~M. Wahl.
\newblock A characterization of quasihomogeneous {G}orenstein surface
  singularities.
\newblock {\em Compositio Math.}, 55(3):269--288, 1985.

\bibitem{Waldi79}
R.~Waldi.
\newblock Deformation von {G}orenstein-{S}ingularit\"{a}ten der {K}odimension
  {$3$}.
\newblock {\em Math. Ann.}, 242(3):201--208, 1979.

\bibitem{Wall81}
C.~T.~C. Wall.
\newblock Finite determinacy of smooth map-germs.
\newblock {\em Bull. London Math. Soc.}, 13(6):481--539, 1981.

\bibitem{Zach17}
M.~Zach.
\newblock {\em Topology of isolated determinantal singularities}.
\newblock Phd-thesis, 2017.

\bibitem{Zach18}
M.~Zach.
\newblock An observation concerning the vanishing topology of certain isolated
  determinantal singularities.
\newblock {\em {M}athematische {Z}eitschrift}, 2018.

\bibitem{Zach18BouquetDecomp}
M.~Zach.
\newblock {B}ouquet decomposition for {D}eterminantal {M}ilnor fibers.
\newblock {\em Journal of Singularities}, 22:190 -- 204, 2020.

\bibitem{Zhang17}
X.~Zhang.
\newblock {L}ocal {E}uler {O}bstruction and {C}hern-{M}ather classes of
  {D}eterminantal {V}arieties.
\newblock 2017.

\end{thebibliography}
\printindex

\end{document}